\documentclass[11pt]{article}
\usepackage[english]{babel}
\usepackage{amsmath,amsthm,amssymb}
\usepackage[all]{xy}
\usepackage{setspace}
\usepackage{color}
\definecolor{grau}{rgb}{0.3,0.3,0.3}
\usepackage[colorlinks, linkcolor=grau, citecolor=grau, urlcolor=grau]{hyperref}
\usepackage{breakurl}
\usepackage{bibspacing}
\usepackage[top=1.3in, bottom=1.6in, left=1.3in, right=1.3in]{geometry}
\usepackage{booktabs}
\usepackage{enumerate}

\newtheorem{theorem}{Theorem}[section]

\newtheorem{corollary}[theorem]{Corollary}
\newtheorem{lemma}[theorem]{Lemma}

\newtheorem{proposition}[theorem]{Proposition}
\newtheorem{Algorithm}[theorem]{Algorithm}

\theoremstyle{definition}
\newtheorem{definition}[theorem]{Definition}

\newtheorem{remark}[theorem]{Remark}

\numberwithin{equation}{section}

\def\sp{\textnormal{Spec}}

\def\P1k{\mathbb P_{k}^{1}}
\def\deg{\textnormal{deg}}

\def\sl2{\textnormal{SL}_2}

\def\sp{\textnormal{Spec}}

\newcommand {\QQ}  {{\mathbb Q}}
\newcommand {\ZZ}  {{\mathbb Z}}
\newcommand {\ord} {{\textnormal{ord}}}
\newcommand {\CC}  {{\mathbb C}}
\newcommand {\RR}  {{\mathbb R}}

\newcommand {\FF}  {{\mathbb F}}

\newcommand{\rad}{\textnormal{rad}}

\newcommand{\ee} {e}
\newcommand{\ellrad}{b}
\newcommand{\indpar}{\nu}

\newcommand{\parshin}{Par\v{s}in}

\newcommand{\wo}{\backslash} 

\newcommand{\floor}[1]{\left\lfloor{#1}\right\rfloor} 
\newcommand{\ceil}[1]{\left\lceil{#1}\right\rceil} 
\newcommand{\incl}{\hookrightarrow} 

\DeclareMathOperator{\st}{: }

\newcommand{\eps}{\varepsilon}

\newcommand{\id}{\textnormal{id}}

\DeclareMathOperator{\divides}{|}

\DeclareMathOperator{\vol}{\textnormal{vol}}

\DeclareMathOperator{\dotcup}{\dot{\cup}}

\newcommand{\abs}[1]{\lvert{}#1{}\rvert}

\renewcommand{\st}{\,;\,}

\date{}

\begin{document}

\author{Rafael von K\"anel and Benjamin Matschke}
\title{{\Large Solving $S$-unit, Mordell,  Thue, Thue--Mahler and generalized Ramanujan--Nagell equations via Shimura--Taniyama conjecture}}
\maketitle


\begin{abstract} {\scriptsize In the first part we construct  algorithms (over $\QQ$) which we apply to solve $S$-unit, Mordell, cubic Thue, cubic Thue--Mahler and generalized Ramanujan--Nagell equations. As a byproduct we obtain alternative practical approaches for various classical Diophantine problems, including the fundamental problem of finding all elliptic curves over $\QQ$ with good reduction outside a given finite set of rational primes. The first type of our algorithms uses modular symbols, and the second type combines explicit height bounds with  efficient sieves. In particular we construct a refined sieve for $S$-unit equations which combines Diophantine approximation techniques of de Weger with new geometric ideas.  To illustrate the utility of our algorithms we determined the solutions of large classes of equations, containing many examples of interest which are out of reach for the known  methods. In addition we used the resulting data to motivate various  conjectures and questions, including Baker's explicit $abc$-conjecture
and a new conjecture on  $S$-integral points of any  hyperbolic genus one curve over $\QQ$.}

{\scriptsize In the second part we establish new results for certain old Diophantine problems (e.g. the difference of squares and cubes) related to Mordell equations, and we prove explicit height bounds for cubic Thue, cubic Thue--Mahler and generalized Ramanujan--Nagell equations. As a byproduct, we obtain here an alternative proof of classical theorems of Baker, Coates and Vinogradov--Sprind{\v{z}}uk. In fact we get refined versions of their theorems, which improve the actual best results in many fundamental cases. We also conduct some effort to work out optimized height bounds for $S$-unit and Mordell equations which are used in our algorithms of the first part. Our results and algorithms all ultimately rely on the method of  Faltings (Arakelov, \parshin, Szpiro) combined with the Shimura--Taniyama conjecture, and they all do not use lower bounds for linear forms in (elliptic) logarithms.}

{\scriptsize In the third part we solve the problem of constructing an efficient sieve for the $S$-integral points of bounded height on any elliptic curve $E$ over $\QQ$ with given Mordell--Weil basis of $E(\QQ)$.  Here we combine a geometric interpretation of the known elliptic logarithm reduction (initiated by Zagier) with several conceptually new ideas. The resulting ``elliptic logarithm sieve" is crucial for some of our algorithms of the first part. Moreover, it considerably extends the class of elliptic Diophantine equations which can be solved in practice: To demonstrate this we solved many notoriously difficult equations by combining our sieve with known height bounds based on the theory of  logarithmic forms.}
\end{abstract}

\newpage

{\scriptsize\tableofcontents}

\onehalfspacing

\newpage

\section{Introduction}

In this paper we combine the method of Faltings~\cite{faltings:finiteness} (Arakelov, \parshin, Szpiro)  with the Shimura--Taniyama conjecture \cite{wiles:modular,taywil:modular,breuil:modular} in order to study various classical Diophantine problems, including $S$-unit equations, Mordell equations, cubic Thue equations, cubic Thue--Mahler equations and generalized Ramanujan--Nagell equations. 

We now begin to discuss the Diophantine equations. Let $S$ be a finite set of rational prime numbers. Write $N_S=1$ if $S$ is empty and $N_S=\prod_{p\in S} p$ otherwise. We denote by $\mathcal O^\times$ the units of
$\mathcal O =\ZZ[1/N_S]$ and we consider the  $S$-unit equation
\begin{equation}\label{eq:sunit}
x+y=1, \ \ \ (x,y)\in\mathcal O^\times\times\mathcal O^\times.
\end{equation}
Many important Diophantine conjectures can be reduced to the study of $S$-unit equations. For example, the $abc$-conjecture of Masser--Oesterl\'e is equivalent to a certain height bound for the solutions of \eqref{eq:sunit}. On using Diophantine approximations in the style of Thue and Siegel, Mahler~\cite{mahler:approx1}  showed that~\eqref{eq:sunit} has only finitely many solutions.  
Furthermore  there already exists a practical method of de Weger~\cite{deweger:phdthesis} which solves $S$-unit equations by using the theory of logarithmic forms~\cite{bawu:logarithmicforms}, see Section~\ref{sec:suheightalgo}. For a detailed discussion of (general) $S$-unit equations we refer to the recent book of Evertse--Gy{\H{o}}ry~\cite{evgy:bookuniteq}. Next we take a nonzero $a\in \mathcal O$ and we consider the Mordell equation 
\begin{equation}\label{eq:mordell}
y^2=x^3+a, \ \ \ (x,y)\in\mathcal O\times\mathcal O.
\end{equation}
This Diophantine equation is a priori more difficult than~\eqref{eq:sunit}. Further, if $\mathcal O=\ZZ$ then resolving~\eqref{eq:mordell}  is equivalent to solving the classical problem, going back at least to Bachet (1621), of finding all perfect squares and perfect cubes with given difference.  In the case $\mathcal O=\ZZ$,  Mordell \cite{mordell:1922,mordell:1923} showed finiteness of~\eqref{eq:mordell} via Diophantine approximation, and Baker--Davenport~\cite{baker:contributions,bada:diophapp} and Masser, Zagier~\cite{masser:ellfunctions,zagier:largeintegralpoints} introduced practical approaches solving \eqref{eq:mordell} via the theory of logarithmic forms (see Section~\ref{sec:malgo}).  Furthermore we shall see that a special class of Mordell equations~\eqref{eq:mordell} covers in particular any generalized Ramanujan--Nagell equation discussed in \eqref{eq:rana} below.  Finally,  we let $m\in\mathcal O$ be nonzero  and we suppose that $f\in \mathcal O[x,y]$ is a homogeneous polynomial of degree three with nonzero discriminant. Consider the cubic Thue equation
\begin{equation}\label{eq:thue}
f(x,y)=m, \ \ \ (x,y)\in \mathcal O\times \mathcal O.
\end{equation} 
Thue (1909) proved that (\ref{eq:thue}) has only finitely many solutions in the case $\mathcal O=\ZZ$. In general, equation \eqref{eq:thue} is essentially equivalent to the  cubic Thue--Mahler equation recalled in~\eqref{eq:thue-mahler} below.  Baker--Davenport~\cite{baker:contributions,bada:diophapp} and Tzanakis--de Weger~\cite{tzde:thue,tzde:thuemahler} obtained practical approaches solving in particular cubic Thue and Thue--Mahler equations via the theory of logarithmic forms~\cite{bawu:logarithmicforms}, see also the discussions in Section~\ref{sec:thuealgo}.

\subsection{Algorithms}\label{sec:ia}

We construct two types of algorithms which we use to solve  $S$-unit equations \eqref{eq:sunit}, Mordell equations \eqref{eq:mordell}, cubic Thue equations~\eqref{eq:thue}, cubic Thue--Mahler equations~\eqref{eq:thue-mahler} and generalized Ramanujan--Nagell equations~\eqref{eq:rana}. Both types do not use the theory of logarithmic forms. Before we discuss our algorithms in more detail, we describe the general strategy. 

\subsubsection{General strategy}\label{sec:igeneralstrategy} As in \cite{rvk:modular} we use the method of Faltings (Arakelov, \parshin, Szpiro) which in our situation is applied as follows: Let $Y(\mathcal O)$ be the set of solutions of any of the above equations.  
Then there is an effective map $\phi$ (\parshin{} construction)  from~$Y(\mathcal O)$
to the set $M(T)$ of isomorphism classes of elliptic curves over a controlled open $T\subset\sp(
\ZZ)$,   $$\phi: Y(\mathcal O)\to M(T).$$ Here effective means that one can compute $\phi^{-1}(E)$ for each $E$ in $M(T)$.  To determine $Y(\mathcal O)$, it thus suffices to compute $M(T)$ (effective Shafarevich theorem). For this purpose we use two types of algorithms: The first type  applies Cremona's algorithm~\cite{cremona:algorithms} involving modular symbols, and the second type combines our optimized height bounds (see Section~\ref{sec:ihb}) with efficient sieves. Both types of algorithms crucially rely on a geometric version of the Shimura--Taniyama conjecture~\cite[Thm A]{breuil:modular}  using inter alia the Tate conjecture~\cite[Thm 4]{faltings:finiteness}, and on isogeny estimates based on the method of Mazur~\cite{mazur:qisogenies,kenku:ellisogenies} or Faltings~\cite{faltings:finiteness,raynaud:abelianisogenies}. In fact the strategy of combining modularity with Faltings' method gives effective finiteness results for considerably more general Diophantine problems, see \cite{rvk:modular,vkkr:intpointsshimura}. However in the present paper we focus on optimizing the strategy for the fundamental Diophantine  equations appearing in the title, and in the future we plan to work out algorithms for other Diophantine problems of interest.

\subsubsection{Algorithms via modular symbols} We next discuss in more detail our first type of algorithms. They crucially rely on Cremona's algorithm~\cite{cremona:algorithms} using modular symbols in order to compute all elliptic curves over $\QQ$ of given conductor. This allows to determine  $M(T)$, since the  curves in $M(T)$ have bounded conductor. Then we compute $Y(\mathcal O)$ by enumerating $\phi^{-1}(E)$ for each $E$ in $\phi(Y(\mathcal O))$. Here we exploit that the maps $\phi$ are effective by classical constructions going back at least to Cayley, Mordell and Frey--Hellegouarch. To illustrate the utility of our first type of algorithms we computed several examples. For instance we solved  the $S$-unit equation~\eqref{eq:sunit} for all  sets $S$ with $N_S\leq 20000$, and  we solved several Mordell equations~\eqref{eq:mordell}. In fact, as already pointed out in \cite[Sect 1.1.2]{rvk:modular}, our first type of algorithms  can in principle solve any Diophantine equation inducing integral points on a moduli scheme of elliptic curves  with an effective \parshin{}  construction $\phi$, see \eqref{eq:mdisplay} and Section~\ref{sec:moduli}. This class of equations contains in particular all equations considered in this paper. Here we mention that the possibility of solving certain cubic Thue--Mahler equations via modular symbols was already discussed in Bennett--Dahmen~\cite[$\mathsection$14]{beda:kleinsuperell},  see also the recent works of Kim~\cite{kimd:modularthuemahler} and Bennett--Billerey~\cite{bebi:sumsofunits}.  Our first type of algorithms is very fast for ``small" parameters, since in this case we can use Cremona's database listing all elliptic curves over $\QQ$ of conductor at most 350000 (as of August 2014).  In particular these algorithms directly benefit from the ongoing extension of such databases. However the approach via modular symbols can usually not compete with the actual most efficient methods solving our equations of interest. Thus we worked out a second type of algorithms. 

\subsubsection{Algorithms via height bounds}\label{sec:ialgoheight} We now give a more detailed description of our second type of algorithms. They rely on an effective Shafarevich theorem in the form of explicit bounds for the Faltings height $h_F$ on~$M(T)$. The space $M(T)$ can be very complicated and it is usually  a difficult task to compute $M(T)$, see the discussions surrounding \eqref{eq:mdisplay}. Hence instead of first computing $M(T)$ and then $Y(\mathcal O)=\phi^{-1}(M(T))$, we often directly work with $Y(\mathcal O)$ by using that the height $\phi^*h_F$ is bounded on~$Y(\mathcal O)$. This has the advantage that we can exploit  extra structures on $Y(\mathcal O)$ in order  to construct efficient sieves for  solutions of bounded height. 

\paragraph{$S$-unit equation.} To solve $S$-unit equations~\eqref{eq:sunit} it is natural to consider the set $\Sigma(S)$  of solutions  of~\eqref{eq:sunit} modulo symmetry. Here solutions $(x,y)$ and $(x',y')$ of~\eqref{eq:sunit} are called symmetric if $x'$ or $y'$ lies in $\{x,\frac{1}{x},\frac{1}{1-x}\}$. One can directly write down all solutions which are symmetric to a given solution and  thus it suffices to determine $\Sigma(S)$ in order to solve \eqref{eq:sunit}. In fact the number of solutions of~\eqref{eq:sunit} is either zero or $6\abs{\Sigma(S)}-3$ where $\abs{\Sigma}$ denotes the cardinality of a set $\Sigma$. For any $n\in\ZZ_{\geq 1}$ we denote by $S(n)$ the set of the $n$ smallest rational primes. Our Algorithm~\ref{algo:suheight}  allows to efficiently solve~\eqref{eq:sunit}, even for sets $S$ with relatively large~$\abs{S}$. To demonstrate this we solved large classes of $S$-unit equations~\eqref{eq:sunit} by using Algorithm~\ref{algo:suheight}. In particular, we obtained the following result.

\vspace{0.4cm}
\noindent{\bf Theorem~A.}
\emph{Suppose that $n\in\{1,2,\ldots,16\}$. Then the cardinality $\#$ of the set $\Sigma(S(n))$ is given in the following table.}
\begin{table}[h]
\begin{center}
{\small
\begin{tabular}{lccccccccccccccc}
$n$ & $1$ & $2$ & $3$ & $4$ & $5$ & $6$ & $7$ & $8$\\
\cmidrule(r){2-9}
\# & $1$ & $4$ & $17$ & $63$ & $190$ & $545$ & $1433$ & $3649$\\

\\

$n$  & $9$ & $10$ & $11$ & $12$ & $13$ & $14$ & $15$ & $16$\\ 
\cmidrule(r){2-9}
\# & $8828$ & $20015$ & $44641$ & $95358$ & $199081$ & $412791$ & $839638$ & $1234567$
\end{tabular}
}
\end{center}
\end{table}
\\
We mention that among all sets $S$ of cardinality $n$ the set $S(n)$ is usually the most difficult case for solving \eqref{eq:sunit}. 
 In \cite[$\mathsection$10]{zagier:polylogs}, Zagier explained that the cardinality of $\Sigma(S)$ plays an important role in certain questions on polylogarithms. In particular he states a table attributed to Gross--Vojta, which for each $n\in\{1,\dotsc,8\}$ lists a lower bound for the cardinality of $\Sigma(S(n))$. Theorem~A proves that the entries of this table are not only lower bounds, but in fact the correct values. Further we point out that the cases $n\in\{1,\dotsc,6\}$ in Theorem~A are not new. They were previously known by the work of de Weger \cite[Thm 5.4]{deweger:lllred}. Our Algorithm~\ref{algo:suheight} substantially improves de Weger's method in~\cite{deweger:lllred} in the following sense: Instead of using inequalities based on the theory of logarithmic forms as done by de Weger, we apply our optimized height bounds  (see Section~\ref{sec:ihb}). These optimized bounds are strong enough such that we can omit de Weger's reduction process. Then to  enumerate all solutions of~\eqref{eq:sunit} of bounded height, we use de Weger's sieve which is efficient as long as $\abs{S}$ is small (e.g. $\abs{S}<6$). To deal efficiently with sets $S$ of larger cardinality, we were forced to introduce new ideas: In Section \ref{sec:dwsieve+} we take into account certain geometric considerations to construct a refined sieve, and in Section~\ref{sec:suenum} we develop a refined enumeration algorithm for solutions of~\eqref{eq:sunit} with very small height.  Our new ideas are crucial to efficiently  solve \eqref{eq:sunit} for sets $S$ with $\abs{S}\geq 6$. Furthermore we prove that our refinements are substantial in the sense that they considerably improve the running time in theory and in practice, see \eqref{eq:rbound} and Section~\ref{sec:suapplications}.
In general we conducted some effort to optimize Algorithm~\ref{algo:suheight}. We refer to Section~\ref{sec:sucomplexity} where we explain and motivate our optimizations. Also we developed a method which (automatically) chooses parameters that are close to optimal in the generic case. This was necessary to obtain our database $\mathcal D_1$ listing the solutions of the $S$-unit equation \eqref{eq:sunit} for many distinct sets $S$, including all sets $S$ with $N_S\leq 10^7$ and all sets $S\subseteq S(16)$. 

\vspace{0.4cm}
\noindent{\bf Theorem~B.}
\emph{For each finite set of rational primes $S$ considered in $\mathcal D_1$, the database $\mathcal D_1$ contains all solutions of the $S$-unit equation \eqref{eq:sunit}.}
\vspace{0.4cm}

\noindent Another useful feature of Algorithm~\ref{algo:suheight} is that it allows to prove properties of $abc$-triples with bounded radical. For example, on using our algorithm we verified Baker's explicit $abc$-conjecture \cite[Conj 4]{baker:abcexperiments} for all $abc$-triples with radical at most $10^7$ or with radical composed of primes in $S(16)$. Furthermore we used our database $\mathcal D_1$ to motivate several new questions.  
 In particular, in view of the construction of the refined sieve, we make the following conjecture describing a property of integral points of $\mathbb P_\ZZ^1-\{0,1,\infty\}$ which is rather unexpected from a general Diophantine geometry perspective.

\vspace{0.3cm}
\noindent{\bf Conjecture 1.}
\emph{There exists $c\in\ZZ$ with the following property: If $n\in\ZZ_{\geq 1}$ then any finite set of rational primes $S$ with  $\abs{S}\leq n$ satisfies $\abs{\Sigma(S)}\leq \abs{\Sigma(S(n))}+c$.}
\vspace{0.3cm}

\noindent Theorem B shows that Conjecture~1 holds with $c=0$ for all sets $S$ in $\mathcal D_1$ and this motivates to ask whether any set of rational primes $S$ with $\abs{S}\leq n$ satisfies $\abs{\Sigma(S)}\leq \abs{\Sigma(S(n))}$?  We remark that additional applications of Algorithm~\ref{algo:suheight} are given in Section~\ref{sec:suapplications}. To conclude the discussion we point out that the geometric main idea (described in Section~\ref{sec:dwsieve+}) underlying our refined sieve is applicable in many other situations where sieves of de Weger type are applied. Here one can mention for example the practical resolution of $S$-unit and Thue--Mahler equations over number fields. We leave this for the future.

\paragraph{Mordell equation.} To solve the Mordell equation~\eqref{eq:mordell} via height bounds, we constructed Algorithm~\ref{algo:mheight}. Generically, this algorithm allows to deal efficiently with huge parameters. 
To illustrate this feature we used Algorithm~\ref{algo:mheight} to create our database $\mathcal D_2$ listing the solutions of \eqref{eq:mordell} for large classes of pairs $(a,S)$ with $a\in\ZZ-\{0\}$, including the classes $\{\abs{a}\leq 10,S\subseteq S(10^5)\}$, $\{\abs{a}\leq 100,S\subseteq S(10^3)\}$ and $\{\abs{a}\leq 10^4,S\subseteq S(300)\}$.

\vspace{0.4cm}
\noindent{\bf Theorem~C.}
\emph{For each pair $(a,S)$ considered in $\mathcal D_2$, the database $\mathcal D_2$ contains all solutions of the Mordell equation \eqref{eq:mordell} defined by $(a,S)$.}
\vspace{0.4cm}

\noindent Here we point out that Gebel--Peth{\H{o}}--Zimmer~\cite{gepezi:mordell} already established the important special case $\{\abs{a}\leq 10^4, S=\emptyset \}$ by using their algorithm~\cite{gepezi:ellintpoints}  based on the elliptic logarithm approach introduced by Masser and Zagier. Algorithm~\ref{algo:mheight} substantially improves the latter approach for \eqref{eq:mordell} in the following sense: Instead of using inequalities based on the theory of logarithmic forms, we apply our optimized height bounds  (see Section~\ref{sec:ihb}). Our bounds are considerably stronger in practice, which leads to significant running time improvements as illustrated in Section~\ref{sec:minitbounds}. Then to enumerate all solutions of~\eqref{eq:mordell} with bounded height, we use the elliptic logarithm sieve constructed in Section~\ref{sec:elllogsieve}. Here our construction combines a geometric interpretation of the known elliptic logarithm  reduction with conceptually new ideas described in Section~\ref{sec:setup}. The elliptic logarithm sieve is very efficient and it considerably improves in all aspects (see Section~\ref{sec:comparisonwithelr}) the known methods enumerating solutions of \eqref{eq:mordell}. However, our sieve requires  an explicit Mordell--Weil basis of the group $E_a(\QQ)$ associated to the elliptic curve $E_a$ defined by \eqref{eq:mordell}. While it is usually possible to determine such a basis in practice, there is so far no general effective method. In fact the dependence on a Mordell--Weil basis is a  disadvantage of Algorithm~\ref{algo:mheight} compared to the classical approach of Baker--Davenport which can be applied in the important case $S=\emptyset$. Their approach is very efficient in solving \eqref{eq:mordell} for varying $a\in\ZZ-\{0\}$ with $\abs{a}$ at most some given bound, see Bennett--Ghadermarzi~\cite{begh:mordell}. On the other hand, an advantage of Algorithm~\ref{algo:mheight} over the known algorithms  is that it can efficiently solve \eqref{eq:mordell} for large sets $S$. This feature allows to study  the  Diophantine problem for hyperbolic curves described in the next paragraph.  Other important features of Algorithm~\ref{algo:mheight} are the following: It can efficiently solve \eqref{eq:mordell} for parameters $a$ with huge height and its underlying correctness proofs are complete (even when $2\in S$ or when $S$ contains bad reduction primes of $E_a$). For example these two features are crucial  to efficiently determine all elliptic curves over $\sp(\ZZ)-S$ by solving certain equations \eqref{eq:mordell}, see below.

\paragraph{Points of hyperbolic curves.}  Suppose that $T$ and $B$ are nonempty open subschemes of $\sp(\ZZ)$, and assume that $T\subseteq B$.  Let $Y\to B$ be  an arbitrary hyperbolic curve   of genus $g$, see for example \cite[p.81]{mochizuki:absanab} for the definition. We denote by $Y(T)$ the set of $T$-points of $Y$ and we now consider the following Diophantine problem.

\vspace{0.3cm}
\noindent{\bf Problem.}
\emph{Describe the set $Y(T)$ in terms of $T$, with $T\subseteq B$ varying.}
\vspace{0.3cm}

\noindent If $g\geq 2$ then a result of Faltings \cite{faltings:finiteness} implies that the cardinality of $Y(T)$ is uniformly bounded in terms of $T$, which in some sense solves the problem for $g\geq 2$. Over the last decades the case $g=0$ was successfully studied  by many authors, including  Bombieri, Erd{\"o}s, Evertse, Gy{\H{o}}ry,  Moree, Silverman, Stewart, Tijdeman \cite{evertse:sunits,erstti:manysol,evgystti:sunitclasses,bomupo:cluster,evmostti:manysol} and more recently Harper, Konyagin, Lagarias, Soundararajan \cite{koso:manysunits,laso:smoothabcsol,harper:manysunits}. However the situation completely changes for $g=1$. In this  case, the problem is essentially not investigated in the literature and is widely open.  
On using Algorithm~\ref{algo:mheight} we study the problem for the families of hyperbolic genus one curves defined by Mordell equations  \eqref{eq:mordell} and cubic Thue equations \eqref{eq:thue}. In particular, motivated by Theorem~C and the construction of the elliptic logarithm sieve, we propose the following conjecture. 

\vspace{0.4cm}
\noindent{\bf Conjecture 2.}
\emph{There are constants $c_a$ and $c_r$, depending only on $a$ and $r$ respectively, such that any finite nonempty set of rational primes $S$ satisfies $\abs{Y_a(\mathcal O)}\leq c_a \abs{S}^{c_r}$.}
\vspace{0.4cm}

\noindent Here $Y_a(\mathcal O)$ denotes the set of solutions of the Mordell equation \eqref{eq:mordell} defined by $(a,S)$ and $r$ is the rank of the free part of the finitely generated abelian group $E_a(\QQ)$. The conjectured  bound is polynomial in terms of $\abs{S}$, while as far as we know all conjectures and results in the literature provide exponential bounds such as in Evertse--Silverman~\cite{evsi:uniformbounds}. We construct an infinite family of sets $S[b]$ which shows that the exponent $c_r$ has to be at least $\tfrac{r}{r+2}$. Furthermore Theorem~C strongly indicates that $c_r=\tfrac{r}{r+2}$ would be still far from optimal for many sets $S$ of interest, including the sets $S(n)$. On taking into account Theorem~C, we ask whether one can replace in Conjecture 2 the quantity $\abs{S}$ by the logarithm of the largest prime in $S$? We  motivate this question by constructing a probabilistic model. Together with a classical Diophantine approximation result of Siegel~(1929) and known estimates for the de Bruijn function, this model predicts a bound for  $\abs{Y_a(\mathcal O)}$ in terms of $S$ which would be optimal in view of the family $S[b]$.

\paragraph{Effective Shafarevich theorem.} Now we take $T=\sp(\ZZ)-S$ and we identify $M(T)$ with the set $M(S)$ of $\QQ$-isomorphism classes of elliptic curves over $\QQ$ with good reduction outside $S$. In the 1960s Shafarevich showed that $M(S)$ is finite: He reduced the problem to Mordell equations \eqref{eq:mordell} and then he applied Diophantine approximations. Coates~\cite{coates:shafarevich} made Shafarevich's proof effective by using the theory of logarithmic forms. In fact there already exist several practical methods which allow to determine the space $M(S)$. We refer to Section~\ref{sec:shaf} for an overview.  On combining Shafarevich's reduction with our Algorithm~\ref{algo:mheight} for Mordell equations~\eqref{eq:mordell}, we obtain Algorithm~\ref{algo:shaf} which allows to compute $M(S)$.  To illustrate the practicality of our approach, we determined the space $M(S)$ for each set $S\in\mathcal S$. Here $\mathcal S$ is a family of sets which contains in particular the set $S(5)$ and  all sets  $S$ with $N_S\leq 10^3$. Motivated by our data, we conjecture that one can replace in Conjecture~1 the moduli scheme $\mathbb P^1_{\ZZ[1/2]}-\{0,1,\infty\}$ of Legendre elliptic curves by the moduli stack $\mathcal M_{1,1}$ of elliptic curves.  In other words for any $n\in\ZZ_{\geq 1}$ our conjecture says that among all sets $S$ with $\abs{S}\leq n$ the cardinality of $M(S)$ is maximal (up to an absolute constant) when $S=S(n)$.  For many sets $S\in\mathcal S$ it seems that computing the space $M(S)$ is out of reach for the known methods, see the discussions in Section~\ref{sec:shaf}. In particular, our approach is significantly more efficient than the method of Cremona--Lingham~\cite{crli:shafarevich}. They use a different reduction to Mordell equations \eqref{eq:mordell} which involves $j$-invariants, and then they solve \eqref{eq:mordell} via the  algorithm of Peth{\H{o}}--Zimmer--Gebel--Herrmann \cite{pezigehe:sintegralpoints} based on the theory of logarithmic forms.   
The input of  Algorithm~\ref{algo:shaf}  requires a Mordell--Weil basis of $E_a(\QQ)$ for $2\cdot 6^{\abs{S}}$ distinct integers $a$. Thus our approach is not practical for large $\abs{S}$. Next, we mention that the problem of explicitly describing the space
\begin{equation}\label{eq:mdisplay}
M(T)=M(S)
\end{equation}
 is of interest for many reasons. For instance,  in \cite{rvk:modular} the moduli formalism was used to reduce many Diophantine problems to the study of $M(T)$.  On combining this strategy with our database listing the set $M(T)=M(S)$ for  $S\in\mathcal S$, we can directly solve any Diophantine problem inducing $T$-points on moduli schemes $Y$ of elliptic curves with effective \parshin{}  construction $\phi:Y(T)\to M(T)$; see Section~\ref{sec:moduli} for details and explicit examples. 
Here it suffices to know the image of $\phi$ in $M(T)$, which is often much smaller than the whole space $M(T)$. Taking this into account, we simplified and optimized the strategy for several classical Diophantine problems. In particular, we worked out the cases of cubic Thue equations \eqref{eq:thue}, cubic Thue--Mahler equations~\eqref{eq:thue-mahler} and generalized Ramanujan--Nagell equations~\eqref{eq:rana}. This led to the following algorithms and results.

\paragraph{Thue equation.} We constructed Algorithm~\ref{algo:theight} which allows to  solve the cubic Thue equation~\eqref{eq:thue}. Our approach is efficient in the generic case and it can deal with large sets $S$.   
To illustrate this we used Algorithm~\ref{algo:theight} in order to compile the database $\mathcal D_3$  containing the solutions of \eqref{eq:thue} for large classes of parameter triples $(f,S,m)$, where $m\in\ZZ$ is nonzero and $f\in\ZZ[x,y]$ is homogeneous of degree three with nonzero discriminant $\Delta$ (see Section~\ref{sec:thueproofs}). 
In particular our database $\mathcal D_3$ covers all  $(f,S,m)$ such that $m=1$ and such that $(\Delta,S)$ lies in  $\{\abs{\Delta}\leq 10^4,S\subseteq S(100)\}$, $\{\abs{\Delta}\leq 100,S\subseteq S(10^3)\}$ or $\{\abs{\Delta}\leq 20,S\subseteq S(10^5)\}$; see   Section~\ref{sec:ttmapp} for more information and additional examples.

\vspace{0.4cm}
\noindent{\bf Theorem~D.}
\emph{For each triple $(f,S,m)$ considered in $\mathcal D_3$, the database $\mathcal D_3$ contains all solutions of the cubic Thue equation~\eqref{eq:thue} defined by $(f,S,m)$.}
\vspace{0.4cm}

\noindent 
This gives in particular a new proof of several results in the literature (see Section~\ref{sec:thuealgo}) which determined the solutions of specific cubic Thue equations~\eqref{eq:thue}. Furthermore Theorem~D motivates new conjectures and questions on the number of solutions of \eqref{eq:thue}, see Section~\ref{sec:ttmapp}. We next describe the main ingredients of Algorithm~\ref{algo:theight}.  As in \cite[Sect 7.4]{rvk:modular} we reduce the problem to Mordell equations: This reduction uses classical invariant theory  which provides an explicit morphism $\varphi:X\to Y$ over $T=\sp(\mathcal O)$, where $X$ and $Y$ are the closed subschemes of $\mathbb A^2_T$  given by the Thue equation \eqref{eq:thue} and by the Mordell equation \eqref{eq:mordell} with $a=432\Delta m^2$ respectively. Then we compute $Y(T)$ using our Algorithm~\ref{algo:mheight} for Mordell equations and  we apply triangular decomposition in order to finally determine  $X(T)=\varphi^{-1}(Y(T))$.  Here we recall that Algorithm~\ref{algo:mheight} requires a Mordell--Weil basis of $E_a(\QQ)$.  Although it turned out that it is usually possible to determine such a basis in practice, the dependence on a Mordell--Weil basis is a disadvantage of our approach compared to the known methods discussed in Section~\ref{sec:thuealgo}. On the other hand, an advantage of our approach is that it can solve \eqref{eq:thue} for huge sets $S$. Here it seems that already sets $S$ with $\abs{S}\geq 10$ are out of reach for the known methods solving \eqref{eq:thue}.

\paragraph{Thue--Mahler equation.} Let $f\in\mathcal O[x,y]$ be a homogeneous polynomial of degree three with nonzero discriminant $\Delta$, and let $m\in\mathcal O$ be nonzero. We constructed Algorithm~\ref{algo:tmheight} which allows in particular to solve the classical cubic Thue--Mahler equation 
\begin{equation}\label{eq:thue-mahler}
f(x,y)=mz,
\end{equation} 
where $x,y,z\in \ZZ$ with $z\in \mathcal O^\times$ and $\gcd(x,y)=1$.  
To demonstrate the practicality of our approach, we used Algorithm~\ref{algo:tmheight} in order to create the database $\mathcal D_4$  listing the solutions of \eqref{eq:thue-mahler} for many triples $(f,S,m)$ with $m=1$ and $f\in\ZZ[x,y]$ as above. In particular $\mathcal D_4$ covers all such triples with $(\Delta,S)$ in  $\{\abs{\Delta}\leq 3000,S\subseteq S(2)\}$, $\{\abs{\Delta}\leq 10^3,S\subseteq S(3)\}$, $\{\abs{\Delta}\leq 100,S\subseteq S(4)\}$ or $\{\abs{\Delta}\leq 16,S\subseteq S(5)\}$; see Section~\ref{sec:ttmapp} for more information. 

\vspace{0.4cm}
\noindent{\bf Theorem~E.}
\emph{For each triple $(f,S,m)$ considered in $\mathcal D_4$, the database $\mathcal D_4$ contains all  solutions of the  cubic Thue--Mahler equation~\eqref{eq:thue-mahler} defined by $(f,S,m)$.}
\vspace{0.4cm}

\noindent 
We mention that $\mathcal D_4$  contains in addition the solutions of  \eqref{eq:thue-mahler} for various other $(f,S,m)$ of interest, including cases with $S=S(6)$. In fact Theorem~E gives in particular a new proof of  several results in the literature (see Section~\ref{sec:thuealgo}) which solved specific equations~\eqref{eq:thue-mahler}.  
We next describe the main ingredients of Algorithm~\ref{algo:tmheight}. On using an elementary standard reduction, we reduce \eqref{eq:thue-mahler} to  $3^{\abs{S}}$ distinct cubic Thue equations~\eqref{eq:thue} and these equations are then solved via Algorithm~\ref{algo:theight}.   Here the applications of Algorithm~\ref{algo:theight} require $3^{\abs{S}}$ distinct Mordell--Weil bases. Hence our approach is not practical when $\abs{S}$ is large.  However for small $\abs{S}$ it turned out  that it is usually possible to determine the required Mordell--Weil bases and then our approach is indeed efficient as illustrated in Section~\ref{sec:ttmapp}. 

\paragraph{Generalized Ramanujan--Nagell equations.} Let now $b$ and $c$ be arbitrary nonzero elements of $\mathcal O$. On using our approach for Mordell equations~\eqref{eq:mordell}, we obtained Algorithm~\ref{algo:ranaheight} which allows to solve the generalized Ramanujan--Nagell equation
\begin{equation}\label{eq:rana}
x^2+b=cy, \ \ \ \ \ (x,y)\in\mathcal O\times \mathcal O^\times.
\end{equation}
There is a vast literature devoted to the study of (special cases of) this Diophantine problem.  See for example the results, discussions and references in  Bugeaud--Shorey~\cite{bush:rana}, 
Bennett--Skinner~\cite[Sect 8]{besk:ternary} and Saradha--Srinivasan~\cite{sasr:rana}. To illustrate the practicality of our approach, we used  Algorithm~\ref{algo:ranaheight} in order to create the database $\mathcal D_5$ listing the solutions of \eqref{eq:rana} for many triples $(b,c,S)$ with $c=1$ and $b\in\ZZ-\{0\}$. In particular our database $\mathcal D_5$ covers all such triples with $(b,S)$ contained in $\{\abs{b}\leq 12,S\subseteq S(5)\}$, $\{\abs{b}\leq 35,S\subseteq S(4)\}$, $\{\abs{b}\leq 250,S\subseteq S(3)\}$ or $\{\abs{b}\leq 10^3,S\subseteq S(2)\}$.

\vspace{0.4cm}
\noindent{\bf Theorem~F.}
\emph{For each triple $(b,c,S)$ considered in $\mathcal D_5$, the database $\mathcal D_5$ contains all solutions of the generalized Ramanujan--Nagell equation~\eqref{eq:rana} defined by $(b,c,S)$.}
\vspace{0.4cm}

\noindent This theorem gives in particular a new proof of many results in the literature (see Section~\ref{sec:ranaalgoheight}) which solved special cases of~\eqref{eq:rana}. If $b\in\ZZ$ is nonzero and $c=1$, then  Peth{\H{o}}--de Weger~\cite{pede:binaryrec1} obtained a practical approach to find all solutions $(x,y)$ of \eqref{eq:rana} with $x,y\in\ZZ_{\geq 0}$. Their method involves binary recurrence sequences and the theory of logarithmic forms. Our approach is completely different: On using an elementary construction, we reduce \eqref{eq:rana} to certain Mordell equations~\eqref{eq:mordell} which we then solve via Algorithm~\ref{algo:mheight}. Here the involved Mordell curves usually have huge height. This is no problem for Algorithm~\ref{algo:mheight} and  it turned out that the bottleneck of our approach is finding the $3^{\abs{S}}$ distinct Mordell--Weil bases required for the applications of Algorithm~\ref{algo:mheight}. In light of this, we worked out a refinement of Algorithm~\ref{algo:ranaheight} in the following special case of \eqref{eq:rana}. For arbitrary nonzero $b,c,d$ in $\ZZ$ with $d\geq 2$, consider the classical Diophantine problem
\begin{equation}\label{eq:rana2}
x^2+b=cd^n, \ \ \ \ \ (x,n)\in\ZZ\times \ZZ.
\end{equation} 
Now, the crucial advantage of our refinement (see Algorithm~\ref{algo:rana2})  is that it only requires three distinct Mordell--Weil bases in order to find all solutions of \eqref{eq:rana2}.  On using Algorithm~\ref{algo:rana2}, we solved \eqref{eq:rana2} for all triples $(7,1,d)$ with $d\leq 888$; we note that here the case $d=2$ corresponds to the classical Ramanujan--Nagell equation. Furthermore, in Section~\ref{sec:ranaapp} we worked out additional applications of Algorithms~\ref{algo:ranaheight} and \ref{algo:rana2}. For example, we apply our approach to the problem of finding all coprime $S$-units $x,y\in\ZZ$ with $x+y$ a square or a cube. Here we solve several new cases of this problem. Also, we show that our approach is a useful tool to study  conjectures of Terai on Pythagorean triples.

\subsection{Diophantine problems related to Mordell equations}\label{sec:im}

We next discuss certain old Diophantine problems which are related to Mordell equations~\eqref{eq:mordell}. After presenting new results for primitive solutions of \eqref{eq:mordell}, we state a corollary on the greatest prime divisor of the difference of coprime squares and cubes.  We also give new height bounds for the solutions of cubic Thue equations~\eqref{eq:thue},  cubic Thue--Mahler equations~\eqref{eq:thue-mahler} and generalized Ramanujan--Nagell equations~\eqref{eq:rana}. As a byproduct, we obtain in this section  alternative proofs of classical theorems of Baker~\cite{baker:contributions}, Coates~\cite{coates:thue1,coates:thue2,coates:shafarevich} and Vinogradov--Sprind{\v{z}}uk~\cite{visp:thuemahler}. 

\subsubsection{Primitive solutions of Mordell equations}\label{sec:introprimsol}

Following Bombieri--Gubler \cite[12.5.2]{bogu:diophantinegeometry}, we say that $(x,y)\in\ZZ\times \ZZ$ is primitive if $\pm 1$ are the only  $n\in\ZZ$ with $n^{6}$ dividing $\gcd(x^3,y^2)$. In particular $(x,y)\in\ZZ\times\ZZ$ is primitive if $x,y$ are coprime. To measure the number $a\in\mathcal O$ and the finite set $S$, we take $$a_S=1728N_S^2\prod p^{\min(2,\ord_p(a))}$$ with the product extended over all rational primes $p\notin S$. Let $h$ be the usual logarithmic Weil height \cite[p.16]{bogu:diophantinegeometry}, with $h(n)=\log\abs{n}$ for $n\in\ZZ-\{0\}$. Building on the arguments of \cite[Cor 7.4]{rvk:modular}, we establish the following result (take $\mu=0$ in Theorem~\ref{thm:m}).

\vspace{0.3cm}
\noindent{\bf Theorem~G.}
\emph{Let $a\in \ZZ$ be nonzero. Assume that $y^2=x^3+a$ has a solution in $\ZZ\times\ZZ$ which is primitive. Then any  $(x,y)\in\mathcal O\times\mathcal O$ with $y^2=x^3+a$ satisfies}
$$\max\bigl(h(x),\tfrac{2}{3}h(y)\bigl)\leq a_S\log a_S.$$

\noindent We now discuss several aspects of this result. A useful feature of Theorem~G is that it does not involve~$\abs{a}$.  To illustrate this we take $n\in \ZZ_{\geq 1}$, we let $\mathcal F_n$ be the infinite family of integers $a$ with radical $\rad(a)$ at most~$n$, and we put $a_{*}=a_\emptyset$. Then it holds $a_*\leq 1728\rad(a)^2$ and Theorem~\ref{thm:m} with $\mu=0$ directly implies the following corollary.

\vspace{0.3cm}
\noindent{\bf Corollary~H.}
\emph{For any integer $n\geq 1$, the set of primitive $(x,y)\in\ZZ\times\ZZ$ with $y^2-x^3\in \mathcal F_n$ is finite and can in principle be determined. Furthermore if $a\in \ZZ$ satisfies $\log\abs{a}\geq a_*\log a_*$, then there are no primitive $(x,y)\in\ZZ\times\ZZ$ with $y^2-x^3=a$.}
\vspace{0.3cm}

\noindent It holds that $(3m^{3n})^2=(2m^{2n})^3+m^{6n}$ for all $m,n\in\ZZ$. Hence we see that one can not remove the assumption in Theorem~G. However, one can weaken the assumption by considering a certain class of (almost primitive) solutions of~\eqref{eq:mordell} which fits into Szpiro's small points philosophy \cite{szpiro:lefschetz}; see Theorem~\ref{thm:m} and the discussions given there.

We also deduce Corollaries~\ref{cor:coates1} and~\ref{cor:coates2} on the difference of perfect squares and perfect cubes. On taking for example $\eps=\frac{1}{10}$ in Corollary~\ref{cor:coates2}, one obtains the following result.

\vspace{0.3cm}
\noindent{\bf Corollary~I.}
\emph{Suppose that $x,y\in\ZZ$ are coprime, and write $X=\max(\abs{x},\abs{y})$. Then the greatest rational prime divisor $p$ of $y^2-x^3$ exceeds
$(1-\frac{1}{10})\log\log X-20.$}
\vspace{0.3cm}

\noindent This improves the  old theorem of Coates  \cite[Thm 2]{coates:shafarevich} in which he established the lower bound $10^{-3}(\log\log X)^{1/4}$. Similarly our Corollary~\ref{cor:coates1} refines \cite[Thm 1]{coates:shafarevich}. 
\noindent \paragraph{Comparison with literature.} We point out that Coates' method, which uses early estimates for logarithmic forms, is completely different to the method applied in this paper. In fact it is possible to prove a weaker version of our Theorem~\ref{thm:m} and to improve Coates' results \cite[Thm 1 and 2]{coates:shafarevich} by using more recent estimates for logarithmic forms. However, it turns out that without introducing new ideas the actual best lower bounds for linear forms in logarithms (see Baker--W\"ustholz~\cite{bawu:logarithmicforms} for an overview) do not give inequalities as strong as those provided by Theorem~\ref{thm:m} and Corollaries~\ref{cor:coates1} and~\ref{cor:coates2}. For example, let us consider our asymptotic version of Corollary~I established in Corollary~\ref{cor:coates2}. For any $\eps>0$ this version gives that the prime $p$ in Corollary~I exceeds
\begin{equation}\label{eq:asymptoticprimelowerbound}
\alpha\log\log X+\beta, \ \ \ \ \ \alpha=1-\eps,
\end{equation}
where $\beta$ denotes an effective constant depending only on $\eps$.  This improves the actual best factor $\alpha=\tfrac{1}{84}-\eps$ contained in the general result of Bugeaud~\cite[Thm 2]{bugeaud:greatestprimefactor}, which was proven by using inter alia a direct and ingenious reduction to lower bounds for logarithmic forms. Here it seems possible that one can slightly improve the factor $\alpha=\tfrac{1}{84}-\eps$ by  updating Bugeaud's approach with the actual best lower bounds for logarithmic forms.  However, the presence of the usual quantity $h(\alpha_1)\cdot\dotsc\cdot h(\alpha_n)$ in these lower bounds shows that 
this approach will always produce a factor $\alpha$ which is smaller than $\tfrac{1}{36}-\eps$. In view of this, our  result \eqref{eq:asymptoticprimelowerbound} seems to be out of reach for the present state of the art in the theory of logarithmic forms. See also the related discussions given at the end of Section~\ref{sec:mordellcoates}.

\paragraph{Idea of proof.} To prove our results for  primitive solutions of \eqref{eq:mordell}, we go into the proof of \cite[Cor 7.4]{rvk:modular} which combines the Shimura--Taniyama conjecture with the method of Faltings (Arakelov, \parshin, Szpiro)  as outlined in Section~\ref{sec:igeneralstrategy}; see also Section~\ref{sec:ihb}. Then we exploit that primitive solutions of \eqref{eq:mordell} induce, via the \parshin{} 
construction $\phi$, elliptic curves  with useful extra properties. For instance, to obtain the factor $\alpha=1-\eps$ in \eqref{eq:asymptoticprimelowerbound}, we use that the \parshin{} construction $\phi$ maps coprime solutions to elliptic curves which have semistable reduction over $\ZZ[1/6]$. The corollaries are then direct consequences of our results for primitive/coprime solutions and of the prime number theorem.

\subsubsection{Height bounds for cubic Thue and Thue--Mahler equations}\label{sec:ithue}

\noindent Baker \cite{baker:contributions} applied his theory of logarithmic forms in order  to establish in particular an effective finiteness result for any cubic Thue equation~\eqref{eq:thue} in the case when $\mathcal O=\ZZ$.
In the general case \eqref{eq:thue} is essentially equivalent to  the cubic Thue--Mahler equation~\eqref{eq:thue-mahler}.
Mahler~\cite{mahler:approx1} showed via Diophantine approximations that \eqref{eq:thue-mahler} has only finitely many solutions. 
Furthermore Coates~\cite{coates:thue1,coates:thue2} and Vinogradov--Sprind{\v{z}}uk~\cite{visp:thuemahler} independently  
proved effective finiteness  via the theory of logarithmic forms. 
We refer to Baker--W\"ustholz \cite{bawu:logarithmicforms} and Evertse--Gy{\H{o}}ry~\cite{evgy:bookuniteq,evgy:bookdiscreq} for an overview on generalizations and improvements of finiteness results for Thue and Thue--Mahler equations. 

\paragraph{Height bounds.} A new effective finiteness proof for any cubic Thue equation \eqref{eq:thue} was obtained in \cite{rvk:modular}. On working out explicitly the arguments of \cite[Sect 7.4]{rvk:modular}, we get explicit height bounds for the solutions of cubic Thue and Thue--Mahler equations. To state our results we denote by $h(f-m)$ the maximum of the logarithmic Weil heights of the coefficients of the polynomial $f-m\in\mathcal O[x,y]$. We put $a=432\Delta m^2$ with $\Delta$ the discriminant of $f$. The next corollary may be viewed as a refinement of \cite[Thm 7.1]{rvk:modular} in the case of moduli schemes (see Section~\ref{sec:moduli}) corresponding to \eqref{eq:thue} or \eqref{eq:thue-mahler}.

\vspace{0.3cm}
\noindent{\bf Corollary~J.}
\emph{The following statements hold.
\begin{itemize}
\item[(i)] Define the number $n$ by putting $n=2$ if $f,m\in\ZZ[x,y]$ and $n=10$ otherwise. Then any solution $(x,y)$ of the cubic Thue equation \eqref{eq:thue} satisfies  $$\max\bigl(h(x),h(y)\bigl)\leq a_S\log a_S+43nh(f-m).$$ 
\item[(ii)] If $(x,y,z)$ is a solution of the cubic Thue--Mahler equation \eqref{eq:thue-mahler} then $$\max\bigl(h(x),h(y),\tfrac{1}{3}h(z)\bigl)\leq 2a_S\log a_S+86nh(f-m).$$
\end{itemize}}
\noindent  
In Corollary~\ref{cor:precthue} we shall establish a more precise version of Corollary~J which provides sharper but more complicated bounds. Furthermore we shall show in Corollary~\ref{cor:precthue} that statement (ii) holds more generally for any primitive solution $(x,y,z)$ of the general cubic Thue--Mahler equation~\eqref{eq:thue-mahler}; the definition of such solutions is given in Definition~\ref{def:primsoltm}.

\paragraph{Comparison with literature.} We now compare Corollary~J with corresponding results in the literature. On using the theory of logarithmic forms, Bugeaud--Gy{\H{o}}ry~\cite[Thm 3 and 4]{bugy:thuemahler}, Bugeaud~\cite[Thm 3]{bugeaud:thuemahler}, Gy{\H{o}}ry--Yu~\cite[Thm 3]{gyyu:sunits} and Juricevic~\cite[$\mathsection$4.2]{juricevic:mordell} obtained  the actual best height bounds\footnote{We point out that these results hold for Diophantine equations which are considerably more general than \eqref{eq:thue} and \eqref{eq:thue-mahler}, and some of these results deal moreover with arbitrary number fields.} for the solutions of \eqref{eq:thue} and \eqref{eq:thue-mahler}.  We do not state these rather complicated bounds, but we mention that each of them has certain advantages and disadvantages. To compare these results with Corollary~J, we may and do assume that $f\in\ZZ[x,y]$ and $m\in\ZZ$. Then it  follows that $$a_S\leq 2^83^5\Delta_2\bigl(\rad(m)N_S\bigl)^2$$ where $\Delta_2=\min\bigl(\rad(\Delta)^2,\abs{\Delta}\bigl)$, and standard height arguments lead to $\abs{\Delta}\leq 3^5H^4$  for $H=\max_i \abs{a_i}$ the maximum of the absolute values of the coefficients $a_i$ of $f$. Therefore Corollary~J gives estimates which are asymptotically of the form $H^4\log H$, improving the actual best bounds $(H\log H)^4$  in terms of  $H$. In particular in the classical case, when $m$ is fixed (usually $m=1$) and $\mathcal O=\ZZ$, our Corollary~J improves the actual best results in all aspects. Furthermore Corollary~J improves the known estimates in terms of $S$ for infinitely many sets $S$, including all sets $S$ with $\abs{S}\leq 3$. On the other hand, our results are worse in terms of $m$ and the bound \cite[(12)]{gyyu:sunits} is significantly better in terms of $S$ for infinitely many sets $S$ including all sets $S=S(n)$ with $n$ large.  Finally we mention that our estimates (see also Corollary~\ref{cor:precthue}) involve small absolute constants and hence they considerably improve the actual best height bounds for all parameters which are not that large.  This might be of interest for the practical resolution of \eqref{eq:thue} and \eqref{eq:thue-mahler}. 

\paragraph{Idea of proof.} Following \cite[Sect 7.4]{rvk:modular}, we deduce our height bounds for cubic Thue equations \eqref{eq:thue} from a result for Mordell equations (Theorem~\ref{thm:m}) discussed above. This deduction uses classical invariant theory which provides an explicit morphism $\varphi:X\to Y$ over $T=\sp(\mathcal O)$, where $X$ and $Y$ are the closed subschemes of $\mathbb A^2_T$  given by the Thue equation \eqref{eq:thue} and by the Mordell equation \eqref{eq:mordell} with $a=432\Delta m^2$ respectively. Then in Proposition~\ref{prop:heightineq} we control the Weil height of any $P\in X(T)$ in terms of the Weil height of $\varphi(P)\in Y(T)$. To prove Proposition~\ref{prop:heightineq} we apply inter alia an effective arithmetic Nullstellensatz over the hypersurface in $\mathbb A_T^3$ given by $f-mz^3$. In fact we use here the Nullstellensatz of D'Andrea--Krick--Sombra~\cite{dakrso:nullstell} which leads to small constants. Finally, we deduce our height bounds for cubic Thue--Mahler equations \eqref{eq:thue-mahler}  by invoking an elementary standard construction which reduces \eqref{eq:thue-mahler} to Thue equations~\eqref{eq:thue}. Alternatively, one can obtain explicit height bounds for the solutions of \eqref{eq:thue} and \eqref{eq:thue-mahler} by directly applying \cite[Thm 7.1]{rvk:modular} with suitable moduli problems; see Sections~\ref{sec:ihb} and \ref{sec:moduli}.

 \subsubsection{Height bounds for generalized Ramanujan--Nagell equations}\label{sec:irana}

An elementary construction reduces the generalized Ramanujan--Nagell equation \eqref{eq:rana}, and the more classical special case \eqref{eq:rana2}, to Mordell equations \eqref{eq:mordell}. In light of this, results of Mordell~(1922) and Mahler~(1933) give finiteness for \eqref{eq:rana2} and \eqref{eq:rana} respectively. Moreover effective finiteness follows from Baker~\cite{baker:mordellequation}  in the case \eqref{eq:rana2} and from Coates~\cite{coates:shafarevich} in the case \eqref{eq:rana}. Our height bounds for Mordell equations lead to the following result which may be viewed as a refinement of \cite[Thm 7.1]{rvk:modular} in the case of moduli schemes (see Section~\ref{sec:moduli}) corresponding to generalized Ramanujan--Nagell equations \eqref{eq:rana}.

\vspace{0.3cm}
\noindent{\bf Corollary~K.}
\emph{If $(x,y)$ satisfies the generalized Ramanujan--Nagell equation \eqref{eq:rana}, then} $$\max\bigl(2h(x),h(y)\bigl)\leq 2a_S+h(a)+3h(c), \ \ \ a=bc^2.$$

\noindent In Corollary~\ref{cor:ranabounds} we shall give a more precise version of this height bound. Furthermore we shall deduce Corollary~\ref{cor:sumsofunits}  which provides explicit height bounds  for ``coprime" $u,v\in\mathcal O^\times$ with $u+v$ a square or cube in $\QQ$.  To discuss an application, we  take arbitrary coprime $m,n\in\ZZ$  and we consider the following simple condition in terms of  $r=\rad(mn)$: 
\begin{itemize}
\item[$(*)$] The natural logarithm of $\abs{m}$ or $\abs{n}$ exceeds $(90r)^2\log (9r)$.
\end{itemize} 
\noindent Now, the height bounds in Corollary~\ref{cor:sumsofunits} imply our next result which shows that condition~$(*)$ is in fact sufficient to rule out that $m+n$ is a perfect square or cube. 

\vspace{0.3cm}
\noindent{\bf Corollary~L.}
\emph{Suppose that $m$ and $n$ are arbitrary coprime rational integers. If condition~$(*)$ holds, then $m+n$ is not a perfect square or cube.}
\vspace{0.3cm}

\noindent One can obtain versions of our corollaries by using height bounds for the solutions of Mordell or Thue--Mahler equations which are based on the actual state of the art in the theory of logarithmic forms. For a comparison of the resulting bounds with our estimates, we refer to the analogous discussions given above and in Section~\ref{sec:simpleboundsm}. We further mention that the strong $abc$-conjecture of Masser--Oesterl\'e  in Remark~\ref{rem:abc} directly implies versions of our corollaries  which are asymptotically considerably better; notice that these implications are (as usual) not compatible with any exponential version of the $abc$-conjecture.

\subsubsection{Proof of the height bounds}\label{sec:ihb}

There is a long tradition of proving effective height bounds for Thue equations, Mordell equations,  Thue--Mahler equations and $S$-unit equations. In fact during the last few decades, one conducted quite some effort to refine the initial effective bounds of Baker~\cite{baker:contributions,baker:mordellequation}, Coates~\cite{coates:thue1,coates:thue2,coates:shafarevich} and  Gy{\H{o}}ry~\cite{gyory:sunitshelvetica}. 
See for example \cite{bawu:logarithmicforms,gyyu:sunits} for an overview on these refinements  which  all\footnote{Except Bombieri's refinement (see  Bombieri--Cohen~\cite{boco:effdioapp2}) of Thue's method  using Diophantine approximations. This method is relatively new and it is essentially self-contained. So far it leads to height bounds which are (slightly) worse compared to those coming from the theory of logarithmic forms.} ultimately rely on the theory of logarithmic forms.  We now discuss the strategy underlying the proofs of our height bounds for the Diophantine problems considered in the present paper. The symbol (ST) refers to the geometric version of the Shimura--Taniyama conjecture~\cite[Thm A]{breuil:modular} which relies inter alia on the Tate conjecture~\cite[Thm 4]{faltings:finiteness}.

\paragraph{$S$-unit equation.} Let $(x,y)$ be a solution of the $S$-unit equation \eqref{eq:sunit}. In the 1990s Frey \cite[p.544]{frey:ternary} (see also Murty~\cite{murty:strongmodular}) remarked that (ST) provides an alternative approach  to bound $h(x)$, and in 2011 it was independently shown  by Murty--Pasten and by the first mentioned author that Frey's ideas  together with (ST) lead to effective  bounds for $h(x)$; see  \cite[Thm 1.1]{mupa:modular} and \cite[Cor 7.2]{rvk:modular}. Here \cite{mupa:modular} works with the coprime Hecke algebra which leads to $h(x)\ll N_S\log N_S$, while \cite{rvk:modular} uses the full Hecke algebra  which in general leads to the (slightly) weaker bound $h(x)\ll N_S(\log N_S)(\log\log N_S)$. Since there are also many situations in which the full Hecke algebra provides the best bounds,  we use in this paper the coprime and the full Hecke algebra approach. To work out the optimized height bounds for Algorithm~\ref{algo:suheight}, we follow closely the arguments of \cite[Thm 1.1]{mupa:modular} and \cite[Cor 7.2]{rvk:modular} and we conduct some effort to refine the involved estimates. 
Asymptotically, the actual best result $h(x)\ll N_S^{1/3}(\log N_S)^3$ is due to Stewart--Yu \cite{styu:abc2}. However, for all sets $S$ with $N_S\leq 2^{100}$ and for many other sets $S$ of practical interest, our optimized height bounds are considerably stronger than those coming from the theory of logarithmic forms; see Section~\ref{sec:simpleboundssu}.  We note that this alternative approach for $S$-unit equations \eqref{eq:sunit}  is in fact (see \cite{rvk:modular}) a special case  of the method of Faltings (Arakelov, \parshin, Szpiro) combined with (ST) as described in Section~\ref{sec:igeneralstrategy} above.

\paragraph{Other Diophantine problems.} Many classical Diophantine problems can be reduced to $S$-unit equations. However in most cases  the known (unconditional) reductions involve number fields larger than $\QQ$. In particular one can not combine these reductions with (unconditional) results for \eqref{eq:sunit} in order to deal with the Diophantine problems considered in the present paper. Instead we use that all these  problems induce integral points on moduli schemes of elliptic curves.  Given this observation, we can apply the  strategy of \cite[Thm 7.1]{rvk:modular} which provides explicit height bounds  for integral points on  moduli schemes of elliptic curves.  This strategy consists of combining  (ST)  with Faltings' method in a way which is similar as described in Section~\ref{sec:igeneralstrategy} above.  Here important ingredients of the proof are the height-conductor inequality \cite[Prop 6.1]{rvk:modular} (proven independently in \cite[Thm 7.1]{mupa:modular}, see Section~\ref{sec:heightcondstatement}) and the moduli formalism.  Besides providing a useful geometrical interpretation of various classical Diophantine problems, the moduli formalism allows to find new explicit applications of the method. Indeed we discovered many results of the present paper by searching for moduli schemes with ``interesting" defining equations:  Given a priori the information that the equation defines a moduli scheme with effective \parshin{} construction $\phi$,  one can explicitly work out the strategy of \cite[Thm 7.1]{rvk:modular} to get effective finiteness results; see \cite[Sect 7]{rvk:modular} and Section~\ref{sec:moduli}. Furthermore, a posteriori one can often  obtain here simpler (but less conceptual) proofs by removing the moduli formalism. In light of the practical purpose of the present article, we conducted some effort to simplify our proofs as much as possible in the case of Mordell, cubic Thue, cubic Thue--Mahler and generalized Ramanujan--Nagell equations. 
For instance on working out the method for the moduli schemes defined by Mordell equations \eqref{eq:mordell}, one obtains  the actual best height bounds \cite[Cor 7.4]{rvk:modular} for the solutions of~\eqref{eq:mordell}. To prove the optimized height bounds for Algorithm~\ref{algo:mheight}, we follow and simplify the arguments of \cite[Cor 7.4]{rvk:modular}. Here we try to optimally estimate the involved quantities. We note that a priori the arguments of \cite[Sect 1-7]{rvk:modular} prove all our explicit simplified height bounds, but of the form $X(\log X)(\log\log X)$. To remove here in addition the factor $\log\log X$, we go into the proof of \cite[Lem 5.1]{rvk:modular} and we now use  an idea of Murty--Pasten~\cite{mupa:modular} involving the coprime Hecke algebra. Their idea was not known to the author of \cite{rvk:modular}. However, for each considered Diophantine problem, there are also many situations in which the full Hecke algebra approach of \cite[Lem 5.1]{rvk:modular} provides the best bounds. Hence  we work out all our height bounds using the coprime and the full Hecke algebra approach.

\subsection{Organization of the paper}

\paragraph{Plan of the paper.} In  Section~\ref{sec:cremonas+st} we  briefly discuss some tools which are crucial for our results and algorithms. In particular, we recall a geometric version of the Shimura--Taniyama conjecture which relies inter alia on the Tate conjecture. 

In Section~\ref{sec:sunitalgo} we present two algorithms for $S$-unit equations. The first approach via modular symbols is worked out in Section~\ref{sec:sucremalgo}. Then in Section~\ref{sec:suheightalgo} we conduct some effort to construct the second algorithm which uses height bounds. Here, after discussing in Section~\ref{sec:dwsieve} a slight variation of de Weger's method, we construct our refined sieve in Section~\ref{sec:dwsieve+} and we develop a refined enumeration in Section~\ref{sec:suenum}. In Section~\ref{sec:suapplications} we present various applications of our algorithm via height bounds. In particular, we discuss our database $\mathcal D_1$ containing the solutions of large classes of $S$-unit equations and we use our data to motivate several Diophantine conjectures related to $S$-unit equations. 

Section~\ref{sec:malgo}
contains two algorithms which allow to solve the Mordell equation. The first approach via modular symbols is worked out in Section~\ref{sec:mcremalgo}.  Then in Section~\ref{sec:mheightalgo} we construct the second algorithm via height bounds. Here, after discussing the main ingredients of our algorithm, including our initial height bounds in Section~\ref{sec:minitbounds} and the elliptic logarithm sieve,  we put everything together in Section~\ref{sec:mordellalgostat}.  We also  present  various applications of our algorithm via height bounds. In Section~\ref{sec:shaf} we apply this algorithm to study the problem of finding all elliptic curves over $\QQ$ with good reduction outside $S$, and in Section~\ref{sec:moduli} we solve certain classes of Diophantine equations by combining our algorithm  with the moduli formalism. Then in Section~\ref{sec:malgoapplications} we discuss  our database $\mathcal D_2$ containing the solutions of large classes of Mordell equations and we motivate new conjectures/questions. In Section~\ref{sec:malgocomparison} we  compare our algorithms with other methods.

In Section~\ref{sec:thuealgo} we present algorithms for cubic Thue and Thue--Mahler equations. Our algorithms use a construction from classical invariant theory which allows to reduce the Diophantine problems to  Mordell equations. After working out some useful properties of this construction in Section~\ref{sec:talgoconst}, we discuss our algorithms via modular symbols in Section~\ref{sec:talgocremona}. Then in Section~\ref{sec:talgoheight} we explain our algorithms via height bounds and we give various applications. In particular, we discuss in Section~\ref{sec:ttmapp} parts of our databases $\mathcal D_3$ and $\mathcal D_4$ containing the solutions of large classes of cubic Thue and Thue--Mahler equations respectively. In Section~\ref{sec:talgocompa} we compare our algorithms with the known methods.

Section~\ref{sec:ranaalgo} contains two algorithms for the generalized Ramanujan--Nagell equation. After presenting our approach via modular symbols in Section~\ref{sec:ranaalgocremona}, we explain the algorithm via height bounds in Section~\ref{sec:ranaalgoheight} and we discuss some applications including our database $\mathcal D_5$. Then we compare our algorithms with the known methods in Section~\ref{sec:ranaalgocomp}.

In Sections~\ref{sec:mordellcoates}, \ref{sec:thueproofs} and \ref{sec:ranaheight} we consider certain classical Diophantine problems related to Mordell equations. In particular, in Section~\ref{sec:mordellcoates} we study properties of (almost) primitive solutions of Mordell equations and we deduce explicit lower bounds for the largest prime divisor of the difference of coprime  squares and cubes. In Sections~\ref{sec:thueproofs} and \ref{sec:ranaheight} we give new explicit height bounds for the solutions of cubic Thue and Thue--Mahler equations and of generalized Ramanujan--Nagell equations. Then in Section~\ref{sec:heightbounds} we prove our results for (almost) primitive solutions and we work out the optimized height bounds for the solutions of Mordell and $S$-unit equations which are used in our algorithms.  

In Section \ref{sec:elllogsieve} we construct the elliptic logarithm sieve. It allows to efficiently find all integral points of bounded height on any elliptic curve $E$ over $\QQ$ with given Mordell--Weil basis of $E(\QQ)$. We refer to the introduction of Section~\ref{sec:elllogsieve} for an overview of the main ideas of our construction. The elliptic logarithm sieve is of independent interest and thus we made the presentation of Section~\ref{sec:elllogsieve} independent of the rest of this paper.

\paragraph{Notation.} We shall use throughout the following (standard) notations and conventions. By $\log$ we mean the principal value of the natural logarithm. We define the product taken over the empty set as~$1$. For any set $M$, we denote by $\lvert M\rvert$ the (possibly infinite) number of distinct elements of~$M$. Let $f_1$ and $f_2$ be real valued functions on $M$. We write $f_1=O(f_2)$ if there is a constant $c$ such that $f_1\leq cf_2$. Further $f_1=O_\eps(f_2^\eps)$ means that for any real number $\eps>0$ there is a constant $c(\eps)$ depending only on $\eps$ such that $f_1\leq c(\eps) f_2^\eps$. For any $n\in\ZZ_{\geq 1}$, we say that $\mathcal E\subset \mathbb R^n$ is an ellipsoid centered at the origin if $\mathcal E=\{x\in\mathbb R^n\st q(x)\leq c\}$ for some positive definite quadratic form $q:\RR^n\to \RR$ and some positive real number $c$.

We denote by $\abs{z}$ the usual complex absolute value of~$z\in\mathbb C$. If $m,n\in\ZZ$ then the symbol $m\mid n$ (resp. $m\nmid n)$ means that $m$ divides $n$ (resp. $m$ does not divide $n$). Further $\gcd(a_1,\dotsc,a_n)$ denotes the greatest common divisor of $a_1,\dotsc,a_n\in\ZZ$. 
The radical of $n\in\ZZ$ is given by $\rad(n)=\prod p$  with the product taken over all rational primes $p$ dividing~$n$.  If $\alpha\in \QQ$ is nonzero and if $p$ is a rational prime, then we write $\ord_p(\alpha)\in\ZZ$ for the order of $p$ in~$\alpha$. We denote by $h(\alpha)$ the usual absolute logarithmic Weil height of $\alpha\in\QQ$, with $h(0)=0$  and $h(\alpha)=\log\max(\abs{m},\abs{n})$ if $\alpha=m/n$ for coprime $m,n\in\ZZ$. Finally for any real number $x\in\mathbb R$, we write $\floor{x}=\max(n\in\ZZ\st n\leq x)$ and $\ceil{x}=\min(n\in\ZZ\st n\geq x)$.

\paragraph{Computer, software and algorithms.} Unless mentioned otherwise, we used a standard personal working computer at the MPI Bonn for our computations. Our algorithms are all implemented in Sage and we shall use functions of the computer algebra systems Pari~\cite{pari:parisystem}, Sage~\cite{sage:sagesystem} and Magma~\cite{magma:magmasystem}. In what follows, we shall sometimes refer by (PSM) to these computer packages in order to simplify the notation.

For each of our algorithms, we conducted some effort to motivate  our  constructions (in theory and in practice), to explain our choice of parameters,  to discuss important complexity aspects,  to give detailed correctness proofs and to  circumvent potential numerical issues. We shall also list the running times of our algorithms for many  examples. The listed times are always upper bounds. In fact some of them were obtained by using older versions of our algorithms, and in many cases the running times would now be significantly better when using the most recent versions (as of February 2016) of our algorithms.

\paragraph{Acknowledgements.}
The research presented in this paper was initiated when we were members at the IAS Princeton (2011/12), it was continued at the IH\'ES (2012/13) and it was completed at the MPIM Bonn (2013-15). We are grateful to these institutions for providing excellent working conditions. The authors were supported by the NSF grant No. DMS-0635607 (2011/12) and by  EPDI fellowships (2012-14). We would like to thank the MPI, in particular Gerd Faltings and Pieter Moree, for support in (2013-15). 

Further, we would like to thank Richard Taylor for motivating and very useful initial discussions. We are  grateful to Yuri Bilu, Enrico Bombieri, Sander Dahmen, Jan-Hendrik Evertse, K\'alm\'an Gy{\H{o}}ry, Pieter Moree, Hector Pasten and Don Zagier for encouraging discussions and/or for informing us about useful literature. Also, we learned a lot from Yuri Bilu, John Cremona, Stephen Donnelly, Nuno Freitas, Steffen M\"uller, Martin Raum and Samir Siksek. We would like to thank all of them for explaining various aspects of computational number theory and/or for answering questions. 

\paragraph{Data and earlier versions.} The present version of the paper (Feb. 2016) extends in particular all results and algorithms presented in our earlier versions of the paper (April 2013 and Nov. 2014). However we removed certain applications and discussions of our algorithms via modular symbols, since they are meanwhile obsolete in view of more recent results. Our data is uploaded on: \url{https://www.math.u-bordeaux.fr/~bmatschke/data/}.

\section{Shimura--Taniyama conjecture}\label{sec:cremonas+st}

A crucial ingredient for all of our results and algorithms is a (geometric) version of the Shimura--Taniyama conjecture which relies inter alia on the Tate conjecture. In the case of our first type of algorithms, another important ingredient is an algorithm of Cremona using modular symbols. In this section we first introduce some notation and then we briefly discuss these ingredients in order to emphasize that they do not depend on results proven by (classical) transcendence or Diophantine approximation techniques. 

Let $N\geq 1$ be an integer. Consider the classical congruence subgroup $\Gamma_0(N)\subset\sl2(\ZZ)$, let $X_0(N)=X(\Gamma_0(N))_\QQ$ be the smooth, projective and geometrically connected model over $\QQ$ of the  modular curve associated to $\Gamma_0(N)$ and denote by $S_2(\Gamma_0(N))$ the complex vector space of cuspforms of weight 2 with respect to~$\Gamma_0(N)$. See for example~\cite{dish:modular} for the definitions.  We denote by $J_0(N)=\textnormal{Pic}^0(X_0(N))$  the Jacobian variety of $X_0(N)$. Let $\mathbb T_\ZZ$ be the subring  of the endomorphism ring of $J_0(N)$, which is generated over $\ZZ$   by 
the usual Hecke operators $T_n$ for all $n\in \ZZ_{\geq 1}$. 
For any $f\in S_2(\Gamma_0(N))$, we denote by $a_n(f)$ the $n$-th Fourier coefficient of $f$ and we say that $f$ is rational if $a_n(f)\in \QQ$ for all $n\in \ZZ_{\geq 1}$.  Further, we say that $f\in S_2(\Gamma_0(N))$ is a newform (of level $N$) if $a_1(f)=1$, if $f$ lies in the new part of $S_2(\Gamma_0(N))$ and if $f$ is an eigenform for all Hecke operators on $S_2(\Gamma_0(N))$. We now suppose that $f\in S_2(\Gamma_0(N))$ is a rational newform. Let $I_f$ be the kernel of the ring homomorphism $\mathbb T_\ZZ\to \ZZ[\{a_n(f)\}]$ which 
is induced 
by $T_n\mapsto a_n(f)$.
It turns out that the image 
$I_fJ_0(N)$ of $J_0(N)$   under $I_f$ is   connected,   and the corresponding quotient
\begin{equation*}
E_f=J_0(N)/I_fJ_0(N)
\end{equation*}
is an elliptic curve over $\QQ$ since $f$ is rational. We next define the modular degree and congruence number of~$f$. On composing the usual embedding $X_0(N)\incl J_0(N)$, which sends the cusp $\infty$ of $X_0(N)$ to the zero element of $J_0(N)$, with the natural projection $J_0(N)\to E_f$, we obtain a finite morphism $\varphi:X_0(N)\to E_f$ of curves over~$\QQ$. The modular degree $m_f$ of $f$ is defined as the degree of~$\varphi$:
\begin{equation}\label{def:mf}
m_f=\deg(\varphi).
\end{equation}
The congruence number $r_f$ of $f$ is defined as the largest integer such that there exists a cusp form $f_c\in S_2(\Gamma_0(N))$ with all Fourier coefficients in $\ZZ$ and 
\begin{equation}\label{def:rf}
(f,f_c)=0 \textnormal{ and } a_n(f)\equiv a_n(f_c)\textnormal{ mod } (r_f) \textnormal{ for all } n\in \ZZ_{\geq 1}.
\end{equation}
Here  $(\ ,\,)$ denotes the usual Petersson inner product. Let $\mathcal E(N)$ be the set of all elliptic curves $E$ over $\QQ$ which are of the form $E=E_f$ for some rational newform $f\in S_2(\Gamma_0(N))$. We say that an elliptic curve $E$ over $\QQ$ is modular if there exists a positive integer $N$ such that the curve $E$ is $\QQ$-isogenous to some elliptic curve in~$\mathcal E(N)$.

 For any given $N\in \ZZ_{\geq 1}$, Cremona's algorithm~\cite{cremona:algorithms} computes in particular the coefficients of minimal Weierstrass equations of all modular elliptic curves over~$\QQ$. A short description of this algorithm may be as follows: One considers~$X_0(N)$, computes its first homology using $M$-symbols (after Manin~\cite{manin:symbols}), computes the action of sufficiently many Hecke operators on it, and determines the one-dimensional eigenspaces with rational eigenvalues. By induction on the  divisors of~$N$, this yields the rational newforms in $S_2(\Gamma_0(N))$, and their period lattices allow then to compute the set~$\mathcal E(N)$. 
Finally on using a theorem of Mazur~\cite{mazur:qisogenies}, one computes all elliptic curves over $\QQ$ which are $\QQ$-isogenous to some curve in $\mathcal E(N)$ and one determines their minimal Weierstrass equations. 

Building on the key breakthroughs by Wiles~\cite{wiles:modular} and by Taylor--Wiles~\cite{taywil:modular}, the Shimura--Taniyama conjecture was finally established by Breuil--Conrad--Diamond--Taylor~\cite{breuil:modular}. This conjecture implies (its geometric version saying) that any elliptic curve over $\QQ$ of conductor $N$ is $\QQ$-isogenous to some curve in~$\mathcal E(N)$.  This implication uses the Tate conjecture~\cite[Thm 4]{faltings:finiteness}. We point out that Faltings'  proof of the Tate's conjecture does not use transcendence theory or classical Diophantine approximations.

\section{Algorithms for $S$-unit equations}\label{sec:sunitalgo}
Let $S$ be a finite set of rational primes, write $N_S=\prod_{p\in S} p$ and denote by $\mathcal O^\times$  the group of units of
$\mathcal O =\ZZ[1/N_S]$. In this section, we are interested to solve the $S$-unit equation
\begin{equation}
x+y=1, \ \ \ (x,y)\in\mathcal O^\times\times \mathcal O^\times. \tag{\ref{eq:sunit}}
\end{equation}
If $\abs{S}\leq 1$ then \eqref{eq:sunit} has either no solutions or  $(2,-1)$, $(-1,2)$ and $(\frac{1}{2},\frac{1}{2})$ are the only solutions. Hence we may and do assume that $\abs{S}\geq 2$ in this section.  
As mentioned in the introduction, there already exists a practical method of de Weger \cite{deweger:lllred} which solves~\eqref{eq:sunit}. De Weger \cite[Thm 5.4]{deweger:lllred} used his method to completely solve~\eqref{eq:sunit} in the case  $S=\{2,3,5,7,11,13\}$. 
Further, Wildanger~\cite{wildanger:unitalgo} and afterwards Smart~\cite{smart:smallsol} generalized the ideas of de Weger and they obtained a practical algorithm which solves~\eqref{eq:sunit} over arbitrary number fields; see also Hajdu~\cite{hajdu:optimalsys} and Evertse--Gy{\H{o}}ry \cite{evgy:bookuniteq}. There is also the recent work of Dan-Cohen--Wewers~\cite{dawe:explicitkim,dawe:sunitalgomotivic}, with the ultimate goal to construct an algorithm \cite{dancohen:sunitalgo3} solving~\eqref{eq:sunit} via ``explicit motivic Chabauty--Kim theory". This method is inspired by Kim's ($p$-adic \'etale) approach \cite{kim:siegel}, see also the discussion in \cite[2.5.2]{dancohen:sunitalgo3} which mentions an additional method of Brown. So far all practical approaches solving \eqref{eq:sunit} crucially rely on the theory of logarithmic forms. 

In the following Sections~\ref{sec:sucremalgo} and~\ref{sec:suheightalgo}, we present and discuss two alternative algorithms which solve $S$-unit equations~\eqref{eq:sunit}. Both of our algorithms do not use  the theory of logarithmic forms. The first algorithm  relies on Cremona's algorithm, 
and we refer to the beginning of Section~\ref{sec:suheightalgo} for a short description of the ingredients of the second algorithm. 
Before we begin to describe our algorithms, we discuss  useful properties of symmetric solutions and we give a lower bound for the complexity of any algorithm solving \eqref{eq:sunit}. 
\begin{definition}\label{def:symmetric}Suppose that $(x,y)$ and $(x',y')$ are solutions of \eqref{eq:sunit}. Then we say that $(x,y)$ and $(x',y')$ are symmetric solutions if $x'$ or $y'$ lies in $\{x,\tfrac{1}{x},\tfrac{1}{1-x}\}$. 
\end{definition}
It turns out that this defines an equivalence relation on the set of solutions of \eqref{eq:sunit}, and hence we can consider the set $\Sigma(S)$ of solutions of \eqref{eq:sunit} up to symmetry. Suppose that $(x,y)$ is a solution of \eqref{eq:sunit}. Then one can directly determine all its symmetric solutions. In fact there are exactly six  solutions of \eqref{eq:sunit} which are symmetric to  $(x,y)$ provided that $(x,y)$ is not equal to $(2,-1)$, $(-1,2)$ or $(\tfrac{1}{2},\tfrac{1}{2})$. In particular, we see that the number of solutions of \eqref{eq:sunit} is either zero or $6\cdot \abs{\Sigma(S)}-3$. The following remark shows that any algorithm solving the $S$-unit equation~\eqref{eq:sunit} can not be too fast in general.  

\begin{remark}[Lower bound for complexity]\label{rem:sulowercomplexity}
The result of Erd{\"o}s--Stewart--Tijdeman \cite[Thm 4]{erstti:manysol}, see also the more recent work of Harper, Konyagin, Lagarias, Soundararajan \cite{koso:manysunits,laso:smoothabcsol,harper:manysunits}, implies the existence of an  effective absolute constant $s_0$ with the following property. For any $s\in\ZZ_{\geq s_0}$ there exists a set $S$ with $\abs{S}=s$ such that the $S$-unit equation~\eqref{eq:sunit} has at least $\exp((s/\log s)^{1/2})$ solutions. 
Furthermore, there are infinitely many sets $S$ such that the $S$-unit equation~\eqref{eq:sunit} has a solution $(x,y)$ with $H(x)\geq N_S$ for $H(x)$ the multiplicative Weil height \cite[p.15]{bogu:diophantinegeometry} of~$x$. 
Therefore we conclude that any algorithm solving~\eqref{eq:sunit} has running time which is in general not better than linear in $\log N_S$ and which is in general not better than $\exp((\abs{S}/\log \abs{S})^{1/2})$. 
\end{remark}

\subsection{Algorithm via modular symbols}\label{sec:sucremalgo}

To assure that our algorithm really computes all solutions of the $S$-unit equation~\eqref{eq:sunit}, we use the following observations. We suppose that $(x,y)$ satisfies~\eqref{eq:sunit}.  Then there exist nonzero  $a,b,c\in\ZZ,$ with $\gcd(a,b,c)=1$ and $\rad(abc)\mid N_S$, such that $x=\frac{a}{c}$, $y=\frac{b}{c}$ and $a+b=c$. In other words, resolving~\eqref{eq:sunit} is  equivalent to resolving 
\begin{equation}
\label{eq:abc}
a+b = c, \quad a,b,c\in\ZZ-\{0\}, \quad \gcd(a,b,c)=1,\quad \rad(abc)\mid N_S.
\end{equation}
Further we observe that to find all solutions of~\eqref{eq:abc} it suffices to consider solutions $(a,b,c)$ of~\eqref{eq:abc} with $a,b,c$ all positive. Indeed if $(a,b,c)$ satisfies~\eqref{eq:abc} then there exists a  solution $(\alpha,\beta,\gamma)$ of~\eqref{eq:abc} such that the sets $\{\alpha,\beta,\gamma\}$ and $\{\abs{a},\abs{b},\abs{c}\}$ coincide.

For any elliptic curve $E$ over~$\QQ$, we denote by $\Delta_E$  the minimal discriminant of $E$ and we write $N_E$ for the conductor of~$E$. A construction of  Frey--Hellegouarch~\cite{frey:curves,hellegouarch:curves} associates to any solution $(a,b,c)$ of~\eqref{eq:abc} an elliptic curve\footnote{We warn the reader that here $E_{abc}$ is not necessarily the usual Frey--Hellegouarch curve with Weierstrass equation $y^2=x(x-a)(x+b)$. See the proof of Lemma~\ref{lem:psu2} in which $E_{abc}$ is denoted by $E$.}   
 $E_{abc}$ over~$\QQ$. On taking the quotient of $E_{abc}$ by the ``subgroup"  generated by a suitable 2-torsion point of~$E_{abc}$, one obtains an elliptic curve $E$ over $\QQ$ with the following properties (see Lemma~\ref{lem:psu2}).

\begin{lemma}\label{lem:psu1}
Suppose that $(a,b,c)$ is a solution of~\eqref{eq:abc}. Then there exists an elliptic curve $E$ over $\QQ$ such that $N_E$ divides $2^4N_S$ and  $\Delta_{E}=2^{8-12m}\abs{ab}c^4$ with $m\in\{0,1,2,3\}$.
\end{lemma}
In fact this lemma may be viewed (see \cite[Prop 3.2]{rvk:modular}) as an explicit \parshin{} construction for integral points on the Legendre moduli scheme $\mathbb P^1_{\ZZ[1/2]}-\{0,1,\infty\}$ of elliptic curves, which is induced by forgetting the Legendre level structure.

\begin{Algorithm}[$S$-unit equation via modular symbols]\label{algo:sucremona} The input is a finite set of rational primes $S$, and the output is the set of solutions $(x,y)$ of the $S$-unit equation~\eqref{eq:sunit}. 

The algorithm\textnormal{:} If $2\notin S$ then output the empty set, and if $2\in S$ then do the following. 

\begin{itemize}
\item[(i)] Use Cremona's algorithm, described in Section~\ref{sec:cremonas+st}, to compute the set $\mathcal T_\Delta$ of minimal discriminants of modular elliptic curves over $\QQ$ of conductor dividing~$2^4N_S$. 
\item[(ii)] Let $\mathcal T=\cup_{m}\mathcal T_m$ be the union of the sets $\mathcal T_m=\{2^{12m-8}\Delta_E\st \Delta_E\in \mathcal T_\Delta\}\cap\ZZ$ where $m=0,1,2,3$. For each $d\in\mathcal T$,
factor $d$ as $d=\prod_{p\in S}p^{n_p}$ with $n_p\in \ZZ$. Then for each disjoint 
partition of the set $S=S_{\alpha}\dot\cup S_{\beta}\dot\cup S_{\gamma}$, with \
$S_{\gamma}\subseteq\{p\in S\st  4\mid n_p\}$, define $\alpha=\prod_{p\in 
S_{\alpha}}p^{n_p}$, $\beta=\prod_{p\in S_{\beta}}p^{n_p}$, and $\gamma=\prod_{p\in 
S_{\gamma}}p^{n_p/4}$.
Then output all $(x,y)$ of the form $(x,y)=(\frac{a}{c},\frac{b}{c})$ with $a,b,c\in\ZZ$ satisfying $a+b=c$ and $\{\abs{a},\abs{b},\abs{c}\}=\{\alpha,\beta,\gamma\}$.
\end{itemize}
\end{Algorithm}
\paragraph{Correctness.}We now verify that this algorithm indeed finds all solutions of any $S$-unit equation~\eqref{eq:sunit}.  Let $(x,y)$ be a solution of~\eqref{eq:sunit}. As explained above~\eqref{eq:abc}, there is a corresponding solution $(a,b,c)$  of~\eqref{eq:abc} with $(x,y)=(\frac{a}{c},\frac{b}{c})$. 
We write $d=\abs{ab}c^4$. Lemma~\ref{lem:psu1} gives an elliptic curve $E$ over $\QQ$ such that $N_E\mid 2^4N_S$ and  $\Delta_E=2^{8-12m}d$ with $m\in\{0,1,2,3\}$. The Shimura--Taniyama conjecture assures that $E$ is modular. This proves that $\Delta_E$ is contained in the set $\mathcal T_\Delta$ computed in step (i) and it follows that $d\in \mathcal T$. Furthermore, 
the number $d=\abs{ab}c^4\in \mathcal T$ factors in step (ii) as 
$d=\prod_{p\in S}p^{n_p}$ with $n_p\in \ZZ$. Thus the disjoint partition $S=S_{\alpha}\dot\cup S_{\beta}\dot\cup S_{\gamma}$, 
with $S_{\alpha}=\{p\st p\mid a\}$, 
$S_{\beta}=\{p\st p\mid b\}$ and 
$S_{\gamma}=S-(S_{\alpha}\cup S_{\beta})$, produces in step (ii) our solution $(x,y)=(\frac{a}{c},\frac{b}{c})$ as desired. Here we used that our coprime $a,b,c\in\ZZ$ satisfy $a+b=c$ and $\{\abs{a},\abs{b},\abs{c}\}=\{\alpha,\beta,\gamma\}$.

\paragraph{Complexity.} In step (i) the algorithm has to compute all (modular) elliptic curves over $\QQ$ of conductor $N$ for precisely\footnote{Notice that $6\cdot 2^{\abs{S}-1}$ is the number of positive rational integers dividing $2^4N_S$, since $2\in S$.} $6\cdot 2^{\abs{S}-1}$ positive integers~$N$.  Here we can exploit the fact that Cremona's algorithm proceeds by induction over the divisors of~$N$, see Section~\ref{sec:cremonas+st}. However at the time of writing it is not clear to us what is the running time of Cremona's algorithm. We also mention that step (i) greatly benefits from the ongoing extension of Cremona's tables which in particular list all (modular) elliptic curves over $\QQ$ of given conductor~$N$. As of August 2014, these tables are complete for all $N\leq 350 000$. 

The complexity of step (ii) crucially depends on the size of~${\mathcal T}$. It is an open (Diophantine) problem to find a simple formula for $\abs{\mathcal T}$ in terms of~$S$. However one can give an upper bound for $\abs{\mathcal T}$ in terms of~$S$. For example, the work of Ellenberg, Helfgott and Venkatesh~\cite{heve:integralpoints,elve:classgroup} implies that $\abs{\mathcal T}\ll N_S^{0.1689}$. 
In step (ii) the algorithm first needs to compute the prime factorization of each $d\in\mathcal T$ and then it needs to compute certain  integers $a,b,c$ for most of the disjoint 3-partitions of the set~$S$. Inequality~\eqref{eq:szpiro} implies that any $d\in\mathcal T$ satisfies $\log d=O(N_S^2)$, and the number of disjoint 3-partitions of $S$ is~$3^{\abs{
S}}$. Hence the above discussions show that the running time of step (ii) is at most polynomial in terms of~$N_S$. This running time estimate can be considerably improved by assuming various conjectures. Indeed for each $d\in\mathcal T$ the $abc$-conjecture ($abc$) recalled in Remark~\ref{rem:abc} would provide that $\log d=O(\log N_S)$, and it follows from Brumer--Silverman \cite[Thm 4]{brsi:number} that the Birch--Swinnerton-Dyer conjecture (BSD) together with the General Riemann Hypothesis (GRH) would give  $\abs{\mathcal T}=O_\eps(N_S^\eps)$. Hence the running time of step (ii) is  $O_\eps(N_S^\eps)$ if the three conjectures (BSD), (GRH) and ($abc$) all hold.

\begin{remark}[Variation]\label{rem:sulegendrealternative}
We now discuss a variation of Algorithm~\ref{algo:sucremona} which in practice provides an improvement for large~$\abs{S}$. The idea is that instead of using in (ii) the curve of Lemma~\ref{lem:psu1}, one can work with Legendre curves  as in \cite[Prop 3.2 (i)]{rvk:modular}. 

\begin{itemize}
\item[(i)'] Use Cremona's algorithm to compute the set $\mathcal T_W$ of minimal Weierstrass models over $\ZZ$ of  modular elliptic curves over $\QQ$ with conductor dividing~$2^4N_S$.
\item[(ii)'] For each $W\in \mathcal T_W$ determine its set $\lambda(W)=\{\lambda,1-\lambda,\lambda^{-1},\dotsc\}$  of Legendre parameters by computing the three roots of $f+f_2^2/4$, where $Y^2+f_2(X)Y=f(X)$ defines~$W$. 
If $\lambda\in \lambda(W)$ and $(\lambda,1-\lambda)$ satisfies~\eqref{eq:sunit} then output $(x,y)=(\lambda,1-\lambda)$.
\end{itemize} This algorithm outputs the set of solutions of~\eqref{eq:sunit}. Indeed for any solution $(x,y)$ of~\eqref{eq:sunit} the proof of Lemma~\ref{lem:psu2} gives $W\in\mathcal T_W$ with $x\in\lambda(W)$. The running times of (i) and (i)' are essentially equal, and it follows for example from \cite[Cor 6.3]{rvk:modular} that the running time of step (ii)' is polynomial in terms of~$N_S$. 
Furthermore one can show that the running time of step (ii)' is in fact~$O_\eps(N_S^\eps)$, provided that all three conjectures (BSD), (GRH) and $(abc)$ hold.  If $\abs{S}$ is large, then step (ii)' is in practice considerably faster than (ii) since the latter  iterates in addition over many disjoint 3-partitions of $S$. However for large $\abs{S}$ the bottleneck of the algorithm is the computation of $\mathcal T_W$ or~$\mathcal T_{\Delta}$, and for small $\abs{S}$ it turns out that (ii)' and (ii) are essentially equally fast. Hence we implemented the algorithm involving (i) and (ii), since the implementation of (i) is simpler compared to (i)'. 
\end{remark}

\paragraph{Applications.} Let $\Sigma(S)$ be the set of solutions  of~\eqref{eq:sunit} up to symmetry. For any $N\in\ZZ_{\geq 1}$, consider the set $\Sigma(N)=\cup_S \Sigma(S)$ with the union taken over all sets $S$ with $N_S\leq N$. On using an implementation in Sage of the above Algorithm~\ref{algo:sucremona}, we computed the sets $\Sigma(N)$ for all $N\leq 20000$. 
This computation was very fast for $N\leq 20000$, since for such $N$  part (i) of our Algorithm~\ref{algo:sucremona} can use Cremona's database which  contains in particular the required data for all elliptic curves of conductor dividing $2^4N< 350 000$. On the other hand, if the required data of the involved elliptic curves is not already known, then our Algorithm~\ref{algo:sucremona} is often not practical anymore. Here the  problem is the application  of Cremona's algorithm (using modular symbols)  in step (i) which requires a huge amount of memory in order to deal with medium sized or large conductors.

\begin{remark}[$abc$-triples]\label{rem:abc}
The $abc$-conjecture $(abc)$  states that for any real number $\eps>0$ there are only a finite number of solutions of~\eqref{eq:abc} with quality larger than $1+\eps$, see for example~\cite{masser:abc}. Here the quality $q=q(a,b,c)$ of a solution $(a,b,c)$ of~\eqref{eq:abc} is defined as $$q=\frac{\log\max(\abs{a},\abs{b},\abs{c})}{\log \rad(abc)}.$$
Among all known high quality $abc$-triples the top two, as of October 2014, can be obtained with our method: $(2,3^{10}109,23^5)$ with $q=1.6299$ due to  Reyssat, and $(11^2,3^25^67^3,2^{21}23)$ with $q=1.6260$ due to de Weger; see \url{abcathome.com}.
\end{remark}

\subsection{Algorithm via height bounds}\label{sec:suheightalgo}

In this section we use the optimized height bounds in Proposition~\ref{prop:algobounds} to construct an algorithm which practically resolves the $S$-unit equation~\eqref{eq:sunit}.
The decomposition of our algorithm is inspired by de Weger's classical method whose main ingredients are as follows: 
\begin{itemize}
\item[(1)]\label{itDeWegerStepBaker} De Weger uses explicit height bounds for the solutions of \eqref{eq:sunit}, based on the theory of logarithmic forms, to rule out the existence of solutions with very ``large" height. See  Section~\ref{sec:heightbounds} for  references and more details.
\item[(2)]\label{itDeWegerStepLLL} Then de Weger tries to further reduce the height bounds in (1) by using the LLL lattice reduction algorithm applied to certain approximation lattices, which are defined using $p$-adic logarithms. If this reduction can be applied, then he can rule out in addition the existence of solutions with ``large" and ``medium sized" height.
\item[(3)]\label{itDeWegerStepSieve}  To find most of the solutions with ``small" height de Weger applies a sieve which we call de Weger's sieve, see Section~\ref{sec:dwsieve}. He repeats this step as many times as required to make sure that the remaining solutions have ``tiny'' height.
\item[(4)] \label{itDeWegerStepEnumerationByHand}
Finally de Weger checks (by hand) all potential solutions of ``tiny" height.
\end{itemize}
In our algorithm we replace the height bounds in (1) by the optimized height bounds which we shall work out in Section~\ref{sec:heightbounds}. These optimized bounds are strong enough such that we can now omit the reduction step (2). In Section~\ref{sec:dwsieve}  we discuss a slight variation of de Weger's sieve described in (3). Then in Section~\ref{sec:dwsieve+} we develop a refined sieve which has considerably improved running time in practice, in particular for all sets $S$ with $\abs{S}>6$. Moreover, in Section~\ref{sec:suenum} we construct an enumeration algorithm which is faster than the standard algorithm in (4). In fact our enumeration is fast enough such that it is now beneficial to go from (3) to (4) at an earlier stage, which leads to  additional running time improvements.  Finally we present and discuss our Algorithm~\ref{algo:suheight}  solving the $S$-unit equation~\eqref{eq:sunit}, and  we also give various applications of our algorithm.

Before we begin to carry out the above program, we introduce some notation which will be used throughout Section~\ref{sec:suheightalgo}. Let $n\in \ZZ_{\geq 1}$ and consider two vectors $l,u\in\mathbb R^n$. We write $l\leq u$ if and only if $l_i\leq u_i$ for all $i\in\{1,\dotsc,n\}$. Further by $l<u$ we mean that $l\leq u$ with $l\neq u$, and we use the symbol $l\not\leq u$ in order to say that $l\leq u$ does not hold. In other words, we use poset (partially ordered set) and not
coordinate-wise notation for $``<"$ and $``\not\leq"$. 
The use of poset notation for $``<"$ and $``\not\leq"$ will  simplify our exposition. Next, we suppose that $(x,y)$ is a solution of the $S$-unit equation~\eqref{eq:sunit}. It will be convenient for the description of the algorithm to use the following quantities associated to~$(x,y)$. For any rational prime $p$ we put $m_p(x,y)=\max (\abs{\ord_p(x)},\abs{\ord_p(y)})$, and then we define
$$m(x,y)= (m_p(x,y))_{p\in S} \ \  \textnormal{ and } \ \ M(x,y)= \max(m_p(x,y)\log p\st p\in S).$$
Proposition~\ref{prop:algobounds} gives an upper bound $M_0$ for~$M(x,y)$. Hence it remains to find the solutions $(x,y)$ of~\eqref{eq:sunit} with $M(x,y)$ between zero and $M_0$, or between $M_{k+1}$ and $M_{k}$ for $k=0,\dotsc,k_0$ and some convenient sequence $0=M_{k_0}<\ldots< M_{1}<M_0$. For this purpose we shall use the following refinements of steps (3) and (4) of de Weger's method.

\subsubsection{De Weger's sieve}
\label{sec:dwsieve}

We now begin to describe the sieve which is used in step (3) of de Weger's method. In fact we shall describe a slight variation of this sieve, see Remark~\ref{rem:teske}.  

For any given $M', M''\in\RR_{\geq 0}$, we want to enumerate all solutions $(x,y)$ of~\eqref{eq:sunit} with $M'<M(x,y)\leq M''$. For this purpose, it suffices by \eqref{def:ulbounds} to solve the following problem: For any two vectors  $l,u\in \ZZ^S$ with $0\leq l\leq u$,  find all solutions $(x,y)$ of~\eqref{eq:sunit} which satisfy
\begin{equation}
\label{prob:mxy}
m(x,y)\not\leq l \ \ \ \textnormal{ and } \ \ \ m(x,y)\leq u.
\end{equation}
If a solution $(x,y)$ of \eqref{eq:sunit} satisfies \eqref{prob:mxy}, then all its symmetric solutions satisfy \eqref{prob:mxy} as well. 
The condition $m(x,y)\not\leq l$ means that there exists at least one ``large" exponent in the prime factorization of $x$ or~$y$. We now exploit this to reduce  (\ref{prob:mxy}) to a Diophantine problem whose solutions can be quickly enumerated. Suppose  that $(x,y)$ is a solution of~\eqref{eq:sunit} which satisfies  (\ref{prob:mxy}). Then there exists $q\in S$ with $m_q(x,y)\geq 1+l_q$.
Further the discussion given above~\eqref{eq:abc} delivers a solution $(a,b,c)$ of~\eqref{eq:abc} with  $x=\frac{a}{c}$ and $y=\frac{b}{c}$.
After possibly replacing $(x,y)$ with a symmetric solution, we may and do assume that $q^{l_q+1}$ divides~$c$. Then it holds that $a+b=0\mod q^{l_q+1}$ and thus  $(a/b)^2 = 1$ in $G$ for $$G=(\ZZ/q^{l_q+1}\ZZ)^\times.$$
Here we squared plainly in order to get rid of the minus sign, and we used that $a,b$ are both coprime to $q$  which provides that $a,b$ are in $G$. 
Next we write $(a/b)^2 =\prod_{p\in {S\wo q}} p^{2\gamma_p}$  with $\gamma_p=\ord_p(a/b)$. Then  we see that any solution $(x,y)$ of~\eqref{eq:sunit} satisfying (\ref{prob:mxy}) induces a solution $\gamma=(\gamma_p)\in\ZZ^{S\wo q}$ of the following problem: If $l,u\in\ZZ^{S}$ are given with $0\leq l\leq u$, then for each $q\in S$ find all  $\gamma\in \ZZ^{{S\wo q}}$ such that $\abs{\gamma_p}\leq u_p$ for all $p\in{S\wo q}$ and such that 
\begin{equation}
\label{eqDeWegerSieve}
\prod_{p\in {S\wo q}} p^{2\gamma_p} = 1 \ \textnormal{ in } \ G.
\end{equation}
Furthermore if $\gamma$ is a solution of~\eqref{eqDeWegerSieve} then one can quickly reconstruct all solutions of~\eqref{eq:sunit} satisfying (\ref{prob:mxy}) which map to $\gamma$ via the above construction.  Indeed one defines $a=\prod_{\gamma_p>0}p^{\gamma_p}$, $b_{+}=\prod_{\gamma_p<0}p^{-\gamma_p}$ and $b_-=-b_+$, and then one checks for each $b\in\{b_+,b_-\}$  whether $\rad(a+b)\divides N_S$ and whether $(x,y)=(\frac{a}{a+b},\frac{b}{a+b})$ satisfies (\ref{prob:mxy}). 
In particular, we see that we can quickly enumerate all solutions of~\eqref{eq:sunit} satisfying (\ref{prob:mxy}) provided that we know all solutions of~\eqref{eqDeWegerSieve}. Finally it remains to enumerate all solutions of~\eqref{eqDeWegerSieve}. For this purpose we observe that the set of vectors $\gamma\in\ZZ^{{S\wo q}}$ satisfying~\eqref{eqDeWegerSieve} is the intersection of a certain lattice $\Gamma\subseteq \ZZ^{S\wo q}$ with the cube $\{\gamma \st |\gamma_p|\leq u_p\}\subset \mathbb R^{S\wo q}$. To determine this intersection we combine the following two algorithms:

\begin{enumerate}
\item[(T)]
The first algorithm is an application of Teske's  ``Minimize''  \cite[Algo 5.1]{teske:minimize}. It takes as the input a set of generators $g_1,\ldots,g_d$ of a finite abelian group\footnote{Here it is not required to know the group $G_0$ explicitly. In fact it suffices  to know a number divisible by $\abs{G_0}$ and to be able to compute the following operations in $G_0$: Multiplying elements, inverting elements and testing whether an element is the neutral element.} $G_0$, 
and it outputs  a basis for the lattice $\Gamma\subseteq\ZZ^d$ formed by those $\gamma\in\ZZ^d$ with $\prod_{i=1}^dg_i^{\gamma_i}=1$.

\item[(FP)]
The second algorithm is a version of the Fincke--Pohst algorithm, see Remark~\ref{rem:fp}. For any $d\in \ZZ_{\geq 1}$ it takes as the input a basis of a lattice $\Gamma\subseteq\ZZ^d$ together with an ellipsoid $\mathcal E\subset\RR^d$ centered at the origin, and
it outputs the inter\-section~$\Gamma\cap \mathcal E$.
\end{enumerate}
More precisely, we apply (T) and (FP) as follows. Let $G_0$ be the subgroup of $G$ which is generated by the squares of the elements in $S\wo q$. 
An application of (T) with $G_0$ gives a basis for the lattice $\Gamma$ underlying \eqref{eqDeWegerSieve}, and then (FP)  computes the intersection of $\Gamma$ with the smallest ellipsoid $\mathcal E\subset \mathbb R^{S\wo q}$ that contains the cube $\{\gamma\st |\gamma_p|\leq u_p\}\subset \mathbb R^{S\wo q}$.

\begin{Algorithm}[De Weger's sieve]
\label{algDeWegersSieve}
The input consists of a finite set of rational primes $S$ together with two vectors $l,u\in\ZZ^S$ such that $0\leq l\leq u$, and the output is the set of solutions $(x,y)$ of the $S$-unit equation~\eqref{eq:sunit} which satisfy~\eqref{prob:mxy}. 

The algorithm\textnormal{:} After possibly shrinking $S$, we may and do assume that all $u_p\geq 1$. If  $2\notin S$ then output the empty set, and if $2\in S$ then  do the following for each $q\in S$.
\begin{itemize}
\item[(i)]
Use the application \textnormal{(T)} of  Teske's algorithm in order to compute a basis for the lattice $\Gamma\subseteq \ZZ^{S\wo q}$ of all $\gamma\in \ZZ^{S\wo q}$ with $\prod_{p\in S\wo q}p^{2\gamma_p}=1$ in $(\ZZ/q^{l_q+1}\ZZ)^\times$.
\item[(ii)] Define the ellipsoid $\mathcal E=\{x\in\RR^{S\wo q}\st \sum_{p\in S\wo q} |x_p/u_p|^2 \leq |S|-1\}$. Then compute the intersection $\Gamma\cap \mathcal E$ using the version \textnormal{(FP)} of the Fincke--Pohst algorithm.
\item[(iii)] For each $\gamma\in \Gamma\cap \mathcal E$ lying in the cube $\{\gamma\st|\gamma_p| \leq u_p\}\subset\mathbb R^{S\wo q}$, define the numbers $a=\prod_{\gamma_p>0}p^{\gamma_p}$,  $b_+=\prod_{\gamma_p<0}p^{-\gamma_p}$ and $b_-=-b_+$. For any $b\in\{b_{+},b_-\}$ such that $c=a+b$ satisfies $\ord_p(c)\leq u_p$ for all $p\in S$,  $q^{l_q+1}\mid c$ and $\prod_{p\in S} p^{\ord_p(c)}=\abs{c}$,  output the solution $(x,y)=(\tfrac{a}{c},\tfrac{b}{c})$ together with all its symmetric solutions.
\end{itemize}
\end{Algorithm}

\paragraph{Correctness.} We now verify that this algorithm indeed finds  all solutions of \eqref{eq:sunit} satisfying \eqref{prob:mxy}. Suppose that $(x,y)$ is such a solution. Then our assumption $m(x,y)\not\leq l$ gives $q\in S$ together with $\gamma\in\ZZ^{S\wo q}$ which is associated to $(x,y)$ via the construction given above \eqref{eqDeWegerSieve}. This $\gamma$ lies in the lattice  $\Gamma$ appearing in (i). Furthermore, our assumption $m(x,y)\leq u$ implies that $\gamma$ is also contained in the ellipsoid $\mathcal E$ from step (ii). In particular $\gamma$ lies in the intersection $\Gamma\cap\mathcal E$ computed in step (ii), and $m(x,y)\leq u$ provides that $\abs{\gamma_p}\leq u_p$ for all $p\in S\wo q$. Therefore on using that $\gamma$ is associated to $(x,y)$ via the construction described above \eqref{eqDeWegerSieve}, we see that step (iii) produces our solution $(x,y)$ as desired.


\paragraph{Complexity.}
Step (i) uses (T). The algorithm (T) is reminiscent of the discrete logarithm problem in the multiplicative group $G$, and the bottleneck for this is the prime factorization of $\abs{G}$. In our situation it holds $\abs{G} = q^{l_q}(q-1)$ and therefore it suffices to factor $q-1$, which is in general much easier than factoring an arbitrary number of size $\abs{G}$. In fact step (i) is in practice not the bottleneck  of Algorithm \ref{algDeWegersSieve}, except in the case when $S$ consists of the prime 2 together with a few large primes. Usually step (ii) is the bottleneck of Algorithm \ref{algDeWegersSieve}, and thus we shall refine this step in  Section~\ref{sec:dwsieve+} below.

\begin{remark}\label{rem:teske}
To find a basis for the lattice $\Gamma$, de Weger used $q$-adic logarithms instead of the application (T) of Teske's algorithm.
In fact (T) was already used in generalizations of de Weger's method to number fields, see Wildanger~\cite{wildanger:unitalgo} and Smart~\cite{smart:smallsol}. 
\end{remark}

\paragraph{Application.} For later use we now mention how we will apply Algorithm \ref{algDeWegersSieve} in order to find  all solutions $(x,y)$ of~\eqref{eq:sunit} satisfying $M'< M(x,y)\leq M''$, for some given $M',M''$ in  $\mathbb R_{\geq 0}$ with $M'<M''$. Consider the two vectors $l=(l_p)$ and $u=(u_p)$ in $\ZZ^S$, with
\begin{equation}\label{def:ulbounds}
l_p=\floor{M'/\log p} \ \ \ \textnormal{ and } \ \ \ u_p=\floor{M''/\log p}.
\end{equation}
We observe that $0\leq l\leq u$. Then an application of  Algorithm \ref{algDeWegersSieve} with $l,u$ finds in particular all solutions $(x,y)$ of~\eqref{eq:sunit} satisfying $M'< M(x,y)\leq M''$, as desired.


\begin{remark}
According to~\eqref{def:m0m1} below, splitting the possible candidates with respect to values of $M(x,y)$ is in particular reasonable in the case when $\log p$ is small compared to $\abs{S}$ for all $p\in S$. In this case, the iteration over all $q\in S$ in Algorithm~\ref{algDeWegersSieve} will take  about equally long for each $q\in S$.
If  $S$ contains primes which are exponentially large in terms of $\abs{S}$, then one should split the initial space  $0\leq M(x,y)\leq M_0$ into more general pieces of the form $\{m(x,y)\not\leq u(k+1) \textnormal{ and } m(x,y)\leq u(k)\}$ for  suitable  $u(k)\in\ZZ^{S}$ with $k=0,1,\dotsc,k_0$. Here the vectors $u(k)\in\ZZ^S$ should satisfy $0=u(k_0)<\ldots<u(1)<u(0)$ and $u(0)_p=\floor{M_0/\log p}$ for $p\in S$. In practice it is reasonable to choose greedily the next $u(k+1)<u(k)$,  such that the subsequent sieving step is as fast as possible.
\end{remark}

\subsubsection{Refined sieve}
\label{sec:dwsieve+}

In this section we continue our notation introduced above. We begin with a short description of the geometric main idea underlying our refinement. Recall that in de Weger's sieve one needs to determine the intersection of a certain lattice $\Gamma\subseteq \ZZ^{S\wo q}$ with the cube $\{x\st \abs{x_p}\leq u_p\}\subset \mathbb R^{S\wo q}$, and for this purpose Algorithm~\ref{algDeWegersSieve}~(ii) first determines $\Gamma\cap \mathcal E$ for $\mathcal E$ the smallest ellipsoid containing the cube. However, in terms of the rank $s-1$ where $s=\abs{S}$,  the volume of $\mathcal E$ is exponentially larger  than the volume of the cube. In our refinement  we essentially truncate from the cube  some regions near the faces of codimension $2,\ldots,\floor{s/3}$. For $s\geq 6$ the resulting geometric object is contained in a notably smaller ellipsoid $\mathcal E'$, which allows  to determine $\Gamma\cap\mathcal E'$ considerably faster than $\Gamma\cap \mathcal E$.

To explain in more detail our refined sieve, we need to introduce additional notation.  Let $(x,y)$ be a solution of the $S$-unit equation~\eqref{eq:sunit}, and let $(a,b,c)$ be a solution of \eqref{eq:abc} with $(x,y)=(\tfrac{a}{c},\tfrac{b}{c})$. We take $j\in\{1,\dotsc,t\}$ for $t=\max(1,\floor{\abs{S}/3})$. This choice of $t$ takes into account that any prime appearing in the prime factorization of $x$ or $y$ divides one of the three coprime integers $a,b,c$. For any $n\in \ZZ$ with $\rad(n)\divides N_S$, we denote by $\mu_j(n)$ the $j$-th largest\footnote{More precisely, $\mu_j(n)$ is the $j$-th largest element of the ordered multi-set of cardinality  $\abs{S}$ obtained by  ordering the $\abs{S}$ non-negative real numbers $\ord_p(n)\log p$, $p\in S$, with respect to their absolute values.} 
 of the  real numbers $\ord_p(n)\log p$, $p\in S$. Then we define
\[
\mu_j(x,y)=\max\bigl(\mu_j(a),\mu_j(b),\mu_j(c)\bigl) \ \ \ \textnormal{ and } \ \ \ \mu(x,y) = \bigl(\mu_1(x,y),\ldots,\mu_t(x,y)\bigl).
\]
We observe that $\mu_1(x,y)=M(x,y)$. However we point out that if $j\in\{2,\dotsc,t\}$ then  $\mu_j(x,y)$ is not necessarily the $j$-th largest  of  the numbers $m_p(x,y)\log p$,  $p\in S$. 
Now we consider the following problem: For any given vectors $\mu', \mu''\in\ZZ^t$ having monotonously decreasing entries such that $0\leq \mu'< \mu''$, find  all solutions $(x,y)$ of \eqref{eq:sunit} that satisfy 
\begin{equation}\label{eqBoundsOnVectorMxy}
\mu(x,y) \not\leq\mu' \ \ \ \textnormal{ and } \ \ \  \mu(x,y) \leq\mu''.
\end{equation}
If a solution $(x,y)$ of \eqref{eq:sunit} satisfies \eqref{eqBoundsOnVectorMxy}, then all its symmetric solutions satisfy \eqref{eqBoundsOnVectorMxy} as well.  Further we note that the condition $\mu(x,y)\not\leq \mu'$ implies, for some $j\in\{1,\dotsc,t\}$, the existence of at least $j$ ``large" exponents in the prime factorization of $x$ or $y$.  In the following algorithm we exploit this in order to work with lattices of rank $\abs{S}-j$.

\begin{Algorithm}[Refined sieve]
\label{algRefinedDeWegerSieve}
The input is a finite set of rational primes $S$, together with two vectors $\mu', \mu''\in\ZZ^t$ having monotonously decreasing entries such that $0\leq \mu'< \mu''$; where $t=\max(1,\floor{\abs{S}/3})$.
The output is the set of  solutions $(x,y)$ of the $S$-unit equation~\eqref{eq:sunit} which satisfy \eqref{eqBoundsOnVectorMxy}. 

The algorithm\textnormal{:} If $2\notin S$ then output the empty set, and if $2\in S$ then do the following for each non-empty subset $T\subseteq S$ of cardinality $\abs{T}\leq t$.
\begin{enumerate}
\item[(i)]
Put $n=\prod_{q\in T} q^{\floor{\mu'_{\abs{T}}/\log q}+1}$, and use the application \textnormal{(T)} of Teske's algorithm to compute a basis of the lattice $\Gamma_T$ of all $\gamma\in  \ZZ^{S\wo T}$ with
$
\prod_{p\in S\wo T} p^{2\gamma_p} = 1$ in $(\ZZ/n\ZZ)^\times.
$
\item[(ii)]
Then use the version \textnormal{(FP)}  of the  Fincke--Pohst algorithm in order to determine  the intersection of the lattice $\Gamma_T\subseteq \ZZ^{S\wo T}$ with the ellipsoid $\mathcal E_T\subset \mathbb R^{S\wo T}$ defined by
\begin{equation}
\label{eqBoundForAlphaForImprovedDeWeger}
\mathcal E_T=\{x\in\mathbb R^{S\wo T}\st\sum_{p\in S\wo T} |x_p \log p|^2 \leq \sum_{i=1}^{|S\wo T|} (\mu''_{\min(\ceil{i/2},t)})^2\}.
\end{equation}
\item[(iii)]
For each  $\gamma\in\Gamma_T\cap\mathcal E_T$, define  $a=\prod_{\gamma_p>0}p^{\gamma_p}$,  $b_+=\prod_{\gamma_p<0}p^{-\gamma_p}$ and $b_-=-b_+$. For any $b\in\{b_{+},b_-\}$ such that $c=a+b$ satisfies $\prod_{p\in S}p^{\ord_p(c)}=\abs{c}$ and such that $(x,y)=(\tfrac{a}{c},\tfrac{b}{c})$ satisfies \eqref{eqBoundsOnVectorMxy},  output $(x,y)$ together with all its symmetric solutions.
\end{enumerate}
\end{Algorithm}

\paragraph{Correctness.} To see that this algorithm works correctly, we suppose that $(x,y)$ is a solution of the $S$-unit equation \eqref{eq:sunit} that satisfies \eqref{eqBoundsOnVectorMxy}. Let $(a,b,c)$ be a solution of \eqref{eq:abc} with $(x,y)=(\tfrac{a}{c},\tfrac{b}{c})$. Our assumption  $\mu(x,y)\not\leq \mu'$ then gives a subset $T\subseteq S$ of cardinality $j=\abs{T}$ in $\{1,\dotsc,t\}$ together with $C\in\{a,b,c\}$ such  that $\ord_q(C)\log q> \mu'_j$ for all $q\in T$. 
Furthermore after replacing $(x,y)$ by a symmetric solution, we may and do assume that $C=c$. It holds that $\ord_q(c)\geq \lfloor \mu'_{\abs{T}}/\log q\rfloor+1$ for all $q\in T$. Hence $n$ divides $c$, and this implies that $a,b$ are both invertible in $\ZZ/n\ZZ$ since $\gcd(a,b,c)=1$. Therefore the equation $a+b=c$ leads to $(a/b)^2=1$ in $(\ZZ/n\ZZ)^\times$. Then on writing $(a/b)^2=\prod_{p\in S\wo T}p^{2\gamma_p}$ with $\gamma_p=\ord_p(a/b)$, we see that the vector $\gamma=(\gamma_p)\in\ZZ^{S\wo T}$  lies  in the lattice $\Gamma_T\subseteq\ZZ^{S\wo T}$ of (i). Moreover,  if $i\in\{1,\dotsc,\abs{S\wo T}\}$ then  the $i$-th largest of the real numbers $\abs{\gamma_p\log p},$ $p\in S\wo T$, is at most $\max(\mu_\iota(a),\mu_\iota(b))$ for $\iota=\min(\ceil{i/2},t)$, and our assumption $\mu(x,y)\leq \mu''$ implies that $\max(\mu_\iota(a),\mu_\iota(b))\leq \mu''_{\iota}$. Hence we deduce that $\gamma$ lies in the ellipsoid $\mathcal E_T$ of (ii) and then we conclude that (iii) produces our solution $(x,y)$ as desired. 

\paragraph{Application.} For any $M\in \ZZ_{\geq 1}$, we would like to find all solutions $(x,y)$ of the $S$-unit equation~\eqref{eq:sunit} with $M(x,y)\leq M$. For this purpose we can use for example  Algorithm~\ref{algRefinedDeWegerSieve}, which we successively apply with $\mu'(n),\mu''(n)\in\ZZ^t$ 
for $n=M+1, M,\dotsc,1$; where
\begin{equation}\label{def:mubounds}
\mu'(n)=\floor{(n-1)\cdot (1,1/2,\dotsc,1/t)} \textnormal{ for } n\in\{1,\dotsc, M+1\},
\end{equation}
$$\mu''(M+1)=M\cdot (1,\dotsc,1), \ \ \ \mu''(n)=\mu'(n+1) \textnormal{ for } n\in\{1,\dotsc, M\}.$$ 
Here $\floor{v}=(\floor{v_i})$ for $v=(v_i)$ a vector with entries $v_i\in\mathbb R$. Suppose that $\abs{S}\geq 6$, $T=\{q\}$ and $\mu''=\mu''(n)$ with $n\in\{1,\dotsc,M\}$. Then the ``radius" $R$ of the ellipsoid in  \eqref{eqBoundForAlphaForImprovedDeWeger} satisfies 
\begin{equation}\label{eq:rbound}
R=\sum_{i=1}^{\abs{S}-1}(\mu''_{\iota(i)})^2\leq n^2\left(\tfrac{t+1}{t^2}+4\sum_{i=1}^{2t}\tfrac{1}{i^2}\right)
\end{equation}
for $\iota(i)=\min(\ceil{i/2},t)$. It follows that $R\leq 7n^2$, and this uniform bound together with $R<(\abs{S}-1)n^2$ shows  that $R$ is considerably smaller (in particular for large $\abs{S}$) than the ``radius" $(\abs{S}-1)n^2$ of the ellipsoid  in Algorithm~\ref{algDeWegersSieve} with $u=(\floor{n/\log p})$. The following complexity discussion provides some motivation for our choice of $t$ and $\mu'(n),\mu''(n)$. 

\paragraph{Complexity.} To discuss the improvements  provided by our refined algorithm for sets $S$ with $\abs{S}\geq 6$, we consider $M\in\ZZ_{\geq 2}$ and we apply Algorithm~\ref{algRefinedDeWegerSieve} with $\mu'(n),\mu''(n)$  for some $n\in\{1,\dotsc, M\}$. Recall that Algorithm~\ref{algRefinedDeWegerSieve} needs to compute in particular $\Gamma_{T}\cap\mathcal E_T$ for all $T\subseteq S$ with $\abs{T}\in\{1,\dotsc,t\}$. Here  the cases $T=\{q\}$ are essentially an application of Algorithm~\ref{algDeWegersSieve} as in \eqref{def:ulbounds}, with $M'=n-1$ and $M''=n$. However  in light of the discussions surrounding~\eqref{eq:rbound}, the crucial difference is that the volume of the ellipsoid $\mathcal E_{\{q\}}$ of our refined Algorithm~\ref{algRefinedDeWegerSieve} is considerably smaller than the volume of the corresponding ellipsoid $\mathcal E$ of Algorithm~\ref{algDeWegersSieve}. In practice, this is the reason for the  significantly improved running time of Algorithm~\ref{algRefinedDeWegerSieve} for $\abs{S}\geq 6$. We note that our refinement needs to iterate in addition  over certain sets $T$ with  $\abs{T}\geq 2$.  However  these additional iterations have little influence on the running time, since the running time for $\abs{T}\geq 2$ is in practice much better  than for $\abs{T}=1$. Indeed if $\abs{T}\geq 2$ then the lattices $\Gamma_T$ have smaller rank $\abs{S}-\abs{T}$, which crucially improves the running time  of (FP) in Algorithm~\ref{algRefinedDeWegerSieve}~(ii).

\begin{remark}[Implementation of Fincke--Pohst]\label{rem:fp}
To avoid numerical issues in Algorithms \ref{algDeWegersSieve} and \ref{algRefinedDeWegerSieve}, we use our own implementation of the version of the Fincke--Pohst algorithm \cite{fipo:algo} described in (FP). Our implementation only uses integer arithmetic, and in particular we do not take square roots.
For this purpose the coordinates and the bound in~\eqref{eqBoundForAlphaForImprovedDeWeger} have been scaled and rounded to integers in such a way, that (FP) will return possibly slightly more candidates $\gamma$ than those fulfilling~\eqref{eqBoundForAlphaForImprovedDeWeger}.
Furthermore  we use an LLL improvement of the original Fincke--Pohst algorithm, such as for example the one in Cohen~\cite{cohen:compant}. This is important since the original implementation~\cite{fipo:algo} becomes in many instances very slow, in fact already for $|S|\geq 10$ it is too slow for our purpose.
\end{remark}

\subsubsection{Refined enumeration}
\label{sec:suenum}
We continue the notation introduced above and we take  $u\in\ZZ^S$ with $u\geq 0$. To find all solutions $(x,y)$ of the $S$-unit equation~\eqref{eq:sunit} satisfying $m(x,y)\leq u$, 
we use the following refined enumeration Algorithm \ref{algEnumerationOfTinySolutions}. In practice, our refined enumeration algorithm is considerably faster than the standard algorithm which is described in the complexity discussion below. For any subset $T\subseteq S$, we define its weight $w(T)=\prod_{p\in T}(1+u_p)$.

\begin{Algorithm}[Refined enumeration]
\label{algEnumerationOfTinySolutions}
The input is a finite set of rational primes $S$, together with a vector $u\in\ZZ^{S}$ such that $u\geq 0$. The output consists of all solutions $(x,y)$ of the $S$-unit equation~\eqref{eq:sunit} that satisfy $m(x,y)\leq u$. 

The algorithm: For each subset $S_a\subseteq S$ with $w(S_a)\geq w(S)^{1/3}$, do the following.
\begin{enumerate}
\item[(i)]
Split the set $S_a=S_{a_1}\dotcup S_{a_2}$ into disjoint subsets $S_{a_1},S_{a_2}$ such that $w(S_{a_1})\leq w(S)^{1/2}$ is fulfilled as tight as possible.
Construct the set $X$, implemented as a hash, of all $a_1\in \ZZ$ such that $\abs{a_1}=\prod_{p\in S_{a_1}} p^{v_p(a_1)}$ with $0\leq v_p(a_1)\leq u_p$ for $v_p=\ord_p$.
\item[(ii)] 
Construct the set $Y=\cup Y(S_b,S_c)$ with the union taken over all pairs $(S_b,S_c)$ such that  $S_a\dotcup S_b\dotcup S_c=S$ and such that either $w(S_b)<w(S_c)$ or $w(S_b)=w(S_c)$ and  $\min S_b< \min S_c$. Here $Y(S_b,S_c)$ is the set of $(b,c,a_2)\in\ZZ^3$ such that $b=\prod_{p\in S_b} p^{v_p(b)}$ with $1\leq v_p(b)\leq u_p$,  $c=\prod_{p\in S_c} p^{v_p(c)}$ with $1\leq v_p(c)\leq u_p$,  $a_2=\prod_{p\in S_{a_2}} p^{v_p(a_2)}$ with $0\leq v_p(a_2)\leq u_p$ and such that $a_2$ divides $b+c$ or $b-c$. 
\item[(iii)]
For each  $(b,c,a_2)\in Y$, check if $a_1:=(b+\eps c)/a_2\in X$ for some $\eps\in\{1,-1\}$.
If so then output all $(x,y)=(\tfrac{\alpha}{\gamma},\tfrac{\beta}{\gamma})$ with $\alpha+\beta=\gamma$ and $\{\abs{\alpha},\abs{\beta},\abs{\gamma}\}=\{\abs{a_1}a_2,b,c\}$.
\end{enumerate}
\end{Algorithm}

In the implementation of the above algorithm, we simultaneously carry out  steps (ii) and (iii) as follows. We iterate over all $(S_b,S_c)$ and over all $(b,c)$: For each  $(b,c)$ we determine the divisors $a_2\in \ZZ_{\geq 1}$ of $b\pm c$ which are only divisible by primes in $S_{a_2}$ and which have exponents bounded by $u$. In the same iteration step,  we also check whether $(b\pm c)/a_2\in X$ and if this is the case then we output the corresponding solutions.

\paragraph{Correctness.} We now show that Algorithm~\ref{algEnumerationOfTinySolutions} indeed finds all solutions $(x,y)$ of~\eqref{eq:sunit} with $m(x,y)\leq u$. Suppose that $(x,y)$ is such a solution, and let $(a,b,c)$ be a solution of~\eqref{eq:abc} with $(x,y)=(\tfrac{a}{c},\tfrac{b}{c})$. After replacing $(x,y)$ with a symmetric solution, we may and do assume that $w(S_a)\geq \max(w(S_b),w(S_c))$ for $S_b=\{p\st p\divides b\}$, $S_c=\{p\st p\divides c\}$ and $S_a=S- (S_b\cup S_c)$. This implies that  $w(S_a)\geq w(S)^{1/3}$, since $w(S)=w(S_a)w(S_b)w(S_c)$. Let $S_{a_1},S_{a_2},X,Y$ be the sets appearing in steps (i) and (ii), and define $a_i=\prod_{p\in S_{a_i}}p^{\ord_p(a)}$ for $i=1,2$. Our assumption $m(x,y)\leq u$ implies that $a_1\in X$ and   $(\abs{b},\abs{c},a_2)\in Y$ since $a_2\mid a=-(b-c)$. Further we observe that $(\abs{b}+\varepsilon \abs{c})/a_2=\pm a/a_2=\pm a_1\in X$ for some $\epsilon\in\{1,-1\}$, and by construction it holds that $a+b=c$ with $\{\abs{a},\abs{b},\abs{c}\}=\{a_1a_2,\abs{b},\abs{c}\}$. Therefore we see that step (iii) produces our solution $(x,y)$ as desired.

\paragraph{Complexity.}   To explain the improvements provided by Algorithm~\ref{algEnumerationOfTinySolutions}, we recall that the standard enumeration of all solutions $(x,y)$ of~\eqref{eq:sunit} with $m(x,y)\leq u$ is as follows: Consider all coprime pairs $(a,b)\in\ZZ\times \ZZ$ with $\rad(ab)\divides N_S$, and with $\max\bigl(\ord_p(a),\ord_p(b)\bigl)\leq u_p$ for all $p\in S$.
If $c=a+b$ satisfies $\rad(c)\divides N_S$ and $\ord_p(c)\leq u_p$ for all $p\in S$, then output the solution $(x,y)=(\tfrac{a}{c},\tfrac{b}{c})$. Now the improved  running time of Algorithm \ref{algEnumerationOfTinySolutions} has the following basic reason.
For fixed $S_a$, $S_b$ and $S_c$, we iterate over all $a_1$ and over all $(b,c)$ in a subsequent way; see  the remark given below Algorithm~\ref{algEnumerationOfTinySolutions}. That is the running time of these two iterations adds, and it does not multiply as in the standard enumeration which iterates over all coprime pairs $(a,b)\in\ZZ\times\ZZ$. Here we note that the splitting of $S_a$ into $S_a=S_{a_1}\dotcup S_{a_2}$  assures that the running time  of (i) does not differ too much from the running time of (ii) together with (iii). Indeed $S_{a_1}$ is chosen such that $w(S_{a_1})$ is approximately $w(S)^{1/2}$ and our assumption $w(S_a)\geq w(S)^{1/3}$ provides that $w(S_b)w(S_c)\leq w(S)^{2/3}$. In particular, for bounded $\abs{S}$ the complexity of our refined enumeration  is asymptotically the $2/3$-th power of the complexity of the standard enumeration algorithm.

\begin{remark}Due to the hash our refinement needs more memory than the standard algorithm.
In practice this is not much compared to the running time,
and in case this becomes an issue we can avoid creating the hash as follows. Iterating over all $(b,c)$ as above, we try to factor $(b\pm c)/a_2$ using only primes in~$S_{a_1}$  and if this succeeds then we output the corresponding solutions provided the exponents are bounded by~$u$.
\end{remark}

\paragraph{Application.}We shall apply Algorithm~\ref{algEnumerationOfTinySolutions} in order to enumerate the solutions $(x,y)$ of \eqref{eq:sunit} with bounded $\mu(x,y)$. More precisely, let $\mu\in\ZZ^t$ with $\mu\geq 0$ and suppose that we want to enumerate all solutions $(x,y)$ of \eqref{eq:sunit} with $\mu(x,y)\leq \mu$. For this purpose it suffices to apply Algorithm~\ref{algEnumerationOfTinySolutions} with $u\in\ZZ^S$ given by $u_p =\floor{\mu_1/\log p}$ for $p\in S$. Indeed any solution $(x,y)$ with $\mu(x,y)\leq \mu$ satisfies $m_p(x,y)\log p\leq\mu_1(x,y)\leq \mu_1$ for all $p\in S$ and thus $m(x,y)\leq u$. This application could output too many solutions $(x,y)$, since not any $(x,y)$ with $m(x,y)\leq u$ satisfies $\mu(x,y)\leq \mu$. To avoid this we can directly  check the latter condition  at each step of the recursions building $X$ and $Y$, and we do this in our implementation. In general an improvement could come from a choice of the weight $w$ making the cardinalities of $X$ and $Y$ even more balanced; we leave this for the future.

\subsubsection{The algorithm} \label{sec:sualgoresults}
We continue the notation introduced above. On putting everything together, we obtain the following algorithm  which solves the $S$-unit equation~\eqref{eq:sunit}.

\begin{Algorithm}[$S$-unit equation via height bounds]\label{algo:suheight}
The input is a finite set of rational primes~$S$, and the output is the set of solutions $(x,y)$ of the $S$-unit equation~\eqref{eq:sunit}. 
\begin{enumerate}
\item[(i)] To rule out solutions with ``large" height, use Proposition~\ref{prop:algobounds} in order to compute $M_0\in\ZZ_{\geq 1}$ such that all solutions $(x,y)$ of \eqref{eq:sunit}  satisfy $M(x,y)\leq M_0$.

\item[(ii)] To find the  solutions of ``medium sized" and ``small" height,  apply de Weger's sieve.

\begin{enumerate}
\item
Find a small $M_1\in\ZZ_{\geq 1}$  such that an application of Algorithm~\ref{algDeWegersSieve}~(ii) as in \eqref{def:ulbounds}, with $M'=M_1$ and $M''=M_0$,  only returns $\gamma=0$.
Here first try $M_1=10$. If this does not work then replace $M_1$ by $\floor{1.3M_1}$, and so on.

\item Having found such an $M_1$,  successively apply  Algorithm~\ref{algDeWegersSieve} as in \eqref{def:ulbounds}, with $M'=M_{k+1}$,  $M''=M_{k}$ and $M_{k+1}=\floor{M_k/1.3}$, for $k=1,2,\dotsc$ until either $M_k=0$ or  Algorithm~\ref{algDeWegersSieve}~(ii) returns more than $10^3$ candidates.
\end{enumerate}

\item[(iii)]
To find the remaining solutions, combine the refined sieve with the refined enumeration. The following two steps run simultaneously, or alternatingly, until the parameters $n,n'$ appearing in these two steps satisfy $n=n'$.
\begin{enumerate}
\item
To enumerate the remaining solutions from above, let $k_0$ be the $k$ where the above step (ii) ended. Then successively apply  Algorithm \ref{algRefinedDeWegerSieve} as in \eqref{def:mubounds}, with $M=M_{k_0}$ and $\mu'(n),\mu''(n)$, for $n=M+1,M,\dotsc$ until $n=n'$.
\item
To enumerate the remaining solutions from below, successively apply the refined enumeration Algorithm~\ref{algEnumerationOfTinySolutions} with $\mu'(n')$ for $n'=1,2,\ldots$ until $n'=n$.
\end{enumerate}
\end{enumerate}
\end{Algorithm}

\paragraph{Correctness.} To verify that this algorithm works correctly, we let $(x,y)$ be a solution of~\eqref{eq:sunit} and we let $M$ be as in step (iii). The construction of $M_0$ in step (i) shows that $M(x,y)\leq M_0$ and hence our solution $(x,y)$ is found in step (ii) if $M(x,y)>M$. 
On the other hand, if $M(x,y)\leq M$ then step (iii) produces our solution $(x,y)$ as desired.

\subsubsection{Complexity}\label{sec:sucomplexity} 
We conducted  some effort to optimize the running time in practice. To explain our optimizations, we now  discuss the above composition of  Algorithm~\ref{algo:suheight} and we motivate our choice of the parameters appearing therein. 
Furthermore, we also consider some additional practical aspects of our algorithm. We continue the notation introduced above.

\paragraph{Step (i).} In this step we need to find a  number $M_0$ with the property that any solution $(x,y)$ of \eqref{eq:sunit} satisfies $M(x,y)\leq M_0$. For this purpose we use Proposition~\ref{prop:algobounds} which requires to compute the number $\alpha(N)$ appearing therein, where $N=2^4N_S$. In case the computation of $\alpha(N)$ takes too long, we can replace   $\alpha(N)$ by the slightly larger number $\bar{\alpha}(N)$ defined below \eqref{def:barb} and this $\bar{\alpha}(N)$ can always be computed very fast. Hence we see that Proposition~\ref{prop:algobounds} allows in any case to quickly determine a relatively small number $M_0$ with the desired property. In practice this $M_0$ is small enough such that we can skip de Weger's first reduction process (2) described at the beginning of Section~\ref{sec:suheightalgo}, and thereby remove an uncertainty in de Weger's original method which crucially relies on (2). Indeed the reduction process (2) is (a priori) an uncertainty, since it is not proved that it always works.  However, we should mention that in practice it is (essentially) always possible to successfully apply (2) by adapting the parameters to the specific situation at hand.

\paragraph{Step (ii).} Here we apply  Algorithm~\ref{algDeWegersSieve} as in \eqref{def:ulbounds}. For this purpose we divided the space $0\leq M(x,y)\leq M_0$ into  subspaces $M'< M(x,y)\leq M''$, with $M'=M_{k+1}$ and $M''=M_k$ for $k=0,\dotsc,k_0$ and for some convenient sequence $M_{k_0}<\ldots< M_{1}<M_0$. 
To explain for which choices of $M',M''$ the application of Algorithm~\ref{algDeWegersSieve} is fast in practice, we use the notation of Section~\ref{sec:dwsieve}. The efficiency of the sieve depends on the number of points in $\Gamma\cap \mathcal E$. In the best case, the sieve \eqref{eqDeWegerSieve} decreases the space of candidates by a factor which is approximately $\abs{G}/2= q^{l_q}(q-1)/2$.
However in the worst case, when the square of each element in $S\wo q$ is $1$ modulo~$q^{l_q+1}$, we obtain $\Gamma=\ZZ^{S\wo q}$ and  the sieve \eqref{eqDeWegerSieve} does not decrease the number of candidates at all.
 In practice we are almost always close to the best case. 
Let $V$ be the euclidean volume of the ball in $\RR^{s-1}$ of square radius $s-1$ for $s=\abs{S}$, and let $\textnormal{covol}(\Gamma)$ denote the covolume of $\Gamma$. 
The  ellipsoid $\mathcal E$ has volume $\vol(\mathcal E)=V\prod_{p\in S\setminus q}u_p$, 
and  in the generic case the cardinality of $\Gamma\cap \mathcal E$ can be approximated  by $\vol(\mathcal E)/\textnormal{covol}(\Gamma)$. 
Thus the sieve is efficient in practice if the ratio  $\vol(\mathcal E)/(\abs{G}/2)$ is ``small". For example this ratio is strictly smaller than $1$ when $M'$ and $M''$ satisfy
\begin{equation}\label{def:m0m1}
M'>(s-1)\log M''+\log V+2\log 2-\sum_{p\in S\setminus q}\log\log p.
\end{equation}
Stirling's approximation leads to a simpler expression for $V$ in terms of $s$. If we now choose $M_1=(s-1)(\log M_0+\tfrac{1}{2}\log(2\pi e))$, then \eqref{def:m0m1} suggests that the sieve \eqref{eqDeWegerSieve} is  strong enough  such that Algorithm~\ref{algDeWegersSieve} (ii) only returns the trivial candidate $\gamma=0$. A relatively small $M_1$ with this property is produced in step~(a) where we start with a very optimistic choice $M_1=10$.  
Step (a) is fast in practice and it improves the running time of (ii) as follows:  Algorithm~\ref{algDeWegersSieve} (iii) is trivial\footnote{The construction of $M_1$ provides that here Algorithm~\ref{algDeWegersSieve} (iii) only needs to output the trivial solutions $(2,-1)$, $(-1,2)$ and $(\frac{1}{2},\frac{1}{2})$ of \eqref{eq:sunit}, which all come from the trivial candidate $\gamma=0$.}  for the large space $M_1<M(x,y)\leq M_0$,  and  for the  space  $M_2<M(x,y)\leq M_1$ the running time of Algorithm~\ref{algDeWegersSieve} (iii) is considerably faster than for $M_2<M(x,y)\leq M_0$ since $M_0$ is much larger than $M_1$. Indeed Algorithm~\ref{algDeWegersSieve}~(ii) applied with $M_0$ can produce much larger candidates $\gamma$, and for large candidates $\gamma$ the reconstruction process in Algorithm~\ref{algDeWegersSieve}~(iii) becomes slow. We next discuss step (b). According to \eqref{def:m0m1}, the choices $M_{k+1}=\floor{M_{k}/1.3}$ for $k\geq 1$ should give a strong sieve in the range where step (b) is applied and this turned out to be true in practice.  We apply step (b)  for $k=1,2,\dotsc$ until Algorithm~\ref{algDeWegersSieve} (ii) returns more than $10^3$ candidates. Here the condition more than $10^3$ candidates means that our refined sieve can find these candidates considerably faster, and thus we switch at this point to step (iii).

\paragraph{Step (iii).} In this step  we are in a situation where we can fully exploit our refinements worked out in the previous sections. To explain this more precisely, we now mention two  points which significantly slow down  Algorithm~\ref{algDeWegersSieve} in the situation of step (iii) where many solutions exist. The first point is the application of (FP) in Algorithm~\ref{algDeWegersSieve}~(ii), which  is the bottleneck of Algorithm~\ref{algDeWegersSieve}. The second point is that Algorithm~\ref{algDeWegersSieve}~(ii) repetitively enumerates the ``same" candidate $\gamma$  in many steps $k$ with $M_k$ small. This is due to the large fraction between the volume of the ellipsoid and the volume of the cube  appearing in Algorithm~\ref{algDeWegersSieve}; indeed  this fraction depends exponentially on the cardinality of $S$. 
To improve these two points we worked out the following refinements:
\begin{enumerate}
\item Concerning the first and second point, we developed our refined Algorithm~\ref{algRefinedDeWegerSieve} which works with smaller ellipsoids. This leads to less repetitions of the candidates $\gamma$, and it considerably improves the running time of the step in which we apply (FP). See also the complexity discussions given in Section \ref{sec:dwsieve+}.
\item Regarding the second point, we constructed the refined enumeration  Algorithm~\ref{algEnumerationOfTinySolutions}. This  enumeration  is fast enough such that  we can now skip the final applications ($M_k=1,2,\dotsc$) of Algorithm~\ref{algDeWegersSieve} which are very slow, and thereby we can in particular circumvent for many candidates $\gamma$ that they get enumerated repetitively. 
\end{enumerate}  
In  (iii) we carry out steps (a) and (b) alternatingly, depending on which step took less time so far. In some cases this leads to a significantly improved running time of  (iii). Further we mention that our choice $\mu''(n)=\mu'(n+1)$ is motivated by \eqref{def:m0m1}. Indeed  according to \eqref{def:m0m1}, one should work with step size $1=(n+1)-n$ in the situation of (iii) where usually $M(x,y)\leq (s-1)(\log s+\log\log s+\tfrac{1}{2}\log(2\pi e))$. See also the discussions surrounding~\eqref{eq:rbound}  which provide additional motivation for our definition of $\mu'(n),\mu''(n)$.

\paragraph{Bottleneck.} Despite our refinements which considerably improve the running time in practice, step (iii) still remains in general the bottleneck of Algorithm~\ref{algo:suheight}. However, for certain special sets $S$ the location of the bottleneck can change.  For example if $S$ consists of the prime 2 together with a few very large primes, then Algorithm~\ref{algo:suheight} often finds all solutions already in step (ii). In this case the bottleneck of Algorithm~\ref{algo:suheight} is located in step (ii) where we apply Algorithm~\ref{algDeWegersSieve}~(i). The main reason is that here  the application of (T) becomes slow since $S$ contains very large primes, and the application of (FP) in Algorithm~\ref{algDeWegersSieve}~(ii) becomes fast since the cardinality of $S$ is small.

\begin{remark}[Parallelization]\label{rem:parallelization}
In our implementation of Algorithm~\ref{algo:suheight} we successfully parallelized essentially everything, except Algorithm~\ref{algDeWegersSieve}~(i) and Algorithm~\ref{algRefinedDeWegerSieve}~(i) which both  involve the application  of Teske's algorithm described in (T).
\end{remark}

\subsubsection{Applications}\label{sec:suapplications}

In this section we give some applications of Algorithm~\ref{algo:suheight}.  In particular we discuss parts of our database $\mathcal D_1$ containing the solutions of the $S$-unit equation \eqref{eq:sunit} for many distinct sets $S$. We also use our database to motivate various Diophantine questions related to \eqref{eq:sunit}, including Baker's explicit $abc$-conjecture and a new conjecture.

We continue the notation introduced above. Recall that $\Sigma(S)$ denotes the set of solutions of the $S$-unit equation~\eqref{eq:sunit} up to symmetry, and for any $n\in\ZZ_{\geq 1}$ we recall that  $S(n)$ denotes the set of the $n$ smallest rational primes. Let $N\in\ZZ_{\geq 1}$ and define $\Sigma(N)=\cup\Sigma(S)$ with the union taken over all finite sets of rational primes $S$ with $N_S\leq N$. 

\paragraph{The sets $\Sigma(S(n))$.} We determined the sets $\Sigma(S(n))$ for all $n\leq 16$. The cardinality of these sets  is given in the table of Theorem~A stated in the introduction. As already mentioned, our Algorithm~\ref{algo:suheight} substantially improves de Weger's original method in \cite{deweger:lllred} which de Weger used to compute the set $\Sigma(S(6))$ in \cite[Thm 5.4]{deweger:lllred}. To illustrate our improvements in practice, we now compare Algorithm~\ref{algo:suheight} with de Weger's original method.  For this purpose we  implemented in Sage de Weger's original method in a slightly improved form (dW), which uses in addition our optimized height bounds. If $S=S(6)$ then (dW)  took 21 seconds, while it took 6 seconds by using Algorithm~\ref{algo:suheight}. For larger $|S|$ our running time improvement significantly increases:
For example if $S=S(10)$ then (dW) takes four days, whereas this decreases to only 25 minutes by using Algorithm~\ref{algo:suheight}. 
Roughly speaking, for large $\abs{S}$ our refinements should save an exponential factor with respect to $|S|$ in comparison to de Weger's original method. Further if  $n>10$ then (dW) becomes too slow and thus  we did not try to use (dW) in order to compute $\Sigma(S(n))$ for $n>10$. To deal with $S=S(n)$ for $10<n\leq 16$ we additionally parallelized (see Remark~\ref{rem:parallelization}) our Algorithm~\ref{algo:suheight}. Then it took 8 days for $n=15$ and 34 days for $n=16$, using 30 CPU's.  

\begin{remark}[Automatically choosing parameters]\label{rem:fullyautomatic}
De Weger's original method  does not specify in general how to choose the involved parameters. For example de Weger's first reduction process requires to make a choice, and to efficiently apply de Weger's sieve one has to choose suitable subspaces   dividing the initial space.
To compute the set $\Sigma(S(6))$, de Weger has chosen by hand the required parameters. Although we can skip de Weger's first reduction process using our optimized height bounds,  we still need to make many choices in Algorithm~\ref{algo:suheight}.
In particular, it is now favourable that the choices required for Algorithm~\ref{algo:suheight} are made by an automatism. 
We implemented such an automatism,  which takes into account \eqref{def:m0m1} and which properly adjusts the parameters during run time. In view of \eqref{def:m0m1}, our automatism chooses parameters for Algorithm~\ref{algo:suheight} such that the running time is considerably less than twice as long as for the optimal parameters. 
However, we do not claim and can not prove that our automatism is optimal.\end{remark}

\paragraph{The sets $\Sigma(N)$.}
We computed the sets $\Sigma(N)$ for all $N\leq 10^7$ in approximately 13 days.  For this computation it was  crucial that  Algorithm~\ref{algo:suheight} automatically chooses all required parameters, as mentioned in Remark~\ref{rem:fullyautomatic}. Indeed to compute the sets $\Sigma(N)$ for all $N\leq 10^7$, we had to apply Algorithm~\ref{algo:suheight} with so many distinct sets $S$ such that it would have been impossible to suitably choose by hand all the involved parameters. 

\paragraph{Explicit $abc$-conjecture.}  Baker \cite[Conj 4]{baker:abcexperiments} proposed the following fully explicit version of the $abc$-conjecture: If $a,b,c\in \ZZ_{\geq 1}$ are coprime with $a+b=c$, then it holds $c\leq\tfrac{6}{5}N(\log N)^{\omega}/\omega!$ for   $\omega$ the number of rational primes dividing the radical $N=\rad(abc)$; here one should exclude the triple $(a,b,c)=(1,1,2)$.  On using our  database $\mathcal D_1$ containing in particular the sets $\Sigma(N)$ for all $N\leq 10^7$, we verified  Baker's explicit $abc$-conjecture  for all coprime $a,b,c\in\ZZ_{\geq 1}$ with  $\rad(abc)\leq 10^7$. Furthermore, Algorithm~\ref{algo:suheight} can be used to verify other properties of all $abc$-triples with bounded radical. In particular all existing $abc$-triples with $\rad(abc)\leq 10^7$ can be directly taken from our database $\mathcal D_1$.

\paragraph{Elliptic curves with full $2$-torsion.}
For any given $N\in\ZZ_{\geq 1}$ we denote by $\mathcal M$ the set of $\QQ$-isomorphism classes of elliptic curves over $\QQ$ of conductor $N$, with all two torsion points defined over $\QQ$. Algorithm~\ref{algo:suheight} allows to efficiently compute the set $\mathcal M$. To see this, we let $S$ be the set of rational primes dividing $2N$. Any elliptic curve in $\mathcal M$ admits a Weierstrass equation $y^2=x(x-a)(x+b)$ such that $a,b\in\ZZ$ have the following properties: $d=\gcd(a,b)$ divides $N_S$ and  $(\frac{a}{a+b},\frac{b}{a+b})$ satisfies the $S$-unit equation~\eqref{eq:sunit}. An application of Algorithm~\ref{algo:suheight} with $S$ determines all possible values of $a/d$ and $b/d$, and this then allows to directly compute the set $\mathcal M$. In particular, we see that  Algorithm~\ref{algo:suheight} provides an alternative way to check the completeness of a small part of Cremona's database \cite{cremona:algorithms}. This application of Algorithm~\ref{algo:suheight} already turned out to be useful in practice.

\paragraph{Conjecture and question.} We next use our database $\mathcal D_1$  to motivate various questions concerning the solutions of $S$-unit equations \eqref{eq:sunit}. First we recall Conjecture~1 which is motivated by our data and by the construction of the refined sieve in Section~\ref{sec:dwsieve+}.

\vspace{0.3cm}
\noindent{\bf Conjecture 1.}
\emph{There exists $c\in\ZZ$ with the following property: If $n\in\ZZ_{\geq 1}$ then any finite set of rational primes $S$ with  $\abs{S}\leq n$ satisfies $\abs{\Sigma(S)}\leq \abs{\Sigma(S(n))}+c$.}
\vspace{0.3cm}

\noindent In other words, if $\mathcal T$ denotes the collection of schemes $T$ such that $T$ can be obtained by removing $n$ closed points of $\sp(\ZZ)$, then Conjecture~1 means the following: Among all $T\in \mathcal T$, the maximal number (up to some constant) of $T$-points of $\mathbb P^1_\ZZ-\{0,1,\infty\}$ is attained  at the scheme in $\mathcal T$ which corresponds to the $n$ closed points of smallest norm. This conjectured property of $\mathbb P^1_\ZZ-\{0,1,\infty\}$ is rather unexpected from a general Diophantine geometry perspective. We further ask whether one can remove the constant.

\vspace{0.3cm}
\noindent{\bf Question 1.1.}
\emph{Does \textnormal{Conjecture}~$1$ hold with $c=0$?}
\vspace{0.3cm}

\noindent Here the main motivation is given by our data. Indeed Question~1.1 has a positive answer for all sets $S$ in our database $\mathcal D_1$.  In view of Theorem A listing the cardinality of $\Sigma(S(n))$ for all $n\in\{1,\dotsc,16\}$, a positive answer to Question 1.1 would give an optimal upper bound for the number of solutions of any $S$-unit equation~\eqref{eq:sunit} with $\abs{S}$ in  $\{1,\dotsc,16\}$.

\subsection{Comparison of algorithms}
To compare Algorithm~\ref{algo:suheight} with our Algorithm \ref{algo:sucremona} using modular symbols, we continue the notation introduced above. Recall that we already computed the sets $\Sigma(N)$ for all $N\leq 20000$ by using Algorithm~\ref{algo:sucremona}, see the examples in Section~\ref{sec:sucremalgo}.  Here it turned out that the output of Algorithm~\ref{algo:sucremona} agrees with the output of Algorithm~\ref{algo:suheight}. 
To determine all solutions of the $S$-unit equation \eqref{eq:sunit},  we recommend to use Algorithm~\ref{algo:sucremona} as long as one already knows the set of elliptic curves over $\QQ$ of conductor dividing $2^4N_S$. If this set is not already known, then it is usually much more efficient to use Algorithm~\ref{algo:suheight}. In fact, as demonstrated in the previous sections, our Algorithm~\ref{algo:suheight} is substantially more efficient in all aspects than the known methods which practically resolve \eqref{eq:sunit}.

\section{Algorithms for Mordell equations}\label{sec:malgo}
Let $S$ be a finite set of rational primes, write $N_S=\prod_{p\in S}p$ and put $\mathcal O=\ZZ[1/N_S]$.
We take a nonzero $a\in \mathcal O$. In this section, we would like to solve the Mordell equation
\begin{equation}
y^2=x^3+a, \ \ \ (x,y)\in\mathcal O\times\mathcal O. \tag{\ref{eq:mordell}}
\end{equation}
This Diophantine problem is a priori more difficult than solving the $S$-unit equation~\eqref{eq:sunit}. Indeed elementary transformations reduce~\eqref{eq:sunit}  to~\eqref{eq:mordell}, while the known reductions of~\eqref{eq:mordell} to~\eqref{eq:sunit} require to solve~\eqref{eq:sunit} in number fields or they require a height bound for the solutions of~\eqref{eq:sunit} which is equivalent to the $abc$-conjecture in Remark~\ref{rem:abc}. 

So far all known practical methods solving~\eqref{eq:mordell}  crucially rely on the theory of logarithmic forms \cite{bawu:logarithmicforms}, see below for an overview. In the following sections we present two alternative algorithms which allow to practically  resolve~\eqref{eq:mordell}. They both do not use the theory of logarithmic forms. In Section~\ref{sec:mcremalgo} we describe the first algorithm which relies on Cremona's algorithm using modular symbols. Then in Section~\ref{sec:mheightalgo} we construct the second algorithm via height bounds.  Here we also give several applications and we discuss various questions motivated by our results, see Sections~\ref{sec:shaf}-\ref{sec:malgoapplications}. Finally in Section~\ref{sec:malgocomparison} we compare our algorithms with the actual best practical methods solving \eqref{eq:mordell}.

\paragraph{Known methods.} We now discuss algorithms and methods in the literature which allow to solve \eqref{eq:mordell}. First we consider the classical case $\mathcal O=\ZZ$. Ellison et al~\cite{ellison:mordell} used the approach of Baker--Davenport~\cite{baker:contributions,bada:diophapp} to solve~\eqref{eq:mordell} for some $a$. Recently the latter approach was refined by Bennett--Ghadermarzi~\cite{begh:mordell} who applied their algorithm to find all solutions of~\eqref{eq:mordell} in $\ZZ\times\ZZ$ for any nonzero $a\in\ZZ$ with $\abs{a}\leq 10^7$; see also the work of Wildanger and J\"atzschmann discussed in Fieker--Ga\'al--Pohst~\cite[p.739]{figapo:mordellcharp}. Alternatively, Masser~\cite{masser:ellfunctions}, Lang \cite{lang:diophantineanalysis}, W\"ustholz~\cite{wustholz:recentprogress} and Zagier \cite{zagier:largeintegralpoints} initiated a practical approach to solve arbitrary elliptic Weierstrass equations $(W)$ over $\ZZ$ via elliptic logarithms. On applying this approach with David's explicit bounds~\cite{david:elllogmemoir}, Stroeker--Tzanakis~\cite{sttz:elllogaa} and Gebel--Peth{\H{o}}--Zimmer~\cite{gepezi:ellintpoints}  obtained independently a practical algorithm solving $(W)$ over $\ZZ$. Gebel--Peth{\H{o}}--Zimmer~\cite{gepezi:mordell} used this algorithm to determine all solutions of \eqref{eq:mordell} in $\ZZ\times\ZZ$ for any nonzero $a\in \ZZ$ with $\abs{a}\leq 10^4$ and for most $a\in\ZZ$ with $\abs{a}\leq 10^5$.  Let $r$ be the Mordell--Weil rank of the group $E(\QQ)$ associated to the elliptic curve $E$ over $\QQ$ defined by $(W)$. In the important special case $r=1$,  there exists in addition a practical approach of Balakrishnan--Besser--M\"uller \cite{babemu:qchabcrelle,babemu:qchabalgo} which is in the spirit of Kim's non-abelian Chabauty program initiated in \cite{kim:siegel}.   

We now discuss practical methods in the literature solving \eqref{eq:mordell} over any ring $\mathcal O$ as above. In practice and in theory, this task is considerably more difficult 
than solving \eqref{eq:mordell} in $\mathcal O=\ZZ$. The method of Bilu~\cite{bilu:superellalgo} and Bilu--Hanrot~\cite{biha:algosuperell} for superelliptic Diophantine equations  allows in particular to solve \eqref{eq:mordell}. 
Further, classical constructions reduce \eqref{eq:mordell} to Thue--Mahler equations which in turn can be solved using the method of Tzanakis--de Weger \cite{tzde:thuemahler}. 
Smart~\cite{smart:sintegralpoints}  extended the above mentioned elliptic logarithm approach to solve  $(W)$ over $\mathcal O$. His algorithm is conditional on explicit lower bounds for linear forms in $p$-adic elliptic logarithms.  R\'emond--Urfels~\cite{reur:padicelllog} proved such bounds\footnote{Hirata-Kohno~\cite{hirata-kohno:p-adicelllogs} recently established the general case. Tzanakis~\cite[Chapt 11]{tzanakis:book} combined her bounds with the elliptic logarithm reduction to solve $(W)$ over $\mathcal O$, see also Hirata-Kohno--Kov\'acs~\cite{hiko:rank3}.} for $r\leq 2$, and these bounds were then applied by  Gebel--Peth{\H{o}}--Zimmer~\cite{gepezi:bordeaux,gepezi:barcelona} to solve \eqref{eq:mordell} for some nonempty $S$.  Furthermore, in a joint work with Herrmann~\cite{pezigehe:sintegralpoints}, they obtained a variation of the elliptic logarithm approach which works also for  $r\geq 3$, see Section \ref{sec:mheightalgo}.  We point out that the elliptic logarithm approach  requires a basis of $E(\QQ)$. There are methods which often can compute such a basis in practice, in particular in our case of elliptic curves defined by \eqref{eq:mordell}. However,  these methods are not (yet) effective in general as discussed in Section~\ref{sec:tor+mwbasis}. For a detailed description of the elliptic logarithm approach, we refer to the excellent book of Tzanakis~\cite{tzanakis:book} which is devoted to this method.

\subsection{Algorithm via modular symbols}\label{sec:mcremalgo}

We continue our notation. For any elliptic curve $E$ over $\QQ$, we  denote by $c_4$ and $c_6$ the usual quantities associated to a minimal Weierstrass model of $E$ over $\ZZ$; see~\cite{tate:aoe}. Write $N_E$ and $\Delta_E$ for the conductor and minimal discriminant of $E$ respectively. We define 
\begin{equation}\label{def:as}
a_S=1728N_S^2\prod p^{\min(\ord_p(a),2)}
\end{equation}
with the product taken over all rational primes $p$ not in~$S$. If $(x,y)$ satisfies the Mordell equation~\eqref{eq:mordell}, then one can consider the elliptic curve $E$ over $\QQ$ which admits the Weierstrass equation $t^2 = s^3 - 27xs - 54y$ with ``indeterminates" $s$ and~$t$.  This construction leads to the following lemma  which will be proven in course of the proof of Lemma~\ref{lem:pm2}.

\begin{lemma}\label{lem:pm1}
Suppose that $(x,y)$ is a solution of the Mordell equation~\eqref{eq:mordell}. Then there exists an elliptic curve $E$ over $\QQ$ such that $N_E\mid a_S$ and such that $c_4 =u^4x$ and $c_6= u^6y$ for some $u\in\QQ$ with $u^{12}=1728\Delta_E\abs{a}^{-1}$.
\end{lemma}

In fact this lemma may be viewed as an explicit \parshin{}  construction for integral points on the moduli scheme of elliptic curves defined by~\eqref{eq:mordell}, see \cite[Prop 3.4]{rvk:modular} for details.

\begin{Algorithm}[Mordell equations via modular symbols]\label{algo:mcremona} The input consists of a finite set of rational primes $S$ together with a nonzero number $a\in\mathcal O$, and the output is the set of solutions $(x,y)$ of the Mordell equation \eqref{eq:mordell}. 

\begin{itemize}
\item[(i)] Define $a_S$ as in~\eqref{def:as}. Then use Cremona's algorithm using modular symbols,   described in Section~\ref{sec:cremonas+st}, to compute the set $\mathcal T\subset \ZZ\times\ZZ$ of quantities $(c_4,c_6)$ which are associated to some modular elliptic curve over $\QQ$ of conductor dividing~$a_S$. 
\item[(ii)]  For each $(c_4,c_6)\in\mathcal T$, write $\abs{(c_4^3-c_6^2)/a}=\frac{m}{n}$ with coprime $m,n\in\ZZ_{\geq 1}$. Compute the subset $\mathcal T_{0}\subseteq \mathcal T$ of $(c_4,c_6)\in\mathcal T$ with  $m=u_1^{12}$ and $n=u_2^{12}$ for $u_1,u_2\in\ZZ$.
\item[(iii)] For any $(c_4,c_6)\in \mathcal T_{0}$ take $u=\frac{u_1}{u_2}$ in $\mathbb \QQ$ with $u^{12}=\abs{(c_4^3-c_6^2)/a}$ and define $x=u^{-4}c_4$ and $y=u^{-6}c_6$. If $x$ and $y$ are both in $\mathcal O$, then output $(x,y)$.
\end{itemize}
\end{Algorithm}
\paragraph{Correctness.}We now verify that this algorithm indeed finds all solutions of any Mordell equation~\eqref{eq:mordell}. Suppose that $(x,y)$ satisfies~\eqref{eq:mordell}. Then Lemma~\ref{lem:pm1} gives an elliptic curve $E$ over $\QQ$ such that $N_E\mid a_S$ and such that $c_4 =u^4x$ and $c_6= u^6y$ for some $u\in\QQ$ with $u^{12}=1728\Delta_E\abs{a}^{-1}$. The Shimura--Taniyama conjecture assures that $E$ is modular. This proves that $(c_4,c_6)$ is contained in the set $\mathcal T$ computed in step (i). Furthermore we obtain that $(c_4,c_6)\in\mathcal T_{0}$, since it holds  $u^{12}=1728\Delta_E\abs{a}^{-1}=\abs{(c_4^3-c_6^2)/a}$ by the definition of the discriminant. Any $u'\in\QQ$ with $u'^{12}=\abs{(c_4^3-c_6^2)/a}$ satisfies $u'^4=u^4$ and $u'^6=u^6$. Therefore we see that step (iii) produces our solution $(x,y)=(u^{-4}c_4,u^{-6}c_6)$ as desired. 

\paragraph{Complexity.} The discussion of the complexity of step (i) is analogous to the discussion of the complexity of Algorithm~\ref{algo:sucremona}~(i) and thus we refer the reader  to Section~\ref{sec:sucremalgo}. For each $(c_4,c_6)\in\mathcal T$, step (ii) needs to check if there is $u=\frac{u_1}{u_2}\in\QQ$ with $u^{12}=\abs{(c_4^3-c_6^2)/a}$ and then step (iii) needs to check if $u_1\in\ZZ$ is only divisible by primes in~$S$. Therefore we see that the complexity discussions of Algorithm~\ref{algo:sucremona} together with the arguments in Remark~\ref{rem:sulegendrealternative} imply the following: The running time of step (ii) and (iii) is at most polynomial in terms of $H(a)N_S$, and is at most $O_\eps((H(a)N_S)^\eps)$ if all three conjectures (BSD), (GRH) and $(abc)$ hold.  Here $H(a)=\exp(h(a))$ denotes the multiplicative  Weil height  of~$a$.

\paragraph{Applications.} We recall from \cite[12.5.2]{bogu:diophantinegeometry} that a solution $(x,y)\in\ZZ\times\ZZ$ of~\eqref{eq:mordell} is called primitive if $\pm 1$ are the only $m\in\ZZ$ with $m^6\mid \gcd(x^3,y^2)$.  For each nonzero $a\in\ZZ$ one can quickly enumerate all $(x,y)\in\ZZ\times\ZZ$ with $y^2=x^3+a$  if one knows the primitive solutions of any Mordell equation $y^2=x^3+a'$ with $a'\in\ZZ$ satisfying $a=a'm^6$ for some $m\in\ZZ$. On using an implementation in Sage of Algorithm~\ref{algo:mcremona}, we computed the set of primitive solutions of the family of Mordell equations~\eqref{eq:mordell} with parameter $a\in\ZZ-\{0\}$ satisfying $r_2(a)\leq 200$ for $r_2(a)=\prod p^{\min(2,\ord_p(a))}$ with the product taken over all rational primes~$p$.   Here the computation of the solutions was very fast, since for each $a\in\ZZ-\{0\}$ with  $r_2(a)\leq 200$  part~(i) of  Algorithm~\ref{algo:mcremona} can use Cremona's database which  contains in particular the required data for all elliptic curves over $\QQ$ of conductor dividing $1728r_2(a)< 350 000$. On the other hand, if the required data of the involved elliptic curves is not already known, then our Algorithm~\ref{algo:mcremona} is usually not practical anymore. Here the problem is the application  of Cremona's algorithm (using modular symbols) in step (i), which requires a huge amount of memory to deal with conductors that are not small.

\subsection{Algorithm via height bounds}\label{sec:mheightalgo}

In this section we use the optimized height bounds of Proposition~\ref{prop:algobounds} to construct Algorithm~\ref{algo:mheight} which allows to solve the Mordell equation~\eqref{eq:mordell}. We continue our notation and we work with the following setup:  We may and do view $y^2=x^3+a$ as a Weierstrass equation of an elliptic curve $E_a$ over $\QQ$, since our $a\in\mathcal O$ is nonzero. A classical result of Mordell gives that the abelian group $E_a(\QQ)$ is finitely generated. Let $P_1,\dotsc,P_r$ be a basis of the free part of $E_a(\QQ)$, and let $E_a(\QQ)_{\textnormal{tor}}$ be the torsion group of $E_a(\QQ)$. We call $r$ the Mordell--Weil rank of $E_a(\QQ)$ and we say that $P_1,\dotsc,P_r$ is a Mordell--Weil basis of $E_a(\QQ)$. Let $(x,y)$ be a solution of \eqref{eq:mordell}. The corresponding point $P$ in $E_a(\QQ)$ takes the form $P=Q+\sum_{i=1}^r n_i P_i$ with $n_i\in\ZZ$ and $Q\in E_a(\QQ)_{\textnormal{tor}}$, and we define $N(x,y)=\max_{i}\abs{n_i}$. 

\paragraph{Decomposition.} Before we describe our Algorithm~\ref{algo:mheight} in detail, we discuss its decomposition which is inspired by the elliptic logarithm approach introduced by Masser--Zagier. A variation of the latter approach was used for example in the algorithm of Peth{\H{o}}--Zimmer--Gebel--Herrmann \cite{pezigehe:sintegralpoints} whose main ingredients  can be described as follows:
\begin{itemize}
\item[(1)]\label{itDeWegerStepBaker} First they try to find a Mordell--Weil basis of $E_a(\QQ)$.
\item[(2)]\label{itDeWegerStepLLL} On using the explicit result of Hajdu--Herendi~\cite{hahe:elliptic} which is based on the theory of logarithmic forms~\cite{bawu:logarithmicforms}, they determine an initial upper bound $N_0$ such that any solution $(x,y)$ of~\eqref{eq:mordell} satisfies $N(x,y)\leq N_0$.
\item[(3)]\label{itDeWegerStepSieve} Following Smart~\cite{smart:sintegralpoints} they apply the elliptic logarithm reduction in order to reduce the initial upper bound $N_0$ to a bound $N_1$ which is usually much smaller.
\item[(4)] \label{itDeWegerStepEnumerationByHand}
Finally they enumerate all solutions $(x,y)$ of \eqref{eq:mordell} with $N(x,y)\leq N_1$. Here in the case $\mathcal O=\ZZ$ the ``inequality trick" usually improves the enumeration process.
\end{itemize}
 Our Algorithm~\ref{algo:mheight} substantially improves in all aspects the elliptic logarithm approach for \eqref{eq:mordell}. More precisely, without introducing new ideas we use the known methods for (1) described in Section~\ref{sec:tor+mwbasis}. To obtain the initial upper bound  for Algorithm~\ref{algo:mheight}, we apply in Proposition~\ref{prop:mwbound} the optimized height bounds  worked out in Section~\ref{sec:heightbounds}. In practice our initial bound is considerably stronger than the initial bound $N_0$ in (2) based on the theory of logarithmic forms, and this leads to  significant  running time improvements of the  reduction process as illustrated in Section~\ref{sec:minitbounds}. Then to enumerate the solutions of bounded height, we apply the elliptic logarithm sieve constructed in Section~\ref{sec:elllogsieve}. This sieve is substantially more efficient than the reduction process (3) together with the subsequent enumeration (4). We now discuss the ingredients of our approach in more detail.

\subsubsection{Torsion group and Mordell--Weil basis}\label{sec:tor+mwbasis}

On using the notation introduced above, we next briefly discuss methods which allow to determine the torsion group $E_a(\QQ)_{\textnormal{tor}}$ of $E_a(\QQ)$ and a Mordell--Weil basis of $E_a(\QQ)$. 

\paragraph{Torsion.} Fueter~\cite{fueter:kubische} completely determined  $E_a(\QQ)_{\textnormal{tor}}$.  To state his result we write $a=e/d$ with coprime $e,d\in\ZZ$. If $d^5e=k\cdot l^6$ for $k,l\in\ZZ$ with $k$ sixth power free, then
\begin{equation}\label{eq:mtorsion}
E_a(\QQ)_{\textnormal{tor}}\cong\begin{cases}
\ZZ/6\ZZ  & \textnormal{if } k=1,\\
\ZZ/3\ZZ & \textnormal{if } k\neq 1 \textnormal{ is a square, or } k=-432,\\
\ZZ/2\ZZ & \textnormal{if } k\neq 1 \textnormal{ is a cube,}\\ 
0 & \textnormal{otherwise.}
\end{cases}
\end{equation} 
This completely determines all solutions of \eqref{eq:mordell} in $\QQ\times\QQ$ in the case when the Mordell--Weil rank $r$ of $E_a(\QQ)$ satisfies $r=0$. Therefore for the purpose of determining all solutions of the Mordell equation \eqref{eq:mordell} we always may assume that $r\geq 1$.

\paragraph{Mordell--Weil basis.} The problem of finding a Mordell--Weil basis of $E_a(\QQ)$ is usually more difficult. If the elliptic curve $E_a$ satisfies the part of the Birch--Swinnerton-Dyer conjecture (BSD) predicting that $r$ coincides with the analytic rank of $E_a$, then a result of Manin \cite{manin:cyclomodcurves} leads to a practical algorithm (see for example Gebel--Zimmer~\cite{gezi:mwalgo}) 
which computes such a basis. On applying this algorithm and an algorithm of Cremona given in~\cite{cremona:algorithms}, Gebel--Peth{\H{o}}--Zimmer~\cite{gepezi:mordell} computed a Mordell--Weil basis of $E_a(\QQ)$  for most $a\in\ZZ$ with $\abs{a}\leq 10^5$. Parts of their database (Mordell$\pm$) are uploaded on their homepage \url{tnt.math.se.tmu.ac.jp/simath/MORDELL}. 
On using  (PSM) we checked their data for all nonzero $a\in \ZZ$ with $\abs{a}\leq 10^4$. Here it turned out that for $a = 7823$ a basis  was missing and that for $a = -7086$ and $a=-6789$ the given bases  were not saturated. In all these three cases we determined a Mordell--Weil basis of $E_a(\QQ)$ using (PSM) and now the (updated)  database  contains a correct Mordell--Weil basis of $E_a(\QQ)$ for any nonzero $a\in \ZZ$ with $\abs{a}\leq 10^4$. We note that in our special case given by the Mordell elliptic curve $E_a$, one can often exploit isogenies to find a Mordell--Weil basis of $E_a(\QQ)$. In fact, it turned out in practice  that   the known techniques implemented in  (PSM) (Generators, two-descent, HeegnerPoint, etc.) usually allow to quickly determine such a basis. 
However we point out that in the case of an arbitrary nonzero $a\in\ZZ$ there is so far no unconditional method which allows to determine a Mordell--Weil basis of $E_a(\QQ)$.

\subsubsection{Elliptic logarithm sieve}

Starting with Zagier \cite{zagier:largeintegralpoints} and de Weger~\cite{deweger:phdthesis}, many authors  developed over the last decades the elliptic logarithm reduction process; see for example Stroeker--Tzanakis \cite{sttz:elllogaa}, Gebel--Peth{\H{o}}--Zimmer~\cite{gepezi:ellintpoints} and Smart~\cite{smart:sintegralpoints}. In practice this process allows to show that the solutions in $\mathcal O\times\mathcal O$ of an arbitrary elliptic Weierstrass equation  have either relatively small or huge height. In Section~\ref{sec:elllogsieve} we constructed the elliptic logarithm sieve which considerably improves the elliptic logarithm reduction and the subsequent enumeration of solutions of small height.
In particular, for any given bound $N\in\ZZ$, the elliptic logarithm sieve solves the problem of efficiently enumerating all solutions $(x,y)$ of \eqref{eq:mordell} with $N(x,y)\leq N$. The sieve combines the core idea of the elliptic logarithm reduction with several conceptually new ideas. We refer to  Section~\ref{sec:elllogsieve} for an overview of the new ideas introduced by the elliptic logarithm sieve and for a detailed discussion of the practical and theoretical improvements provided by these ideas.

\subsubsection{Initial bounds}\label{sec:minitbounds}

We continue the notation introduced above.  In this section we give an initial upper bound for various  heights attached to the solutions of the Mordell equation \eqref{eq:mordell}. We also compare our bound with the actual best results in the literature and we explain how our  result improves the  running time of the reduction process in the elliptic logarithm sieve.

 We recall that $P_1,\dotsc,P_r$ denotes a Mordell--Weil basis of $E_a(\QQ)$, and for any solution $(x,y)$ of \eqref{eq:mordell} we write as above $N(x,y)=\max \abs{n_i}$ for the ``infinity norm" of the corresponding point $P=Q+\sum n_i P_i$ in $E_a(\QQ)$. Let $\hat{h}$ be the canonical N\'eron-Tate height on $E_a(\QQ)$. Here we work with the natural normalization of $\hat{h}$ which divides by the degree of the involved rational function, see for example \cite[p.248]{silverman:aoes}.

\paragraph{Initial bounds.} Let $(x,y)$ be a solution of \eqref{eq:mordell}, with corresponding point $P\in E_a(\QQ)$. To deduce an upper bound for $N(x,y)$, we recall standard properties of  $\hat{h}$. One can control $\hat{h}(P)$ in terms of the usual logarithmic Weil height $h$  as follows 
\begin{equation}\label{eq:canonheightcomp}
\hat{h}(P)\leq \tfrac{1}{2}h(x)+\tfrac{m}{6}h(a)+1.58, \ \ \ \ m =\begin{cases}
 1 & \textnormal{if } a\in\ZZ,\\
12 & \textnormal{otherwise.}
\end{cases}
\end{equation}
Indeed in the integral case $a\in\ZZ$ the displayed inequality directly follows for example from Silverman~\cite[Thm 1.1]{silverman:heightcomparison}. 
To deal with any nonzero $a\in\mathcal O$, we write $a=c/d$ with coprime $c,d\in\ZZ$ and we consider  $E_{b}$ for $b=d^6a\in\ZZ$. There is an isomorphism $\varphi:E_a\to E_{b}$ induced by $x\mapsto d^2x$. Then an application of the already verified integral case of \eqref{eq:canonheightcomp}, with $E_b$ and $b\in\ZZ$,  gives an upper bound for $\hat{h}(P)=\hat{h}(\varphi(P))$ in terms of $h(d^2x)$ and $h(b)$ which proves \eqref{eq:canonheightcomp} as desired. 
Further it is known that $\hat{h}$ defines a positive definite quadratic form on the real vector space $E(\QQ)\otimes_\ZZ\RR$. On using the basis $P_1,\dotsc,P_r$ we identify $E(\QQ)\otimes_\ZZ\RR$ with $\RR^r$ and then we denote by $\lambda$ the smallest eigenvalue of the matrix defining the binary form associated to the quadratic form  $\hat{h}$ on $\RR^r$. 

We now use the optimized height bound in Proposition~\ref{prop:algobounds}, or the $abc$-conjecture of Masser-Oesterl\'e~\cite{masser:abc} stated in Remark~\ref{rem:abc}, in order to obtain initial bounds. 

\begin{proposition}\label{prop:mwbound}
Suppose that $(x,y)$ is a solution of \eqref{eq:mordell}, and denote by $P$ the corresponding point in $E_a(\QQ)$. Then the following statements hold.
\begin{itemize}
\item[(i)] Let $\alpha=\alpha(a_S)$  be the number from Proposition~\ref{prop:algobounds}, and recall that $m=1$ if $a\in\ZZ$ and $m=12$ otherwise. It holds $\lambda N(x,y)^2\leq\hat{h}(P) \leq M_0$ for some $M_0\in\ZZ$ with $$M_0\leq \tfrac{m+1}{6}h(a)+2\alpha+\log(\alpha+16.52)+52.12.$$
\item[(ii)] Suppose that $a\in\ZZ$, and assume that the $abc$-conjecture holds. Then for any real number $\epsilon>0$ there exists a constant $c_\epsilon$ depending only on $\epsilon$ such that $$\hat{h}(P)\leq (1+\epsilon)(\log N_S+\tfrac{7}{6}h(a))+c_\epsilon.$$
\end{itemize}
\end{proposition}

In case the computation of the number $\alpha$ from Proposition~\ref{prop:algobounds} takes too long, one can replace $\alpha$ by the (slightly) larger number $\bar{\alpha}=\bar{\alpha}(a_S)$ which is defined below \eqref{def:barb}.  The number $\bar{\alpha}$ has the advantage that  it can be quickly computed in all cases. We further mention that the bound in Proposition~\ref{prop:mwbound}~(ii)  is a direct consequence of a result in Bombieri--Gubler~\cite[Thm 12.5.12]{bogu:diophantinegeometry}; this bound is optimal in terms of $N_S$. 
\begin{proof}[Proof of Proposition~\ref{prop:mwbound}]
We first prove (i).  Linear algebra leads to  $\lambda N(x,y)^2\leq\hat{h}(P)$, since $\hat{h}$ is a quadratic form on  $E(\QQ)\otimes_\ZZ\RR$. Further, the estimate for $M_0$ in (i) follows by combining \eqref{eq:canonheightcomp}  with the upper bound for $h(x)$ given in Proposition~\ref{prop:algobounds}. 

To show assertion (ii) we assume that $a\in \ZZ$ and we write $x=x_1/d^2$ and $y=y_1/d^3$ with $x_1,y_1,d\in\ZZ$ satisfying $\gcd(d,x_1y_1)=1$ 
and $d>0$. Let $n$ be the largest element in $\ZZ$ with $n^6\mid \gcd(x_1^3,y_1^2)$. On dividing the equation $x_1^3-y_1^2=-ad^6$ by $n^6$, we obtain a new equation of the form $u^3-v^2=w$ with $u,v,w\in\ZZ$ and $\gcd(u^3,v^2)=1$. Now, on assuming the $abc$-conjecture, we see that \cite[Thm 12.5.12]{bogu:diophantinegeometry} gives estimates for $\abs{u},\abs{v}$. These estimates lead to an upper bound for $h(x)$ which together with \eqref{eq:canonheightcomp} proves  (ii). 
\end{proof}

\paragraph{Comparison with literature.} There are several explicit bounds for $M_0$ and $N(x,y)$ in the literature. They are all based on the theory of logarithmic forms. In fact this theory allows to effectively solve Diophantine equations which are considerably more general than Mordell equations \eqref{eq:mordell}, see for example \cite{bawu:logarithmicforms}. To compare Proposition~\ref{prop:mwbound}~(i) with the actual best bounds for \eqref{eq:mordell} in the literature, we use a simpler but weaker version of our bound. On replacing in the proof of Proposition~\ref{prop:mwbound}~(i) our optimized height bounds by the simplified height bounds in Proposition~\ref{prop:m}, we obtain 
\begin{equation}\label{eq:simplmwbound}
\lambda N(x,y)^2\leq M_0\leq \tfrac{m+1}{6}h(a)+\tfrac{1}{2}a_S\log a_S.
\end{equation}
For the purpose of the following discussion, we recall that in the case  $a\in\ZZ$ it holds that $a_S\leq 1728\abs{a}N_S^2$ and $m=1$. 
The actual best effective upper bound for $N(x,y)$ and $M_0$ was established by Peth{\H{o}}--Zimmer--Gebel--Herrmann~\cite[Thm]{pezigehe:sintegralpoints}. 
Their result is based on the work of Hajdu--Herendi~\cite{hahe:elliptic} which in turn relies on the theory of logarithmic forms. To state the rather complicated bound for $N(x,y)$ provided by Peth{\H{o}} et al, we need to introduce some notation.
As in~\cite[Thm]{pezigehe:sintegralpoints} we define the constants 
$$k_3=\frac{32}{3}\Delta_0^{\frac{1}{2}}(8+\frac{1}{2}\log \Delta_0)^4, \ \ \ k_4=10^4\cdot 256\cdot \Delta_0^{\frac{2}{3}}, \ \ \ \Delta_0=27\lvert a\rvert^2.$$  Further, we write $s=\lvert S\rvert$ and $q=\max S$ (with $q=1$ if $S=\emptyset$). Then we define
$$\kappa_1=\tfrac{7}{2}\cdot 10^{38s+87}(s+1)^{20s+35}q^{24}\max(1,\log q)^{4s+2} k_3(\log k_3)^2(k_3+20sk_3+\log (ek_4)).$$
We mention that the result in \cite[Thm]{pezigehe:sintegralpoints} is stated under the assumption that the given Weierstrass equation over $\ZZ$ is minimal at all primes in $S$. However, on looking at the proof one sees that this minimality assumption is not necessary for the portion of the theorem which provides an upper bound for $N(x,y)$. We conclude that \cite[Thm]{pezigehe:sintegralpoints} provides in general that any solution $(x,y)$ of (\ref{eq:mordell}) with $a\in\ZZ-\{0\}$ satisfies
\begin{equation}\label{eq:pzghmwbound}
\lambda N(x,y)^2\leq M_0\leq \kappa_1+\tfrac{1}{3}\log \abs{4\cdot 6^3a}.
\end{equation}
In our simplified bound \eqref{eq:simplmwbound} the dependence on $a\in \ZZ$  is of the form $\lvert  a\rvert \log \lvert a\rvert$, while in \eqref{eq:pzghmwbound} it is of the weaker form $\lvert a \rvert^2(\log \lvert a\rvert)^{10}$. Further we see that \eqref{eq:simplmwbound} considerably improves \eqref{eq:pzghmwbound}  for essentially all  sets $S$ of practical interest, in particular for all sets $S$ with $N_S\leq 2^{1200}$ or $s\leq 12$ and for all sets $S$ of the form $S=S(n)$ where $S(n)$ denotes the set of the first $n$ primes for some $n\in\ZZ_{\geq 1}$. 
 We now choose parameters $\mathcal A$ and $\mathcal S$ as follows: The set $\mathcal S$ is given by $\{\emptyset,S(1),S(10)\}$, and the set  $\mathcal A$ consists of 24 distinct nonzero $a\in\ZZ$  such that for each $r\in\{1,\dotsc,12\}$ there are precisely two $a$ in $\mathcal A$ with $E_a(\QQ)$ of rank $r$; here we tried\footnote{For $r\leq 6$ we found the ``smallest" possible $a$.  Further, we note that for our purpose of illustrating the running time improvements of the reduction process it suffices to work with $r$ independent points of $E_a(\QQ)$; for $r\geq 9$ we could not prove (unconditionally) that our $r$ independent points form a basis.} to choose these elements $a\in\ZZ$ with $\abs{a}$ as small as possible. To illustrate that our bound leads to significant running time improvements,  we computed for all parameter pairs $(a,S)\in\mathcal A\times\mathcal S$ the running times $\rho$ and $\rho^*$ of the elliptic logarithm reduction in Algorithm~\ref{algo:elllogsieve}~(ii) using Proposition~\ref{prop:mwbound}~(i) and \eqref{eq:pzghmwbound} respectively.  In the case $S=\emptyset$, it turned out that we obtain a running time improvement by a factor $\rho^*/\rho$ which is approximately $2$ for small/medium $\abs{a}$ and which is close to $4$ for large $\abs{a}$. The running time improvements become more significant in the case $S=S(1)$. Here the factor $\rho^*/\rho$ is approximately $30$ for small/medium $\abs{a}$ and it is approximately $60$ for large $\abs{a}$.
Finally we achieve big running time improvements when  $S=S(10)$. In this case the factor $\rho^*/\rho$ is approximately $300$ for small $\abs{a}$, it lies between $500$ and $10^3$ for medium sized $\abs{a}$ and it varies between $10^3$ and $10^4$ for large $\abs{a}$. For example if $a=-2520963512$ ($r=8$) and $S=S(10)$ then $\rho$ is less than 23 seconds while $\rho^*$ exceeds 2 days.

In the classical case $\mathcal O=\ZZ$, there is also a fully explicit estimate $N_0\geq N(x,y)$ which was independently established by Stroeker--Tzanakis~\cite{sttz:elllogaa} and Gebel--Peth{\H{o}}--Zimmer~\cite{gepezi:ellintpoints}. This estimate is based on lower bounds for linear forms in elliptic logarithms (see Masser~\cite{masser:ellfunctions}, W\"ustholz~\cite{wustholz:recentprogress}, Hirata-Kohno~\cite{hirata-kohno:ellloginvent} and David~\cite{david:elllogmemoir}).  We do not state here $N_0$ in its precise form, since $N_0$ is even more complicated than the bound in \eqref{eq:pzghmwbound}. To see that our result improves $N_0$ for essentially all $a\in\ZZ$ of practical interest, it suffices to consider the following simpler lower bound 
\begin{equation}\label{eq:tzmwbound}
N_0\geq \lambda^{-1/2}10^{3(r+2)}4^{(r+1)^2}(r+2)^{(r^2+13r+23.3)/2} \prod_{i=1}^{r} \max\bigl(\hat{h}(P_i)^{1/2},\log (4\abs{a})\bigl).
\end{equation}
This lower bound follows for example from \cite[p.386-387]{pezigehe:sintegralpoints}, see also the recent book of Tzanakis~\cite{tzanakis:book}. In \eqref{eq:tzmwbound} we may and do assume that $r\geq 1$ by Fueter's result \eqref{eq:mtorsion}. Then we see that our simplified bound \eqref{eq:simplmwbound}  improves \eqref{eq:tzmwbound} for all nonzero $a\in\ZZ$ with $\abs{a}\leq 10^{40}$.  In the case of arbitrary $\mathcal O\neq \ZZ$ and $r=2$, one can deduce an explicit estimate $N_0\geq N(x,y)$  by using  lower bounds of  David \cite{david:elllogmemoir} and R\'emond--Urfels~\cite{reur:padicelllog}; see also the recent work of Hirata-Kohno--Kov\'acs~\cite{hirata-kohno:p-adicelllogs,hiko:rank3} removing the assumption $r=2$. Here the quantity $N_0$ is very complicated and its dependence  on $S$ is quite involved. In any case $N_0$ is larger than the lower  bound in \eqref{eq:tzmwbound} and thus our result is better than the estimate $N_0\geq N(x,y)$ for all pairs $(a,S)$ with $a_S\leq 10^{40}$. On the other hand, for large $a_S$ it is rather difficult (when not impossible in general) to compare Proposition~\ref{prop:mwbound}~(i) with the corresponding results based on lower bounds for linear forms in elliptic logarithms. The reason is that the involved quantities are quite different. However, the dependence of \eqref{eq:tzmwbound} on the rank $r$ means that our result leads to significant running time improvements  in the notoriously difficult case when $r$ is not small. To illustrate this we computed for all parameter pairs $(a,S)\in\mathcal A\times \mathcal S$ the running times $\rho$ and $\rho'$ of the elliptic logarithm reduction in Algorithm~\ref{algo:elllogsieve}~(ii) using  Proposition~\ref{prop:mwbound}~(i) and $N_0$ respectively. We note that instead of implementing the very complicated estimate $N_0$ in its precise form, we used here the simpler lower bound in \eqref{eq:tzmwbound}. In other words the running time $\rho'$ is slightly too good, which means that our running time improvements are slightly better than illustrated by the numbers appearing in the following discussion.  In the case $S=\emptyset$, we obtain a running time improvement by a factor $\rho'/\rho$ which is approximately $2$ when $2\leq r\leq 4$ and which lies between $3$ and $10$ in the range $5\leq r\leq 12$. The running time improvements become more significant in the case $S=S(1)$. Here the factor $\rho'/\rho$ lies between $2$ and $20$ when $2\leq r\leq 4$, it varies between $50$ and $100$  in the range $5\leq r\leq 8$, and it lies between $500$ and $10^3$ for $9\leq r\leq 12$. 
Finally we obtain big running time improvements when $S=S(10)$. In this case the factor $\rho'/\rho$ varies between $2$ and $10$ in the range $2\leq r\leq 4$, it lies between $30$ and $500$ for $r\leq 5\leq 8$, and it varies between $700$ and $3000$ in the range $r\leq 9\leq 12$. For example, if $S=S(10)$ then there is an $a\in\mathcal A$ with  $r=12$ such that our running time $\rho$ is less than 3 minutes while $\rho'$ is approximately 5 days.

\subsubsection{The algorithm}\label{sec:mordellalgostat}

We continue the notation introduced above. On combining the ingredients of the previous sections, we obtain an algorithm which allows to solve the Mordell equation \eqref{eq:mordell}. Here we point out that our algorithm requires an explicitly given  Mordell--Weil basis of $E_a(\QQ)$. While it is usually possible   to determine such a basis in practice (see Section~\ref{sec:tor+mwbasis}), there is so far no unconditional method which in principle works for an arbitrary nonzero $a\in\ZZ$. In view of this we included a Mordell--Weil basis of $E_a(\QQ)$ in the input.

\begin{Algorithm}[Mordell equations via height bounds]\label{algo:mheight} The inputs are  a finite set of rational primes $S$, a nonzero number $a\in\mathcal O$  and the coordinates of a Mordell--Weil basis of $E_a(\QQ)$. The output is the set of solutions $(x,y)$ of  \eqref{eq:mordell}. 
\begin{itemize}
\item[(i)]  Use Proposition~\ref{prop:mwbound}~(i) to compute an initial bound $M_0$ such that for any solution $(x,y)$ of \eqref{eq:mordell} the corresponding point $P\in E_a(\QQ)$ satisfies $\hat{h}(P)\leq M_0$.  
\item[(ii)] Apply the elliptic logarithm sieve in Algorithm~\ref{algo:elllogsieve} in order to find all solutions $(x,y)$ of \eqref{eq:mordell} with corresponding point $P\in E_a(\QQ)$ satisfying $\hat{h}(P)\leq M_0$.
\end{itemize}
\end{Algorithm}
\paragraph{Correctness.} If $P\in E_a(\QQ)$ corresponds to a solution $(x,y)$ of \eqref{eq:mordell}, then Proposition \ref{prop:mwbound}~(i) gives that $\hat{h}(P)\leq M_0$. Thus the application of the elliptic logarithm sieve in step (ii) produces all solutions of \eqref{eq:mordell} as desired, see Remark~\ref{rem:elllogsievegen} when $a\notin \ZZ$.

\paragraph{Complexity.} We now discuss various aspects which significantly influence the running time of Algorithm~\ref{algo:mheight}. In view of the remark given below Proposition~\ref{prop:mwbound}~(i), the computation of the initial upper bound $M_0$ in step (i) is always very fast. The running time of step (ii) crucially depends on the size of $M_0$, the height $h(a)$, the rank $r$ and the cardinality of $S$.  For a complexity discussion of the elliptic logarithm sieve used in step (ii) we refer to Section~\ref{sec:elllogsieve}. Therein we explain in detail various complexity aspects  and we also discuss in detail the influence of the parameters $M_0,r,h(a)$ and $\abs{S}$ on the running time in practice (and in theory). See also Section~\ref{sec:minitbounds} where we illustrated the improvements provided by the sharpened initial bound $M_0$ obtained in Proposition~\ref{prop:mwbound}~(i). 

\begin{remark}[Generalizations]
Algorithm~\ref{algo:mheight} allows to solve more general Diophantine equations associated to a Mordell  curve, that is an elliptic curve with vanishing $j$-invariant. Assume that we are given the coefficients $a_1\dotsc,a_6\in\QQ$ of a Weierstrass equation
\begin{equation}\label{eq:mweieq}
y^2+a_1xy+a_3y=x^3+a_2x^2+a_4x+a_6
\end{equation}
of an elliptic curve $E$ with $j$-invariant $j=0$, and suppose that we know a basis of the free part of $E(\QQ)$. Then Algorithm~\ref{algo:mheight} allows to find all solutions $(x,y)$  of \eqref{eq:mweieq} in $\mathcal O\times\mathcal O$. Indeed there is an explicit isomorphism which transforms any such solution of \eqref{eq:mweieq}  into a solution of \eqref{eq:mordell} for some explicit $a\in\ZZ-\{0\}$, and hence an application of Algorithm~\ref{algo:mheight} with this $a$ produces the set of solutions of \eqref{eq:mordell} and then of \eqref{eq:mweieq}. In fact we implemented this slightly more general version of Algorithm~\ref{algo:mheight}. To conclude  we mention that   further generalizations are possible by using the arguments in \cite[Chapt 8]{tzanakis:book}.
\end{remark}

\subsubsection{Elliptic curves with good reduction outside a given set of primes}\label{sec:shaf}

We continue the notation introduced above. Let $M(S)$ be the set of $\QQ$-isomorphism classes of elliptic curves over $\QQ$ with good reduction outside a given finite set of rational primes $S$. In this section we apply Algorithm~\ref{algo:mheight} in order to compute the set $M(S)$.

\paragraph{Known methods.} There are already several practical methods in the literature which allow to determine $M(S)$.  Agrawal--Coates--Hunt--van der Poorten~\cite{agcohuva:conductor11} computed the semi-stable locus of $M(\{11\})$ by using an approach via Thue--Mahler equations which ultimately relies on the theory of logarithmic forms.  Their work builds on Coates' effective proof \cite{coates:shafarevich} of Shafarevich's theorem  mentioned in the introduction.
Alternatively one can compute $M(S)$ by using the Shimura--Taniyama conjecture and modular symbols, see Cremona~\cite{cremona:algorithms}.  There are also two more recent approaches which ultimately rely on the theory of logarithmic forms: The method of Cremona--Lingham~\cite{crli:shafarevich} discussed in the introduction, and the very recent approach of Koutsianas~\cite{koutsianas:shafarevich} via $S$-unit equations over number fields. Furthermore, very recently Bennett--Rechnitzer~\cite{bere:compell,bere:compelloneprime}  substantially refined (in particular for $\abs{S}=1$) the above mentioned classical Thue--Mahler approach: In the irreducible case they use the Thue--Mahler algorithm of Tzanakis--de Weger~\cite{tzde:thuemahler} and in the rational two torsion case they apply the algorithm of de Weger~\cite{deweger:phdthesis,deweger:sumsofunits} for sums of units being a square.
Finally, several authors used ingenious ad hoc methods to determine $M(S)$ for specific sets $S$. For an overview, see for example the discussions and references in  \cite[Sect 1]{cremona:condhistory} and \cite[Sect 1]{bere:compelloneprime}.

\paragraph{The algorithm.} On combining Shafarevich's classical reduction to Mordell equations \eqref{eq:mordell} with Algorithm~\ref{algo:mheight}, we obtain an alternative approach to determine $M(S)$. Here we do not use modular symbols or lower bounds for linear forms in logarithms. To state our algorithm, we introduce some terminology. We may and do identify  any $[E]$ in $M(S)$ with the pair $(c_4,c_6)$ associated by Tate \cite[p.180]{tate:aoe} to a minimal  Weierstrass model  of  $E$ over $\ZZ$.   Further, for arbitrary $s,t\in\QQ$ we say that an elliptic curve over $\QQ$ is given by $(s,t)$ if the equation $y^2=x^3-27sx-54t$ defines an affine model of the curve.  If $s,t$ are in $\QQ$ with $s^3-t^2$ nonzero and if $E$ denotes an elliptic curve  over $\QQ$ given by $(s,t)$, then  Tate's algorithm~\cite{tate:algo} allows to  compute the pair $(c_4,c_6)$ associated to a minimal  Weierstrass model  of $E$ over $\ZZ$ and it allows to check whether $[E]$ lies in $M(S)$.

\begin{Algorithm}\label{algo:shaf} The inputs are  a finite set of rational primes $S$  and a Mordell--Weil basis of $E_a(\QQ)$ for all $a=1728w$ with $w\in\ZZ$ dividing $N_S^5$. The output is the set $M(S)$. 

The algorithm: For each $a=1728w$ with $w\in\ZZ$ dividing $N_S^5$, do the following.
\begin{itemize}
\item[(i)] Apply Algorithm~\ref{algo:mheight} in order to determine the set $Y_a(\mathcal O)$ formed by the solutions of the Mordell equation \eqref{eq:mordell} defined by the parameter pair $(a,S)$. 
\item[(ii)] For each $(x,y)\in Y_a(\mathcal O)$ and for any $d\in\ZZ_{\geq 1}$ dividing $N_S$, let $E$ be the elliptic curve over $\QQ$ given by  $(d^2x,d^3y)$ and output the pair $(c_4,c_6)$ associated to a minimal Weierstrass model of $E$ over $\ZZ$ provided that $[E]$ lies in $M(S)$.
\end{itemize}
\end{Algorithm}

\paragraph{Correctness.} We take $[E]=(c_4,c_6)$ in $M(S)$. The minimimal discriminant $\Delta$ of $E$ lies in $\mathcal O^\times$. Hence there are integers $u,w,d\in\mathcal O^\times$, with $d\in\ZZ_{\geq 1}$ dividing $N_S$ and $w$ dividing $N_S^5$, such that $\Delta=-wd^6u^{12}$. We define $x=\tfrac{c_4}{u^4d^2}$ and $y=\tfrac{c_6}{u^6d^3}$.  The formula $1728\Delta=c_4^3-c_6^2$ shows that  $(x,y)$ lies in the set $Y_a(\mathcal O)$  computed in step (i) for $a=1728w$, and the elliptic curve $E$ is given by $(d^2x,d^3y)$. Hence we see that step (ii) produces our $[E]$ as desired.

\paragraph{Complexity.} The running time of Algorithm~\ref{algo:shaf} is essentially determined by step~(i). Therein we compute the sets $Y_a(\mathcal O)$ for all $a=1728w$ with $w\in\ZZ$ dividing $N_S^5$ and for this purpose we need to apply Algorithm~\ref{algo:mheight} with $2\cdot 6^{\abs{S}}$ distinct parameters $a$. Hence the running time of Algorithm~\ref{algo:shaf} crucially depends on $\abs{S}$ and on the complexity of Algorithm~\ref{algo:mheight} which we already discussed in the previous section. In step (ii), it might be possible that one can omit to check  whether $[E]$ lies in $M(S)$. In any case this check is always very quick and  it has no influence on the running time in practice. 

\paragraph{Input obstruction and the family $\mathcal S$.} The input of Algorithm~\ref{algo:shaf} requires $2\cdot 6^{\abs{S}}$ distinct Mordell--Weil bases. In fact, for large $\abs{S}$, it usually happens that one can not determine unconditionally all bases and then our Algorithm~\ref{algo:shaf} can not be used to compute $M(S)$. However, for small $\abs{S}$ it turned out in practice that one can often efficiently compute the required bases by using the known techniques in (PSM). For example, without introducing crucial new ideas, we  computed the required bases for each set $S$ in  $\mathcal S$.  Here $\mathcal S$ is a family of sets which contains in particular the set $S(5)$ and  all sets $S$ with $N_S\leq 10^3$. We observe that any elliptic curve over $\QQ$ with good reduction outside $S$ has conductor dividing  $N_S^{\textnormal{cond}}=\prod_{p\in S} p^{f_p}$, where $(f_2,f_3)=(8,5)$ and  $f_p=2$ if $p\geq 5$. It holds that $N_S^2$ divides $N_S^{\textnormal{\textnormal{cond}}}$, and thus $\mathcal S$ contains in particular all sets $S$ with $N_S^{\textnormal{\textnormal{cond}}}\leq 10^6$.

\paragraph{Applications.} On using Algorithm~\ref{algo:shaf}, we determined the sets $M(S)$ for all $S\in \mathcal S$. This took less than $2.5$ hours for $S=S(5)$,  and on average it took   approximately 30 seconds for sets $S\in\mathcal S$ with $\abs{S}=2$, roughly 2.5 minutes for sets $S\in\mathcal S$ with $\abs{S}=3$ and approximately 8 minutes for sets $S\in\mathcal S$ with $\abs{S}=4$. Here we did not take into account the time required to determine the Mordell--Weil bases for the input. In fact if the bases for the input are not already known, then their computation is usually the bottleneck of our approach to determine $M(S)$ via Algorithm~\ref{algo:shaf}. 
Let $T$ be a nonempty open subscheme of $\sp(\ZZ)$. Inspired by our Conjecture~1 on $T$-points of $\mathbb P^1_\ZZ-\{0,1,\infty\}$, we propose the following analogous conjecture on $T$-points of the moduli stack $\mathcal M_{1,1}$ of elliptic curves.

\vspace{0.3cm}
\noindent{\bf Conjecture 1 for $\mathcal M_{1,1}$.}
\emph{Does there exist $c\in\ZZ$ with the following property: If $n\in\ZZ_{\geq 1}$ then any  set of rational primes $S$ with $\abs{S}\leq n$ satisfies $\abs{M(S)}\leq \abs{M(S(n))}+c?$}
\vspace{0.3cm}

\noindent Our database listing the sets $M(S)$ for all $S\in\mathcal S$ shows the following: For any $n\in\ZZ_{\geq 1}$ and for each $S\in\mathcal S$ with $\abs{S}\leq n$, it holds that $\abs{M(S)}$ is at most $\abs{M(S(n))}$.  In light of this we ask whether the above conjecture is true with the optimal constant $c=0$?

We point out that for certain sets $S\in\mathcal S$ one can compute the spaces $M(S)$ by using different methods.  For example Cremona--Lingham~\cite{crli:shafarevich}  determined $M(S)$ for $S=\{2,p\}$ with $p\leq 23$,  and Koutsianas~\cite{koutsianas:shafarevich} moreover computed $M(S)$ for $S=\{2,3,23\}$ and $S=\{2,p\}$ with $p\leq 127$. Further Cremona's database  \cite{cremona:algorithms} allows to directly determine the space $M(S)$ for all sets $S$ with  $N_S^{\textnormal{cond}}\leq 380000$ (as of February 2016). We also mention that the case $\abs{S}=1$ was studied by Edixhoven--Groot--Top in \cite{edgrto:primecond}. In particular they showed  that $M(\{p\})$ is empty for many rational primes $p$, see \cite[Cor 1]{edgrto:primecond} which explicitly lists such primes $p$.  Furthermore very recently Bennett--Rechnitzer~\cite{bere:compelloneprime}  determined $M(\{p\})$ for all primes $p< 2\cdot 10^9$, and for all $p< 10^{12}$ conditional on an explicit version of Hall's conjecture with ``Hall ratio" $10^{14}$. Here to prove their unconditional results, they  exploit that $\abs{S}=1$ in order to reduce to Thue equations which can be solved much more efficiently than Thue--Mahler equations. Their ingenious reduction uses in particular a result of Mestre--Oesterl\'e~\cite[Thm 2]{meoe:weilcurvediscriminants} which in turn relies (inter alia) on the geometric version of the Shimura--Taniyama conjecture.

\paragraph{Comparison.} We now briefly discuss advantages and disadvantages of the different methods which allow to compute $M(S)$ in practice. Our Algorithm~\ref{algo:shaf} significantly improves the method (CL) of Cremona--Lingham~\cite{crli:shafarevich}.  Indeed our Algorithm~\ref{algo:mheight} is considerably more efficient in solving \eqref{eq:mordell} than the algorithm \cite{pezigehe:sintegralpoints}  used in (CL).  To  illustrate that our improvements are significant, we used (CL) to determine $M(S)$ for $S=S(3)$.  This took more than 35 minutes\footnote{This is a lower bound for the time required for (CL) to solve the involved equations \eqref{eq:mordell}. Here we used the official Sage implementation of \cite{pezigehe:sintegralpoints} which works with an ``absolute" reduction process. This means that the running times of (CL) are in fact  larger than the numbers we listed.}, while it took Algorithm~\ref{algo:shaf} less than 2 minutes. Furthermore there are several sets $S\in\mathcal S$ which seem to be out of reach for (CL). For example in the case $S=S(4)$ it took Algorithm~\ref{algo:shaf} less than 20 minutes to compute $M(S)$, while (CL) did not terminate  within 2 months. In comparison with the other practical methods which allow to compute $M(S)$, the main disadvantage of our approach and of (CL) is that they both require $2\cdot 6^{\abs{S}}$ distinct Mordell--Weil bases.  The  modular symbols method (see Cremona \cite{cremona:algorithms}) can efficiently compute the curves in $M(S)$ with small conductor, while the curves of large conductor cause memory problems. We note that even for relatively small sets $S$ the maximal conductor $N_S^{\textnormal{cond}}$ can be large. For example if $S$ contains $\{2,3,p\}$ for some $p\geq 13$ then it holds that $N_S^{\textnormal{cond}}\geq 10^8$ and thus the practical computation of $M(S)$ seems to be out of reach for the modular symbol method.  On the other hand, the modular symbol method deals much more efficiently with the important related problem of compiling a database which lists all elliptic curves over $\QQ$ of given conductor $N\leq 380000$. The efficiency of the approach of Koutsianas~\cite{koutsianas:shafarevich} strongly depends on the size of $\abs{S}$ and on the involved number fields (quadratic, cubic, or $S_3$-extension) in which one has to solve the unit equations. We point out that (CL) and the method of Koutsianas \cite{koutsianas:shafarevich} both allow to deal with more general number fields $K$, while our approach currently only works in the considerably simpler case $K=\QQ$.  As already mentioned, Bennett--Rechnitzer~\cite{bere:compell,bere:compelloneprime} substantially refined the classical Thue--Mahler approach in order to compute $M(\{p\})$ for all primes $p<2\cdot 10^9$.  This computation is unfavorable for our method, since finding unconditionally all the required Mordell--Weil bases would (when possible) take a long time with the known techniques.  In general, \cite{bere:compelloneprime} crucially depends on the algorithms \cite{tzde:thuemahler,deweger:sumsofunits} for which we are  not aware of a complexity analysis. In particular, if $\abs{S}\neq 1$ then it is not clear to us how efficient is the Thue--Mahler  approach of \cite{bere:compell,bere:compelloneprime}. 

To compare some data, we computed the space $M(S)$ for all sets $S$ considered in the papers of Cremona--Lingham~\cite{crli:shafarevich} and Koutsianas~\cite{koutsianas:shafarevich} and for all sets $S$ which can be covered by Cremona's database \cite{cremona:algorithms} (as of February 2016). In all cases it turned out that our Algorithm~\ref{algo:shaf} produced exactly the same number of curves.

\subsubsection{Integral points on moduli schemes}\label{sec:moduli}
We continue our notation. Many Diophantine equations can be reduced via the moduli formalism to the study of $M(S)$. To explain this more precisely, we use the notation and terminology of \cite[Sect 3]{rvk:modular}.  The set $M(S)$ identifies with the set $M(T)$ of isomorphism classes of elliptic curves over the open subscheme $T$ of $\sp(\ZZ)$ given by $T=\sp(\ZZ)-S$.  Let $Y$ be a $T$-scheme and suppose that $Y=M_\mathcal P$ is a moduli scheme of elliptic curves.  We further assume that  the \parshin{} construction $\phi:Y(T)\to M(T)$, induced by forgetting the level structure $\mathcal P$, is effective in the sense that for each $E\in M(T)$ one can determine the set $\mathcal P(E)$.   Then  \cite[Thm 7.1]{rvk:modular} and the discussions in \cite[Sect 3]{rvk:modular} show that one can in principle determine $Y(T)$.  Furthermore, if $S\in \mathcal S$ then one can indeed determine $Y(T)$ by applying our explicit results for $M(T)=M(S)$.  This strategy allows to efficiently solve various classical Diophantine problems, including the following  equations.
\begin{itemize}
\item[(i)]  One can directly solve the $S$-unit equation \eqref{eq:sunit} for any set $S\in\mathcal S$. Here one works with the  moduli problem $\mathcal P=[Legendre]$ as in the proof of \cite[Prop 3.2]{rvk:modular}.
\item[(ii)] We can directly solve any Mordell equation \eqref{eq:mordell} defined by  $(a,S)$ such that  $6a$ is invertible in $\ZZ[1/N_{S'}]$ for some  $S'\in \mathcal S$ with $S\subseteq S'$. Here one works with the moduli problem $\mathcal P_b=[\Delta=b]$ as in the proof of \cite[Prop 3.4]{rvk:modular}, where $1728b=-a$. 
\item[(iii)] One can directly solve any cubic Thue equation~\eqref{eq:thue} defined by  $(f,S,m)$ such that $6\Delta m$ is invertible in $\ZZ[1/N_{S'}]$ for some $S'\in\mathcal S$ with $S\subseteq S'$, where $\Delta$ denotes the discriminant of $f$. Here one works with the moduli problem obtained by pulling back the problem $\mathcal P_b$, with $4b=-\Delta m^2$, along the morphism $\varphi$ given in \eqref{eq:thuemap}.
\item[(iv)] We can directly solve any cubic Thue--Mahler equation~\eqref{eq:thue-mahler} defined by  $(f,S,m)$ such that $6\Delta m$ is invertible in $\ZZ[1/N_{S'}]$ for some $S'\in\mathcal S$ with $S\subseteq S'$. Here we work with the moduli problem $\mathcal P$ represented by the elliptic curve $E$ over the moduli scheme $Y=\sp\bigl(\mathcal O[x,y,\tfrac{1}{d}]\bigl)$ for $d=6\Delta f^2$, where $E$ is given by the closed subscheme of $\mathbb P^2_Y$ defined by $v^2w=u^3+3\mathcal Huw^2+Jw^3$ with $\mathcal H$ and $J$ the covariants (Hessian and Jacobian) of the cubic form $f$ normalized as in \eqref{def:covarpol2} and \eqref{def:covarpol3}.
\item[(v)] We can directly solve any generalized Ramanujan--Nagell equation~\eqref{eq:rana} defined by  $(b,c,S)$ such that $2bc$ is invertible in $\ZZ[1/N_{S'}]$ for some $S'\in\mathcal S$ with $S\subseteq S'$. Here we work with the moduli problem $\mathcal P$ represented by the elliptic curve $E$ over the moduli scheme $Y=\sp\bigl(\ZZ[a_2,a_4,\tfrac{1}{\delta}]\bigl)$ for $\delta=16a_4^2(a_2^2-4a_4)$, where $E$ is given by the closed subscheme of $\mathbb P^2_Y$ defined by $v^2w=u(u^2+a_2uw+a_4w^2)$. Note that $\mathcal P$ is related to the classical moduli problem $[\Gamma_1(2)]$, see for example \cite{kama:moduli}. 
\end{itemize}
In (iv) and (v), any elliptic curve $E$ over $T$ has either no or infinitely many level $\mathcal P$-structures and Tate's formulas \cite[p.181]{tate:aoe} allow to explicitly determine the set of level $\mathcal P$-structures $\mathcal P(E)$ of $E$.  Here, in (iv) one can proceed  similarly as in (R1) of Section~\ref{sec:talgoconst} and in (v) we exploit that one can compute the two torsion of the group $E(T)$.  
In fact for each moduli problem used in (i)-(v), one can quickly compute the preimage of the involved \parshin{} construction $\phi:Y(T)\to M(T)$ and one can directly determine whether a given point in $Y(T)$ corresponds to a solution of the considered Diophantine problem.
Hence, if $M(T)$ is known for some $T$ then one can directly solve the equations in (i)-(v) defined by parameters satisfying the mentioned conditions with respect to $T$; for example the parameters need to be invertible in $\mathcal O_T(T)$.  On the other hand, if $M(T)$ is not already known, then our algorithms via height bounds are more efficient than first computing $M(T)$ and afterwards the preimage of $\phi$. Here the main reason is that these algorithms only need to compute the image of $\phi$ inside $M(T)$ and this image is usually much smaller than the whole space $M(T)$. We conclude by mentioning that in (i) 
we do not use modular symbols as in Algorithm~\ref{algo:sucremona} or de Weger's sieve as in Algorithm~\ref{algo:suheight}.

\subsubsection{Applications}\label{sec:malgoapplications}

In this section  we present other applications of Algorithm~\ref{algo:mheight}. We first discuss parts of our database $\mathcal D_2$ containing the solutions of large classes of Mordell equations~\eqref{eq:mordell}.  Then we use $\mathcal D_2$ to motivate a conjecture and two questions on the number of solutions of \eqref{eq:mordell}. Here we also construct a probabilistic model providing additional motivation.

We continue the notation introduced above. Let $Y_a(\mathcal O)$ be the set of solutions of \eqref{eq:mordell} and recall that $S(n)$ denotes the set of the first $n$ rational primes. 
To determine a Mordell--Weil basis of $E_a(\QQ)$ which is required in the input of Algorithm~\ref{algo:mheight}, we used the methods discussed in Section~\ref{sec:tor+mwbasis}. Further we mention that among all sets $S$ of cardinality $n$ the set $S(n)$ is usually the most difficult case to determine $Y_a(\mathcal O)$. In particular the following running times of Algorithm~\ref{algo:mheight} would be considerably better if $S(n)$ is replaced by any set $S$ of $n$ large rational primes. The reason is that the elliptic logarithm sieve becomes considerably stronger for large primes. In fact one would already obtain significant running time improvements  by removing from $S(n)$ the notoriously difficult prime $2$.

\paragraph{The case $\abs{a}\leq 10^4$.} 
We solved the Mordell equation \eqref{eq:mordell} for all pairs $(a,S)$ such that $S\subseteq S(300)$ and such that $a\in\ZZ$ is nonzero with $\abs{a}\leq 10^4$. Here the important special case  $S=\emptyset$ was already established by Gebel--Peth{\H{o}}--Zimmer~\cite{gepezi:mordell} using a different algorithm. Further we mention that for many $a\in\ZZ$ with $\abs{a}\leq 10^4$ we determined $Y_a(\mathcal O)$ for sets $S$ which are considerably  larger than $S(300)$. For example, in the ranges $\abs{a}\leq 10$ and $\abs{a}\leq 100$ we computed $Y_a(\mathcal O)$ for all $S\subseteq S(10^5)$ and all $S\subseteq S(10^3)$ respectively.  

\paragraph{Huge $a$ and $S$.} In practice the most common (nontrivial) case is when the Mordell--Weil rank of $E_a(\QQ)$ is one, and in this case our algorithm allows to deal efficiently with huge parameters $a$ and $S$. To illustrate this feature, we have randomly chosen 100 distinct rank one curves $E_a$ with $\abs{a}\geq 10^{10}$ and  for each of these curves we then determined the sets $Y_a(\mathcal O)$ for all $S\subseteq S(10^5)$. On average it took Algorithm~\ref{algo:mheight} approximately 0.15 seconds, 6 seconds and 5 hours for $S=\emptyset$, $S=S(100)$ and $S=S(10^5)$ respectively.

\paragraph{Small rank.} The efficiency of Algorithm~\ref{algo:mheight} crucially depends on the Mordell--Weil rank $r$ of $E_a(\QQ)$. We recall that Fueter's result \eqref{eq:mtorsion} completely determines the set $Y_a(\mathcal O)$ when $r=0$. Thus we assume that $r\geq 1$ in the following discussion.  In the generic case, the Mordell curve $E_a$ has small rank $r$ and then our algorithm is very fast.

(Rank $1$). As already mentioned, in this situation our algorithm can deal efficiently with huge sets $S$. In particular for each rank one curve $E_a$ with $\abs{a}\leq 10^4$ we computed the set $Y_a(\mathcal O)$ for all $S\subseteq S(10^4)$. There are 9546 such rank one curves and on average it took Algorithm~\ref{algo:mheight} approximately 20 minutes to determine $Y_a(\mathcal O)$ for $S=S(10^4)$.

(Rank $2$ and $3$). These cases also  appear quite often in practice. For example in the range $\abs{a}\leq 10^4$ there are 3426 curves $E_a$ of rank two and 478 curves $E_a$ of rank three. For these curves, we computed the set $Y_a(\mathcal O)$ for all $S\subseteq S(300)$ and on average it took less than 5 hours and 7 hours in the case of a curve of rank two and three respectively. 
 
\paragraph{Large rank.} The situation $r\geq 4$ is rather uncommon in practice. However the notoriously  difficult case of large rank $r$ is of particular interest, since it is the most challenging for the known  methods computing $Y_a(\mathcal O)$ inside the Mordell--Weil group $E_a(\QQ)$. We mention that in the present case $r\geq 4$ the following running times can be significantly improved  by parallelizing the elliptic logarithm sieve which is used in Algorithm~\ref{algo:mheight}.

(Rank $4$, $5$ and $6$). We computed $Y_a(\mathcal O)$ for $18$  rank four curves with $S=S(300)$, for $12$ rank five curves with $S= S(100)$ and for $2$ rank six curves with $S=S(50)$. On average the corresponding running time was roughly $4$ days, $2$ days and $19$ hours in the case of a curve of rank four, five and six respectively.  The running times considerably increased for larger $S$. For example on enlarging $S(100)$ to $S(150)$ and $S(50)$ to $S(75)$, the running time  was on average $6$ days and $5$ days in the case of a curve of rank five and six respectively.

(Rank $7$ and $8$). We determined the set $Y_a(\mathcal O)$ for $2$  rank seven curves  with $S=S(40)$ and for $4$ rank eight curves with $S=S(30)$. On average the corresponding running time was less than $3$ days and $5$ days in the case of a curve of rank seven and eight respectively. Here again, the running times significantly increased for larger sets $S$. For instance on enlarging $S(40)$ to $S(50)$ and $S(30)$ to $S(40)$, the running time  was on average approximately $5$ days and $14$ days in the case of a curve of rank seven and eight respectively.

(Rank at least $9$). This situation is extremely rare. However there exist Mordell curves $E_a$ with $E_a(\QQ)$ of rank at least nine. Unfortunately we could not find such a curve  for which we were able to determine a Mordell--Weil basis of $E_a(\QQ)$; here we usually could only prove that our candidate ``basis" generates a subgroup of $E_a(\QQ)$ which has full rank.

\paragraph{Conjecture and questions.} 
We next use our database $\mathcal D_2$ to motivate various questions on the cardinality of the set $Y_a(\mathcal O)$ of solutions of \eqref{eq:mordell}. First we recall Conjecture~2 which is motivated by our data and by the construction of the elliptic logarithm sieve; see also the discussion at the end of this paragraph for additional motivation.

\vspace{0.3cm}
\noindent{\bf Conjecture~2.}
\emph{There are constants $c_a$ and $c_r$, depending only on $a$ and $r$ respectively, such that any nonempty finite set of rational primes $S$ satisfies} $$\abs{Y_a(\mathcal O)}\leq c_a \abs{S}^{c_r}.$$

\noindent  
We now discuss the exponent $c_r$ in this conjecture. For any $b\in\ZZ_{\geq 1}$ we denote by $S[b]$ the smallest set of rational primes such that for any nonzero $P\in E_a(\QQ)$ with $\hat{h}(P)\leq b$ the corresponding solution $(x,y)$ of $y^2=x^3+a$ lies in $Y_a(S[b])$. The N\'eron--Tate height $\hat{h}$ defines a positive definite quadratic form on $E_a(\QQ)\otimes_\ZZ\RR\cong\RR^r$. Therefore we obtain that $b^{r/2}=O(\abs{Y_a(S[b])})$ and we deduce that $\abs{S[b]}=O(b^{r/2}\cdot b)$ since  all nonzero $P\in E_a(\QQ)$ satisfy $\abs{\tfrac{1}{2}h(x)-\hat{h}(P)}=O(1)$; here the $O$ constants depend only on $a$.  It follows that the exponent $c_r$ has to be at least $\tfrac{r}{r+2}$ and this leads us to the following question.

\vspace{0.3cm}
\noindent{\bf Question 2.1.}
\emph{What is the optimal exponent $c_r$ in Conjecture~$2$?}
\vspace{0.3cm}

\noindent  
In addition our database $\mathcal D_2$ strongly indicates that the exponent $c_r=\tfrac{r}{r+2}$ is still far from optimal for many families of sets $S$ of interest, including the family $S(n)$ with $n\in\ZZ_{\geq 1}$. More precisely, together with the bound \eqref{refquestbound}, our database $\mathcal D_2$  motivates the following question concerning the dependence on $q=\max S$.

\vspace{0.3cm}
\noindent{\bf Question 2.2.}
\emph{Are there constants $c_a$ and $c_r$, depending only on $a$ and $r$ respectively, such that any nonempty finite set of rational primes $S$ with $q=\max S$ satisfies} $$\abs{Y_a(\mathcal O)}\leq c_a(\log q)^{c_r} \, \textnormal{?}$$

\noindent  In the case $S=S(n)$ with $n\geq 2$, one can replace here $q$ by $n\log n$ without changing the question. However the above discussion of  Conjecture~2 shows that Question~2.2 has in general a negative answer when $q$ is replaced by any power of $\max(2,\abs{S})$. Further, on considering again the family $S[b]$, we see that  the exponent $c_r$ of  Question~2.2 has to be at least $r/2$. Now we ask whether Question~2.2 has a positive answer for the exponent 
\begin{equation}\label{refquest}
c_r=r/2 \, \textnormal{?}
\end{equation}
To motivate this refined question,  we may and do assume that $a\in\ZZ$. Recall that $S$ is nonempty with $q=\max S$. Mahler's result (1933) gives that $\abs{Y_a(S')}$ is bounded for all sets of rational primes $S'$ with $\max S'\leq q$. Thus we may and do assume in addition that $q$ is large.  Now we take a nonzero point $P\in E_a(\QQ)$ and we denote by $(x,y)$ the corresponding solution of $y^2=x^3+a$.  We write $x=x_1/d^2$ and $y=y_1/d^3$ with $x_1,y_1,d\in\ZZ$ satisfying $\gcd(d,x_1y_1)=1$ 
and $d>0$.  Further we define $\rho(0)=0$ and $\rho(P)=\tfrac{n(P)}{d(P)}$, where $d(P)=d$ and $n(P)$ is the number of positive integers  $n\in \mathcal O^\times$ with $n\leq d(P)$.  It holds that $d(P)\in\mathcal O^\times$ if and only if $(x,y)$ lies in  $Y_a(\mathcal O)$. In light of this we would like to interpret $\rho(P)$ as the probability of the event that  $P\in E_a(\QQ)$ corresponds to some $(x,y)\in Y_a(\mathcal O)$. More precisely, putting $\mu_P(\{1\})=\rho(P)$ defines a probability measure $\mu_P$ on the space $\Omega_P=(\{0,1\},\mathcal P)$ for $\mathcal P$ the power set of $\{0,1\}$. Consider the associated product probability space $\Omega=(\prod\Omega_P,\prod \mu_P)$ with the product taken over all nonzero points $P\in E_a(\QQ)$. 
It follows that the random variable $\abs{\tilde{Y}_a(S)}=\sum \omega_P$ on $\Omega$ has expected value 
$$\mathbb E\bigl(\abs{\tilde{Y}_a(S)}\bigl)=\sum_{P\in E_a(\QQ)} \rho(P)$$ 
where $\omega_P:\Omega\to \Omega_P$ denotes the coordinate function.
We next estimate this expected value.  For each $n\in\ZZ_{\geq 1}$ we denote by $\Psi(n,q)$ the de Bruijn function, that is the number of $q$-smooth numbers which are at most $n$. We observe that $\rho(P)\leq \Psi(d(P),q)/d(P)$ and de Bruijn (1951) gives absolute constants $c_1,c_2\in\RR_{>0}$ such that $\tfrac{1}{n}\Psi(n,q)\leq c_1 n^{-c_2/\log q}$. 
Further, for each $\varepsilon>0$  a classical Diophantine approximation result of Siegel (1929) implies  that $\hat{h}(P)\leq (1+\varepsilon)d(P)+c_3$ with a constant $c_3$ depending only on $a$ and $\varepsilon$. We also recall that $E_a(\QQ)_{\textnormal{tor}}$ has bounded cardinality and that $\hat{h}$ defines a positive definite quadratic form on $E_a(\QQ)\otimes_\ZZ \RR\cong \RR^r$. Therefore, on combining the above observations, we see that elementary analysis gives a constant $c_a$ depending only on $a$ such that
\begin{equation}\label{refquestbound}
\mathbb E\bigl(\abs{\tilde{Y}_a(S)}\bigl)\leq c_a(\log q)^{r/2}.
\end{equation}
This motivates  Question~2.2 and its refinement in \eqref{refquest}. Moreover, the above arguments allow to describe explicitly the constant $c_a$ of \eqref{refquestbound} in terms of $r$, $a$, the regulator of $E_a(\QQ)$, the cardinality of $E_a(\QQ)_{\textnormal{tor}}$ and a  constant given by an effectively computable integral involving the Dickman function.  To control here the constant $c_3$ in terms of $a$, one can use Baker's explicit abc-conjecture stated in Section~\ref{sec:suapplications}. We point out that all constructions of this paragraph do not use that $E_a$ is a Mordell curve. In fact they can be  directly applied  to motivate the corresponding conjecture and questions for any hyperbolic genus one curve over $\sp(\ZZ)-S$. We refer to Section~\ref{sec:elllogsieveapp} for details.

\subsection{Comparison of algorithms}\label{sec:malgocomparison} In this section we discuss advantages and disadvantages of Algorithms~\ref{algo:mcremona} and \ref{algo:mheight}. We also compare our approach to the actual best methods solving \eqref{eq:mordell}. 

\paragraph{Advantages and disadvantages.} Our Algorithm~\ref{algo:mcremona} via modular symbols is very fast for all parameters $S$ and~$a$ which are small enough such that the image of the \parshin{} construction $\phi$  is contained (see Sections~\ref{sec:ia} and \ref{sec:mcremalgo}) in a database listing all elliptic curves over $\QQ$ of given conductor. Unfortunately this image is usually not contained in the actual largest known database (due to Cremona) and then the computation of the required elliptic curves via modular symbols is not efficient; here the main problem is the memory.  Thus in most cases Algorithm~\ref{algo:mcremona} can presently not compete with other approaches. 

In the generic case, our Algorithm~\ref{algo:mheight} considerably improves the actual best methods resolving \eqref{eq:mordell}. In particular it is significantly faster than the known algorithms using the elliptic logarithm approach. Indeed our optimized height bounds are sharper in practice and our elliptic logarithm sieve substantially improves in all aspects the known enumerations. Furthermore,  an important feature of Algorithm~\ref{algo:mheight} is that it allows to efficiently solve \eqref{eq:mordell} for large sets $S$. This seems to be out of reach for approaches via logarithmic forms which usually reduce to Thue(--Mahler) equations or to $S$-unit equations over number fields. Here we point out that  in the important special case $S=\emptyset$ and varying $a\in\ZZ-0$ with $\abs{a}\leq A$ for some given $A\in \ZZ_{\geq 1}$, the classical Baker--Davenport approach via logarithmic forms is very efficient. As already mentioned, Bennett--Ghadermarzi~\cite{begh:mordell} refined this approach and computed the solutions of \eqref{eq:mordell} in $\ZZ\times\ZZ$ for all nonzero $a\in\ZZ$ with $\abs{a}\leq 10^7$.  This computation involves many distinct parameters $a$, which is unfavorable for our approach since finding unconditionally all the required Mordell--Weil bases would (when possible) take a long time with the known techniques. In particular this highlights the disadvantage of Algorithm~\ref{algo:mheight} which is its dependence on an explicitly given Mordell-Weil basis. On   the other hand, for $a\in\ZZ-0$ fixed one can usually determine a basis in practice and then Algorithm~\ref{algo:mheight} is very fast even when $\abs{a}$ is huge; see Section~\ref{sec:malgoapplications}.

\paragraph{Comparison of data.} Some parts of our database $\mathcal D_2$ containing the solutions of large classes of Mordell equations \eqref{eq:mordell}  were already computed by other authors using different methods; see the work of Gebel--Peth{\H{o}}--Zimmer~\cite{gepezi:bordeaux,gepezi:barcelona,gepezi:mordell} and Bennett--Ghadermarzi~\cite{begh:mordell}. On comparing the data in the overlapping cases, it turned out that our Algorithm~\ref{algo:mheight} never produced less solutions. In particular, for all parameters in the class $\{\abs{a}\leq 10^4, S=\emptyset\}$ one verifies that our data coincides with the corresponding results in the database obtained by Bennett--Ghadermarzi~\cite{begh:mordell}.

\section{Algorithms for Thue and Thue--Mahler equations}\label{sec:thuealgo}

In \cite[Sect 7.4]{rvk:modular} an effective finiteness proof (see Section~\ref{sec:thueproofs}) for arbitrary cubic Thue equations was obtained by using inter alia an explicit reduction to a specific Mordell equation. In the present section we combine the same strategy with our algorithms for Mordell equations in order to solve  cubic Thue and Thue--Mahler equations.

We continue the  notation introduced in the previous sections. In particular we denote by $S$ an arbitrary finite set of rational prime numbers and we write $\mathcal O=\ZZ[1/N_S]$ for $N_S=\prod_{p\in S} p$.  Let $f\in\mathcal O[x,y]$ be a homogeneous polynomial of degree 3 with nonzero discriminant and let $m\in\mathcal O$ be nonzero. We recall the cubic Thue equation  
\begin{equation}
f(x,y)=m, \ \ \  (x,y)\in\mathcal O\times\mathcal O. \tag{\ref{eq:thue}}
\end{equation}
In theory, the problem of solving cubic Thue equations \eqref{eq:thue} is equivalent to the problem of finding all primitive solutions of general cubic Thue--Mahler equations~\eqref{eq:thue-mahler}.  

\begin{definition}\label{def:primsoltm}
We say that $(x,y,z)$ is a primitive solution of the general cubic Thue--Mahler equation~\eqref{eq:thue-mahler} if  $x,y,z\in\ZZ$ satisfy the equation $f(x,y)=mz$ with $z\in\mathcal O^\times$ and if $\pm 1$ are the only $d\in\ZZ$ with the property that  $d\mid\gcd(x,y)$ and $d^3\mid z$.
\end{definition}
If $(x,y,z)$ is a solution of the cubic Thue--Mahler equation \eqref{eq:thue-mahler} discussed in the introduction, then $(x,y,z)$ is in particular a primitive solution in the sense of Definition~\ref{def:primsoltm}. In fact one can directly write down all solutions of the  Thue--Mahler equation \eqref{eq:thue-mahler} if one knows all primitive solutions of the general cubic Thue--Mahler equation~\eqref{eq:thue-mahler}. 

\paragraph{Known methods.}  Baker--Davenport~\cite{baker:contributions,bada:diophapp} obtained a practical approach (see e.g. Ellison et al \cite{ellison:mordell}) to solve the cubic Thue equation \eqref{eq:thue} in $\ZZ\times\ZZ$. See also the variation of Peth{\H{o}}--Schulenberg~\cite{pesc:thue} which uses in addition the $L^3$ algorithm. Moreover, Tzanakis--de Weger~\cite{tzde:thue,tzde:thuemahler} and  Bilu--Hanrot~\cite{biha:thuehighdeg,biha:thuecomposite}  constructed practical algorithms solving Thue and Thue--Mahler equations of arbitrary degree  by applying the theory of logarithmic forms~\cite{bawu:logarithmicforms}.  We further remark that the classical $p$-adic method of Skolem often allows to find all solutions of the cubic Thue equation \eqref{eq:thue} in $\ZZ\times\ZZ$. In fact several authors used this method to practically resolve specific Thue equations. See for instance Stroeker--Tzanakis~\cite{sttz:skolem} and the references therein. There is also a recent algorithm for \eqref{eq:thue-mahler} due to Kim~\cite{kimd:modularthuemahler}, which we shall  discuss in Section~\ref{sec:talgocompa}.

\subsection{Preliminary constructions}\label{sec:talgoconst}

 In this section we discuss various constructions which shall be used in our algorithms for Thue and Thue--Mahler equations. We continue the notation introduced above.

\paragraph{Invariant theory.} To reduce our given cubic Thue equation \eqref{eq:thue} to some specific Mordell equation~\eqref{eq:mordell}, we use classical invariant theory for cubic binary forms going back at least to Cayley. We write  $\Delta$ for the discriminant of $f$ and we denote by $\mathcal H$ and $J$ the covariant polynomials of $f$ of degree two and three respectively; see Section~\ref{sec:ans+covariants} for the definitions and for our normalizations. Classical invariant theory gives that $u=-4\mathcal H$ and $v=4J$ satisfy the relation 
$v^2=u^3+432\Delta f^2$ in $\mathcal O[x,y]$. This induces a morphism 
\begin{equation}\label{eq:thuemap}
\varphi: X\to Y
\end{equation}
of $\mathcal O$-schemes, where $X$ and $Y$ are the closed subschemes of $\mathbb A^2_\mathcal O$ associated to the Thue equation \eqref{eq:thue} and to the Mordell equation \eqref{eq:mordell} with  $a=432\Delta m^2$  respectively. The solution sets of \eqref{eq:thue} and \eqref{eq:mordell} identify with the sets of sections $X(\mathcal O)$ and $Y(\mathcal O)$ of the $\mathcal O$-schemes $X$ and $Y$ respectively. Further, the projective closure inside $\mathbb P^2_\QQ$ of the generic fiber of $Y$ coincides with the elliptic curve $E_a$ over $\QQ$ appearing in previous sections.

\paragraph{The preimage of $\varphi$.} The above morphism $\varphi:X\to Y$ is effective in the following sense:  
For any given $Q\in Y(\bar{\QQ})$, one can determine all $P\in X(\bar{\QQ})$ with $\varphi(P)=Q$. Indeed this follows for example directly from the explicit height inequality in Proposition~\ref{prop:heightineq}. 
Alternatively,  for any given point $Q\in Y(\QQ)$ one can efficiently determine all $P\in X(\QQ)$ with $\varphi(P)=Q$ by using triangular decomposition.  In particular if we are given all points in $Y(\mathcal O)$, then we can efficiently reconstruct the set $X(\mathcal O)$ as follows:
\begin{itemize}
\item[(R1)]
For any given $Q\in Y(\mathcal O)$ do the following: First determine the set $Z(\QQ)$ by  applying a function in Sage (Singular)  based on triangular decomposition, where $Z$ is the spectrum of $\QQ[x,y]/I$ for $I=\bigl(4\mathcal H+u,v-4J\bigl)$ with $(u,v)$ the solution of \eqref{eq:mordell} corresponding to $Q$. Then output the points of $Z(\QQ)$ which are in $X(\mathcal O)$. 
\end{itemize}
Here one can apply triangular decomposition with the affine scheme $Z$, since it has dimension zero. Indeed it turns out (see Section~\ref{sec:thueproofs}) that $\varphi$ induces  a finite morphism $\bar{X}\to E_a$ of degree 3, where $\bar{X}$ is the projective closure inside $\mathbb P^2_\QQ$  of the generic fiber of $X$.
 
\paragraph{Reduction to Thue equations.} To find all primitive solutions of the general cubic  Thue--Mahler equation \eqref{eq:thue-mahler}, it suffices to solve certain cubic Thue equations \eqref{eq:thue}.  We now consider an elementary standard reduction: For any $w\in \ZZ_{\geq 1}$ dividing $N_S^2$, we denote by $X_w$ the closed subscheme of $\mathbb A^2_\mathcal O$ given by  $f=mw$. Suppose that $(x,y,z)$ is a primitive solution of the general cubic Thue--Mahler equation~\eqref{eq:thue-mahler}. On using that the integer $z$ lies in $\mathcal O^\times$, we may and do write $z=w\epsilon^3$ with an integer $\epsilon\in \mathcal O^\times$ and $w\in\ZZ_{\geq 1}$ dividing $N_S^2$. Then $u=x/\epsilon$ and $v=y/\epsilon$ are elements in $\mathcal O$ which satisfy the Thue equation $f(u,v)=mw$. In other words $(u,v)$ lies in  $X_w(\mathcal O)$. This motivates to consider the following  reconstruction: 
\begin{itemize}
\item[(R2)] For each $w\in\ZZ_{\geq 1}$ dividing $N_S^2$ and for any point $(u,v)$ in $X_w(\mathcal O)$, define $x=lu$, $y=lv$ and $z=l^3w$ for $l\in\ZZ_{\geq 1}$ the least common multiple of the denominators of $u$ and $v$ and output the two primitive solutions  $\pm(x,y,z)$.
\end{itemize}
Here one verifies that $\pm(x,y,z)$ are indeed primitive solutions by using that the integer $w\mid N_S^2$ is cube free.
Suppose now that we are given the sets $X_w(\mathcal O)$ for all $w\in \ZZ_{\geq 1}$ dividing $N_S^2$. Then an application of (R2) produces all primitive solutions of the general cubic Thue--Mahler equation~\eqref{eq:thue-mahler}. To prove this statement, we assume that $(x,y,z)$ is such a primitive solution. Then the construction described above (R2) gives $w\in\ZZ_{\geq 1}$ dividing $N_S^2$ and $(u,v)\in X_w(\mathcal O)$. If $x',y',z'$ are the integers in (R2) associated to $w$ and $(u,v)$, then there exists $\delta\in\mathcal O^\times$ such that $(x,y,z)=(\delta x',\delta y',\delta^3z')$. We deduce that $\delta=\pm 1$, since the triples are primitive. Hence (R2) produces all primitive solutions as desired.

\subsection{Algorithms via modular symbols}\label{sec:talgocremona}
We continue the above notation. Further we denote by $\mathcal I(S,f,m)$  the  data consisting of a finite set of rational primes $S$, the coefficients of a homogeneous polynomial $f\in\mathcal O[x,y]$ of degree three with nonzero discriminant $\Delta$ and a nonzero number $m\in\mathcal O$.

\begin{Algorithm}[Thue equation via modular symbols]\label{algo:tcremona} The input is the data $\mathcal I(S,f,m)$ and  the output is the set of solutions $(x,y)$ of the Thue equation \eqref{eq:thue}. 

The algorithm: First use Algorithm~\ref{algo:mcremona} in order to compute the set $Y(\mathcal O)$ and then apply the reconstruction algorithm described in \textnormal{(R1)}. 
\end{Algorithm}

\begin{Algorithm}[Thue--Mahler equation via modular symbols]\label{algo:tmcremona} The input consists of the data $\mathcal I(S,f,m)$ and  the output is the set formed by the primitive solutions $(x,y,z)$ of the general cubic Thue--Mahler equation \eqref{eq:thue-mahler}. 

The algorithm: First use Algorithm~\ref{algo:tcremona} in order to determine the sets $X_w(\mathcal O)$ for all $w\in\ZZ_{\geq 1}$ dividing $N_S^2$ and then apply the reconstruction described in \textnormal{(R2)}. 
\end{Algorithm}

\paragraph{Correctness.} The discussions surrounding the reconstruction (R1) imply that Algorithm~\ref{algo:tcremona} finds all solutions of the cubic Thue equation \eqref{eq:thue} as desired. Furthermore, in view of the arguments given below the reconstruction (R2), we see that Algorithm~\ref{algo:tmcremona} indeed produces all primitive solutions of the general cubic Thue--Mahler equation \eqref{eq:thue-mahler}.

\paragraph{Complexity.} The set $Y(\mathcal O)$ appearing in Algorithm~\ref{algo:tcremona} contains very few elements in practice and then the reconstruction (R1) is always very efficient. In fact the bottleneck of Algorithm~\ref{algo:tcremona} is usually the application of Algorithm~\ref{algo:mcremona} whose complexity is discussed in Section~\ref{sec:mcremalgo}. We further mention that the running time of Algorithm~\ref{algo:tmcremona} is essentially determined by the computation of the sets $X_w(\mathcal O)$ for all $w\in\ZZ_{\geq 1}$ dividing $N_S^2$. 

\paragraph{Applications.} To discuss practical applications,  we define $a=432\Delta m^2$ and we let $a_S$ be as in \eqref{def:as}. In the case $a_S\leq 350 000$,   Algorithm~\ref{algo:tcremona} efficiently solves the cubic Thue equation \eqref{eq:thue} and Algorithm~\ref{algo:tmcremona} quickly finds all primitive solutions of the general cubic Thue--Mahler equation \eqref{eq:thue-mahler}.  Indeed in this case the applications of Algorithm~\ref{algo:mcremona} are very efficient, since the involved elliptic curves are given in Cremona's database (see  Section~\ref{sec:mcremalgo}).  On the other hand, if the required data of the involved elliptic curves is not already known, then our Algorithms~\ref{algo:tcremona} and \ref{algo:tmcremona} are often not practical anymore. Here the problem is Cremona's algorithm involving modular symbols, which is used in Algorithm~\ref{algo:mcremona} and which requires a huge amount of memory for large parameters.

\subsection{Algorithms via height bounds}\label{sec:talgoheight}

We continue the above notation. In view of the discussions at the beginning of Section~\ref{sec:mordellalgostat}, we included the required Mordell--Weil bases in the input of the following algorithms. We refer to Section~\ref{sec:tor+mwbasis} for methods  computing such a basis in practice.

\begin{Algorithm}[Thue equation via height bounds]\label{algo:theight} The input is the data $\mathcal I(S,f,m)$  together with the coordinates of a Mordell--Weil basis of $E_a(\QQ)$ for $a=432\Delta m^2$. The output is the set of solutions $(x,y)$ of the cubic Thue equation \eqref{eq:thue}. 

The algorithm: First use Algorithm~\ref{algo:mheight} in order to compute the set $Y(\mathcal O)$ and then apply the reconstruction algorithm described in \textnormal{(R1)}.
\end{Algorithm}

\begin{Algorithm}[Thue--Mahler equation via height bounds]\label{algo:tmheight} The input consists of the data $\mathcal I(S,f,m)$ together with the coordinates of a Mordell--Weil basis of $E_a(\QQ)$ for all parameters $a=432\Delta (m w)^2$ with $w\in\ZZ_{\geq 1}$ dividing $N_S^2$.  The output is the set of primitive solutions $(x,y,z)$ of the general cubic Thue--Mahler equation \eqref{eq:thue-mahler}. 

The algorithm: First use Algorithm~\ref{algo:theight} in order to determine the sets $X_w(\mathcal O)$ for all $w\in\ZZ_{\geq 1}$ dividing $N_S^2$ and then apply the reconstruction described in \textnormal{(R2)}.
\end{Algorithm}

\paragraph{Correctness.} On using the arguments appearing in the correctness proof of Algorithms~\ref{algo:tcremona} and \ref{algo:tmcremona}, we see that Algorithms~\ref{algo:theight} and \ref{algo:tmheight} work correctly. 

\paragraph{Complexity.} We first discuss aspects influencing the running time of Algorithm~\ref{algo:theight} in practice. In this algorithm the reconstruction \textnormal{(R1)} is always very fast, while the running time of the computation of $Y(\mathcal O)$ is determined by  the efficiency of the application of  Algorithm~\ref{algo:mheight} with $a=432\Delta m^2$. Here  the efficiency crucially depends on $\abs{S}$ and on the size of the Mordell--Weil rank of $E_a(\QQ)$, see the complexity discussions in Section~\ref{sec:mordellalgostat}. The computation of  $Y(\mathcal O)$ is usually the bottleneck of Algorithm~\ref{algo:theight}. We next discuss Algorithm~\ref{algo:tmheight}. The running time of this algorithm is essentially determined by the computation of the sets $X_w(\mathcal O)$ for all $w\in\ZZ_{\geq 1}$ dividing $N_S^2$. For this computation we need to apply Algorithm~\ref{algo:theight} with $3^{\abs{S}}$ distinct inputs $\mathcal I(S,f,m')$, where $m'$ is of the form $m'=mw$ with $w\in\ZZ_{\geq 1}$ dividing $N_S^2$. In particular the running time of Algorithm~\ref{algo:tmheight} crucially depends on $\abs{S}$ and on the aspects influencing the complexity of Algorithm~\ref{algo:theight} as discussed above.

\paragraph{Input obstruction.}

The inputs of the above algorithms require a Mordell--Weil basis of $E_a(\QQ)$ for certain parameters $a$. In the case of Algorithm~\ref{algo:theight}, one needs to determine such a basis for only one parameter $a$ and this is usually possible in practice (see Section~\ref{sec:ttmapp}) by using the known techniques implemented in (PSM).  On the other hand, the input of Algorithm~\ref{algo:tmheight} requires a Mordell--Weil basis of $E_a(\QQ)$ for $3^{\abs{S}}$ distinct parameters $a$. Here, for large $\abs{S}$, it often happens in practice that one can not determine unconditionally all required bases in an efficient way and then our Algorithm~\ref{algo:tmheight} can not be applied to find all primitive solutions of the general cubic Thue--Mahler equation~\eqref{eq:thue-mahler}. However for small $\abs{S}$ it turned out in practice that the  known techniques  are usually efficient enough to  determine unconditionally the required bases, see Section~\ref{sec:ttmapp}.

\subsubsection{Cubic forms of given discriminant}\label{sec:redcubicforms}
There are infinitely many cubic Thue and Thue--Mahler equations of some given nonzero discriminant. However, to solve all these equations, it essentially suffices to consider the equations up to the equivalence relation induced by the action of $\textnormal{GL}_2(\ZZ)$. In this section we discuss certain aspects of this equivalence relation and we explain how to efficiently determine an explicit equation in each equivalence class.  We continue our notation.

\paragraph{Equivalence classes.} We say that a polynomial in $\QQ[x,y]$ is a cubic form if it is homogeneous of degree three with nonzero discriminant.  The group $G=\textnormal{GL}_2(\ZZ)$ acts on the set of cubic forms in the usual way. If $f,f'\in\mathcal O[x,y]$ are cubic forms with $f'=g\cdot f$ for some $g\in G$, then their discriminants  coincide and there is an explicit isomorphism between $X(\mathcal O)$ and $X'(\mathcal O)$ induced by $g$;
here $X$ and $X'$ are the closed subschemes of $\mathbb A^2_\mathcal O$ given by  $f-m$ and $f'-m$ respectively. To determine the set of solutions $X(\mathcal O)$ of the cubic Thue equation~\eqref{eq:thue}, it now suffices to know $g$ together with the set $X'(\mathcal O)$. Similarly if one is given $g\in G$ together with the set of primitive solutions of the general cubic Thue--Mahler equation~\eqref{eq:thue-mahler} defined by $(f',S,m)$ with $f'=g\cdot f$, then one can directly write down all primitive solutions of the general cubic Thue--Mahler equation~\eqref{eq:thue-mahler} defined by $(f,S,m)$.

\paragraph{Reduced cubic forms.} The reduction theory of binary forms over $\ZZ$ is well-developed. See for example the recent book of Evertse--Gy{\H{o}}ry~\cite[Sect 13.1]{evgy:bookdiscreq}. Let  $f\in \ZZ[x,y]$ be a cubic form. We next discuss how to obtain a cubic form in $G\cdot f$ which is reduced in some sense. The notion of a  reduced cubic form varies a lot in the literature and therefore we now explain in detail the notion which we shall use in this paper.

We first consider the case when $f$ is irreducible in $\QQ[x,y]$. In this case, Belabas showed in \cite[Cor 3.3 and Lem 4.3]{belabas:cubicfields} that the orbit $G\cdot f$  contains a unique cubic form $f'\in\ZZ[x,y]$ which is reduced in the sense of \cite[Def 3.2 and 4.1]{belabas:cubicfields}; this notion of a reduced form is inspired by the work of Hermite (1848/1859) if the roots of $f(x,1)$ are all real  and of Mathews (1912) otherwise. Furthermore the arguments given in Belabas~\cite[Sect 3 and 4]{belabas:cubicfields} can be transformed into a simple algorithm which allows to efficiently determine the reduced form $f'$ together with $g\in G$ satisfying $f'=g\cdot f$.

Suppose now that $f$ is reducible in $\QQ[x,y]$. In this case we work with a notion of a reduced form which is very simple and which is convenient in the sense that one can trivially determine such a form in each equivalence class.    
More precisely, the orbit  $G\cdot f$ contains a cubic form $\sum a_i x^{3-i}y^i$ in $\ZZ[x,y]$  which is reduced in our following sense: 
\begin{equation}\label{def:reducedcubicforms}
0=a_3\leq a_1\leq a_2.
\end{equation}
To prove this statement we first observe that we may assume that $f$ is primitive, that is the greatest common divisor of the coefficients $a_i$ of our cubic form $f=\sum a_i x^{3-i}y^i\in\ZZ[x,y]$ is one. Hence we assume that $f$ is primitive. We next show that one can obtain that $a_3=0$. Suppose that $a_3$ is nonzero. After possibly exchanging $x$ and $y$, we can assure that $a_0$ is nonzero. Then $f(x,1)$ is reducible in $\QQ[x]$, and  thus it is reducible in $\ZZ[x]$ by Gauss' lemma and by our assumption that $f$ is primitive. Hence, on exploiting again that $f$ is primitive, we see that the extended Euclidean algorithm provides a transformation in $G$ which makes $a_3=0$ as desired. Furthermore, after possibly replacing $x$ by $-x$, we can assure that $a_2\geq 0$. It follows that $a_2\geq 1$,  since the discriminant of $f$ is nonzero and since $a_3=0$. Then on using that $a_2\in\ZZ_{\geq 1}$, we find $\alpha\in\ZZ$ depending on $a_1,a_2$ such that  $(x,y)\mapsto (x,y+\alpha x)$  leads to $-a_2<a_1\leq a_2$. Finally,  after possibly replacing $y$ by $-y$ we can assure that $a_1\geq 0$. We conclude that the orbit $G\cdot f$  indeed contains a cubic form in $\ZZ[x,y]$ which is reduced in the sense of \eqref{def:reducedcubicforms}. Here the reduced form may not be unique in its $G$-orbit, which is no disadvantage for our purpose of solving equations. 

\paragraph{Given discriminant.} For any given nonzero $\Delta\in\ZZ$, we now explain how we determine all reduced cubic forms in $\ZZ[x,y]$ of discriminant $\Delta$. First we apply the results \cite[Lem 3.5 and 4.4]{belabas:cubicfields} of Belabas in order to list all desired forms which are irreducible in $\QQ[x,y]$. Then to find the remaining cubic forms we proceed as follows: If $\sum a_ix^{3-i}y^i$  is a reduced cubic form in $\ZZ[x,y]$ which is reducible in $\QQ[x,y]$, then the property $a_3=0$ assures that $a_2^2\mid \Delta$. Hence we can directly write down all possible values for $a_2$, which together with $0\leq a_1\leq a_2$ allows to list all possible values for $a_1$. Finally we find all possible values for $a_0$ by using an explicit expression for  $a_0$ in terms of $\Delta,a_1,a_2$. Here the explicit expression for $a_0$ can be obtained by inserting $a_3=0$ in the discriminant equation.

\subsubsection{Applications}\label{sec:ttmapp}

In this section we discuss applications of Algorithms~\ref{algo:theight} and \ref{algo:tmheight}.  After explaining the database $\mathcal D_3$  containing the solutions of large classes of cubic Thue equations~\eqref{eq:thue}, we motivate new conjectures and questions on the number of solutions of \eqref{eq:thue}.  Then we discuss the database $\mathcal D_4$ listing the primitive solutions of many general cubic Thue--Mahler equations~\eqref{eq:thue-mahler} and we consider generalized superelliptic equations studied by Darmon--Granville~\cite{dagr:superell} and Bennett--Dahmen~\cite{beda:kleinsuperell}. We continue the above notation.

\paragraph{Preliminaries.}  In our databases $\mathcal D_3$ and $\mathcal D_4$ we use the set $\mathcal F_\Delta$ of reduced cubic forms in $\ZZ[x,y]$ of given nonzero discriminant $\Delta$.  This is sufficient to cover the general case  of an arbitrary cubic form $f\in\ZZ[x,y]$ of discriminant $\Delta$.  Indeed the arguments of the previous section allow to quickly find a reduced cubic form $f'\in\ZZ[x,y]$ and $g\in \textnormal{GL}_2(\ZZ)$ with $f'=g\cdot f$, and then one can directly write down the solutions with respect to $f$ using the solutions in $\mathcal D_3$ and $\mathcal D_4$. On applying the techniques described in Section~\ref{sec:redcubicforms}, we computed in 3 seconds the sets $\mathcal F_\Delta$ for all nonzero $\Delta\in\ZZ$ with $\abs{\Delta}\leq 10^4$. The database $\mathcal F$ containing the 17044 reduced forms is uploaded on our homepage: There are $2683$ distinct $\textnormal{GL}_2(\ZZ)$-orbits of cubic forms in $\ZZ[x,y]$ which are irreducible in $\QQ[x,y]$ and $\mathcal F$ lists the reduced form of each such orbit. In addition $\mathcal F$ contains the 14361 reduced cubic forms in $\ZZ[x,y]$ which are reducible in $\QQ[x,y]$. Further, we determined the Mordell--Weil bases required in the inputs of Algorithms~\ref{algo:theight} and \ref{algo:tmheight} by using the known techniques implemented in (PSM) without introducing new ideas; see also Section~\ref{sec:tor+mwbasis}.

\paragraph{Thue equation.} For any $n\in\ZZ_{\geq 1}$ we recall that $S(n)$ denotes the set of the first $n$ rational primes. Our database $\mathcal D_3$ contains in particular the solutions of the cubic Thue equation \eqref{eq:thue} for all parameter triples $(f,S,m)$ such that $f\in\mathcal F_\Delta$ with $1\leq \abs{\Delta}\leq d$, $S\subseteq S(n)$ and  $m\in\ZZ-0$ with $\abs{m}\leq \mu$, where $(d,n,\mu)$ is as in the following discussion.

In the case  $(d,n,\mu)=(10^4,300,1)$, 
we could quickly compute almost all of the required Mordell--Weil bases: If the rank was not one then this took a few seconds (in rare cases a few minutes), and also for most rank one curves we could instantly  determine a generator. However there were a few rank one curves of large regulator for which it took several hours to compute a generator by using methods in (PSM) (2, 4 and 8-descent, Heegner points).  On average it then took Algorithm~\ref{algo:theight} approximately 5 seconds and 5 minutes in order to solve \eqref{eq:thue} for $S$ the empty set and $S=S(100)$ respectively. 

In certain situations we can make $S$ huge. For example in the case  $(d,n,\mu)=(100,10^3,1)$, 
we could compute the required bases in less than 1 minute and on average  it then took approximately 1.3 hours and 12 hours in order to solve \eqref{eq:thue} for  $S=S(500)$ and $S=S(10^3)$ respectively.  
Furthermore, in the case  $(d,n,\mu)=(20,10^5,1)$, 
we instantly found the required Mordell--Weil bases and on average we then solved  \eqref{eq:thue} for  $S=S(10^4)$ and $S=S(10^5)$ in less than 20 minutes  and 5 hours respectively.  
 Finally for the classical form $f=x^3+y^3$ we solved \eqref{eq:thue} for all $(S,m)$ as follows\footnote{For varying $m\in \ZZ$ with $\abs{m}$ bounded, it suffices to consider the case $m\geq 1$. Indeed the polynomial $f$ is homogeneous of odd degree and thus the equation $f(x,y)=m$ is equivalent to $f(-x,-y)=-m$.}. In the case $(S(10^5),m)$ with  $m\in\ZZ_{\geq 1}$ satisfying $m\leq 15$, we computed the required bases in less than 1 second and on average it then took roughly 2 hours to solve \eqref{eq:thue}.  Further in the case $(S(10^3),m)$ with $m\in\ZZ_{\geq 1}$ satisfying $m\leq 100$, we computed the required bases in less than 1 second and on average it then took approximately 6 hours to solve \eqref{eq:thue}.

\begin{remark}[Dependence on rank]
We created our database $\mathcal D_3$ with $\Delta$, $S$ and $m$ in a given range. These parameters directly influence the efficiency of the known methods solving  \eqref{eq:thue}. However for our approach the crucial parameter is the involved Mordell--Weil rank $r$, which does not depend on $S$ and which is usually small even for huge $\Delta$, $m$. Hence the discussions in Section~\ref{sec:malgoapplications}, containing running times for any given $r$, might be more meaningful than the above running times.  Finally we mention that in the generic situation where $r\leq 2$, our Algorithm~\ref{algo:theight} is very fast even for huge parameters $\Delta$, $S$, $m$. 
\end{remark}

\paragraph{Conjectures and questions.} The morphism $\varphi:X\to Y$ in \eqref{eq:thuemap} induces a finite morphism $\bar{X}\to \bar{Y}$ of degree 3, where $\bar{X}$ and $\bar{Y}$ are the projective closures inside $\mathbb P^2_\QQ$ of the generic fibers of $X$ and $Y$ respectively.  It follows that $\abs{X(\mathcal O)}\leq 3\abs{Y(\mathcal O)}$.
 Hence on applying our conjectures and questions in Section~\ref{sec:malgoapplications} with $Y=Y_a$ for $a=432\Delta m^2$, we directly obtain the analogous conjectures and questions on upper bounds for the number of solutions of the cubic Thue equation~\eqref{eq:thue} in terms of $S$ and the Mordell--Weil rank $r$ of $\textnormal{Pic}^0(\bar{X})(\QQ)$.  Our database $\mathcal D_3$ motivates these analogous conjectures and questions for cubic Thue equations~\eqref{eq:thue}. Furthermore, it might be possible to obtain more precise conjectures for \eqref{eq:thue}  by analyzing in addition the fibers of $\varphi$. We leave this for the future.

\paragraph{Thue--Mahler equation.} We next discuss our database $\mathcal D_4$.  In what follows, by  solving \eqref{eq:thue-mahler} for $(f,S)$ we mean finding all primitive solutions of the general cubic Thue--Mahler equation~\eqref{eq:thue-mahler} defined by $f$, $S$ and $m=1$;  note that  any such primitive solution $(x,y,z)$ satisfies $\gcd(x,y)=1$ provided that $f\in\ZZ[x,y]$.  We solved \eqref{eq:thue-mahler} for all $(f,S)$  such that  $f\in\mathcal F_\Delta$ with $1\leq\abs{\Delta}\leq d$ and $S\subseteq S(n)$, where $(d,n)$ is of the form $(3000,2)$, $(1000,3)$, $(100,4)$ or $(16,5)$. 
Here again we could quickly determine almost all of the required Mordell--Weil bases. However for increasing $\abs{S}$ and $\abs{\Delta}$ there were more and more rank one curves of large regulator, and finding a generator for these rank one curves was often the bottleneck of our approach. Given the input, Algorithm~\ref{algo:tmheight} was fast in all cases. 

To give the reader an idea of our running times, we now discuss some equations appearing in the literature.  We solved the  equation of Tzanakis--de Weger~\cite{tzde:thuemahlercomp}, and we determined all solutions of the equation of Agraval--Coates--Hunt--van der Poorten~\cite{agcohuva:conductor11}. 
Here our total running times were 3 minutes and 3 seconds, which includes  the 2.5 minutes and 1.5 seconds that were required to compute the involved bases. In addition we solved the equation of Tzanakis--de Weger~\cite{tzde:thuemahler} which they used to illustrate the practicality of their method. Here we determined the required  bases in less than a day and then it took Algorithm~\ref{algo:tmheight} approximately 1 minute to solve the Thue--Mahler equation.

Further, we also solved \eqref{eq:thue-mahler} for all $(f,S(6))$ with $f\in\mathcal F_\Delta$ and $\Delta\in\{-9,-1,3,27\}$. Here the case $\Delta=27$ covers in particular  $f=x^3+y^3$, and the case $\Delta=-1$ corresponds to the $S$-unit equation \eqref{eq:sunit} which means that our Theorem~B solves in particular \eqref{eq:thue-mahler} for all $(f,S)$ with $f\in\mathcal F_{-1}$ and with $S$ satisfying $S\subseteq S(16)$ or $N_S\leq 10^7$.

To conclude we mention that our Algorithm~\ref{algo:tmheight} can be used to study properties of certain generalized superelliptic equations. More precisely, let $f\in\ZZ[x,y]$ be a cubic form with nonzero discriminant $\Delta$ and take $l\in\ZZ$ with $l\geq 4$.   Darmon--Granville~\cite{dagr:superell} deduced from the Mordell conjecture \cite{faltings:finiteness}  that the generalized superelliptic equation
\begin{equation}\label{eq:superelliptic}
f(x,y)=z^l, \ \ \ (x,y,z)\in\ZZ^3
\end{equation}
with $\gcd(x,y)=1$ has at most finitely many solutions.  Moreover on using inter alia  modularity of certain Galois representations, level lowering, classical invariant theory and properties of elliptic curves with isomorphic mod-$n$ Galois representations, Bennett--Dahmen~\cite[Thm 1.1]{beda:kleinsuperell} proved: The equation $f(x,y)=z^l$ has only finitely many solutions  $(x,y,z,l)\in\ZZ^4$ with $l\geq 4$ and $\gcd(x,y)=1$ if the following condition $(*)$ holds. 
\begin{itemize}
\item[$(*)$] The polynomial $f$ is irreducible 
and  there are no solutions of the Thue--Mahler equation \eqref{eq:thue-mahler} defined by $f$, $S=\{p\st p\mid 2\Delta\}$ and $m=1$.
\end{itemize}
Bennett--Dahmen  explicitly constructed in \cite[Thm 1.2]{beda:kleinsuperell} an infinite family of polynomials satisfying condition $(*)$ and they explained in \cite[Sect 12]{beda:kleinsuperell} a heuristic indicating that ``almost all" cubic forms should satisfy $(*)$. Now, for any given cubic form $f\in\ZZ[x,y]$ with nonzero discriminant $\Delta$, our Algorithm~\ref{algo:tmheight} allows to verify in practice whether condition $(*)$ holds. In other words, one can check condition $(*)$  without using algorithms which ultimately rely on the theory of logarithmic forms. For example, we used Algorithm~\ref{algo:tmheight} to verify that $3x^3+2x^2y+5xy^2+3y^3$ satisfies condition $(*)$; note that according to \cite[p.174]{beda:kleinsuperell} this is the cubic form of minimal $\abs{\Delta}$ which satisfies condition $(*)$.

\subsection{Comparison of algorithms}\label{sec:talgocompa}
In this section we compare our algorithms for cubic Thue equations~\eqref{eq:thue} and cubic Thue--Mahler equations~\eqref{eq:thue-mahler} with the actual best practical methods in the literature. 

\paragraph{Advantages and disadvantages.} We begin by discussing Algorithms~\ref{algo:tcremona} and \ref{algo:tmcremona} for \eqref{eq:thue} and \eqref{eq:thue-mahler} using modular symbols (Cremona's algorithm).  In the recent work \cite{kimd:modularthuemahler}, Kim  independently constructed an algorithm for cubic Thue--Mahler equations \eqref{eq:thue-mahler} using modular symbols and the Shimura--Taniyama conjecture. Kim's method differs from our strategy in the sense that he is not using the route via Thue and Mordell equations, but directly associates to each solution  of \eqref{eq:thue-mahler} a certain elliptic curve. It turns out that his method is more efficient in terms of $S$ and our strategy is more efficient in terms of $\Delta$. In fact both approaches  are very fast for all parameters such that the involved elliptic curves are already known. However, usually these curves are not already known and computing these curves via modular symbols is currently not efficient for large parameters; here the main problem is the memory. Thus in most cases the algorithms via modular symbols can presently not compete with approaches solving \eqref{eq:thue} and \eqref{eq:thue-mahler} via height bounds.

We next compare our Algorithm~\ref{algo:theight} with the actual best  methods  in the literature solving cubic Thue equations~\eqref{eq:thue} using height bounds.   Our algorithm requires a Mordell--Weil basis in its input. In practice this basis can usually be computed and then our approach is very efficient. In the important special case when $S$ is empty, the already mentioned method of Tzanakis--de Weger~\cite{tzde:thue} works very well in practice and it usually allows to efficiently solve \eqref{eq:thue}.  Their method has the advantage of not requiring a Mordell--Weil basis in the input. On the other hand, an advantage of Algorithm~\ref{algo:theight} is that it efficiently deals with large sets $S$. For example it seems that already sets $S$ with $\abs{S}\geq 10$ are out of reach for the known methods solving \eqref{eq:thue}, while in the generic case Algorithm~\ref{algo:theight} allows to efficiently solve \eqref{eq:thue} for essentially all sets $S$  with $\abs{S}\leq 10^3$.

It remains to discuss our Algorithm~\ref{algo:tmheight} for cubic Thue--Mahler equations \eqref{eq:thue-mahler}. Its input requires $3^{\abs{S}}$ distinct Mordell--Weil bases  and thus our approach is not practical when $\abs{S}$ is large. However for small sets $S$ it turned out in practice that one can usually determine the required bases  and then our approach is efficient as illustrated in Section~\ref{sec:ttmapp}. If $f\in\ZZ[x,y]$ is irreducible, then the above mentioned method (TW) of Tzanakis--de Weger~\cite{tzde:thuemahler} solves in particular any cubic Thue--Mahler equation \eqref{eq:thue-mahler}. We are not aware of a complexity analysis of (TW) and thus we restrict ourselves to the following comments. There are several results in the literature which resolved specific equations \eqref{eq:thue-mahler} using (TW).  As far as we know, these results all involve small sets $S$ with $\abs{S}\leq 4$ and  (TW) is quite practical for such small sets.  On the other hand, sets $S$ of large cardinality are also problematic for (TW) since this method needs to enumerate points in lattices of rank at least $\abs{S}$.

\paragraph{Comparison of data.} We are not aware of any database in the literature which contains the solutions of large classes of cubic Thue equations \eqref{eq:thue} or cubic Thue--Mahler equations \eqref{eq:thue-mahler}. To compare at least some data, we solved the equations of \cite{agcohuva:conductor11,tzde:thuemahlercomp,tzde:thuemahler}  and in all cases it turned out that we found the same set of solutions.

\section{Algorithms for generalized Ramanujan--Nagell equations}\label{sec:ranaalgo}
In the present section we use our approaches for Mordell equations~\eqref{eq:mordell} in order to construct algorithms for the generalized Ramanujan--Nagell equation~\eqref{eq:rana}.

We continue the  notation of the previous sections. In particular we denote by $S$ an arbitrary finite set of rational prime numbers and we let $\mathcal O^\times$ be the group of units of $\mathcal O=\ZZ[1/N_S]$ for $N_S=\prod_{p\in S} p$.  Further we suppose that $b$ and $c$ are arbitrary nonzero elements of $\mathcal O$. Now we recall the generalized Ramanujan--Nagell equation 
\begin{equation}
x^2+b=cy, \ \ \ \ \ (x,y)\in\mathcal O\times \mathcal O^\times. \tag{\ref{eq:rana}}
\end{equation}
We observe that this Diophantine problem is equivalent to the a priori more general Diophantine problem obtained by replacing in \eqref{eq:rana} the polynomial $x^2+b$ by any given polynomial $f\in\mathcal O[x]$ of degree two with nonzero discriminant. 

\paragraph{Known methods.} As mentioned in the introduction, if $b\in\ZZ$ is nonzero and $c=1$ then Peth{\H{o}}--de Weger~\cite{pede:binaryrec1}  already obtained a practical method to find all solutions $(x,y)$ of \eqref{eq:rana} with $x,y\in\ZZ_{\geq 0}$. They use inter alia the theory of logarithmic forms and binary recurrence sequences; see also de Weger~\cite{deweger:phdthesis,deweger:sumsofunits}. In addition, Kim~\cite[Sect 8]{kimd:modularthuemahler} and Bennett--Billerey~\cite[Sect 5]{bebi:sumsofunits} recently obtained other practical approaches for \eqref{eq:rana}  which are briefly discussed in Sections~\ref{sec:talgocompa} and \ref{sec:ranaapp} respectively.

\subsection{Algorithm via modular symbols}\label{sec:ranaalgocremona}

We continue the notation introduced above. For any nonzero $a\in\mathcal O$, we denote by $Y_a(\mathcal O)$ the set of solutions of the Mordell equation \eqref{eq:mordell} defined by $(a,S)$. The following algorithm is a direct application of our Algorithm~\ref{algo:mcremona} for Mordell equations~\eqref{eq:mordell}.

\begin{Algorithm}[Ramanujan--Nagell equation via modular symbols]\label{algo:ranacremona} The input consists of a finite set of rational primes $S$ together with nonzero $b,c\in\mathcal O$. The output is the set of solutions $(x,y)$ of the generalized Ramanujan--Nagell equation \eqref{eq:rana}. 

The algorithm: For each  $\epsilon\in\ZZ_{\geq 1}$ dividing $N_S^2$,  use Algorithm~\ref{algo:mcremona}  to determine $Y_a(\mathcal O)$ with $a=-b(\epsilon c)^2$  and for any $(u,v)\in Y_a(\mathcal O)$ output $(\tfrac{v}{\epsilon c},\tfrac{u^3}{\epsilon^2 c^3})$ if it satisfies \eqref{eq:rana}.
\end{Algorithm}

\paragraph{Correctness.} To show that this algorithm works correctly, we take a solution $(x,y)$ of \eqref{eq:rana}. We write $y=\epsilon y'^3$ with $y'\in\mathcal O^\times$ and  $\epsilon\in\ZZ_{\geq 1}$  dividing $N_S^2$, and we define $a=-b(c\epsilon)^2$. Further we put $u=\epsilon c y'$ and $v=\epsilon c x$. It follows that $(u,v)$ lies in $Y_a(\mathcal O)$ and thus we see that the above Algorithm~\ref{algo:ranacremona} indeed finds all solutions of \eqref{eq:rana} as desired.

\paragraph{Complexity.} The running time of Algorithm~\ref{algo:ranacremona} is essentially determined by the applications of Algorithm~\ref{algo:mcremona} whose complexity is discussed in Section~\ref{sec:mcremalgo}. 

\paragraph{Applications.} To discuss practical applications of Algorithm~\ref{algo:ranacremona}, we define $a=bc^2$ and we let $a_S$ be as in \eqref{def:as}. In the case $a_S\leq 350 000$,  our Algorithm~\ref{algo:ranacremona} allows to efficiently determine all solutions of the generalized Ramanujan--Nagell equation \eqref{eq:rana}. Indeed in this case the applications of Algorithm~\ref{algo:mcremona} are very efficient, since the involved elliptic curves can be found in Cremona's database (see Section~\ref{sec:mcremalgo}). For example, one can quickly resolve the classical Ramanujan--Nagell equation: $x^2+7=2^n$ with $x,n\in\ZZ_{\geq 1}$. 
Resolving this Diophantine equation is equivalent to the problem of finding all triangular Mersenne numbers, and any solution lies in the set  $\{(1,3),(3,4),(5,5),(11,7),(181,15)\}$.  The latter assertion was conjectured by Ramanujan (1913) and was  proven by Nagell (1948).  One obtains an alternative proof of Nagell's result by using  Algorithm~\ref{algo:ranacremona}.  To conclude the discussion we mention that our Algorithm~\ref{algo:ranacremona} is often not practical anymore if the elliptic curves induced by the solutions of \eqref{eq:rana} need to be computed via Cremona's algorithm involving modular symbols (see Algorithm~\ref{algo:mcremona}). Here the problem is that Cremona's algorithm requires a huge amount of memory for all parameters which are not small. 

\subsection{Algorithm via height bounds}\label{sec:ranaalgoheight}

We continue the above notation. In the next algorithm we apply Algorithm~\ref{algo:mheight} several times. These applications require certain Mordell--Weil bases which we included  in the input. See Section~\ref{sec:tor+mwbasis} for  methods computing such bases in practice.

\begin{Algorithm}[Ramanujan--Nagell equation via height bounds]\label{algo:ranaheight} The inputs are  nonzero $b,c\in\mathcal O$, a finite set of rational primes $S$ and the coordinates of a Mordell--Weil basis of $E_a(\QQ)$ for all parameters $a=-b(\epsilon c)^2$ with $\epsilon\in\ZZ_{\geq 1}$ dividing $N_S^2$.  The output is the set of solutions $(x,y)$ of the generalized Ramanujan--Nagell equation \eqref{eq:rana}. 

The algorithm: For each  $\epsilon\in\ZZ_{\geq 1}$ dividing $N_S^2$,  use Algorithm~\ref{algo:mheight}  to determine $Y_a(\mathcal O)$ with $a=-b(\epsilon c)^2$  and for any $(u,v)\in Y_a(\mathcal O)$ output $(\tfrac{v}{\epsilon c},\tfrac{u^3}{\epsilon^2 c^3})$ if it satisfies \eqref{eq:rana}.
\end{Algorithm}

\paragraph{Correctness.} The arguments given in the correctness proof of Algorithm~\ref{algo:ranacremona} show  that the above Algorithm~\ref{algo:ranaheight} indeed finds all solutions of \eqref{eq:rana} as desired.

\paragraph{Complexity.} To compute in Algorithm~\ref{algo:ranaheight} the sets $Y_a(\mathcal O)$, we need to apply Algorithm~\ref{algo:mheight} with $3^{\abs{S}}$ distinct parameters $a$. In particular the running time of Algorithm~\ref{algo:ranaheight} crucially depends on $\abs{S}$ and on the complexity of Algorithm~\ref{algo:mheight} discussed in Section~\ref{sec:mordellalgostat}. Here we mention that the involved Mordell--Weil ranks are usually small in practice and then our Algorithm~\ref{algo:ranaheight} is very fast even for parameters $b,c$ with huge height. 

\paragraph{Refinement.} The input of Algorithm~\ref{algo:ranaheight} requires $3^{\abs{S}}$ distinct Mordell--Weil bases.  In practice it turned out that the computation of these bases is currently the bottleneck of our approach solving \eqref{eq:rana} via height bounds.  Consider arbitrary nonzero $b,c,d$ in $\ZZ$ with $d\geq 2$.  We now work out a refinement of Algorithm~\ref{algo:ranaheight} which only requires three distinct Mordell--Weil bases to  find all solutions of the classical Diophantine problem
\begin{equation}
x^2+b=c d^n, \ \ \ \ \ (x,n)\in\ZZ\times \ZZ \tag{\ref{eq:rana2}}.
\end{equation}
This is a special case of \eqref{eq:rana} with $y=d^n$ and  $S$ given by $S_d=\{p\st p\mid d\}$. In fact many authors refer by ``generalized Ramanujan--Nagell equation"  to (special cases of) the Diophantine problem \eqref{eq:rana2}. We obtain the following algorithm for \eqref{eq:rana2}. 

\begin{Algorithm}[Refinement]\label{algo:rana2}
The input consists of nonzero $b,c,d$ in $\ZZ$ such that $d\geq 2$ together with the coordinates of a Mordell--Weil basis of $E_a(\QQ)$ for all $a=-b(\epsilon c)^2$ with $\epsilon\in\{1,d,d^2\}$. The output is the set of solutions $(x,n)$ of \eqref{eq:rana2}.

The algorithm: For each $\epsilon\in\{1,d,d^2\}$,  use Algorithm~\ref{algo:mheight} to find $Y_a(\mathcal O)$ with $(a,S)=(-b(\epsilon c)^2,S_d)$ and for any $(u,v)\in Y_a(\mathcal O)$ output $\bigl(\tfrac{v}{\epsilon c},\log_d(\tfrac{u^3}{\epsilon^2 c^3})\bigl)$ if it satisfies \eqref{eq:rana2}.
\end{Algorithm}

Here for any $z\in\RR$ we define $\log_d(z)=(\log z)/\log d$ if $z>0$ and $\log_d(z)=-\infty$ otherwise. To prove that the above algorithm indeed finds all solutions of \eqref{eq:rana2}, we suppose that $(x,n)$ is such a solution. We write $d^n=\epsilon d^{3m}$ with $\epsilon\in\{1,d,d^2\}$ and $m\in\ZZ$. Further we define $v=\epsilon c x$ and $u= \epsilon cd^m$. Then we observe that $(u,v)$ lies in $Y_a(\mathcal O)$ for $(a,S)=(-b(\epsilon c)^2,S_d)$ and thus we see that Algorithm~\ref{algo:rana2} finds all solutions of \eqref{eq:rana2} as desired.

\subsubsection{Applications}\label{sec:ranaapp}
In this section we give some applications of Algorithms~\ref{algo:ranaheight} and \ref{algo:rana2}. In particular we discuss the database $\mathcal D_5$  containing the solutions of many generalized Ramanujan--Nagell equations~\eqref{eq:rana} and of many equations which are of the more classical form \eqref{eq:rana2}. We also explain how to apply our approach in order to study $S$-units $m,n\in\ZZ$ with $m+n$  a square or cube, and we provide some motivation for Terai's conjectures on Pythagorean numbers.

\paragraph{Preliminaries.} We continue the above notation. Further for any $n\in\ZZ_{\geq 1}$ we denote by $S(n)$ the set of the first $n$ rational primes.  In what follows in this section, the running time $(t_1,t_2)$ of our approach via (Al) is given by the time $t_1$ which was required to compute via (PSM) the Mordell--Weil bases for the input of (Al) and the time $t_2$ which was required to solve via (Al) the discussed equation. Here (Al) is either Algorithm~\ref{algo:ranaheight} or \ref{algo:rana2}.

\paragraph{Generalized Ramanujan--Nagell equation.}  Our database $\mathcal D_5$ contains in particular the solutions of the generalized Ramanujan--Nagell equation \eqref{eq:rana} for all parameter triples $(b,c,S)$ such that $b\in\ZZ$ is nonzero with $\abs{b}\leq B$, $c=1$ and $S\subseteq S(n)$, where $(B,n)$ is of the form $(6,5)$, $(35,4)$, $(250,3)$ or $(10^3,2)$.  For our algorithms via height bounds, the running time to solve the equation~\eqref{eq:rana} defined by $(b,c,S)$ is essentially the same as the one to solve the Thue--Mahler equation~\eqref{eq:thue-mahler} defined by $(f,S,m)$ with $m=c$ and $f\in\mathcal O[x,y]$ of discriminant $\Delta=b$. 
 Hence we refer to the discussions in Section~\ref{sec:ttmapp} which contain in particular our running times for many distinct Thue--Mahler equations \eqref{eq:thue-mahler}.

We next discuss a problem inspired by the original Ramanujan--Nagell equation. Recall that the original equation is $x^2+7=y$ with $(x,y)\in\ZZ\times\ZZ$ and $y=2^m$ for some $m\in\ZZ$. Now we put $b=7$, $c=1$ and $S=S(n)$ with $n\in\ZZ_{\geq 1}$ and we consider the problem of finding all solutions $(x,y)$ of \eqref{eq:rana} with $x,y\in\ZZ$.  Here the assumption $x,y\in\ZZ$ considerably simplifies the problem in practice.  Indeed we can remove all primes $p\in S$ with $-7$ not a square modulo $p$, and we know that $\ord_7(y)$ is either zero or one. Hence to solve the problem for $n=8$, it suffices to find all solutions of \eqref{eq:rana} for $(b,c,S)=(7,c,\{2,11\})$  with $c=1$ and $c=7$. These solutions were computed by Peth{\H{o}}--de Weger~\cite[Thm 5.1]{pede:binaryrec1}.   We obtained an alternative proof of their theorem by using Algorithm~\ref{algo:ranaheight}.  Indeed it took our approach less than (5 sec,\,4 sec) and (3 sec,\,4 sec) to find all solutions in the case $c=1$ and $c=7$ respectively.  In particular, we solved the case $n=8$ of the problem in less than 16 seconds.  To settle in addition the open case $n=9$, we need to find all solutions of \eqref{eq:rana} defined by $(b,c,S)=(7,c,\{2,11,23\})$  with $c=1$ and $c=7$. Here it took our approach less than (15 sec,\,15 sec) and (12 sec,\,15 sec) respectively, which means that we solved the case $n=9$ of the problem in less than 1 minute.  However in the cases $n\geq 10$ the running times $t_1$ become significantly larger, since for increasing $n$ we have to deal with more and more rank one curves of large regulator. For example to establish the case $n=11$ of the problem, it took our approach approximately (5 hours, 2 minutes) in total. Here we needed to solve \eqref{eq:rana} for $(b,c,S)=(7,c,\{2,11,23,29\})$  with $c=1$ and $c=7$.

\paragraph{Classical case.} Our database $\mathcal D_5$ contains in addition the solutions of the more classical equation \eqref{eq:rana2} for all $(b,c,d)$ of the form $(7,1,d)$ with $d\leq 888$.  To give the reader an idea of the running times of our approach via Algorithm~\ref{algo:rana2},  we determined the solutions of various equations \eqref{eq:rana2} of interest which were already solved by different methods.  For instance, we found in (1 sec,\,1 sec) and (3 sec,\,4 sec) all solutions of the equations appearing in the title of the papers of Leu--Li~\cite{leli:rana} and Stiller~\cite{stiller:rana} respectively.   The Diophantine problem \eqref{eq:rana2} was intensively studied in the literature when $d$ is prime, see for example \cite[Sect 2]{sasr:rana} and \cite[Sect 8.2]{besk:ternary} for an overview.  Here we solved in (8 sec,\,8 sec) all four exceptional equations  \cite[(2.7)]{sasr:rana} which appear in the classification of Le  initiated in \cite{le:firstrana}. Further,  it took less than (30 sec,\,20 sec) in total to solve all nine exceptional equations appearing in the classification of Bugeaud--Shorey~\cite[Thm 2]{bush:rana}. In particular, this includes the two exceptional equations $x^2+19=55^n$ and $x^2+341=377^n$ of \cite[Thm 2]{bush:rana} which we solved in (3 sec,\,4 sec)  and (19 sec,\,4 sec) respectively.

\paragraph{Sums of units being a square or cube.}  We now consider the problem of finding all integers $m,n\in\mathcal O^\times$ with $\gcd(m,n)$ square-free and $m+n$ a perfect square. On using inter alia the theory of logarithmic forms and generalized recurrences, de Weger~\cite{deweger:phdthesis,deweger:sumsofunits} obtained a practical approach for this problem which he used to settle the case $S=S(4)$. Suppose that $l\in\{2,3\}$. In a recent work, Bennett--Billerey~\cite{bebi:sumsofunits} show in particular how to practically solve the following problem  (in which $l=2$ is the original problem)
\begin{equation}\label{eq:powersumsofunits}
m+n=z^l, \ \ \ (m,n,z)\in\ZZ^3,
\end{equation}  
where  $m,n\in\mathcal O^\times$ have $l$-th power free $\gcd(m,n)$. They use a different method which combines the Shimura--Taniyama conjecture, modular symbols (Cremona's algorithm) and Frey--Hellegouarch curves. On using in addition congruence arguments and Cremona's database of elliptic curves of given conductor, they solved \eqref{eq:powersumsofunits} for $S=S(4)$  and $S=\{2,3,p\}$ with $p\leq 100$. Here for various sets $S=\{2,3,p\}$ they moreover applied the archimedean elliptic logarithm approach in the form of \cite{sttz:elllogaa,gepezi:ellintpoints}.

In the case $l=2$, we directly obtain an alternative  approach for \eqref{eq:powersumsofunits} by applying  Algorithm~\ref{algo:ranaheight} with $(-b,1,S)$ for all $b\in\ZZ_{\geq 1}$ dividing $N_S$.  Indeed for any solution $(m,n,z)$ of \eqref{eq:powersumsofunits} we may and do write $\max(m,n)=bu^2$ with $b,u$ positive integers in $\mathcal O^\times$  such that  $b\mid N_S$ and then we see that   $x=z/u$ and $y=\min(m,n)/u^2$ satisfy the generalized Ramanujan--Nagell equation \eqref{eq:rana} defined by $(-b,1,S)$.  Here $z$ and $u$ are coprime, since $\gcd(m,n)$ is square-free. 
Similarly in the case $l=3$, we directly obtain an alternative  approach for \eqref{eq:powersumsofunits} by combining Algorithm~\ref{algo:mheight} with an elementary reduction to Mordell equations \eqref{eq:mordell} with parameter $a\in\ZZ$ dividing $N_S^5$; see the proof of Corollary~\ref{cor:sumsofunits}~(ii). 
To illustrate the practicality of our approach, we solved \eqref{eq:powersumsofunits} for all $S\subseteq S(5)$ and all $S$ with $N_S\leq 10^3$. Given the required bases, this took less than 1 day in total. Here we could use our database which already contained the required bases. For $l=2$ (resp. $l=3$) we need in general $6^{\abs{S}}$ (resp. $2\cdot 6^{\abs{S}}$) distinct Mordell--Weil bases to solve \eqref{eq:powersumsofunits}, which means that our approach is not practical when $\abs{S}$ is large. On the other hand, for small $\abs{S}$ it turned out that one can usually determine the required bases in practice and then our approach is efficient. We compared our data with the known  results, obtained by  de Weger~\cite{deweger:phdthesis,deweger:sumsofunits}  ($S=S(4)$, $l=2$) and  Bennett--Billerey~\cite{bebi:sumsofunits} ($S=S(4)$ and $S=\{2,3,p\}$ with $p\leq 100$). In all cases it turned out that we found  the same solutions.

We briefly discuss advantages and disadvantages of the different approaches.  De Weger's method for $l=2$ is quite involved and we are not aware of its strengths and weaknesses. The strategy of Bennett--Billerey via modular symbols\footnote{In fact on replacing in our approach the involved algorithms via height bounds by the corresponding Algorithms~\ref{algo:ranacremona} and \ref{algo:mcremona} via modular symbols, we would directly obtain an alternative approach for \eqref{eq:powersumsofunits} via modular symbols. However we did not include this, since the arguments of Bennett--Billerey (using a careful analysis of conductors of Frey--Hellegouarch curves) are  more direct and more efficient.} has the usual advantages and disadvantages of an effective method involving  modular symbols (Cremona's algorithm), see the analogous discussions in Section~\ref{sec:malgocomparison}. A weakness of our approach is its dependence on  many Mordell--Weil bases, and a strength is its efficiency in the case when these bases can be determined. We also mention that Bennett--Billerey used additional tools (such as for example level lowering and the theory of logarithmic forms) to moreover prove explicit finiteness results for \eqref{eq:powersumsofunits} for all $l\geq 4$, see \cite[Sect 7]{bebi:sumsofunits}. Without introducing crucial new ideas,  such results are out of reach for our approach.  

\paragraph{Pythagorean numbers.} We next illustrate that Algorithm~\ref{algo:ranaheight} is a useful tool to study certain classical Diophantine problem on Pythagorean numbers which appear in the literature.  To state the first Diophantine equation we take coprime  $a,b,c\in\ZZ_{\geq 1}$ with $a$ even, and we assume that $a^2+b^2=c^2$. Inspired by the works of Je\'smanovic and Sierpi\'nski  published in 1956, Terai~\cite{terai:conjecture} conjectured that $(a,2,2)$ is the unique solution of 
\begin{equation}\label{eq:teraiconj}
x^2+b^m=c^n, \ \ \ (x,m,n)\in\ZZ^3,
\end{equation} 
with $x,m,n$ all positive.  Several authors settled special cases of \eqref{eq:teraiconj}, see for example \cite[p.21]{terai:conjaust} for an overview. We observe that equation \eqref{eq:teraiconj} is a special case of \eqref{eq:powersumsofunits} with $l=2$ and hence the above described approach via Algorithm~\ref{algo:ranaheight} allows to  solve \eqref{eq:teraiconj} for any given Pythagorean triple $(a,b,c)$. 
For example we verified Terai's conjecture for all triples $(a,b,c)$ with  $c\leq 85$.  Given the required Mordell--Weil bases, this took  less than 1 minute in total. In fact to deal with \eqref{eq:teraiconj} we used a modified version of Algorithm~\ref{algo:ranaheight}, which exploits the special shape of \eqref{eq:teraiconj} in order to reduce the number of required Mordell--Weil bases to 18. More recently, Terai~\cite{terai:conjaust} studied the following variation of \eqref{eq:teraiconj}. Let $d\geq 2$ be a rational integer and consider the Diophantine equation 
\begin{equation}\label{eq:teraiconj2}
x^2+(2d-1)^m=d^n, \ \ \ (x,m,n)\in\ZZ^3,
\end{equation} 
with $x,m,n$ all positive.  Terai conjectured in \cite[Conj 3.1]{terai:conjaust} that $(d-1,1,2)$ is the unique solution of \eqref{eq:teraiconj2}.  He verified his conjecture for certain values of $d$, including all $d\leq 30$  except the two cases  $d=12,24$ which were both settled independently by Deng~\cite{deng:teraiconj} and Bennett--Billerey~\cite[Prop 5.5]{bebi:sumsofunits}.  As above, we observe that our approach via Algorithm~\ref{algo:ranaheight} allows to solve \eqref{eq:teraiconj2} for any given $d$. To illustrate the utility of this strategy, we also verified \cite[Conj 3.1]{terai:conjaust} for all $d\leq 30$.  Here it was no problem to compute the required bases   and then it took the modified version of Algorithm~\ref{algo:ranaheight} less than 1 minute to solve \eqref{eq:teraiconj2} for all $d\leq 30$. We can also prove new cases of Terai's conjecture concerning \eqref{eq:teraiconj2}. However the running times to compute the required bases explode for larger $d$. For example, in the range $30< d< 35$  it took several days to compute the bases  and then we solved all equations \eqref{eq:teraiconj2}  in roughly 30 seconds by using the modified version of Algorithm~\ref{algo:ranaheight}. It turned out that Terai's conjecture holds in the range  $30< d< 35$.

\subsection{Comparison of algorithms}\label{sec:ranaalgocomp}
In this section we compare our algorithms for the generalized Ramanujan--Nagell equation~\eqref{eq:rana} and its more classical form \eqref{eq:rana} with the known practical methods.

\paragraph{Advantages and disadvantages.}Algorithms solving \eqref{eq:rana} via modular symbols (Cremona's algorithm) have the usual strengths and weaknesses, see the analogous discussion in Section~\ref{sec:talgocompa}.  We further mention that our Algorithm~\ref{algo:ranacremona} uses a reduction to Mordell equations~\eqref{eq:mordell}, while Kim's approach \cite[Sect 8]{kimd:modularthuemahler} for \eqref{eq:rana} works with a reduction to cubic Thue--Mahler equations~\eqref{eq:thue-mahler}. In fact Kim studies the  special case $x^2+7=y$ with $x,y\in\ZZ$ and $y\in\mathcal O^\times$. Similar as in the case of \eqref{eq:thue-mahler},  Kim's method is more efficient in terms of $S$ and our Algorithm~\ref{algo:ranacremona} is more efficient in terms of $b$. 

We next discuss approaches via height bounds. The discussion of advantages and disadvantages of Algorithm~\ref{algo:ranaheight} (resp. of Algorithm~\ref{algo:rana2}) is analogous to the corresponding discussion of Algorithm~\ref{algo:tmheight}  for cubic Thue--Mahler equations \eqref{eq:thue-mahler} (resp. of Algorithm~\ref{algo:theight} for cubic Thue equations \eqref{eq:thue}). Hence we refer to Section~\ref{sec:talgocompa}. The approach of Peth{\H{o}}--de Weger~\cite{pede:binaryrec1}, using inter alia the theory of logarithmic forms and binary recurrence sequences, is quite involved and we are not aware of its strengths and weaknesses.

\paragraph{Comparison of data.} In the cases where our results were already known (see above), we compared the data. In all cases it turned out that we obtained the same solutions.

\section{Mordell equations and almost primitive solutions}\label{sec:mordellcoates}

In this section we state and discuss Theorem~\ref{thm:m} which gives results for almost primitive solutions of Mordell equations. We also deduce Corollaries~\ref{cor:coates1} and~\ref{cor:coates2} on the difference of perfect squares and perfect cubes, and we discuss how these corollaries improve old theorems of Coates~\cite{coates:shafarevich}  which are based on the theory of logarithmic forms.

Let $S$ be a finite set of rational primes. We write $\mathcal O=\ZZ[1/N_S]$ for $N_S=\prod_{p\in S} p$ and we denote by $h$ the logarithmic Weil height. Let $a\in \mathcal O$ be nonzero and define 
\begin{equation}\label{def:asr2}
a_S=1728N_S^2r_2(a), \ \ \ \ r_2(a)=\prod p^{\min(2,\ord_p(a))}
\end{equation}
with the product taken over all rational primes $p$ not in $S$. It holds that $\log r_2(a)\leq h(a)$, and if $a\in\ZZ-\{0\}$ then we observe that  $r_2(a)\leq \min(\abs{a},\rad(a)^2)$. In view of the definition of Bombieri--Gubler \cite[12.5.2]{bogu:diophantinegeometry}, we say that $(x,y)\in\ZZ\times\ZZ$ is primitive if $\pm 1$ are the only $n\in\ZZ$ with $n^{6}$ dividing $\gcd(x^3,y^2)$. We recall the Mordell equation  
\begin{equation}
y^2=x^3+a, \ \ \  (x,y)\in\mathcal O\times\mathcal O. \tag{\ref{eq:mordell}}
\end{equation}
Inspired by Szpiro's small points philosophy (see e.g. \cite[Sect 2]{szpiro:lefschetz}) discussed below, we consider a certain class of solutions of~\eqref{eq:mordell} which contains all primitive solutions of~\eqref{eq:mordell}. For want of a
better name we call\footnote{We note that~\eqref{eq:mordell} defines naturally a moduli scheme of elliptic curves. Then one observes that the notions minimal and almost minimal solutions may be more appropriate (from a geometric point of view) than primitive and almost primitive solutions respectively.} these solutions ``almost primitive".

\begin{definition}\label{def:ap}
Let $\mu:\ZZ_{\geq 1}\to \mathbb R_{\geq 0}$ be an arbitrary function, and let $(x,y)\in\mathcal O\times\mathcal O$. We define $u_{x,y}=u=u_1/u_2$ with $u_1\in\ZZ_{\geq 1}$ minimal such that $u_1x,u_1y$ are in $\ZZ$ and with $u_2\in\ZZ$ maximal such that $u_2^6$ divides $\gcd\bigl((u_1^2x)^3,(u_1^3y)^2\bigl)$. Then $(u^2x,u^3y)\in\ZZ\times\ZZ$ is primitive and we say that $(x,y)$ is almost primitive  with respect to $\mu$ if $h(u)\leq \mu(a_S)$.
\end{definition} 
We notice that  $(x,y)\in\ZZ\times\ZZ$ is primitive if and only if it is almost primitive  with respect to all functions $\mu:\ZZ_{\geq 1}\to\mathbb R_{\geq 0}$. Intuitively, one may view almost primitive elements of $\mathcal O\times\mathcal O$ as those elements  with ``non-primitive part" bounded in terms of $a_S$.


Building on the arguments of \cite[Cor 7.4]{rvk:modular}, we obtain the following result  which depends on the quantity $a_S$ but which does not involve the height $h(a)$ of~$a$.

\begin{theorem}\label{thm:m}
The following statements hold.
\begin{itemize} 
\item[(i)] Let $\mu:\mathbb \ZZ_{\geq 1}\to\mathbb R_{\geq 0}$ be an arbitrary function.  Assume that~\eqref{eq:mordell} has a solution which is almost primitive with respect to~$\mu$. Then any solution $(x,y)$ of~\eqref{eq:mordell} satisfies $$\max\bigl(h(x),\tfrac{2}{3}h(y)\bigl)\leq \tfrac{2}{3}a_S\log a_S+\tfrac{1}{4}a_S\log\log\log a_S+\tfrac{3}{5}a_S+2\mu(a_S).$$
Moreover, if $(x,y)$ is in addition almost primitive with respect to $\mu$ then 
$$\max\bigl(h(x),\tfrac{2}{3}h(y)\bigl)\leq \tfrac{2}{9}a_S\log a_S+\tfrac{1}{12}a_S\log\log\log a_S+\tfrac{1}{4}a_S+2\mu(a_S).$$
\item[(ii)] Suppose that $a\in\ZZ$ with $\abs{a}\to\infty$ and define $a_*=1728\prod_{p\mid a} p^{\min(2,\ord_p(a))}$. Then any primitive solution $(x,y)\in\ZZ\times\ZZ$ of the Mordell equation~\eqref{eq:mordell}  satisfies 
$$\max\bigl(h(x),\tfrac{2}{3}h(y)\bigl)\leq  \tfrac{1}{6}a_*\log a_* + \frac{(\frac{2}{9}\log 2+o(1))}{\log\log a_*}a_*\log a_*.$$
\end{itemize}
\end{theorem}

We now make some remarks which complement our discussions of Theorem~\ref{thm:m} given in the introduction: Any solution $(x,y)$ of~\eqref{eq:mordell} is by definition almost primitive with respect to the constant function $\mu=h(u_{x,y})$ on $\ZZ_{\geq 1}$, and therefore we see that Theorem~\ref{thm:m}~(i) provides in particular an explicit upper bound in terms of $a_S$ and~$h(u_{x,y})$.

Theorem~\ref{thm:m} and its proof fit into Szpiro's small points philosophy  \cite[Sect 2]{szpiro:lefschetz} for hyperbolic curves $X$ of genus at least two defined over  number fields. This philosophy says, roughly speaking, that the rational points of $X$ have small\footnote{Here small means that the Arakelov height is effectively  bounded from above only in terms of the bad reduction places of~$X$, the genus of $X$ and the given base field over which $X$ is defined.}  Arakelov height  defined with the minimal regular model of~$X$. In our case of the hyperbolic genus one curve~\eqref{eq:mordell}, the existence of an almost primitive solution assures that~\eqref{eq:mordell} is sufficiently minimal, and then the height bound in Theorem~\ref{thm:m}~(i) shows that all solutions of~\eqref{eq:mordell} are small\footnote{A solution of~\eqref{eq:mordell} is small if its Weil height, given on the subvariety of $\mathbb A^2$ defined by~\eqref{eq:mordell}, is effectively bounded from above only in terms of the bad reduction places of the affine curve~\eqref{eq:mordell}. The Weil height is more suitable in our case of the affine genus one curve~\eqref{eq:mordell}, since Szpiro's result \cite[Thm 2]{szpiro:grothendieck} implies that the Arakelov height  is in fact constant on the rational points of any elliptic curve.}. Furthermore our proof of Theorem~\ref{thm:m} uses inter alia~\eqref{eq:szpiro} which is in fact a version\footnote{It follows for example from \cite[Thm 2]{szpiro:grothendieck} that~\eqref{eq:szpiro} is indeed a version of Szpiro's small points conjecture for elliptic curves over~$\QQ$.} of Szpiro's small points conjecture \cite[p.101]{szpiro:faltings} for elliptic curves over~$\QQ$. 

\subsection{The difference of perfect squares and perfect cubes}

Let $a\in\ZZ-\{0\}$  and let $S$ be a finite set of rational primes. Following Coates~\cite{coates:shafarevich}, we denote by $f\in\ZZ$ the largest divisor of $a$ which is only divisible by primes in~$S$. Then $\abs{a/f}$ is the largest divisor of $a$ which is coprime to $N_S=\prod_{p\in S}p$.  We define
$$\alpha_S=\frac{a_S}{N_S}=1728N_Sr_2(a)$$
for  $a_S$ and $r_2(a)$ the quantities given in (\ref{def:asr2}). Now we can state the following corollary.

\begin{corollary}\label{cor:coates1}
If $(x,y)\in\ZZ\times\ZZ$ with $y^2-x^3=a$ and $\gcd(x,y,N_S)=1$, then $$\log \max(\abs{x},\abs{y})\leq \tfrac{1}{2}\log \abs{a/f}+2\alpha_S\log\alpha_S+\tfrac{3}{4}\alpha_S\log\log\log \alpha_S+6\alpha_S.$$
\end{corollary}
\begin{proof}
It follows that $(x,y)$ is almost primitive with respect to $\mu=\frac{1}{6}\log\abs{a/f}$. Then Theorem~\ref{thm:m}~(i) proves the desired bound, but with $a_S$ in place of $\alpha_S$. To obtain the bound involving $\alpha_S$ claimed by Corollary~\ref{cor:coates1}, one uses a version of Theorem~\ref{thm:m}~(i) which takes into account that $\gcd(x,y,N_S)=1$. See Section~\ref{sec:mproofs} for details.
\end{proof}

We write $s=\abs{S}$ and we define $P=\max(S\cup\{2\})$. On using a completely different method, based on the theory of logarithmic forms, Coates \cite[Thm 1]{coates:shafarevich} obtained that any $(x,y)\in\ZZ\times\ZZ$ with $y^2-x^3=a$ and $\gcd(x,y,N_S)=1$ satisfies
\begin{equation}\label{eq:coates}
\log\max(\abs{x},\abs{y})\leq 2^{10^7(s+1)^4}P^{10^9(s+1)^3}\abs{a/f}^{10^6(s+1)^2}.
\end{equation}
We observe that $\alpha_S\leq 1728N_S\abs{a/f}$ and it holds that $N_S\leq P^s$. Therefore we see that Corollary~\ref{cor:coates1} improves Coates' result \cite[Thm 1]{coates:shafarevich}  stated in~\eqref{eq:coates}. 

We also obtain the following corollary on the size of the greatest rational prime divisor of the difference of (coprime) perfect squares and perfect cubes.

\begin{corollary}\label{cor:coates2} 
For any real number $\eps>0$ there is an effective constant~$c(\eps)$, depending only on $\eps$, with the following property: Suppose that $x$ and $y$ are coprime rational integers, and write $X=\max(\abs{x},\abs{y})$. Then  the greatest rational prime factor of $y^2-x^3$ exceeds $$ (1-\eps)\log\log X+c(\eps).$$
For example, if $\eps=\frac{1}{10}$ then one can take here the constant $c(\eps)=-20$.
\end{corollary}
\begin{proof}
Let $S$ be the set of rational primes dividing $a=y^2-x^3$ and write $q=\max(S)$. The explicit version of the prime number theorem given in \cite[Thm 4]{rosc:formulas} shows that $\log N_S\leq \sum_{p\leq q}\log p\leq q\bigl(1+\frac{1}{2\log q}\bigl)$. Thus Corollary~\ref{cor:coates1} implies Corollary~\ref{cor:coates2}.
\end{proof}

We conclude with several remarks. Coates obtained in \cite[Thm 2]{coates:shafarevich} the weaker lower bound $10^{-3}(\log\log X)^{1/4}$ by using his result \cite[Thm 1]{coates:shafarevich} displayed in~\eqref{eq:coates}. 

The proofs of Corollaries~\ref{cor:coates1} and~\ref{cor:coates2}  show in addition that  one can  weaken the assumptions $\gcd(x,y,N_S)=1$ and $\gcd(x,y)=1$ in Corollaries~\ref{cor:coates1} and~\ref{cor:coates2} by slightly changing the bounds. 
Further on using the link $(abc)\Rightarrow\eqref{eq:mordell}$ in \cite[p.429]{bogu:diophantinegeometry}, between the $abc$-conjecture $(abc)$ and height bounds for the solutions of~\eqref{eq:mordell}, one can show that $(abc)$ gives grosso modo our inequalities with the logarithmic Weil height $h$ replaced by~$\exp(h)$. 

Finally we point out that the known links $(abc)\Rightarrow\eqref{eq:mordell}$ are not compatible with exponential versions of~$(abc)$. In particular one can not combine the exponential version of $(abc)$ given in Stewart--Yu \cite[Thm 1]{styu:abc2} with the known links $(abc)\Rightarrow\eqref{eq:mordell}$ in order to improve our results in this Section~\ref{sec:mordellcoates} or Proposition~\ref{prop:m} below.

\section{Height bounds for Thue and Thue--Mahler equations}\label{sec:thueproofs}

To deduce Corollary~J from our results for  Mordell equations, we work out explicitly the arguments of \cite[Sect 7.4]{rvk:modular}. In fact in the present section we shall establish a more precise version of Corollary~J. This version provides optimized and sharper height bounds for the solutions of any cubic Thue and Thue--Mahler equation, see Corollary~\ref{cor:precthue}. We continue the terminology and notation introduced in Section~\ref{sec:thuealgo}.

\subsection{Reduction to Mordell equations}

We begin by recalling that $X$ is given by the Thue equation $f-m=0$, where $m\in\mathcal O$ is  nonzero and $f\in \mathcal O[x,y]$ is a homogeneous polynomial of degree 3 with nonzero discriminant $\Delta$. Further we recall that $Y$ is given by the Mordell equation $y^2-(x^3+a)=0$ with parameter  $a=432\Delta m^2$. We shall work with the following morphism $$\varphi:X\to Y$$ of $\mathcal O$-schemes, which was  constructed in \eqref{eq:thuemap} using classical invariant theory. For any polynomial $g$ with rational coefficients $a_\alpha$, we define $h(g)=\max_\alpha h(a_\alpha)$.   On recalling that $X$ and $Y$ are affine, we see that the usual logarithmic Weil height $h$ (see \cite[p.16]{bogu:diophantinegeometry}) naturally defines a height function $h$ on $X(\bar{\QQ})$ and $Y(\bar{\QQ})$. We obtain the following result. 

\begin{proposition}\label{prop:heightineq}
Assume that $f,m\in\ZZ[x,y]$. Then any $P\in X(\bar{\QQ})$ satisfies  $$h(P)\leq \tfrac{1}{3}h(\varphi(P))+12\bigl(h(f-m)+6h(f)+186\bigl).$$
\end{proposition} 
We point out that the first term $\tfrac{1}{3}h(\varphi(P))$ is optimal here.  The coefficients in the second term come from a recent version of the arithmetic Nullstellensatz discussed in the next paragraph. These coefficients can be slightly improved in certain cases. For example,  the proof of Proposition~\ref{prop:heightineq} shows in addition that one can replace the coefficient 12 by the smaller number 6 in the case when $P\in X(\ZZ)$. We also mention that one can directly remove the assumption $f,m\in\ZZ[x,y]$ in Proposition~\ref{prop:heightineq} by replacing in the bound the quantity $h(f-m)+6h(f)$ by the  larger expression $7h(f-m)+28h(f)$.

\subsubsection{Arithmetic Nullstellensatz and covariants of cubic forms}\label{sec:ans+covariants}

We continue our notation. In this section we collect some results which shall be used in the proof of Proposition~\ref{prop:heightineq}. In particular we discuss effective versions of the arithmetic Nullstellensatz and we recall classical properties of covariants of cubic forms. 

\paragraph{Arithmetic Nullstellensatz.} An important ingredient for our proof of Proposition~\ref{prop:heightineq} is an  effective version of the Nullstellensatz. Masser--W\"ustholz~\cite[Thm IV]{mawu:ellfunc} worked out a fully explicit version of the strong arithmetic Nullstellensatz. Their result would be sufficient for our purpose in the sense that it would give a version of Proposition~\ref{prop:heightineq} in which the coefficient $12$ is replaced by a number exceeding $24^{15}$. To obtain the considerably smaller number $12$, we shall apply a recent result of D'Andrea--Krick--Sombra~\cite[Thm 2]{dakrso:nullstell} providing an explicit version of the strong arithmetic Nullstellensatz over an affine variety of pure dimension. More precisely we shall work over the affine hypersurface $V\subset \mathbb A^3_\QQ$ defined by the polynomial $g=f-mz^3$. On using our assumptions that $m$ and the discriminant of $f$  are both nonzero, we see that $g$ is geometrically irreducible  
and therefore the affine variety  $V_{\bar{\QQ}}$ is  of pure dimension two. 
Write $h(V)$ for the projective logarithmic Weil height $h$ (see \cite[p.15]{bogu:diophantinegeometry}) of the point in projective space determined by the coefficient vector of the polynomial $g$, and let $\hat{h}(V)$ be the canonical height (see \cite[p.589]{dakrso:nullstell}) of the projective closure $\bar{V}$ of $V$ inside $\mathbb P^3_\QQ$. It holds
\begin{equation}\label{eq:canheightest}
\hat{h}(V)\leq h(V)+3\log 5.
\end{equation}
To verify this inequality,  we temporarily write  $R=\ZZ[x_1,\dotsc,x_4]$ and we let $D$ be  the effective Cartier/Weil divisor of $\mathbb P^3_\QQ$ given by the irreducible hypersurface $\bar{V}\subset \mathbb P^3_\QQ$. 
On using the terminology of \cite[$\mathsection$1.1]{dakrso:nullstell}, we denote by $f_{D}\in R$ a primitive polynomial determined by the Cartier divisor $D$.
The homogeneous polynomial $g$ is  irreducible in $R_\QQ=R\otimes_\ZZ\QQ$ and it holds that $\bar{V}\cong\textnormal{Proj}\bigl(R_\QQ/(g)\bigl)$. It follows that there exists $\epsilon\in\QQ^\times$ such that $\epsilon f_{D}$ coincides with $g$ in $R_\QQ$, which implies that $h(f_D)=h(V)$ since $f_D$ is primitive. 
Let $m(f_D)$ be the  Mahler measure of $f_D$ defined in \cite[$\mathsection$2.2]{dakrso:nullstell}. An application of \cite[Prop 2.39]{dakrso:nullstell} with  $D$ gives that $\hat{h}(V)=m(f_D)$, and we deduce from \cite[Lem 2.30]{dakrso:nullstell} that  $m(f_D)\leq h(f_{D})+3\log 5$.  Hence the equality $h(f_D)=h(V)$ proves \eqref{eq:canheightest}.

\paragraph{Covariants of cubic forms.} We next recall some properties of cubic forms which shall be used in the proof of Proposition~\ref{prop:heightineq}. Write $f(x,y)=ax^3+bx^2y+cxy^2+dy^3$ with $a,b,c,d\in\QQ$ and denote by $\mathcal H$ the covariant of $f$ of degree two. The form $\mathcal H$ is the Hessian of $f$ in the sense of \cite[p.175]{salmon:book}; we shall work with the  normalization
\begin{equation}\label{def:covarpol2}
\mathcal H(x,y)=Ax^2+Bxy+Cy^2, \ \ \ A=3ac-b^2, \ B=9ad-bc, \ C=3bd-c^2.
\end{equation}
To avoid a collision of notation, we write throughout this section $k=432\Delta m^2$ for the parameter $a$ appearing in the   equation \eqref{eq:mordell} which defines $Y$. Here $\Delta$ denotes the discriminant  of $f$ defined in \cite[p.175]{salmon:book}; note the sign convention used in \cite{salmon:book}.
Further we denote by $J$ the covariant of $f$ of degree three. It is a binary form
 \begin{equation}\label{def:covarpol3}
J(x,y)=\sum a_i x^{3-i}y^{i}
\end{equation}
given by the Jacobian of the forms $f$ and $\mathcal H$ in the sense of \cite{salmon:book}; for our purpose it will be convenient to normalize the coefficients $a_i$ of the polynomial $J$ as follows:
\begin{align*}
a_0&= 27a^2d-9abc+2b^3, &  a_1&=3(b^2c+9abd-6ac^2),\\  
a_2&= 3(6b^2d-bc^2-9acd), &  a_3&=9bcd-27ad^2-2c^3.
\end{align*}
To determine the set $Z(f,J)\subset\bar{\QQ}^2$ of common zeroes of $f$ and $J$, we suppose that $(x,y)\in Z(f,J)$. A formula going back at least to Cayley (see \cite[p.177]{salmon:book}) shows that the polynomials $u=-4\mathcal H$ and $v=4J$ satisfy the relation $v^2=u^3+432\Delta f^2$. This implies that $(x,y)$ is a zero of $\mathcal H$. In the case $a=0=d$, we then see that our assumption  $\Delta\neq 0$ together with \eqref{def:covarpol2}  implies that $(x,y)=0$ and therefore we obtain
\begin{equation}\label{eq:zfjcomp}
Z(f,J)=0.
\end{equation}
To prove that \eqref{eq:zfjcomp} holds in general, we now may and do assume that the coefficient $a$ of $f$ is nonzero. Then the polynomial $f(t,1)=a\prod (t-\gamma_j)$ in $\QQ[t]$ has degree 3. Furthermore the roots $\gamma_j$ of $f(t,1)$ are distinct since $\Delta\neq 0$. Thus $\mathcal H(\gamma_j,1)$ is nonzero by the formula for the Hessian given in \cite[p.175]{salmon:book}. It follows that $(x,y)=0$,  since $f,\mathcal H$ are homogeneous and since $a\neq 0$. Hence the set $Z(f,J)$ is trivial.  Alternatively, one can use here invariant theory providing that the resultant of $f$ and $H$ is a power of the invariant $\Delta$ up to a sign.

\subsubsection{Proof of Proposition~\ref{prop:heightineq}} We continue our notation. In this section we use the above results to prove Proposition~\ref{prop:heightineq}.
\begin{proof}[Proof of Proposition~\ref{prop:heightineq}] 
To obtain the desired height inequality, we work in the projective space. Let $\bar{X}$ and $\bar{Y}$ be the projective closures  in $\mathbb P^2_\QQ$ of the affine curves  $X_\QQ$ and $Y_\QQ$ respectively. It follows from \eqref{eq:zfjcomp} that $\varphi:X\to Y$ induces a finite morphism  $$\bar{\varphi}:\bar{X}\to\bar{Y}$$ of degree three, which is given by $\varphi_1=-4x_3\mathcal H(x_1,x_2)$, $\varphi_2=4J(x_1,x_2)$ and $\varphi_3=x_3^3$ in terms of coordinates $x_i$ on $\mathbb P^2_\QQ$.  To simplify the exposition, we write $R=\QQ[x_1,x_2,x_3]$ and we shall identify (when convenient) a polynomial in $R$ with its image in $\mathcal O_V(V)$. 

We next apply the Nullstellensatz to express $x_i$ in terms of $\varphi_j$. Let $Z\subset V(\bar{\QQ})$ be the set of common zeroes  of the two functions $\varphi_1,\varphi_2\in\mathcal O_V(V)$. Suppose that $(x,y,z)\in Z$. Then it holds that $z=0$ or $\mathcal H(x,y)=0$ and thus the identity $432\Delta f^2=\varphi_2^2-(\varphi_1/x_3)^3$  together with $\Delta\neq 0$ implies $(x,y)\in Z(f,J)$. We conclude that $(x,y)=0$, since $Z(f,J)$ is trivial by  \eqref{eq:zfjcomp}. It follows that the functions $x_i$ vanish on $Z$ for $i=1,2$. Furthermore, our assumptions $f,m\in\ZZ[x,y]$ together with \eqref{def:covarpol2} and \eqref{def:covarpol3} show that  $\varphi_1,\varphi_2$ have coefficients in $\ZZ$. Therefore applications of the strong arithmetic Nullstellensatz over $V$, with $\varphi_1, \varphi_2$ and $x_i$,  give  $\alpha_i,e_i\in\ZZ_{\geq 1}$ and $\rho_{ij}\in R$  such that the two functions $\alpha_i x_i^{e_i}$ and $\sum \rho_{ij}\varphi_j$ coincide on $V(\bar{\QQ})$ for $i=1,2$. In other words the difference of these two functions vanishes on $V(\bar{\QQ})$, 
and then $V=\sp\bigl(R/(g)\bigl)$ implies that this difference is divisible in $R$ by the geometrically irreducible polynomial $g$. Thus there exist $\rho_{i0}\in R$ such that 
\begin{equation}\label{eq:nssatz}
x_i^{e_i}=\sum (\rho_{ij}/\alpha_i)\varphi_j, \ \ \ i=1,2,
\end{equation}
where $\varphi_0=g$. On multiplying here both sides with $x_i^{e-e_i}$ for $e=\max e_i$, we may and do assume that  $e_i=e$. Furthermore we may and do assume that  all $\rho_{ij}$ are homogeneous of degree $e-3$, since the polynomials $\varphi_j$ are all homogeneous of degree 3. 
 
We now control the right hand side of \eqref{eq:nssatz} in terms of $\varphi$ and $g$. Put $h(\rho/\alpha)=h(Q)$ for $Q$  the point in projective space whose coordinates are given by the coefficients of the four polynomials $\rho_{ij}/\alpha_i\in R$ with $i,j\in\{1,2\}$.  Each polynomial $\rho_{ij}$ has at most $\tbinom{e-1}{2}$ nonzero coefficients, and the function $\varphi_0$ vanishes on the $\bar{\QQ}$-points of $\bar{X}=\textnormal{Proj}\bigl(R/(\varphi_0)\bigl)$. Therefore on combining \eqref{eq:nssatz} with the above observations,  we see that standard height arguments lead to the following statement: Any $\bar{\QQ}$-point $P$ of $\bar{X}$ satisfies
\begin{equation}\label{eq:projheightineq}
3h(P)\leq  h(\bar{\varphi}(P))+h(\rho/\alpha)+\log \bigl(2\tbinom{e-1}{2}\bigl).
\end{equation}  
In particular any $P\in X(\bar{\QQ})$ satisfies  \eqref{eq:projheightineq} with $\bar{\varphi}$ replaced by $\varphi$. Now  we see that the explicit version \cite[Thm 2]{dakrso:nullstell} of the strong arithmetic  Nullstellensatz  gives the following: In \eqref{eq:nssatz} one can choose  $\alpha_i$, $e_i$ and $\rho_{ij}$, with $e_i\leq 54$ and $h(\rho/\alpha)$ explicitly bounded in terms of $\hat{h}(V)$ and $h(\varphi_j)$, such that  \eqref{eq:projheightineq}, \eqref{eq:canheightest}, \eqref{def:covarpol2} and \eqref{def:covarpol3} lead to the desired height inequality. This completes the proof of Proposition~\ref{prop:heightineq}.\end{proof} 

The above proof gives moreover a version of Proposition~\ref{prop:heightineq} for projective closures inside $\mathbb P^2_\QQ$. Indeed it follows from \eqref{eq:projheightineq} that any $\bar{\QQ}$-point $P$ of $\bar{X}$ satisfies the height inequality in Proposition~\ref{prop:heightineq} with $\varphi$ replaced by the morphism $\bar{\varphi}:\bar{X}\to \bar{Y}.$

\subsection{Optimized height bounds}

We continue the notation introduced above. In this section we give Corollary~\ref{cor:precthue} which provides optimized height bounds for the solutions of cubic Thue and Thue--Mahler equations. To state our height bounds, we recall that $k=432\Delta m^2$ and we denote by $\Omega_{\textnormal{sim}}$ the simplified height height bound given in Proposition~\ref{prop:m} with $a=k$. It holds  $$\Omega_{\textnormal{opt}}\leq\Omega_{\textnormal{sim}}=\tfrac{1}{3}h(k)+\tfrac{4}{9}k_S\log k_S+\tfrac{1}{6}k_S\log\log\log k_S+\tfrac{2}{5}k_S$$   
for $k_S$ as in \eqref{def:as} and $\Omega_{\textnormal{opt}}$ the optimized height bound provided by Proposition~\ref{prop:algobounds} with $a=k$.
The next result is a sharper (but more complicated) version of Corollary~J.
\begin{corollary}\label{cor:precthue}Assume that $f,m\in\ZZ[x,y]$. Then the following statements hold.
\begin{itemize}
\item[(i)] Suppose that $(x,y)$ is a solution of the cubic Thue equation \eqref{eq:thue}, and put $n=1$ if $(x,y)\in\ZZ\times\ZZ$ and $n=2$ otherwise. Then  $h(x)$ and $h(y)$ are at most $$\tfrac{n}{2}\Omega_{\textnormal{opt}}+6n\bigl(h(f-m)+6h(f)+186\bigl).$$
\item[(ii)]Suppose that $(x,y,z)$ is a primitive solution of the general cubic Thue--Mahler equation \eqref{eq:thue-mahler}. Then  $h(x)$, $h(y)$ and $\tfrac{1}{3}h(z)$ are at most  $$2\Omega_{\textnormal{opt}}+51\log N_S+24\bigl(h(f-m)+6h(f)+186\bigl).$$
\end{itemize}
\end{corollary} 
One can directly remove the extra assumption $f,m\in\ZZ[x,y]$ by multiplying the equation $f(x,y)=m$ with  a suitable integer, see \eqref{eq:precthue} for the resulting bounds. 

\begin{proof}[Proof of Corollary~\ref{cor:precthue}]
To prove (i) we take $P\in X(\mathcal O)$. Then  $(u,v)=\varphi(P)$ satisfies the Mordell equation \eqref{eq:mordell} with parameter $a=k$.  The number $k=432\Delta m^2$ is  nonzero,  since $m\Delta\neq 0$ by assumption. Hence an application of Propositions~\ref{prop:m}~and~\ref{prop:algobounds} with $(u,v)$ gives an upper bound for $h(\varphi(P))$ which together with Proposition~\ref{prop:heightineq} implies (i). Here we used that one can replace in  Proposition~\ref{prop:heightineq} the coefficient 12 by 6 if $P\in X(\ZZ)$.

To show (ii) we assume that $(x,y,z)$ is a primitive solution of the general Thue--Mahler equation \eqref{eq:thue-mahler}. We write $z=z_0\epsilon^3$ with $z_0,\epsilon\in\ZZ$ such that $\pm 1$ are the only  $l\in\ZZ$ with $l^{3}$ dividing $z_0$. Then $u=x/\epsilon$ and $v=y/\epsilon$ satisfy the Thue equation $f(u,v)=m'$ with $m'=mz_0$.  On exploiting our assumption that $(x,y,z)$ is primitive,  one controls the absolute values of $x,y,z$ in terms of the Weil heights of $u,v$ and then (ii) follows from the height bound for $(u,v)$ obtained in (i). Here we used that $k'=432\Delta m'^2$ satisfies $k_S=k'_S$ and that $h(k')$ is at most $h(k)+4\log N_S$.  This completes the proof of Corollary~\ref{cor:precthue}.
\end{proof}

To remove the extra assumption $f,m\in\ZZ[x,y]$ in Corollary~\ref{cor:precthue}, we define $f^*=lf$ and $m^*=lm$ for $l$ the least common multiple of the denominators of the coefficients of the polynomial $f-m\in\mathcal O[x,y]$.  Any solution $(x,y)$ of the Thue equation \eqref{eq:thue} satisfies $f^*(x,y)=m^*$, and any primitive solution $(x,y,z)$ of the general Thue--Mahler equation~\eqref{eq:thue-mahler} satisfies $f^*(x,y)=m^*z$. Therefore an application of Corollary~\ref{cor:precthue} with $f^*,m^*\in\ZZ[x,y]$ implies the following: Statements (i) and (ii) of Corollary~\ref{cor:precthue} hold without the extra assumption $f,m\in\ZZ[x,y]$ if the bounds in (i) and (ii) are replaced by
\begin{equation}\label{eq:precthue}
B_1=\Omega_{\textnormal{opt}}+86\bigl(4h(f)+h(m)+26\bigl) \ \ \  \textnormal{ and }  \ \ \ 2B_1+51\log N_S
\end{equation}
respectively.  Here we used that $h(f^*),h(m^*)$ are at most $4h(f)+h(m)$ and we exploited  that $k^*_S=k_S$,  where  $k^*=432\Delta^* (m^*)^2$ for $\Delta^*$  the discriminant of $f^*$. Note that $k^*_S=k_S$ follows from $k^*=l^6k$ and from our assumptions $f,m\in\mathcal O[x,y]$ which assure that $l\in \mathcal O^\times$. Finally we deduce Corollary~J by simplifying the bounds in  Corollary~\ref{cor:precthue}.
\begin{proof}[Proof of Corollary~J]
In the case $f,m\in\ZZ[x,y]$, we see that Corollary~\ref{cor:precthue} together with $h(\Delta)\leq 4h(f)+5\log 3$ implies  (i) and (ii). In general, on following the proof of  \eqref{eq:precthue} we reduce to the case $f^*,m^*\in\ZZ[x,y]$  and then we apply the already established case of  Corollary~J with $f^*,m^*\in\ZZ[x,y]$. This completes the proof of Corollary~J.
\end{proof}

\section{Height bounds for Ramanujan--Nagell equations}\label{sec:ranaheight}

In this section we give explicit height bounds for the solutions of the generalized Ramanujan--Nagell equation. We also study pairs of units whose sum is a square or cube.

As in the previous sections let $S$ be a finite set of rational primes, write $h$ for the usual logarithmic Weil height and denote by $\mathcal O^\times$ the units of $\mathcal O=\ZZ[1/N_S]$ for $N_S=\prod_{p\in S} p$. Further let $b,c\in\mathcal O$ be nonzero and recall the generalized Ramanujan--Nagell equation
\begin{equation}
x^2+b=cy, \ \ \ \ \ (x,y)\in\mathcal O\times \mathcal O^\times. \tag{\ref{eq:rana}}
\end{equation}
To state our height bounds for the solutions of \eqref{eq:rana}, we define $a=bc^2$ and we denote by $\Omega_{\textnormal{sim}}=\Omega_{\textnormal{sim}}(a,S)$ the simplified height height bound given in Proposition~\ref{prop:m}. It holds  $$\Omega_{\textnormal{opt}}\leq\Omega_{\textnormal{sim}}=\tfrac{1}{3}h(a)+\tfrac{4}{9}a_S\log a_S+\tfrac{1}{6}a_S\log\log\log a_S+\tfrac{2}{5}a_S$$   
for $a_S$ as in \eqref{def:as} and $\Omega_{\textnormal{opt}}=\Omega_{\textnormal{opt}}(a,S)$ the optimized height bound in Proposition~\ref{prop:algobounds}.
The next result is a direct consequence of our height bounds for Mordell equations.

\begin{corollary}\label{cor:ranabounds}
If $(x,y)$ satisfies \eqref{eq:rana} then $h(x^2),h(y)\leq 3\Omega_{\textnormal{opt}}+3h(c)+8\log N_S.$
\end{corollary}
\begin{proof}
We write $y=\epsilon y'^3$ with $y'\in\mathcal O^\times$ and $\epsilon\in\ZZ_{\geq 1}$ dividing $N_S^2$. Then $u=\epsilon c y'$ and $v=\epsilon c x$ satisfy the Mordell equation $v^2=u^3+a'$ with $a'=-b(\epsilon c)^2$. It holds that $a'\neq 0$ and $a'_S=a_S$, since  $bc\neq 0$ and $\epsilon\in\mathcal O^\times$. Thus Proposition~\ref{prop:algobounds} implies Corollary~\ref{cor:ranabounds}. 
\end{proof}

It holds that $\Omega_{\textnormal{opt}}\leq \Omega_{\textnormal{sim}}$ and we observe that $3\Omega_{\textnormal{sim}}+8\log N_S$ is at most $2a_S+h(a)$.  Therefore we see that Corollary~\ref{cor:ranabounds} proves Corollary~K stated in the introduction. 

\begin{remark}[Generalization]Suppose that $f\in\mathcal O[x]$ is a polynomial of degree two, with nonzero discriminant $\Delta$. Then we claim that Corollary~\ref{cor:ranabounds} gives an explicit 
height bound for any solution $(x,y)\in\mathcal O\times\mathcal O^\times$ of the more general Diophantine equation $f(x)=cy.$
To prove this claim we suppose that $(x,y)\in\mathcal O\times\mathcal O^\times$ satisfies $f(x)=cy$.   We write $f(x)=a_1x^2+a_2x+a_3$ with $a_i\in\mathcal O$ and we put $x'=2a_1x+a_2$. It follows that $(x',y)$ is a solution of \eqref{eq:rana} with parameters $b=-\Delta$ and $c=4a_1c$. Hence an application of Corollary~\ref{cor:ranabounds} with $(x',y)$ gives an explicit upper bound for $h(x)$ and $h(y)$ as claimed. 
\end{remark}

We next consider the problem of finding all  $(u,v)\in\mathcal O^\times\times\mathcal O^\times$  with $u+v$ a square or cube in $\QQ$. To obtain here finiteness statements, one has to work modulo the actions of $\mathcal O^\times$ arising naturally in this context. The above Corollary~\ref{cor:ranabounds} leads to the following fully explicit results which involve the quantity $\Omega=3\Omega_{\textnormal{opt}}(1,S)+9\log N_S$. 

\begin{corollary}\label{cor:sumsofunits}
Assume that $u,v$ are in $\mathcal O^\times$. Then the following statements hold.
\begin{itemize}
\item[(i)] If $u+v$ is a square in $\QQ$, then there is $\epsilon\in\mathcal O^\times$ such that $h(\epsilon^2u),h(\epsilon^2v)\leq \Omega$.
\item[(ii)] If $u+v$ is a cube in $\QQ$, then there is $\delta\in\mathcal O^\times$ such that $h(\delta^3u),h(\delta^3v)\leq \Omega$.
\end{itemize}
\end{corollary}
\begin{proof} To prove (i) we assume that $u+v$ is a square in $\QQ$. Then there is $\epsilon\in\mathcal O^\times$ such that $\epsilon^2u$ and $\epsilon^2v$ are in $\ZZ$, with $\epsilon^2u+\epsilon^2v$ a perfect square and $\gcd(\epsilon^2u,\epsilon^2v)$ square-free. If $m,n$ are in $\ZZ\cap\mathcal O^\times$ with $m+n$ a perfect square  and $\gcd(m,n)$ square-free, then we claim that
\begin{equation}\label{claimsquare}
h(n)\leq \Omega. 
\end{equation}
To prove this inequality we take $l\in\ZZ$ with $l^2=m+n$ and we write $m=m'^2m_0$ with $m',m_0\in\ZZ$ such that $m_0\mid N_S$. Further we define $x=l/m'$ and $y=n/m'^2$.  Then $(x,y)$  satisfies \eqref{eq:rana} with $b=-m_0$ and $c=1$. Thus an application of Corollary~\ref{cor:ranabounds} with $(x,y)$ implies  \eqref{claimsquare} as claimed. Here we used that $n$ and $m'$ are coprime which follows from  our assumption that $\gcd(m,n)$ is square-free. Now, an application of \eqref{claimsquare} with $m=\epsilon^2u$ and $n=\epsilon^2v$ shows that $\epsilon$ has the desired properties. This  proves assertion (i).

The following proof of (ii) uses the arguments of (i) with some modifications. We assume that $u+v$ is a cube in $\QQ$. Then there is $\delta\in\mathcal O^\times$ such that $\delta^3u$ and $\delta^3v$ are  in $\ZZ$, with $\delta^3u+\delta^3v$ a perfect cube and  $\gcd(\delta^3u,\delta^3v)$ cube-free.  We now consider the following claim: If $m,n$ are in $\ZZ\cap\mathcal O^\times$ with $m+n$ a perfect cube  and $\gcd(m,n)$ cube-free, then 
\begin{equation}\label{claimcube}
h(n)\leq \Omega. 
\end{equation}
To prove this claim we take $l\in\ZZ$ with $l^3=m+n$ and we write $m=m'^3m_0$ with $m',m_0\in\ZZ$ such that $m_0\mid N_S^2$.  Further we define $x=l/m'$  and $y=n/m'^3$. Then $(x,y)$ lies in $\mathcal O\times\mathcal O^\times$ and satisfies $x^3-m_0=y$.  On writing $y=w y'^2$ with $y'\in\mathcal O^\times$ and $w\in\ZZ$ dividing $N_S$, we see that $(xw,y'w^2)$ is a solution of the Mordell equation \eqref{eq:mordell} with parameter $a=-m_0w^3$.  Here $a\in \ZZ$ is nonzero. Thus on using that $\gcd(m,n)$ is cube-free, we see that Proposition~\ref{prop:algobounds} implies our claim in \eqref{claimcube}. Finally we deduce  (ii) by applying \eqref{claimcube} with $m=\delta^3u$ and $n=\delta^3v$. This completes the proof of Corollary~\ref{cor:sumsofunits}.
\end{proof}

We now suppose that $m,n\in\ZZ$ are coprime. Then $\gcd(m,n)$ is in particular square-free and cube-free, and $m,n$ are in $\mathcal O^\times$ for $S=\{p\st p\mid mn\}$ with $N_S=\rad(mn)$. Therefore  \eqref{claimsquare} and \eqref{claimcube} imply  Corollary~L stated in the introduction.

\section{Height bounds for Mordell and $S$-unit equations}\label{sec:heightbounds}

In this section we  prove the results of Section~\ref{sec:mordellcoates}, and we establish in Proposition~\ref{prop:algobounds} the height bounds for  Mordell and $S$-unit equations which are used in the algorithms of Sections~\ref{sec:suheightalgo} and~\ref{sec:mheightalgo}. The plan of Section~\ref{sec:heightbounds} is as follows: In Section~\ref{sec:simpleversions} we give simplified versions of the height bounds. Then we prove the results for Mordell and $S$-unit equations in Sections~\ref{sec:mproofs} and~\ref{sec:suproofs} respectively. Finally, in Section~\ref{sec:hc}, we work out the height conductor inequalities for elliptic curves over $\QQ$ which are used in our proofs.

\subsection{Simplified versions}\label{sec:simpleversions}

The precise form of our height bounds in Proposition~\ref{prop:algobounds} is fairly complicated. 
To give the reader an idea of the size of the used height bounds, we worked out simplified versions of our bounds for Mordell equations (see Proposition~\ref{prop:m}) and for $S$-unit equations (see Propositions~\ref{prop:su} and~\ref{prop:abc}). These simplified versions  slightly improve several estimates in the literature and they will allow us (up to some extent) to compare our optimized height bounds with results based on the theory of logarithmic forms.

As in the previous sections we let $S$ be a finite set of rational primes, we write $\mathcal O=\ZZ[1/N_S]$ for $N_S=\prod_{p\in S} p$, and we denote by $h$ the logarithmic Weil height.

\subsubsection{Mordell equation}\label{sec:simpleboundsm}

Let $a\in\mathcal O$ be nonzero and let $a_S=1728N_S^2\prod_{p\notin S}p^{\min(2,\ord_p(a))}$ be as in (\ref{def:asr2}). On using and refining the arguments of \cite[Cor 7.4]{rvk:modular}, we obtain the following result. 

\begin{proposition}\label{prop:m}
If $x,y\in\mathcal O$ satisfy $y^2=x^3+a$, then 
$$
\max\bigl(h(x),\tfrac{2}{3}h(y)\bigl)\leq \tfrac{1}{3}h(a)+\tfrac{4}{9}a_S\log a_S+\tfrac{1}{6}a_S\log\log\log a_S+\tfrac{2}{5}a_S.
$$
\end{proposition}

To compare this result with the actual best height bounds for Mordell equations~\eqref{eq:mordell}, we suppose that $x,y\in\mathcal O$ satisfy $y^2=x^3+a$ and we observe that Proposition~\ref{prop:m} implies 
\begin{equation}\label{eq:simplemordell}
\max\bigl(h(x),h(y)\bigl)\leq \tfrac{1}{2}h(a)+a_S\log a_S.
\end{equation}
We first consider the classical case $\mathcal O=\ZZ$. The discussions in \cite[Sect 7.3]{rvk:modular} show that~\eqref{eq:simplemordell} improves the results of Baker~\cite{baker:mordellequation}, Stark~\cite{stark:mordell} and Juricevic~\cite{juricevic:mordell} which are all based on the theory of logarithmic forms~\cite{bawu:logarithmicforms}.
In fact Proposition~\ref{eq:simplemordell} provides the actual best height bound for Mordell equations~\eqref{eq:mordell}, since it updates the inequality $\max\bigl(h(x),h(y)\bigl)\leq h(a)+4\cdot 36a_S\log (36a_S)^2$ in \cite[Cor 7.4]{rvk:modular}.  

We now discuss the case $\mathcal O\supsetneq \ZZ$. On using the theory of logarithmic forms, Hajdu--Herendi \cite[Thm 2]{hahe:elliptic} obtained explicit height bounds for the solutions in $\mathcal O$ of arbitrary elliptic Diophantine equations. The discussion in \cite[Sect 7.3.3]{rvk:modular} shows that Proposition~\ref{prop:m}  improves the case of  Mordell equations~\eqref{eq:mordell} in \cite[Thm 2]{hahe:elliptic} for ``small'' sets $S$, in particular for 
any set $S$ with $N_S\leq 2^{1200}$ or 
 $\abs{S}\leq 12$. 
If $N_S\gg\lvert a\rvert$ then there are sets $S$  for which Proposition~\ref{prop:m} is better than \cite[Thm 2]{hahe:elliptic}, and vice versa. To conclude our comparisons we mention that the  theory of logarithmic forms allows to deal with more general Diophantine equations over arbitrary number fields; see~\cite{bawu:logarithmicforms}.

\subsubsection{$S$-unit equations}\label{sec:simpleboundssu}

On using and refining the arguments of \cite[Thm 1.1]{mupa:modular} and \cite[Cor 7.2]{rvk:modular},  which were discovered independently in 2011 by Murty--Pasten and by the first mentioned author, we obtain the following update of the explicit height bounds in \cite[Thm 1.1]{mupa:modular} and \cite[Cor 7.2]{rvk:modular}; see also Frey's remark in  \cite[p.544]{frey:ternary} and Proposition~\ref{prop:abc}.

\begin{proposition}
\label{prop:su}
Any solution $(x,y)$ of the $S$-unit equation~\eqref{eq:sunit} satisfies   $$\max\bigl(h(x),h(y)\bigl) \leq \tfrac{5}{2}N_S\log N_S+9N_S \ \textnormal{ and }$$
$$\max\bigl(h(x),h(y)\bigl) \leq \tfrac{12}{5}N_S\log N_S+\tfrac{9}{10}N_S\log\log\log (16N_S)+8.26N_S+28.
$$
\end{proposition}
To compare this result with the actual best height bounds for $S$-unit equations~\eqref{eq:sunit}, we suppose that $(x,y)$ satisfies~\eqref{eq:sunit} and we write $h=\max(h(x),h(y))$. Starting with Gy{\H{o}}ry~\cite{gyory:sunitshelvetica} several authors proved explicit bounds for $h$ by using the theory of logarithmic forms~\cite{bawu:logarithmicforms}; see the references in~\cite{gyyu:sunits}. In particular Gy{\H{o}}ry--Yu \cite[Thm 2]{gyyu:sunits} obtained that 
$h\leq 2^{10s+22}s^4q\prod\log p$
with the product taken over all rational primes $p\in S-\{q\}$ for $q=\max S$ and $s=\abs{S}$. Further  Proposition~\ref{prop:su} updates the inequalities $h\leq 4.8N_S\log N_S+13N_S+25$ in Murty--Pasten \cite[Thm 1.1]{mupa:modular} and $h\leq  3\cdot 2^6N_S\log(2^7N_S)^2+65$ in \cite[Cor 7.2]{rvk:modular}.
Hence Proposition~\ref{prop:su} establishes the actual best height bound for all sets $S$ with ``small"~$N_S$, in particular for all sets $S$ with  $N_S\leq 2^{90}.$ Further we see that there are sets $S$ with arbitrarily large $N_S$ for which Proposition~\ref{prop:su} is sharper than \cite[Thm 2]{gyyu:sunits}, and vice versa. If $N_S\to \infty$ then our bounds (see also~\eqref{eq:asymptoticsu}) are worse than $h\leq O(N_S^{1/3}(\log N_S)^3)$ in Stewart--Yu \cite[Thm 1]{styu:abc2}. To conclude our comparison, we mention that the theory of logarithmic forms and Bombieri's refinement of the Thue--Siegel method  both give effective height bounds for the solutions of $S$-unit equations in arbitrary number fields; see for example~\cite{gyyu:sunits,boco:effdioapp2}.

\subsection{Notation}\label{sec:notations}

In the remaining of Section~\ref{sec:heightbounds} we shall use throughout the following notation. Let $S$ be a finite set of rational primes. We write $\mathcal O=\ZZ[1/N_S]$ for $N_S=\prod_{p\in S}p$ and we denote by $h$ the usual logarithmic Weil height. For any elliptic curve $E$ over~$\QQ$, we denote by $c_4$ and $c_6$ the usual quantities (see for example~\cite{tate:aoe}) associated to a minimal Weierstrass equation of $E$ over~$\ZZ$, and we write $N_E$ and $\Delta_E$ for the conductor and minimal discriminant of $E$ respectively. We denote by $h(E)$ the relative Faltings height of~$E$, defined for example in \cite[Sect 2]{rvk:modular} using Faltings' original normalization \cite[p.354]{faltings:finiteness} of the metric.

\subsection{Mordell equation}\label{sec:mproofs}

We use the notation introduced in Section~\ref{sec:notations} above. Let $a\in \mathcal O-\{0\}$ and recall that $a_S=1728N_S^2r_2(a)$ for $r_2(a)=\prod p^{\min({2,\ord_p(a)})}$ with the product taken over all rational primes $p$ not in~$S$. Further we recall the Mordell equation 
\begin{equation}
y^2=x^3+a, \ \ \ (x,y)\in\mathcal O\times\mathcal O.\tag{\ref{eq:mordell}}
\end{equation}
The following lemma gives an explicit \parshin{} construction for the solutions of this equation.

\begin{lemma}\label{lem:pm2}
Suppose that $(x,y)$ satisfies the Mordell equation~\eqref{eq:mordell}. Then there exists an elliptic curve $E$ over $\QQ$ with the following properties. It holds that $c_4=u^4x$ and $c_6=u^6y$ for some $u\in \QQ$ with $u^{12}=1728\Delta_E\abs{a}^{-1}$, the conductor $N_E$ divides $a_S$ and  $$\max\bigl(h(x),\tfrac{2}{3}h(y)\bigl)\leq\tfrac{1}{3}h(a)+8h(E)+2\log\max(1,h(E))+36.$$
\end{lemma}
We conducted here some effort to refine the weaker result  $N_E\mid 2^23^2a_S$ which follows directly from \cite[Prop 3.4]{rvk:modular}. In fact Lemma~\ref{lem:pm2} is based on \cite[Prop 3.4]{rvk:modular},  and we mention that the height inequality in \cite[Prop 3.4]{rvk:modular} relies inter alia on explicit versions of certain results of Faltings~\cite{faltings:finiteness,faltings:arithmeticsurfaces} and Silverman~\cite{silverman:arithgeo}. 

\begin{proof}[Proof of Lemma~\ref{lem:pm2}]
Suppose that $(x,y)$ is a solution of~\eqref{eq:mordell}. Then an application of \cite[Prop 3.4]{rvk:modular} with $S=\sp(\mathcal O)$ and $T=\sp(\mathcal O[1/(6a)])$ gives an elliptic curve over~$T$, with generic fiber $E=E_\QQ$, such that $N_E\mid 2^{2}3^2a_S$ and such that $h(x)$ is bounded in terms of $h(E)$ and $h(a)$ as desired. Let  $c_4$ and $c_6$ be the quantities associated to~$E$. The proof of \cite[Prop 3.4]{rvk:modular} shows in addition that
 $c_4=u^4x$ and $c_6=u^6y$ for some $u\in \QQ$ with $u^{12}=1728\Delta_E\abs{a}^{-1}$.  Thus on combining \cite[Lem 3.3]{rvk:modular} with \cite[Lem 3.5]{rvk:modular}, we deduce an upper bound for $h(y)=\frac{1}{2}h(u^{-12}c_6^2)$ in terms of $h(E)$ and $h(a)$ as desired. 
Furthermore, it was shown in the proof of \cite[Prop 3.4]{rvk:modular}  that  $E$ admits a Weierstrass model 
$W$ over $\mathcal O$ defined by the Weierstrass equation $t^2 = s^3 - 27xs - 54y$
with ``indeterminates'' $s$ and~$t$. This Weierstrass equation has discriminant $\Delta'=-2^63^9a$. 

It remains to show that $N_E\mid a_S$. For this purpose, we observe that $2^83^5=2^23^2\cdot 1728$ and we recall that $N_E\mid 2^23^2a_S$. Hence it suffices to consider the exponents $f_2$ and $f_3$ of $N_E$ at the primes 2 and 3 respectively. They always satisfy $f_2\leq 8$ and $f_3\leq 5$, see for example  \cite[Thm 6.2]{brkr:conductor}. Thus we obtain that $f_2\leq 8= \ord_2(a_S)$ if $2\in S$ and  $f_3\leq 5= \ord_3(a_S)$ if $3\in S$, and then we see that the desired result $N_E\mid a_S$ holds if $2\in S$ or $3\in S$. Hence we are left to treat the case where 2 and 3 are both not in $S$, or equivalently that $2,3\in \sp(\mathcal O)$. In this case, after localizing at 2 and 3,  we may and do view the model $W$ over $\mathcal O$ as an integral model of $E$ over the local rings at 2 and 3 respectively. 

We first consider~$f_2$. It holds that $\ord_2(\Delta_E)\leq \ord_2(\Delta')$, since $W$ is a Weierstrass model of $E$ over the local ring at 2 and since $\Delta_E$ is minimal at~2. Thus the formula of Ogg--Saito \cite[Cor 2]{saito:conductor} implies that $f_2\leq \ord_2(\Delta')=6+\ord_2(a)$, and then the inequality $f_2\leq 8$ proves that $f_2\leq 6+\min(2,\ord_2(a))=\ord_2(a_S)$ as desired.

We now consider~$f_3$. If $W$ is not the minimal Weierstrass model of $E$ over the local ring at~$3$, then we get that $\ord_3(\Delta')\geq 12$ and hence $\ord_3(a)\geq 3$. This together with $f_3\leq 5$ shows the desired inequality $f_3\leq 3+\min(2,\ord_3(a))=\ord_3(a_S)$ when $W$ is not minimal at~$3$. Hence we may and do assume in addition that $W$ is minimal at~$3$. Then we claim that  the ``reduction" of $E$ at $3$ is of Kodaira type IV$^*$, III$^*$, or II$^*$ and that $f_3 = \ord_3(\Delta')-6$, $f_3=\ord_3(\Delta')-7$, or $f_3=\ord_3(\Delta')-8$  respectively. It follows that $f_3\leq 3+\ord_3(a)$ and this together with the  inequality $f_3\leq 5$ shows that $f_3\leq 3+\min(2,\ord_3(a))=\ord_3(a_S)$ as desired. It remains to verify our claim. 
For this purpose we use Tate's algorithm~\cite{tate:algo} and we let $a_1,\ldots,a_6$ and $b_2,\ldots,b_8$ denote the usual quantities associated to the Weierstrass equation $t^2 = s^3 - 27xs - 54y$ which is minimal at 3 by assumption. 
We observe that
$3\divides\Delta'$,
$3\divides a_1=0$,
$3\divides a_2=0$, 
$3^3\divides a_3=0$,
$3^2\divides a_4=-27 x$,
$3^3\divides a_6=-54y$,
$3\divides b_2=0$,
$3^3\divides b_6=-216y$, and
$3^3\divides b_8=-729x^2$. 
This brings us into case 6) of Tate's algorithm.
Here one considers the polynomial $P(T):=T^3+a_{2,1}T^2+a_{4,2}T+a_{6,3}=T^3-3xT-2y$. In~$\FF_3$, it reduces to $T^3-2y$ which is purely inseparable with $3$-fold root~$2y$. Therefore 8), 9) and 10) of Tate's algorithm prove our claim. This completes the proof of Lemma~\ref{lem:pm2}.\end{proof}

In what follows we need to control the function $\nu(n)=n\nu^*(n)$ on $\ZZ_{\geq 1}$, where $\nu^*$ is defined in Proposition~\ref{prop:explbounds}~(ii). If $m,n$ are in $\ZZ_{\geq 1}$ with $m$ dividing $n$, then it holds
\begin{equation}\label{eq:nuineq}
\nu(m)\leq \nu(n). 
\end{equation}
Indeed this inequality directly follows  by unwinding the definitions and by taking into account rational primes $p$ with $\ord_p(m)=1$ and $\ord_p(n)\geq 2$. We now combine the above Lemma~\ref{lem:pm2} with Proposition~\ref{prop:explbounds} in order to prove Proposition~\ref{prop:m}.

\begin{proof}[Proof of Proposition~\ref{prop:m}]
Suppose that $(x,y)$ is a solution of Mordell's equation~\eqref{eq:mordell}. Then  Lemma~\ref{lem:pm2}  provides an elliptic curve $E$ over $\QQ$ such that $N_E\mid a_S$ and such that $$\max\bigl(h(x),\tfrac{2}{3}h(y)\bigl)\leq\tfrac{1}{3}h(a)+8h(E)+2\log\max(1,h(E))+36.$$
It follows from~\eqref{eq:nuineq} that $\nu(N_E)\leq \frac{2}{3}a_S$.
Therefore on combining the displayed inequality with the explicit bound for $h(E)$ in terms of $\nu(N_E)$ given in Proposition~\ref{prop:explbounds}~(ii), we deduce an estimate for $\max\bigl(h(x),\frac{2}{3}h(y)\bigl)$ in terms of $a_S$ and $h(a)$ as stated in Proposition~\ref{prop:m}.
\end{proof}

We recall that a solution $(x,y)\in\ZZ\times\ZZ$ of~\eqref{eq:mordell} is  primitive if $\pm 1$ are the only  $n\in\ZZ$ with $n^{6}\mid \gcd(x^3,y^2)$. The following result refines Lemma~\ref{lem:pm2} for primitive solutions.

\begin{lemma}\label{lem:pm3}
Suppose that $(x,y)\in\ZZ\times\ZZ$  is a primitive solution of~\eqref{eq:mordell}. If Lemma~\ref{lem:pm2} associates $(x,y)$ to the elliptic curve $E$, then $E$ has in addition the following properties. 
\begin{itemize}
\item[(i)] It holds that  $\Delta_E=2^m3^n\lvert a\rvert$ with $m\in\{-6,6\}$ and $n\in\{-3,9\}$.
\item[(ii)] The curve $E$ is semi-stable at each rational prime $p\geq 5$ with $\gcd(x,y,p)=1$.
\item[(iii)] There is the refined height inequality $$\max\bigl(h(x),\tfrac{2}{3}h(y)\bigl)\leq 4h(E)+2\log\max(1,h(E))+28.$$
\end{itemize}
\end{lemma}
\begin{proof}
Let $(x,y)\in\ZZ\times\ZZ$ be a primitive solution  of~\eqref{eq:mordell}. Suppose that Lemma~\ref{lem:pm2} associates $(x,y)$ to the elliptic curve $E$ over $\QQ$. To prove (i) we recall from the proof of Lemma~\ref{lem:pm2} that $t^2 = s^3 - 27xs - 54y$ is a Weierstrass equation of $E$
with ``indeterminates'' $s$ and $t$ and with discriminant $\Delta'=-2^63^9a$. This Weierstrass equation has associated quantities $c_4'=6^4x$ and $c_6'=6^6y$. On using that $x,y$ are both in $\ZZ$ and that $\Delta_E$ is the discriminant of a minimal Weierstrass model of $E$ over $\ZZ$, we obtain a nonzero $u\in\ZZ$ with 
\begin{equation}\label{eq:mc4c6}
u^4c_4=6^4x,  \ \ \ u^6c_6=6^6y, \ \ \ u^{12}\Delta_E=2^63^9\abs{a}.
\end{equation}
It follows that $u^{12}$ divides $6^{12}\gcd(x^3,y^2)$ and hence we deduce that $\abs {u}\in\{1,2,3,6\}$ since $(x,y)$ is primitive by assumption. Therefore the equality $u^{12}\Delta_E=2^63^9\abs{a}$ implies that $2^m3^n\lvert a\rvert=\Delta_E$ with $m\in\{-6,6\}$ and $n\in\{-3,9\}$. This proves assertion~(i).

To show the semi-stable properties of $E$ claimed in~(ii), we take a rational prime $p$. If $p\nmid \Delta_E$ then $E$ has good (and thus semi-stable) reduction at $p$. We now assume that $p\mid \Delta_E$, that  $p\geq 5$ and that $\gcd(x,y,p)=1$.
Then~\eqref{eq:mc4c6} together with $\abs{u}\in \{1,2,3,6\}$ implies that $\gcd(c_4,c_6,p)=1$. Hence we obtain that $p\nmid c_4$ since $p\geq 5$ divides $12^3\Delta_E=\abs{c_4^3-c_6^2}$ by assumption. Therefore the semi-stable criterion in \cite[p.196]{silverman:aoes} shows that $E$ has semi-stable reduction at each $p\geq 5$ with $\gcd(x,y,p)=1$. This proves assertion~(ii).

It remains to show the refined height inequality in~(iii). The identities in~\eqref{eq:mc4c6} together with \cite[Lem 3.5]{rvk:modular} and $\abs{u}\leq 6$ lead to an explicit upper bound for $\max\bigl(h(x),\frac{2}{3}h(y)\bigl)$ in terms of $h(E)$ as claimed in~(iii). This completes the proof of Lemma~\ref{lem:pm3}. 
\end{proof}

The above proof shows in addition that Lemma~\ref{lem:pm3} holds more generally for all solutions $(x,y)\in\ZZ\times\ZZ$ of~\eqref{eq:mordell} such that $\pm 1$ are the only $n\in\ZZ$ with $n^{12}\mid \gcd(x^3,y^2)$. We now use Lemma~\ref{lem:pm3} and Proposition~\ref{prop:m} to prove Theorem~\ref{thm:m}.

\begin{proof}[Proof of Theorem~\ref{thm:m}]
We begin to prove the first part of assertion~(i). Let $\mu:\ZZ_{\geq 1}\to\mathbb R_{\geq 0}$ be an arbitrary function, and suppose that $(x,y)$ is a solution of~\eqref{eq:mordell} which is almost primitive with respect to~$\mu$. Let $u=u_{x,y}$ be as in Definition~\ref{def:ap}, and define $x'=u^2x$ and $y'=u^3y$. Then it follows that $(x',y')\in\ZZ\times\ZZ$ is a primitive solution of the Mordell equation $y'^2=x'^3+a'$ for $a'=u^6a\in\ZZ$.  Therefore Lemmas~\ref{lem:pm2} and~\ref{lem:pm3}  give an elliptic curve $E$ over $\QQ$ together with integers $m\in \{-6,6\}$ and $n\in\{-3,9\}$ such that 
\begin{equation}\label{eq:deltaea}
\Delta_E=2^{m}3^n\abs{a'} \ \ \ \textnormal{ and } \ \ \ N_E\mid a'_{S}.
\end{equation}
The construction of the number $u=u_1/u_2$ in Definition~\ref{def:ap} shows that any rational prime $p$ with $\ord_p(a')>\ord_p(a)$ satisfies $p\mid u_1$ and thus $p\in S$ since $x,y\in \mathcal O$. It follows that $\ord_p(a')\leq\ord_p(a)$ for all rational primes $p$ not in $S$ and therefore we find that $a'_S$ divides $a_S.$ 
Further it holds that $h(a)\leq 6h(u)+\log\abs{a'}$ and $h(u)\leq \mu(a_S)$ since $(x,y)$ is almost primitive with respect to~$\mu$. Then on combining~\eqref{eq:deltaea} and $N_E\mid a'_S\mid a_{S}$ with  the estimate for $\log \Delta_E$ in terms of $N_E$ given in~\eqref{eq:szpiro}, we derive an explicit upper bound for $h(a)$ in terms of $a_S$ and $\mu(a_S)$ which together with Proposition~\ref{prop:m} implies the first part of~(i).

We now show the second part of~(i). Since $(x,y)$ is almost primitive with respect to $\mu$ we obtain that $h(x)\leq 2\mu(a_S)+h(x')$ and $\frac{2}{3}h(y)\leq 2\mu(a_S)+\frac{2}{3}h(y')$. Hence an application of Lemma~\ref{lem:pm3}~(iii) with the primitive $(x',y')\in\ZZ\times\ZZ$ leads to  
\begin{equation}\label{eq:primitivesolestimate}
\max\bigl(h(x),\tfrac{2}{3}h(y)\bigl)\leq 2\mu(a_S)+4h(E)+2\log\max(1,h(E))+28.
\end{equation}
Then on combining~\eqref{eq:deltaea} and $a'_S\mid a_{S}$ with the explicit estimate for $h(E)$ in terms of $N_E$ given in Proposition~\ref{prop:explbounds}, we deduce  an inequality as claimed in the second part of~(i).

To prove~(ii) we may and do assume that $\abs{a}\to \infty$ and that $S$ is empty. Let $(x,y)\in\ZZ\times\ZZ$ be a primitive solution of~\eqref{eq:mordell}, and  take $(x',y')=(x,y)$ and $a'=a$ in the proof of~(i). Then~\eqref{eq:deltaea} and~\eqref{eq:szpiro} show that $\log\abs{a}\leq O( N_E^2)$  and thus $\abs{a}\to \infty$ forces $N_E\to\infty$. Hence on using~\eqref{eq:primitivesolestimate} with $\mu=0$  and on applying the asymptotic bound for $h(E)$ in terms of $N_E$ given in Proposition~\ref{prop:explbounds}, we see that~\eqref{eq:deltaea} together with $a_S=a_*$ leads to the asymptotic estimate claimed in~(ii). This completes the proof of Theorem~\ref{thm:m}.
\end{proof}

Next we give a proof of the inequality claimed in Corollary~\ref{cor:coates1} by using a version of Theorem~\ref{thm:m}~(i) which takes into account Lemma~\ref{lem:pm3}~(ii).

\begin{proof}[Proof of Corollary~\ref{cor:coates1}]
Recall that $a\in\ZZ$ is nonzero, $S$ is an arbitrary finite set of rational primes and $f\in\ZZ$ is the largest divisor of $a$ which is only divisible by primes in~$S$. We take $(x,y)\in\ZZ\times\ZZ$ with $y^2-x^3=a$ and $\gcd(x,y,N_S)=1$ as in the statement.   

We first show that $(x,y)$ is almost primitive with respect to $\mu=\frac{1}{6}\log\abs{a/f}$. Let $m=u_2\in\ZZ$ be maximal such that $m^{6}\mid\gcd(x^3,y^2)$, and  define $x'=m^{-2}x$ and $y'=m^{-3}y$. Then $(x',y')\in\ZZ\times\ZZ$ is a primitive solution of the Mordell equation $y'^2=x'^3+a'$ for $a'=m^{-6}a\in\ZZ$. Further, $\gcd(x,y,N_S)=1$ together with $\rad(f)\mid N_S$ implies that $\gcd(x,y,f)=1$, and hence  $\gcd(m,f)=1$  since $m$ divides $x$ and~$y$. Thus on using that $m^6\mid a$, we see that $m^{6}$ divides~$a/f$. This shows that the number $u_{x,y}=1/m$ in Definition~\ref{def:ap} satisfies $h(u_{x,y})\leq \mu$ as desired. 
We now apply (a version of) Theorem~\ref{thm:m}~(i). Suppose that Lemma~\ref{lem:pm2} associates the primitive $(x',y')\in\ZZ\times\ZZ$ to the elliptic curve $E$ over $\QQ$.  Lemma~\ref{lem:pm3}~(ii) gives that $E$ is semi-stable at each $p\in S$ with $p\geq 5$, since $\gcd(x',y',N_S)=1$. Then (the proof of) Lemma~\ref{lem:pm2} shows that $N_E$ divides $6\cdot 1728N_Sr_2(a')$, and hence we obtain that $N_E\mid 6\alpha_S$ since $a'\mid a$. 
Thus, on replacing $a_S$ by $6\alpha_S=6a_S/N_S$ in the proof of Theorem~\ref{thm:m}~(i) and on taking $\mu=\frac{1}{6}\log\abs{a/f}$, we deduce Corollary~\ref{cor:coates1}.
\end{proof}

\subsection{$S$-unit equations}\label{sec:suproofs}

We use the notation introduced in Section~\ref{sec:notations} above. The discussions in Section~\ref{sec:sucremalgo} show that bounding the Weil height of the solutions of the $S$-unit equation~\eqref{eq:sunit} is equivalent to estimating the absolute value of the integers satisfying the Diophantine equation 
\begin{equation}
a+b = c, \quad a,b,c\in\ZZ-\{0\}, \quad \gcd(a,b,c)=1,\quad \rad(abc)\mid N_S. \tag{\ref{eq:abc}}
\end{equation}
We first prove Lemma~\ref{lem:psu2} which is used in the proof of Proposition~\ref{prop:su} given below. 
\begin{lemma}\label{lem:psu2}
Suppose that $(a,b,c)$ is a solution of~\eqref{eq:abc}. Then there exist $\QQ$-isogenous elliptic curves $E$ and $E'$ over $\QQ$ such that $N_E$ divides $2^4N_S$, $\Delta_E=2^n(abc)^2$ with $n\in\{4,-8\}$, and  $\Delta_{E'}=2^{8-12m}\abs{ab}c^4$ with $m\in\{0,1,2,3\}$.
\end{lemma}

The proof of Lemma~\ref{lem:psu2} uses inter alia a formula of Diamond--Kramer~\cite{dikr:modularity} and classical results which are given for example in~\cite{silverman:aoes}. In fact the proof consists of explicit computations with Frey--Hellegouarch elliptic curves. A more conceptual proof of (parts of) Lemma~\ref{lem:psu2} can be given by using the point of view in \cite[Sect 3]{rvk:modular}.
\begin{proof}[Proof of Lemma~\ref{lem:psu2}]
We now use in particular several (known) reductions. To make sure that these reductions are compatible with each other we prefer to give all details. 

We suppose that $(a,b,c)$ is a solution of~\eqref{eq:abc}. Then exactly one of the numbers $a, b, c$ is even, since they are coprime and satisfy $a+b=c$. We denote this number by $B'$.  Now we define $A'=a$ if $B'=b$, $A'=b$ if $B'=a$, and $A'=-a$ if $B'=c$.  Further we define  $(A,B,C)=(A',B',A'+B')$ if $A'=-1\pmod{4}$, and $(A,B,C)=(-A',-B',-A'-B')$ otherwise. Let $E$ be the elliptic curve over $\QQ$, defined by the Weierstrass equation
$$y^2=x(x-A)(x+B)$$
with discriminant $\Delta=2^4(ABC)^2$.
We observe that $(A,B,C)$ has the following properties: $A,B,C$ are coprime, $A=-1\pmod{4}$,  $B$ is even, $A+B=C$ and $(ABC)^2=(abc)^2$. Thus  \cite[Lem 1 and Lem 2]{dikr:modularity} give that the conductor $N_E$ of $E$ divides $2^4\rad(abc)$. Our assumption, that $(a,b,c)$ satisfies~\eqref{eq:abc}, provides that $\rad(abc)\mid N_S$ and hence we obtain that $N_E$ divides~$2^4 N_S$.
It follows from \cite[Lem 2]{dikr:modularity} and \cite[p.257]{silverman:aoes} that the minimal discriminant $\Delta_E$ of $E$ satisfies $\Delta_E=\Delta$ if $\ord_2(abc)\leq 3$, and $\Delta_E=2^{-8}(abc)^2$ if $\ord_2(abc)\geq 4$. We conclude that $\Delta_E=2^n(abc)^2$ with $n\in\{4,-8\}$, and since $N_E$ divides $2^4N_S$ we see that the elliptic curve $E$ has all the desired properties. 

To construct the elliptic curve $E'$, we notice that $A,B,C$ are in $\{\pm a,\pm b,\pm c\}$. We define $(\alpha,\beta,\gamma,x')=(A,B,C,x)$ if $C=\pm c$, $(\alpha,\beta,\gamma,x')=(C,-B,A,x+B)$ if $A=\pm c$, and $(\alpha,\beta,\gamma,x')=(-A,C,B,x-A)$ if $B=\pm c$. It follows that $\gamma=\pm c$, that $\alpha+\beta=\gamma$, that $(\alpha,\beta)=1$ and that $E$  admits the Weierstrass equation $y^2=x'(x'-\alpha)(x'+\beta).$ 
Then we obtain that $E$ is $\QQ$-isogenous to the elliptic curve $E'$ over $\QQ$, which is defined by the Weierstrass equation (see for example \cite[p.70]{silverman:aoes}) $$w^2=z^3-2(\beta-\alpha)z^2 +\gamma^2z$$
with discriminant $\Delta'=-2^8\alpha\beta\gamma^4$. This Weierstrass equation has associated quantities $c_4=16(\alpha^2-14\alpha\beta+\beta^2)$ and $c_6=64(\alpha^3+33\alpha^2\beta-33\alpha\beta^2-\beta^3)$. Further we observe that $\Delta'=-2^8abc^4$. To determine the minimal discriminant $\Delta_{E'}$ of $E'$ we use a strategy inspired by \cite[p.258]{silverman:aoes}. An application of the Euclidean algorithm gives the identities
\begin{align*}
4(395\alpha^2-430\alpha\beta-13\beta^2 )c_4 + (181\alpha - 13\beta)c_6 &= 2^{12} 3^2\alpha^4,\\
4(-13\alpha^2 - 430\alpha\beta + 395\beta^2)c_4 + (13\alpha - 181\beta)c_6 &= 2^{12} 3^2\beta^4.
\end{align*}
On using that $\Delta',c_4,c_6\in\ZZ$ and that $\Delta_{E'}$ is the absolute value of the discriminant of a minimal Weierstrass model of $E'$ over $\ZZ$, we obtain $u\in\ZZ$ such that $u^{12}\Delta_{E'}=\pm\Delta'$, $u^4\mid c_4$ and  $u^6\mid c_6$. 
Thus the displayed identities together with $(\alpha,\beta)=1$ imply that $u^4\mid 2^{12} 3^2$ and
hence $\pm u\in\{1,2,4,8\}$. We deduce that  $\Delta_{E'}=2^{8-12m}\abs{ab}c^4$ with $m\in\{0,1,2,3\}$, and therefore  $E'$ has the desired properties. This completes the proof of Lemma~\ref{lem:psu2}.
\end{proof}

We remark that the application of \cite[Lem 2]{dikr:modularity} in the above proof  shows in addition the following: Suppose that $(a,b,c)$ is a solution of~\eqref{eq:abc}. If  Lemma~\ref{lem:psu2}  associates $(a,b,c)$ to the elliptic curve $E$ over $\QQ$, then the conductor $N_E$ of $E$ satisfies
\begin{equation}\label{eq:refinedcondbound}
\ord_2(N_E)=\ee+1, \ \ \ (\ee,\lambda)=\begin{cases}
(4,12) & \textnormal{if }\ord_2(abc) = 1,\\
(2,3) & \textnormal{if }\ord_2(abc) = 2, 3,\\
(-1,\frac{1}{2}) & \textnormal{if }\ord_2(abc) = 4,\\
(0,1) & \textnormal{if }\ord_2(abc) \geq 5.\\
\end{cases}
\end{equation}
We note that it always holds that $\ord_2(abc)\geq 1$, since $a+b=c$ and  $(a,b,c)=1$. The following result is in fact a (slightly) more precise version of Proposition~\ref{prop:su}. 
\begin{proposition}\label{prop:abc}
Suppose that $(a,b,c)$ is a solution of~\eqref{eq:abc}, and let $(\ee,\lambda)$ be the numbers in~\eqref{eq:refinedcondbound} associated to $(a,b,c)$. Then it holds $$\log\max\left(\abs{a},\abs{b},\abs{c}\right)\leq \tfrac{\lambda}{5}N_S\log(2^{\ee}N_S)+\tfrac{3\lambda}{40} N_S\log\log\log(2^{\ee}N_S)+\tfrac{2\lambda}{15}N_S+28.$$
\end{proposition}

\begin{proof}
We begin by noting that the elliptic curves appearing in Lemma~\ref{lem:psu2} have the same conductor since they are $\QQ$-isogenous. Thus an application of Lemma~\ref{lem:psu2} with the solution $(a,b,c)$ of~\eqref{eq:abc} gives an elliptic curve $E'$ over~$\QQ$, with conductor $N_{E'}$ and minimal discriminant~$\Delta_{E'}$,  such that $\abs{ab}c^4\leq 2^{28}\Delta_{E'}$ and such that $N_{E'}\mid 2^{\ee}N_S$ for $\ee$ as in~\eqref{eq:refinedcondbound}. The equality $a+b=c$ proves that $\abs{c}-1\leq \abs{ab}$ 
and then we deduce $$(\abs{c}-1)^5\leq (\abs{c}-1)c^4\leq \abs{ab}c^4\leq 2^{28}\Delta_{E'}.$$
Further, 
$N_{E'}\mid 2^{\ee}N_S$ with $\ord_2(N_{E'})=\ee+1$ implies that $\nu(N_{E'})\leq\lambda N_S$ 
for $\lambda$ as in~\eqref{eq:refinedcondbound} and for $\nu$ the function on $\ZZ_{\geq 1}$ defined in Proposition~\ref{prop:explbounds}~(ii). Then we see that~\eqref{eq:szpiro} leads to Proposition~\ref{prop:abc} provided that $H=\abs{c}$ for $H=\max(\abs{a},\abs{b},\abs{c})$. To deal with the remaining cases $H=\abs{a}$ and $H=\abs{b}$, we notice that $(b,-c,-a)$ and $(a,-c,-b)$ are solutions of~\eqref{eq:abc} as well. Thus applications of the above arguments with $(b,-c,-a)$ and $(a,-c,-b)$ prove  Proposition~\ref{prop:abc} in the cases $H=\abs{a}$ and $H=\abs{b}$ respectively. 
\end{proof}

Let $(a,b,c)$ be a solution of~\eqref{eq:abc}. Lemma~\ref{lem:psu2} associates to $(a,b,c)$ an elliptic $E$ over $\QQ$ with conductor $N_E$  that satisfies $N_E\to \infty $ if $\rad(abc)\to \infty$. Therefore the arguments of Proposition~\ref{prop:abc} together with  the asymptotic bound obtained below~\eqref{eq:szpiro}  show that any solution $(a,b,c)$ of~\eqref{eq:abc} with $\rad(abc)=r$ satisfies
\begin{equation}\label{eq:asymptoticsu}
\log \max(|a|,|b|,|c|)\leq \tfrac{9}{5}r\log r + O\Big(\frac{r\log r}{\log\log r}\Big) \ \  \textnormal{ for } r\to \infty.
\end{equation}
This  (slightly) improves the bound $4r\log r+O(r\log\log r)$ obtained by Murty--Pasten in \cite[Thm 1.1]{mupa:modular}. However, our  estimate displayed in~\eqref{eq:asymptoticsu} is still worse than the actual best asymptotic bound $O(r^{1/3}(\log r)^3)$ of Stewart--Yu \cite[Thm 1]{styu:abc2}. 

\begin{proof}[Proof of Proposition~\ref{prop:su}] We suppose that $(x,y)$ satisfies the $S$-unit equation~\eqref{eq:sunit}. Then there exists a solution $(a,b,c)$ of~\eqref{eq:abc} with $(x,y)=(\frac{a}{c},\frac{b}{c})$. The number $\max(h(x),h(y))$ equals $\log\max(\abs{a},\abs{b},\abs{c})$ and thus Proposition~\ref{prop:abc} implies Proposition~\ref{prop:su}.
\end{proof}

\subsection{Optimized height bounds and height conductor inequalities}\label{sec:hc}

We use the notation of Sections~\ref{sec:cremonas+st} and~\ref{sec:notations}. Let $N\geq 1$ be an integer. We now define constants $\alpha$, $\beta$ and $\beta^*$ depending on $N$, which will appear in the optimized height bounds. 

\subsubsection{The constants $\alpha$, $\beta$, $\beta^*$ and optimized height bounds}\label{sec:optimizedbounds} To define $\beta$ and $\beta^*$ we let $m$ be the number of newforms of level dividing $N$, and we let $g=g(N)$ be the genus of $X_0(N)$. We write $l=\lfloor\frac{N}{6}\prod (p+1)\rfloor$ and $l^*=\lfloor\frac{N}{6}\prod (1+1/p)\rfloor$ with both products taken over all rational primes $p$ dividing $N$, where $\lfloor r\rfloor=\max(n\in\ZZ; \, n\leq r)$ for any $r\in\mathbb R$. Let $\tau(n)$ be the number of divisors of any $n\in\ZZ_{\geq 1}$. We define
\begin{equation}\label{def:bb*}
\beta=\tfrac{1}{2}m\log m+\max_J \sum_{j\in J}\log(\tau(j)j^{1/2}) \  \textnormal{ and }  \ \beta^*=\tfrac{1}{2}g\log g+\max_J \sum_{j\in J}\log(\tau(j)j^{1/2})
\end{equation}
with the first  maximum  taken over all subsets $J\subseteq\{1,\dotsc,l\}$ of cardinality $m$ and with the second maximum taken over all subsets $J\subseteq\{1,\dotsc,l^*\}$ of cardinality~$g$. On comparing $\beta$ with $\beta^*$ we notice that $\beta$ involves the smaller number $m\leq g$ at the expense of depending on the larger parameter $l\geq l^*$. It turns out that $\beta\leq O(N\log N)$, while such an upper bound can not hold for $\beta^*$ since there is a constant $r\in\mathbb R$ such that  infinitely many  $n\in\ZZ_{\geq 1}$ satisfy $n\log\log n\leq rg(n)$. 
On the other hand, it holds that $\beta^*< \beta$ for infinitely many $N$ and thus we shall work in our algorithms with the quantity $$\alpha=\min(\beta,\beta^*).$$  
We define $\kappa=4\pi+\log(163/\pi)$ and we mention that in the statement of the following Proposition~\ref{prop:algobounds} we use the notation of Propositions~\ref{prop:m}  and~\ref{prop:abc}.

\begin{proposition}\label{prop:algobounds}
Proposition~\ref{prop:abc} holds with the bound $\frac{6}{5}\alpha(2^{\ee}N_S)+28,$ and Proposition~\ref{prop:m} holds with the bound $\frac{1}{3}h(a)+4\alpha(a_S)+2\log(\alpha(a_S)+\kappa)+35+4\kappa$.  
\end{proposition}
\begin{proof}
We observe that $\alpha(m)\leq \alpha(n)$ for all $m,n$ in $\ZZ_{\geq 1}$ with $m$ dividing $n$. Therefore Proposition~\ref{prop:algobounds} follows directly by using in the above proofs of Propositions~\ref{prop:m} and~\ref{prop:abc} the optimized bound given in Proposition~\ref{prop:explbounds}~(i).
\end{proof}
We remark that for any given $N$, one can practically compute $\alpha$ by using the formulas for $m$ and $g$ in  \cite[Thm 4]{martin:dimension} and \cite[p.107]{dish:modular} respectively. However, if $N$ becomes large then the precise computation of $\alpha$ becomes slow and in this case we shall use  
\begin{equation}\label{def:barb}
\bar{\beta}=\tfrac{1}{2}m\log m+\tfrac{5}{8}m(18+\log l) \ \textnormal{ and } \ \bar{\beta^*}=\tfrac{1}{2}g\log (gl^*)+\tfrac{1}{2}l^*\log(4+4\log l^*).
\end{equation}
Remark~\ref{rem:tau} gives that $\tau(j)\leq 45197j^{1/8}$ for all $j\in\ZZ_{\geq 1}$ and hence $\beta\leq\bar{\beta}$.
To prove that $\beta^*\leq \bar{\beta^*}$ we may and do assume that $g\geq 1$.  Let $J\subseteq\{1,\dotsc,l^*\}$ be a subset of cardinality $g$ and observe that $g\leq n=\lfloor l^*/2\rfloor$. 
Then the elementary inequalities $l^*\leq 4n$, 
$\prod_{j\in J} \tau(j)\leq \bigl(\frac{1}{n}\sum_{j=1}^{l^*}\tau(j)\bigl)^n$ 
and $\frac{1}{l^*}\sum_{j=1}^{l^*}\tau(j)\leq 1+\log l^*$ show that $\sum_{j\in J}\log\tau(j)\leq n\log(4+4\log l^*)$. This implies that $\beta^*\leq\bar{\beta^*}$ as desired. 
It follows that $\alpha\leq \bar{\alpha}$ for
$$\bar{\alpha}=\min(\bar{\beta},\bar{\beta^*}).$$
 We take here the minimum 
since there are infinitely many $N$ for which $\bar{\beta^*}<\bar{\beta}$ and vice versa. Finally we point out that the computation of $\bar{\alpha}$ is very fast, even for large~$N$.

It remains to work  out  the explicit height conductor inequalities for elliptic curves over $\QQ$ which are used in our proofs and this will be done in the next section. 
\subsubsection{Height and conductor of elliptic curves over $\QQ$}\label{sec:heightcondstatement}

The geometric version of the Shimura--Taniyama conjecture  gives that any elliptic curve $E$ over $\QQ$ of conductor $N$ is $\QQ$-isogenous to $E_f$ for some rational newform $f\in S_2(\Gamma_0(N))$, see  Section~\ref{sec:cremonas+st}.  We say that $f$ is the newform associated to $E$. Let $m_f$  be the modular degree of $f$ and let $r_f$ be the congruence number of $f$, defined in (\ref{def:mf}) and (\ref{def:rf}) respectively. 

 On using and refining the arguments\footnote{These arguments use an approach of Frey which involves the modular degree $m_f$ and the geometric version of the Shimura--Taniyama conjecture, see for example Frey~\cite[p.544]{frey:ternary}.} of \cite[Thm 7.1]{mupa:modular} and \cite[Prop 6.1]{rvk:modular}, which were discovered independently in 2011 by Murty--Pasten and by the first mentioned author, we obtain the following update of several results (see below) in~\cite{mupa:modular,rvk:modular}.
\begin{proposition}\label{prop:explbounds}
Let $\beta$ and $\beta^*$ be as in \textnormal{(\ref{def:bb*})}. Suppose that $E$ is an elliptic curve over $\QQ$ of conductor~$N$, with associated newform~$f$.  Then the following statements hold.
\begin{itemize}
\item[(i)] Define $\kappa=4\pi+\log(163/\pi)$.  There are inequalities
$$
2h(E)-\kappa\leq \log m_f \leq \log r_f\leq \alpha=\min(\beta,\beta^*).$$
\item[(ii)] Let $\nu^*$ be the multiplicative function on $\ZZ_{\geq 1}$ defined by $\nu^*(p)=1$ for $p$ a rational prime  and $\nu^*(p^k)=1-1/p^{2}$ for $k\in\ZZ_{\geq 2}$, and put $\nu=N\nu^*(N)$. It holds $$\beta\leq \tfrac{1}{6}\nu\log N + \tfrac{1}{16}\nu\log\log\log N + \tfrac{1}{9}\nu, \ \ \ \nu\leq N.$$
\item[(iii)]  If $N\to \infty$ then $\beta\leq \frac{1}{8}\nu\log N + \frac{(\frac{1}{6}\log 2+o(1))}{\log\log N}\nu\log N.$
\end{itemize}
\end{proposition}
It is known (see e.g.~\cite{frey:linksulm,silverman:arithgeo}) that $h(E)$ is related to~$m_f$, and a result attributed to Ribet provides that~$m_f\mid r_f$. We now compare Proposition~\ref{prop:explbounds} with the literature. It improves $h(E)\leq (2N)^{40^2}$ in\footnote{The result \cite[Thm 3.6]{rvk:thesis} bounds a ``naive" height of $E$. However, Silverman's arguments in \cite[Sect 2]{silverman:arithgeo} lead to an explicit upper bound for $h(E)$ in terms of the ``naive"  height of $E$ used in~\cite{rvk:thesis}.} \cite[Thm 3.6]{rvk:thesis} and $h(E)\leq (25N)^{162}$ in \cite[Thm 3.1]{rvk:height} which were established by different methods:  \cite{rvk:thesis} uses the effective reduction theory of Evertse--Gy{\H{o}}ry  and~\cite{rvk:height} combines Legendre level structure with the theory of logarithmic forms. Furthermore, Proposition~\ref{prop:explbounds} updates $2h(E)\leq \frac{1}{5}N\log N+22$ in Murty--Pasten \cite[Thm 7.1]{mupa:modular} and $2h(E)\leq \frac{1}{2}N(\log N)^2+18$ in \cite[Prop 6.1]{rvk:modular}. It also updates the asymptotic bound $2h(E)\leq \frac{1}{6}N\log N+O(N\log\log N)$  in Murty--Pasten \cite[Thm 7.1]{mupa:modular} and the bounds for $m_f$ and $r_f$  in  \cite[Thm 4.3]{mupa:modular} and  \cite[Lem 5.1 (i)]{rvk:modular}. 

We proved Proposition~\ref{prop:explbounds} in our unpublished 2014 preprint  ``Solving S-unit and Mordell equations via Shimura--Taniyama conjecture". Hector Pasten told us that in Taipei (May 20, 2016) he will announce the following results which he obtained independently.

\begin{remark}[Recent results of Hector Pasten]  Let $S$ be a finite set of rational primes. For all elliptic curves $E$ over $\QQ$ that are semistable outside S, we have $$h(E) \ll_S \phi(N)\log N.$$ Furthermore, if we assume GRH then we have $h(E) \ll_S \phi(N)\log \log N.$ Here $\phi$ denotes the Euler totient function, and the constants  $\ll_S$ are effective, small and  explicit. Further,  it holds that $\log m_f\leq (o(1)+\tfrac{1}{24})N\log N$ as $N\to \infty$, and if $E$ is semistable then $\log m_f\leq (o(1)+\tfrac{1}{24})\phi(N)\log N $ as $N\to\infty$; both estimates are effective and can be made explicit with good constants.  We thank Hector Pasten for sending us the statements of his results.
\end{remark}

\subsubsection{Proof of Proposition~\ref{prop:explbounds}}

As already mentioned, the proof of Proposition~\ref{prop:explbounds} uses the ideas of \cite[Thm 7.1]{mupa:modular} and \cite[Prop 6.1]{rvk:modular}; see also Frey \cite[p.544]{frey:ternary}. For example the arguments of \cite[Lem 5.1, Prop 6.1]{rvk:modular} directly lead to $2h(E)-\kappa\leq \log m_f \leq \log r_f\leq  \beta^*.$
Then to replace here $\beta^*$ with $\beta$ we go through the proof of \cite[Lem 5.1]{rvk:modular} by using the ``coprime" matrix constructed in Murty--Pasten~\cite{mupa:modular}. To show (ii) we combine a formula of Martin~\cite{martin:dimension} with classical analytic number theory. The latter is used to explicitly bound the quantities $\prod_{p\mid n}(1+1/p)$ and $\tau(n)$ in terms of $n\in\ZZ_{\geq 1}$. Finally we deduce (iii) by replacing in the proof of (ii) our explicit estimate for $\tau(n)$ by Wigert's asymptotic bound.

\begin{proof}[Proof of Proposition~\ref{prop:explbounds}]
We first prove (i). To show that $\log r_f\leq \beta$ we denote by $\mathcal B=\{f_1,\ldots,f_g\}$  the Atkin--Lehner basis \cite[Thm 5]{atle:hecke} for $S_2(\Gamma_0(N))$, indexed in such a way such that $f_1=f$ and such that $f_1,\ldots,f_m$ are the newforms of level dividing $N$.
We write $I=\{1,\dotsc,m\}$ and $J'=\{j\in \{1,\dotsc,l\}; (j,N)=1\}$. On using Atkin--Lehner theory, Murty--Pasten proved in \cite[Prop 3.2]{mupa:modular} that the matrix $F'=(F_{ij})$  has full rank $m$ for $F_{ij}=a_j(f_i)$ with $i\in I$ and $j\in J'$.
Hence there is a subset $J\subseteq J'$ of cardinality $m$ such that the matrix $F=(F_{ij})$, with $i\in I$ and $j\in J$, is invertible. The Ramanujan--Petersson bounds for Fourier coefficients imply that  $\lvert F_{ij}\rvert\leq\tau(j)j^{1/2}$ for all $i\in I$ and $j\in J$. Thus Hadamard's determinant inequality shows that $\log \det(F)\leq \beta$ and then the claim $r_f^2\mid \det(F)^2\in \ZZ$ proves that $\log r_f\leq \beta$ as desired. To verify the claim $r_f^2\mid \det(F)^2\in \ZZ$ we take $f_c\in S_2(\Gamma_0(N))$ as in (\ref{def:rf}). Then there exists  $y=(y_j)\in\ZZ^J$ such that 
\begin{equation}
\label{eq:coeffOfFandFc}
a_j(f_c)=a_j(f)+y_jr_f \  \textnormal{ for all } j\in J.
\end{equation}
We write $f_c=\sum_{i=1}^g k_if_i$ with $(k_i)\in\CC^g$.
It holds that $k_1=0$, since $\mathcal B$ is an orthogonal basis and since $(f,f_c)=0$ by (\ref{def:rf}).
Therefore on comparing Fourier coefficients we deduce
\begin{equation}
\label{eq:FcInTermsOfFis}
a_j(f_c)=\sum_{i=2}^{g} k_ia_j(f_i)=\sum_{i=2}^m k_ia_j(f_i) \  \textnormal{ for all } j\in J.
\end{equation}
Here we used that $a_j(f_i)=0$ for all $i\geq m+1$ and $j\in J$. To see that $a_j(f_i)=0$ for all $i\geq m+1$ and $j\in J$, one recalls that each $j\in J\subseteq J'$ is coprime to $N$ and  any $f_i\in \mathcal B$ with $i\geq m+1$ is of the following form: $f_i(\tau)=f^*(n\tau)$ 
with  $f^*\in S_2(\Gamma_0(M))$ a  newform, $M$ a proper divisor of~$N$, and $n\geq 2$ a divisor of~$N/M$. 
Now on using exactly the same arguments as in the proof of \cite[Lem 5.1 (ii)]{rvk:modular}, we see that the formulas~\eqref{eq:coeffOfFandFc} and~\eqref{eq:FcInTermsOfFis} imply the claim $r_f^2\mid \det(F)^2$. It follows that $\log r_f\leq \min(\beta,\beta^*)$ since \cite[Lem 5.1 (ii)]{rvk:modular} gives the upper bound $\log r_f\leq \beta^*$. Further \cite[Thm 2.1]{agrist:congruence} implies that $\log m_f\leq \log r_f$.  Finally,  the desired lower bound for $\log m_f$ follows for example from the explicit inequality $h(E)\leq \frac{1}{2}\log m_f+2\pi+\frac{1}{2}\log(163/\pi)$ which was obtained in course of the proof of \cite[Prop 6.1]{rvk:modular}. This completes the proof of assertion~(i).

We now prove statement~(ii).
Martin obtained in \cite[Thm 4]{martin:dimension} an explicit formula for $m$ in terms of $N$ and~$\nu$. This formula implies the estimate
\begin{equation}
\label{eq:martinbound}
m\leq \frac{\nu}{12}-\frac{1}{2}+\frac{1}{3}+\frac{1}{4}=\frac{\nu+1}{12}.
\end{equation}
We next work out an  explicit upper bound for $l=\lfloor\frac{N}{6}\prod_{p\mid N}p(1+1/p)\rfloor$ with the product taken over all rational primes $p$ dividing~$N$.
Let $\gamma=0.577\dotsc$ denote the Euler--Mascheroni constant. In a first step we show that any $n\in\ZZ_{\geq 3}$ satisfies
\begin{equation}
\label{eq:keyexplbound}
\prod_{p\divides n} (1+1/p) \leq \frac{6 e^{\gamma}}{\pi^2}\left(\log\log n + \frac{2}{\log\log n}\right)
\end{equation}
with the product taken over all rational primes $p\mid n$.
To prove~\eqref{eq:keyexplbound} we may and do assume that $n$ is of the form $n=\prod_{p\leq x}p=e^{\vartheta(x)}$ with the product taken over all rational primes $p$ at most $x=x(n)\in\ZZ_{\geq 2}$. Indeed this follows by observing that the function $\log\log n+2/\log\log n$ is monotonously increasing for $n\geq 62$, by considering special cases such as for example the case when $n$ is a prime power and by checking (for example with Sage) all $n\leq 62$. 
On writing $1+1/p = (1-1/p^2)/(1-1/p)$, we obtain
$$\prod_{p\leq x}(1+1/p)=\left(\prod_{p}(1-\frac{1}{p^2})\cdot \bigl(\prod_{p>x}(1-\frac{1}{p^2})\bigl)^{-1}\right)/\prod_{p\leq x}(1-1/p)$$
with the products taken over all rational primes $p$ satisfying the specified conditions. The effective version of Merten's theorem in \cite[Thm 7]{rosc:formulas} provides
\[
\prod_{p\leq x}(1-1/p) > \frac{e^{-\gamma}}{\log x}\Big(1-\dfrac{1}{2(\log x)^2}\Big) \ \textnormal{ if } x\geq 285.
\]
Euler's product formula gives that $\prod_p (1-1/p^2)=\zeta(2)^{-1}=6/\pi^2$ with the product taken over all rational primes $p$. Further, we deduce the inequalities
\[
\log\prod_{p>x}\Big(1-\frac{1}{p^2}\Big) \geq -\sum_{p>x}\frac{1}{p^2}\Big(1+\frac{1}{2p^2}\Big)\geq -\Big(1+\frac{1}{2x^2}\Big)\int_{x}^{\infty}\frac{1}{t^2}\,dt = -\frac{1}{x}\Big(1+\frac{1}{2x^2}\Big).
\]
On combining the results collected above with the  effective prime number theorem of the following form (see \cite[Thm 4]{rosc:formulas})
\[
x-\frac{x}{2\log x} < \vartheta(x)=\log n < x+\frac{x}{2\log x} \ \textnormal{ if } x\geq 563,
\]
we see that the claimed inequality~\eqref{eq:keyexplbound} holds for all $n=e^{\vartheta(x)}$ with $x>10^4$.
Finally, one checks (for example with Sage) that~\eqref{eq:keyexplbound} holds in addition in the remaining cases $n=e^{\vartheta(x)}$ for $2\leq x\leq 10^4$. The conductor $N$ is at least $11$ and hence~\eqref{eq:keyexplbound} gives
\begin{equation}
\label{eq:lbound}
l\leq \frac{e^\gamma}{\pi^2}N^2\bigl(\log\log N + \frac{2}{\log\log N}\bigl).
\end{equation}
To estimate $\tau(n)$ we consider the real valued function $u(n)=\tau(n)/n^{1/4}$ on $\ZZ_{\geq 1}$. This function is multiplicative and  satisfies  $u(n)= \prod_{p}(n_p+1)p^{-n_p/4}$ with the product taken over all rational primes $p$, where $n_p=\ord_p(n)$ denotes the order of $p$ in $n\in\ZZ_{\geq 1}$. To find the maximum of $u$ we look at each factor separately.
It holds that $(n_p+1)p^{-n_p/4}=1$ when $n_p=0$, and if $n_p\geq 1$ then we observe that $(n_p+1)p^{-n_p/4}< 1$ provided that $p\geq 17$ or $n_p\geq 17$.
Thus after checking (for example with Sage) the remaining cases we find that $\sup(u)=8.44\dotsc$; in fact this supremum is attained at  $n=2^5\cdot 3^3\cdot 5^2\cdot 7\cdot 11\cdot 13$.
Therefore we conclude that any $n\in \ZZ_{\geq 1}$ satisfies the inequality
\begin{equation}
\label{eq:tau}
\tau(n)\leq 8.5\, n^{1/4}.
\end{equation}
To put everything together we recall that in $\beta=\frac{m}{2}\log m+\max_J \sum_{j\in J}\log(\tau(j)j^{1/2})$ the maximum is taken over all subsets $J\subseteq \{1,\dotsc,l\}$ of cardinality~$m$. Hence~\eqref{eq:tau} gives
\begin{equation*}
\beta\leq \beta'=m\left(\tfrac{1}{2}\log m+\tfrac{3}{4}\log l+\log 8.5\right).
\end{equation*}
Further Euler's product formula shows that $\nu\geq 6N/\pi^2$ and then~\eqref{eq:martinbound} together with~\eqref{eq:lbound} leads to $\beta'\leq \frac{\nu}{6}\log N+\frac{\nu}{16}\log\log\log N+\frac{\nu}{9}$ for all $N\geq 23$. 
Moreover, one checks (for example with Sage) that this upper bound holds in addition for all $N$ with $11\leq N< 23$.  Therefore the displayed inequality $\beta\leq \beta'$ proves (ii).

To show~(iii) we observe that $N\to \infty$ implies $l\to \infty$ and $m\to \infty$. Hence we may and do assume that $l\to \infty$ and $m\to \infty$. If $n\in \ZZ_{\geq 1}$ with $n\to \infty$, then Wigert's bound gives  $$\log \tau(n)\leq(\log 2+o(1))\frac{\log n}{\log\log n}.
$$
This implies that  $\sum_{j\in J}\log(\tau(j))\leq m(\log 2+o(1))\frac{\log l}{\log\log l}$ for any subset $J\subseteq \{1,\dotsc,l\}$ of cardinality~$m$. 
Therefore on recalling the definition of $\beta$ we obtain
\begin{equation*}
\beta\leq m\left(\tfrac{1}{2}\log m+\tfrac{1}{2}\log l+(\log 2+o(1))\frac{\log l}{\log\log l}\right).
\end{equation*}
Then we see that the above estimates for $m$ and $l$, given in~\eqref{eq:martinbound} and~\eqref{eq:lbound} respectively, lead to statement~(iii). This completes the proof of  Proposition~\ref{prop:explbounds}.
\end{proof}

\begin{remark}\label{rem:tau}
To prove (ii) we used the explicit bound $\tau(n)\leq 8.5n^{1/4}$ obtained in~\eqref{eq:tau}, since the constant $8.5$ is reasonably small. On enlarging the constant $8.5$, one could replace the exponent $1/4$ by any other positive real number. However, for small exponents the corresponding (effective) constants become quite large. For example, if we take the exponent $1/8$ then the constant $8.5$ needs to be replaced by $45196.8$.\end{remark}

We mention that Proposition~\ref{prop:explbounds} allows to update several bounds in the literature. For example, Proposition~\ref{prop:explbounds} together with \cite[Lem 3.3]{rvk:modular} directly implies that any elliptic curve $E$ over $\QQ$ of conductor $N$ and minimal discriminant $\Delta_E$ satisfies 
\begin{equation}\label{eq:szpiro}
\log\Delta_E\leq \nu\log N+\tfrac{3}{8}\nu\log\log\log N+\tfrac{2}{3}\nu+115.1
\end{equation}
and if $N\to\infty$  then $\log \Delta_E\leq \frac{3}{4}\nu\log N+\frac{\log 2+o(1)}{\log\log N}\nu\log N$.  Here $\nu$ is as in Proposition~\ref{prop:explbounds}~(ii).  These inequalities update the discriminant conductor inequalities in Murty--Pasten \cite[Thm 7.1]{mupa:modular} and   in \cite[Cor 6.2]{rvk:modular}; see also \cite[Thm 3.3]{rvk:szpiro} for more general (but weaker) discriminant bounds  based on the theory of logarithmic forms.

\section{The elliptic logarithm sieve}\label{sec:elllogsieve}

\subsection{Introduction}\label{sec:setup}

In this section we solve the problem of constructing an efficient sieve for the $S$-integral points of bounded height on any elliptic curve $E$ over $\QQ$ with given Mordell--Weil basis of $E(\QQ)$. Our construction combines a geometric interpretation of the known elliptic logarithm reduction (initiated by Zagier~\cite{zagier:largeintegralpoints}) with several conceptually new ideas. The resulting ``elliptic logarithm sieve" considerably extends the class of elliptic Diophantine equations which can be solved in practice. To illustrate this we solved many notoriously difficult equations by applying our sieve. 
We also used the resulting data and our sieve to motivate new conjectures and questions on the number of $S$-integral points of $E$. 

The precise construction of the elliptic logarithm sieve is rather technical. We now begin to describe  the underlying ideas using a geometric point of view: After fixing the setup and briefly discussing the known elliptic logarithm approach, we explain the main ingredients of our sieve and we describe the improvements provided by our new ideas. 

\paragraph{Setup.}Throughout this section we shall work with the following setup. Let $E$ be an elliptic curve over $\QQ$. We suppose that we are given an arbitrary Weierstrass equation of $E$, with coefficients $a_1,\dotsc,a_6$ in $\ZZ$, of the following form
\begin{equation}\label{eq:weieq}
y^2+a_1xy+a_3y=x^3+a_2x^2+a_4x+a_6.
\end{equation}
For any field $K$ containing $\QQ$, we often identify a nonzero point in $E(K)$ with the corresponding solution of \eqref{eq:weieq} and vice versa. We further assume that we are given a basis $P_1,\dotsc,P_r$ of the free part of the finitely generated abelian group $E(\QQ)$, see Section~\ref{sec:compmwbasis}.
Let $S$ be a finite set of rational primes and let $\Sigma(S)$ be the set of $(x,y)$ in $\mathcal O\times\mathcal O$ satisfying \eqref{eq:weieq}, where $\mathcal O=\ZZ[1/N_S]$ and $N_S=\prod_{p\in S} p$. 
Finally, we suppose that we are given an initial bound $M_0$, that is $M_0\in\ZZ_{\geq 1}$ such that any $P\in \Sigma(S)$ satisfies $\hat{h}(P)\leq M_0$ for $\hat{h}$ the canonical N\'eron--Tate height of $E$. In fact everything works equally well with an initial bound for the usual Weil height or for the infinity norm $\|\cdot\|_\infty$, see Remark~\ref{rem:elllogsievegen}~(ii).

\subsubsection{Elliptic logarithm approach}\label{sec:elllintroelr}

Starting with Masser~\cite{masser:ellfunctions}, Lang~\cite{lang:diophantineanalysis} and W\"ustholz~\cite{wustholz:recentprogress}, many authors developed a practical approach to determine $\Sigma(S)$ using elliptic logarithms. Here a fundamental ingredient is a technique introduced by Zagier~\cite{zagier:largeintegralpoints} which we call the elliptic logarithm reduction. 
In practice this technique allows to considerably reduce the initial bound $N_0$ coming from transcendence theory;  $\|P\|_\infty\leq N_0$ for all $P\in \Sigma(S)$. More precisely on combining Zagier's arguments  with de Weger's approach via $L^3$~\cite{lelelo:lll},
Stroeker--Tzanakis~\cite{sttz:elllogaa} and  Gebel--Peth{\H{o}}--Zimmer~\cite{gepezi:ellintpoints,gepezi:mordell} showed independently the following when $\mathcal O=\ZZ$: The elliptic logarithm reduction produces  a relatively small number $N_1<N_0$ such that any point $P\in \Sigma(S)$ with $N_1< \|P\|_\infty \leq N_0$ has to be exceptional (Definition~\ref{def:exceptpoint}).  Smart~\cite{smart:sintegralpoints} extended 
the method to general $\mathcal O$, see also \cite{pezigehe:sintegralpoints} and the recent book of Tzanakis~\cite{tzanakis:book}  devoted to the elliptic logarithm approach. 
If $N_1^r$ is small enough then $\Sigma(S)$ can be enumerated by checking all remaining candidates $P$ with $\|P\|_\infty\leq N_1$ and by finding the exceptional points. However there are many situations of interest in which  $N_1^r$ is usually too large to determine $\Sigma(S)$ via the known methods. In view of this, resolving the following problem would be of fundamental importance. 

\vspace{0.3cm}
\noindent{\bf Problem.}
\emph{Construct an efficient sieve for the points $\Sigma(S)\subseteq E(\QQ)$ of bounded height.}
\vspace{0.3cm}

\noindent  Here the known sieves are often useless in practice. For example, working in the finite groups obtained by ``reducing the curve $E$ mod $p$" for suitable primes $p$ is usually not efficient (see Smart \cite[p.398]{smart:sintegralpoints}).  In fact, since an efficient sieve for $\Sigma(S)$ inside $E(\QQ)$ was not available, various authors conducted some effort to  develop other techniques to enumerate $\Sigma(S)$ in certain cases when $N_1^r$ is not small enough; see Section~\ref{sec:comparisonwithelr}.

\subsubsection{The elliptic logarithm sieve}

Building on the core  idea of the elliptic logarithm reduction, we construct the elliptic logarithm sieve which resolves in particular the above problem. Here we introduce several conceptually new ideas. They all rely on a geometric point of view and to explain our ideas we now give a geometric interpretation of the known elliptic logarithm reduction:  For each $v$ in  $S^*=S\cup \{\infty\}$ one uses elliptic logarithms to construct a  lattice $\Gamma_v\subset\ZZ^d$ of rank $d$ such that any non-exceptional $P\in \Sigma(S)$ with $\|P\|_\infty>N_1$ is essentially determined  by a nonzero point in some $\Gamma_v$.  Then 
one tries to show via $L^3$  that a certain cube $Q_v\subset\RR^{d}$ satisfies $\Gamma_v\cap Q_v=0$ proving that all $P\in \Sigma(S)$ with $\|P\|_\infty>N_1$ are exceptional.  Here $d$ equals $r$ or $r+1$, and for any $v\in S$ the cube $Q_v\subset\RR^d$ is given by $\{\|P\|_\infty\leq N_0\}$ inside $E(\QQ)\otimes_\ZZ\RR\cong \RR^d$.  Further, increasing $N_1$ enlarges the co-volume of each $\Gamma_v$ and hence $\Gamma_v\cap Q_v$ is usually trivial for sufficiently large $N_1$. In fact the cube $Q_v$ always contains a certain ellipsoid $\mathcal E_v\subset \RR^d$ arising from $\hat{h}$. Now, our new ideas can be described as follows:

\paragraph{Global sieves.} We use $\Gamma_v\cap \mathcal E_v$ to construct various global sieves for $\Sigma(S)$ inside $E(\QQ)$. Here, instead of computing a lower bound for the length of the shortest nonzero vector in $\Gamma_v$, we actually determine the points in $\Gamma_v\cap \mathcal E_v$ using Fincke--Pohst~\cite{lelelo:lll,fipo:algo} and we check if these points come from $\Sigma(S)$. This has the following advantages: 
\begin{itemize}
\item[(i)] We can further reduce $N_1$ in the crucial situation where the usual reduction is not working anymore (e.g. the shortest nonzero vector of $\Gamma_v$ actually lies in $Q_v$).
\item[(ii)] On ``covering" the set $\Sigma(S)$ by the local sieves $\Gamma_v\cap\mathcal E_v$,  $v\in S^*$, we obtain a global sieve for $\Sigma(S)$ inside $E(\QQ)$ which is more efficient than the standard enumeration.
\end{itemize}
\paragraph{Refined coverings.}  On using the geometric point of view, we construct in Proposition~\ref{prop:refinedcov} refined ``coverings" of $\Sigma(S)$ in order to improve the global-local passage in (ii). This leads to a refined sieve which enhances our sieve in (ii) and which allows to reduce $N_1$ even further. Here the construction is inspired by our refined sieve for $S$-unit equations in Section \ref{sec:dwsieve+}. However, in the present case of elliptic curves, the technicalities arising from $v$-adic elliptic logarithms at $v=2$ and $v=\infty$ are more involved.

\paragraph{Height-logarithm sieve.} We construct a  sieve for $\Sigma(S)$ inside $E(\QQ)$ by exploiting that for any non-exceptional point $P\in \Sigma(S)$ the height $\hat{h}(P)$ is essentially determined by the local $v$-adic elliptic logarithms with $v\in S^*$.  The height-logarithm sieve is a crucial ingredient of our global sieves discussed above. For many involved points $P$, it allows to avoid the slow process of testing whether the coordinates of $P$ are in fact $S$-integers.  

\paragraph{Ellipsoids.} Instead of using the infinity norm $\|\cdot\|_\infty$ as done by all other authors, we work directly with the canonical height $\hat{h}$.  Our approach using $\hat{h}$ is more efficient than the known improvements of the elliptic logarithm reduction (see Section~\ref{sec:comparisonwithelr}), since the cubes $Q_v$ arising from $\|\cdot\|_\infty$ always contain our ellipsoids $\mathcal E_v$ determined by $\hat{h}$. In fact working here with ellipsoids is optimal from a geometric point of view and it is crucial for the construction of our sieves. To circumvent issues with the real valued function $\hat{h}$, we constructed a suitable rational approximation of the quadratic form $\hat{h}$ on $E(\QQ)\otimes_\ZZ\RR$.

\paragraph{Exceptional points.} We conducted some effort to avoid as much as possible working with the coordinates of the points. For example to deal with exceptional points, we prove the crucial Proposition~\ref{prop:refinedcov} which allows here to work entirely in the finitely generated group $E(\QQ)$. This considerably improves the ``extra search" for exceptional points. In fact in most cases Proposition~\ref{prop:refinedcov} completely removes the ``extra search".

\paragraph{Generic situation.}  The case $r\leq 1$ is of particular importance, since it represents the most common situation in practice; see also Katz--Sarnak~\cite{kasa:random} and Bhargava--Shankar~\cite{bhsh:avranklessthan1}. Furthermore, one can efficiently verify in practice whether $r\leq 1$ by using for example the work of Kolyvagin~\cite{kolyvagin:bsd} and Gross--Zagier--Zhang~\cite{zhang:gengz}. 
Also one can directly determine $\Sigma(S)$ when $r=0$. In view of this we tried to further improve our sieves for $r=1$. On exploiting that $\Gamma_v$ has rank $r=1$ for $v\in S$, we  optimized the reduction process at $v\in S$ and we enhanced our height-logarithm sieve for huge sets $S$.

\subsubsection{Discussion}

We shall motivate (using geometry) the ideas and constructions underlying the elliptic logarithm sieve. Further, for each of our algorithms, we conducted some effort to discuss important complexity aspects, to motivate our choice of parameters, to circumvent  potential numerical issues, and to give detailed correctness proofs. In particular we prove in detail that our constructions involving $v$-adic elliptic logarithms have the required properties at the problematic places $v$ of $\QQ$, that is $v=\infty$, $v=2$ and bad reduction $v$.

\paragraph{Improvements.} The elliptic logarithm sieve improves in all aspects the known elliptic logarithm reduction and its subsequent enumeration. In Sections~\ref{sec:globalsieve} and \ref{sec:comparisonwithelr}, we shall demonstrate (in theory and in practice) that our improvements are substantial. In fact we obtain running time improvements by a factor which is exponential in terms of the rank $r$, and which is exponential in terms of $\abs{S}$ when $\max\abs{a_i}$ is large.  Furthermore in the case of a generic Weierstrass equation \eqref{eq:weieq} we can efficiently determine all $S$-integral solutions for huge sets $S$. For example sets $S$ with $\abs{S}=10^5$ are usually no problem here.  Also, if $\abs{S}$ is very small then our sieve allows to deal efficiently with large ranks $r$ such as  $r=14,\dotsc,19$. In particular, in the case when $\mathcal O=\ZZ$, the sieve is practical even for huge ranks such as $r=28$.  The elliptic logarithm sieve considerably extends the class of elliptic Diophantine equations which can be solved  in practice.  We shall demonstrate this by solving several notoriously difficult Diophantine problems which appear to be completely out of reach for the known methods, see Section~\ref{sec:elllogsieveapp} for explicit examples.

\subsubsection{Applications}\label{sec:elllintroapp}

We solved large classes of elliptic Diophantine equations by applying our sieve. 
In particular, we efficiently solved several Diophantine problems in which the involved rank $r$ is large.  Further, for each globally minimal Weierstrass equation \eqref{eq:weieq} of any elliptic curve over $\QQ$ of conductor at most $100$ (resp. $1000$), we determined its set of $S$-integral solutions with $S$ given by the first $10^4$ (resp. $20$) primes. See Section~\ref{sec:elllogsieveapp} for more information. 

\paragraph{Conjecture and questions.} We used our data to motivate various questions on points of hyperbolic  curves $Y=(X,D)$ of genus one.  More precisely, let  $B$ be a nonempty open subscheme of $\sp(\ZZ)$ and let $X\to B$ be a smooth, proper and geometrically connected genus one curve.  Let $Y\hookrightarrow X$ be an open immersion onto the complement $X-D$ of a nontrivial relative Cartier divisor $D\subset X$ which is finite \'etale over $B$.  We now state the following conjecture involving the rank $r$ of the group formed by the $\QQ$-points of $\textnormal{Pic}^0(X_\QQ)$.

\vspace{0.3cm}
\noindent{\bf Conjecture.}
\emph{There are constants $c_Y$ and $c_r$, depending only on $Y$ and $r$ respectively, such that any nonempty finite set of rational primes $S$ with $T=\sp(\ZZ)-S$ satisfies} $$\abs{Y(T)}\leq c_Y \abs{S}^{c_r}.$$

\noindent Our initial motivation for making this conjecture is explained in Section~\ref{sec:ialgoheight}. Further, the above conjecture generalizes our conjecture for Mordell equations which we  discussed and motivated in Section~\ref{sec:malgoapplications}. In fact our  discussion and motivation given there, including the construction of our probabilistic model, can be applied in exactly the same way in the case when $Y$ is a Weierstrass curve.  Here $Y$ is a Weierstrass curve if the Cartier divisor $D$ is given by the image of a section of $X\to B$.  We shall also discuss and motivate various questions related to  the above conjecture. For example, we ask  whether the above conjecture holds  with $\abs{S}$ replaced by the logarithm of the largest prime in $S$?

\subsubsection{Organization of the section}

\paragraph{Plan.} In Section~\ref{sec:heightselllog} we discuss a  suitable rational approximation of the N\'eron--Tate height on $E(\QQ)$. The subsequent Sections~\ref{sec:archisieve} and \ref{sec:nonarchisieve} contain our construction of the local sieves at the archimedean place and the non-archimedean places. In Sections~\ref{sec:heightlogsieve} and \ref{sec:refinedenumell} we explain the height-logarithm sieve and the refined enumeration. Then we construct the refined sieve in Section~\ref{sec:refinedsieveell}. In Sections~\ref{sec:globalsieve} and \ref{sec:elllogsievealgo}, we put everything together to obtain the elliptic logarithm sieve. Here we also compare our sieve with the known approach. Then, after recalling in Section~\ref{sec:input} results and methods which allow to compute the required input data,  we discuss applications of our sieve in Section~\ref{sec:elllogsieveapp}. Finally,  we explain in Section~\ref{sec:compuaspects}  computational aspects of our constructions.

\paragraph{Notation.} Throughout this section we shall use the following  conventions.   By $\log$ we mean the principal value of the natural logarithm. Unless mentioned otherwise, $\abs{z}$ denotes the usual complex absolute value of $z\in\CC$ and the product taken over the empty set is~$1$.  Further  $\textnormal{lcm}(a_1,\dotsc,a_n)$ denotes the least common multiple of $a_1,\dotsc,a_n\in\ZZ$. For any real number $x\in\mathbb R$, we write $\floor{x}=\max(n\in\ZZ\st n\leq x)$ and $\ceil{x}=\min(n\in\ZZ\st n\geq x)$.  We denote by $h(\alpha)$ the usual absolute logarithmic Weil height of $\alpha\in\QQ$, with $h(0)=0$  and $h(\alpha)=\log\max(\abs{m},\abs{n})$ if $\alpha=m/n$ for coprime $m,n\in\ZZ$. If $\alpha\in \QQ$ is nonzero and if $p$ is a rational prime, then we write $\ord_p(\alpha)\in\ZZ$ for the order of $p$ in~$\alpha$. For any set $M$, we denote by $\lvert M\rvert$ the (possibly infinite) number of distinct elements of~$M$. Finally, for any $n\in\ZZ_{\geq 1}$, we say that $\mathcal E\subset \mathbb R^n$ is an ellipsoid if $\mathcal E=\{x\in\mathbb R^n\st q(x)\leq c\}$ for some positive definite quadratic form $q:\RR^n\to \RR$ and some positive real number $c$.

\paragraph{Computer, software and algorithms.} We implemented all our algorithms in Sage. Here a significant part of our program code is devoted to assure that the numerical aspects of the algorithms are all correct, see Section~\ref{sec:compuaspects} for certain important numerical aspects.   We point out that we shall use various  functions of Sage~\cite{sage:sagesystem} which in fact are direct applications of the corresponding functions of Pari~\cite{pari:parisystem}. Further, to compute the Mordell--Weil bases required for the input of the elliptic logarithm sieve, we used the techniques implemented in the computer packages Pari, Sage and Magma~\cite{magma:magmasystem}. For all our computations, we used a standard personal working computer at the MPIM Bonn. We shall list the running times of our algorithms for many  examples. In fact the listed times are always upper bounds and some of them were obtained by using older versions of our algorithms. Here in many cases the running times would now be significantly better when using the most recent versions (as of February 2016) of our algorithms.

\paragraph{Acknowledgements.} Our construction of the  elliptic logarithm sieve crucially builds on  ideas and techniques of the authors  who developed, generalized and/or refined the elliptic logarithm reduction  over the last 30 years (see \cite{tzanakis:book} for an overview).

\subsection{Heights}\label{sec:heightselllog}

In this section we first recall useful results for heights of rational points on our given elliptic curve $E$ over $\QQ$ with Weierstrass equation \eqref{eq:weieq}. Then we construct a suitable rational approximation of the canonical N\'eron--Tate height on $E(\QQ)$, and we fix some terminology.

\paragraph{Canonical height.}Let $\hat{h}$ be the canonical N\'eron--Tate height on $E(\QQ)$. Here we use the natural normalization which divides by the degree of the involved rational function, see for example \cite[p.248]{silverman:aoes}. For any nonzero $P\in E(\QQ)$, it is known that the logarithmic Weil height $h(x)$ of the corresponding solution $(x,y)$ of \eqref{eq:weieq} can be explicitly compared with $\hat{h}(P)$. For instance  Silverman~\cite[Thm 1.1]{silverman:heightcomparison} used an approach of Lang to obtain an explicit constant $\mu(E)$, depending only on the coefficients $a_i$ of \eqref{eq:weieq}, such that 
\begin{equation}\label{eq:nthcompa}
-\tfrac{1}{24}h(j_E)-\mu\leq \hat{h}(P)-\tfrac{1}{2}h(x)\leq \mu
\end{equation}
for $\mu=\mu(E)+1.07$ and $j_E$ the $j$-invariant of $E$. Further it is known that $\hat{h}$ defines a positive semi-definite quadratic form on the geometric points of $E$, and for any point $P\in E(\QQ)$ it holds $\hat{h}(P)=0$ if and only if $P$ lies in the torsion subgroup $E(\QQ)_{\textnormal{tor}}$ of $E(\QQ)$. Therefore on identifying  the real vector space $E(\QQ)\otimes_\ZZ \RR$ with $\RR^r$ via our given basis $P_1,\dotsc,P_r$ of the free part of $E(\QQ)$, we obtain that $\hat{h}$ extends to a positive definite quadratic form on $\RR^r$. 
Let $\lambda$ be the smallest eigenvalue of the matrix $(\hat{h}_{ij})$ in $\RR^{r\times r}$ defining the bilinear form associated to $\hat{h}$. Linear algebra gives that any point $P$ in $E(\QQ)$  satisfies 
\begin{equation}\label{eq:nthlowerbound}
\lambda \|P\|^2_\infty\leq \hat{h}(P)\leq r\lambda'\|P\|^2_\infty 
\end{equation}
for $\lambda'$ the largest eigenvalue of $(\hat{h}_{ij})$. Here the infinity norm is defined by $\|P\|_\infty=\max \abs{n_i}$, where $P=Q+\sum n_i P_i$ with $n_i\in\ZZ$ and $Q\in E(\QQ)_{\textnormal{tor}}$. 
 We point out that in practice it is  always possible to quickly determine the points in $E(\QQ)_{\textnormal{tor}}$. For what follows we therefore always may assume that the rank $r\geq 1$.
To avoid numerical problems with the real valued function $\hat{h}$, we shall work with a rational approximation of $\hat{h}$.

\paragraph{Rational approximation.} 
We next explain our construction of a suitable rational approximation of $\hat{h}$. Let $k\in\ZZ_{\geq 1}$ and  define the norm $\|\hat{h}_{ij}\|$ of $(\hat{h}_{ij})$ by $\|\hat{h}_{ij}\|=\max \abs{\hat{h}_{ij}}$. On using continued fractions we obtain $f\in\QQ$ which approximates the real number $2^{k}/\|\hat{h}_{ij}\|$ up to any required precision, see Section~\ref{sec:compuaspects}. We identify the rational vector space $E(\QQ)\otimes_\ZZ \QQ$ with $\QQ^r$ via the basis $P_1,\dotsc,P_r$ and we consider the quadratic form
\begin{equation}
\hat{h}_k:E(\QQ)\otimes_\ZZ\QQ\to \QQ
\end{equation}
 associated to $\tfrac{1}{f}([f\hat{h}]-r\cdot \id)\in \QQ^{r\times r}$. Here $[f\hat{h}]$ denotes the symmetric matrix in $\ZZ^{r\times r}$ with $ij$-th entry given by $[f\hat{h}_{ij}]$ for $[\cdot]$ the rounding ``function" defined in Section~\ref{sec:compuaspects}. The following lemma  compares $\hat{h}$ with the natural extension of $\hat{h}_k$ to $E(\QQ)\otimes_\ZZ\RR=\RR^r$.

\begin{lemma}\label{lem:hk}
If $\|\cdot\|_2$ denotes the euclidean norm on $\RR^r$, then any $x\in \RR^r$ satisfies $$\hat{h}(x)-\tfrac{2r}{f}\|x\|_2^2\leq \hat{h}_k(x)\leq \hat{h}(x).$$
\end{lemma}
\begin{proof}
To simplify notation we write $V=\RR^r$ and we denote by $q$ the quadratic form on $V$ which is associated to $(\delta_{ij})=(\hat{h}_{ij})-\tfrac{1}{f}[f\hat{h}]$. We take $x\in V$ and we deduce
$$\hat{h}_k(x)=\hat{h}(x)-q(x)-\tfrac{r}{f}\|x\|_2^2.$$
It holds that $f\abs{\delta_{ij}}\leq 1$ and the Cauchy--Schwarz inequality implies that $\sum_{ij}\abs{x_ix_j}\leq r\|x\|_2^2$. Therefore we obtain  $\abs{q(x)}\leq \tfrac{r}{f}\|x\|_2^2$ and then we see that the displayed formula leads to the claimed inequality. This completes the proof of the lemma. 
\end{proof}
If $k$ is sufficiently large then the above lemma implies  that the quadratic form $\hat{h}_k$ is positive definite and is close to $\hat{h}$. 
For what follows we fix an element $k\in\ZZ$ such that $\hat{h}_k$ is positive definite and is close to $\hat{h}$, see also the discussions in Section~\ref{sec:compuaspects}.  

\paragraph{Terminology.} To introduce some terminology, we take $\sigma\in \RR_{>0}$ and we consider a place $v$ of $\QQ$. The $v$-adic elliptic logarithm is of local nature, while $\hat{h}$ and $\hat{h}_k$ are global height functions. In our global sieve we shall need to measure the ``weight" of the $v$-adic norm of the $v$-adic elliptic logarithm inside $\hat{h}_k$. For this purpose, we shall work with the set
\begin{equation}\label{def:sigmavsigma}
\Sigma(v,\sigma) 
\end{equation}
formed by the nonzero points $P\in E(\QQ)$ whose corresponding solution $(x,y)$ of \eqref{eq:weieq} satisfies $\tfrac{1}{2}\log \abs{x}_v\geq \tfrac{1}{\sigma}(\hat{h}_k(P)-\mu).$
Here we write $\abs{x}_v=p^{-\ord_p(x)}$ if $v$ is a finite place given by the rational prime $p$, and if $v=\infty$ then  $\abs{x}_v$ is defined by $\abs{x}_v=\abs{x}$ for $\abs{\cdot}$ the usual complex absolute value. 
We note that one can define the set $\Sigma(v,\sigma)$ more intrinsically using (Arakelov) intersection theory. However, it is not clear to us if this provides a significant advantage in practice and thus we work with \eqref{def:sigmavsigma} using \eqref{eq:nthcompa}. For any rational integer $n\geq r$, we say that $P\in E(\QQ)$ is determined modulo torsion by  $\gamma\in \ZZ^n$ if there exists $Q\in E(\QQ)_{\textnormal{tor}}$ such that $P=Q+\sum \gamma_i P_i$. Further we denote by $\Gamma_E=\ZZ^r$ the lattice inside $\RR^r$ given by the image of $E(\QQ)$ inside $E(\QQ)\otimes_\ZZ\RR=\RR^r$ using the identification via $P_i$.

\subsection{Archimedean sieve}\label{sec:archisieve}

 Building on ideas of Zagier~\cite{zagier:largeintegralpoints}, de Weger~\cite{deweger:phdthesis}, Stroeker--Tzanakis~\cite{sttz:elllogaa} and Gebel--Peth{\H{o}}--Zimmer~\cite{gepezi:ellintpoints}, we construct in this section our archimedean sieve. We shall use this sieve to improve inter alia the known reduction process at infinity of the elliptic logarithm method, see the discussions in Sections~\ref{sec:globalsieve} and \ref{sec:comparisonwithelr}.  

Throughout this Section~\ref{sec:archisieve} we use the setup of Section~\ref{sec:setup} and we continue the notation introduced above. Further throughout this section we write $\abs{\cdot}=\abs{\cdot}_\infty$.

\paragraph{Real elliptic logarithm.}We shall work with the following normalization of the elliptic logarithm on the identity component $E^0(\RR)$ of the real Lie group $E(\RR)$. 
First we recall that the uniformization theorem for complex elliptic curves gives a lattice $\Lambda=\omega_1\ZZ+\omega_2\ZZ$ inside $\CC$ with $\omega_1\in\RR_{>0}$ and an isomorphism  $\CC/\Lambda\xrightarrow{\sim} E(\CC)$ whose inverse we denote by $$\log: E(\CC)\xrightarrow{\sim} \CC/\Lambda.$$ 
To describe more explicitly the restriction to $E^0(\RR)$ of the displayed morphism,  we write  $x=x'-\tfrac{1}{12}b_2$ with $b_2=a_1^2+4a_2$  and $y=\tfrac{1}{2}(y'-a_1x-a_3)$ and we transform \eqref{eq:weieq} into the Weierstrass equation $y'^2=4x'^3-g_2x'-g_3$ whose complex solutions we identify with the nonzero points of $E(\CC)$.
We may and do assume
that $g_i=g_i(\Lambda)$ is associated to $\Lambda$ as in \cite[p.169]{silverman:aoes} and then the isomorphism $\CC/\Lambda\xrightarrow{\sim} E(\CC)$ is given outside zero by $z\mapsto (\wp(z),\wp'(z))$ for $\wp=\wp(\Lambda)$  the Weierstrass $\wp$-function and  $\wp'$ its derivative.  Hence 
we deduce for example from \cite[p.174]{silverman:aoes} that the restriction of $\log: E(\CC)\xrightarrow{\sim} \CC/\Lambda$ to $E^0(\RR)$ is of the form $E^0(\RR)\xrightarrow{\sim} \RR/(\omega_1\ZZ)$, 
 which in turn induces a bijective map
\begin{equation}\label{def:reallog}
\log: E^0(\RR)\to \{z\in\RR\st 0\leq z<\omega_1\}.
\end{equation}
Explicitly if $P\in E^0(\RR)$ corresponds to a real solution  $(x,y)$ of \eqref{eq:weieq} then it holds that $\log(P)=\tfrac{y'}{\abs{y'}} \int^{x'}_{\infty}\tfrac{dz}{f(z)^{1/2}}$ mod $(\omega_1\ZZ)$ for $f(z)=4z^3-g_2z-g_3$. 
Here one can compute the real number $\log(P)$ up to any required precision, see for example Zagier~\cite[p.430]{zagier:largeintegralpoints}. 
Further, we denote by $e_t$ the exponent of the finite group  $E(\QQ)_{\textnormal{tor}}$ and we define 
\begin{equation}\label{def:marchi}
m=\textnormal{lcm}(e_t,
\iota )
\end{equation}
for  $\iota $  the index of $E^0(\RR)$ inside $E(\RR)$.  
It holds that $\iota \in\{1,2\}$, 
since $E(\RR)$ is either connected  or isomorphic to $E^0(\RR)\times (\ZZ/2\ZZ)$. 
We recall that the points $P_1,\dotsc,P_r$ form a basis of the free part of $E(\QQ)$. Now any  $P=Q+\sum n_i P_i$ in $E(\QQ)$, with $Q\in E(\QQ)_{\textnormal{tor}}$ and $n_i\in\ZZ$, satisfies $mP=\sum n_i (mP_i)\in E^0(\RR)$.
Next we take $\kappa\in\ZZ_{\geq 1}$ and we define 
\begin{equation}\label{def:x0}
x_0(\kappa)=(\kappa+1)(\abs{b_2}/12+\max\abs{\xi_i})
\end{equation}
 for $\{\xi_i\}$ the set of roots of $f(z)=4z^3-g_2z-g_3$. 
If  $P\in E^0(\RR)$ corresponds to a solution $(x,y)$ of \eqref{eq:weieq}, then the next lemma allows  to control $\log(P)$ in terms of $\abs{x}$.
\begin{lemma}\label{lem:archiest}
The following  statements hold. 
\begin{itemize}
\item[(i)] Suppose that $P\in E(\RR)$ corresponds to a real solution $(x,y)$ of \eqref{eq:weieq} with $\abs{x}\geq x_0(\kappa)$. Then $P$ lies in $E^{0}(\RR)$ and there is $\epsilon\in\{0,-1\}$ such that any $n\in\ZZ$ satisfies 
\begin{equation*}
\abs{n\log(P)+n\epsilon\omega_1}\leq \abs{n}\bigl(1+\tfrac{1}{\kappa}\bigl)^2\abs{x}^{-1/2}.
\end{equation*}
\item[(ii)]If $P=Q+\sum n_i P_i$ lies in $E^0(\RR)$ with $Q\in E(\QQ)_{\textnormal{tor}}$ and $n_i\in\ZZ$, then there exists $l\in\ZZ$ with $\abs{l}\leq m+ \sum \abs{n_i}$ such that $m\log(P)=\sum n_i\log(mP_i)+l\omega_1$.
\end{itemize}
\end{lemma}
\begin{proof}
We first prove assertion (i). Our assumption implies that $x'=x+\tfrac{1}{12}b_2$ is positive and that $x'$ strictly exceeds the largest real root of $f(z)=4z^3-g_2z-g_3$. 
Hence we conclude that $P\in E^0(\RR)$.
To verify the second statement of (i) we observe that any $z\in\RR$ with $z\geq (\kappa+1)\max\abs{\xi_i}$ satisfies $f(z)\geq 4\bigl(\tfrac{\kappa}{\kappa+1}\bigl)^3z^3$. It follows that $\abs{\int^{x'}_\infty \tfrac{dz}{f(z)^{1/2}}}^2$ is at most $\bigl(\tfrac{\kappa+1}{\kappa}\bigl)^{3}\abs{x'}^{-1}$, since our assumption gives  $\abs{x'}\geq (\kappa+1)\max\abs{\xi_i}$. Furthermore our assumption 
provides that $\abs{x'}\geq \tfrac{\kappa}{\kappa+1}\abs{x}$, and then we see that there exists $\epsilon\in\{0,-1\}$ such that the claimed inequality holds for $n=1$ and thus for all $n\in\ZZ$.

It remains to show (ii). The points $mP_i$ are all in $E^0(\RR)$ since $\iota $ divides $m$, and the point $P$ is in  $E^0(\RR)$ by assumption. Thus, on exploiting that the real elliptic logarithm is induced by a group isomorphism $E^0(\RR)\xrightarrow{\sim} \RR/(\omega_1\ZZ)$, we find $l',l''\in\ZZ$ with $m\log(P)=\log(mP)+l'\omega_1$ and $\log(mP)=\sum n_i\log(mP_i)+l''\omega_1$.
Then on using that $\log(P)$, $\log(mP)$ and $\log(mP_i)$ are in the interval $[0,\omega_1[$, we deduce that $\abs{l'}\leq m$ and $\abs{l''}\leq \sum \abs{n_i}$. 
Hence the integer  $l=l'+l''$ has the desired property. This completes the proof of the lemma.
\end{proof}

\paragraph{Construction of $\Gamma$ and $\mathcal E$.}

Let $\sigma>0$ be a real number,  let $\mu$ be as in \eqref{eq:nthcompa} and write $\Sigma$ for the set $\Sigma(\infty,\sigma)$ defined in \eqref{def:sigmavsigma}. Suppose that $M',M\in\ZZ$ with $\mu\leq M'<M$ and let $\kappa\in\ZZ_{\geq 1}$. We would like to construct a lattice $\Gamma\subset \ZZ^{r+1}$ and an ellipsoid $\mathcal E\subset \RR^{r+1}$ such that any  $P\in \Sigma$ with $M'<\hat{h}_k(P)\leq M$ is determined modulo torsion by a  point in $\Gamma\cap\mathcal E$. The following construction depends on a suitable choice of a parameter $c\in \ZZ_{\geq 1}$, which we shall explain below \eqref{archivolumecomp}. We write $\alpha_i=\log(mP_i)$ for $i\in\{1,\dotsc,r\}$ and we denote by $$\Gamma\subset\ZZ^{r+1}$$ the lattice formed by the elements $\gamma\in\ZZ^{r+1}$ such that $\gamma_{r+1}=l[c\omega_1]+\sum \gamma_i [c\alpha_i]$ for some $l\in\ZZ$.   
Next we choose a positive number $\delta\in\QQ$ as explained in the discussion surrounding \eqref{def:delta12} and  we denote by $q$ the positive definite quadratic form on $\RR^{r+1}$ which is  given by $q(z)=\hat{h}_k(z_1,\dotsc,z_r)+(M/\delta^2)z_{r+1}^2$ for any $z\in \RR^{r+1}$. Now we define the ellipsoid
$$\mathcal E=\{z\in\mathbb R^{r+1}\st q(z)\leq 2M\}.$$
Let $x_0=x_0(\kappa)$ be as in \eqref{def:x0} and let $\Sigma(x_0)$ be the set of points $P$ in $\Sigma$ with $\abs{x}>x_0$, where $(x,y)$ is the solution of \eqref{eq:weieq} corresponding to $P$. We obtain the following lemma.

\begin{lemma}\label{lem:archicov}
Suppose that $P\in \Sigma(x_0)$ satisfies $M'<\hat{h}_k(P)\leq M$. Then the point $P$ is determined modulo torsion by some lattice point $\gamma$ in $\Gamma\cap \mathcal E$.
\end{lemma}

\begin{proof}
 Let $(x,y)$ be the  solution of \eqref{eq:weieq} corresponding to $P$, and write $P=Q+\sum n_i P_i$ with $Q\in E(\QQ)_{\textnormal{tor}}$ and $n_i\in\ZZ$.  
It holds that $\abs{x}\geq x_0$ since $P\in \Sigma(x_0)$ and therefore Lemma~\ref{lem:archiest}~(i) shows that $P\in E^0(\RR)$. Hence we see that Lemma~\ref{lem:archiest}~(ii) gives $l_0\in\ZZ$ with $\abs{l_0}\leq m+\sum\abs{n_i}$ such that $m\log(P)=\sum n_i\alpha_i+l_0\omega_1$. On inserting this into the inequality in Lemma~\ref{lem:archiest}~(i) with $n=m$, we obtain $l\in\ZZ$ with $\abs{l}\leq 2m+\sum\abs{n_i}$ 
such that
\begin{equation}\label{eq:archifundineq}
\abs{\sum n_i\alpha_i+l\omega_1}\leq m\bigl(1+\tfrac{1}{\kappa}\bigl)^2e^{-\tfrac{1}{\sigma}(\hat{h}_k(P)-\mu)}.
\end{equation}
Here we used our assumption $P\in \Sigma(x_0)$, which  provides that $\tfrac{1}{2}\log\abs{x}\geq\tfrac{1}{\sigma}(\hat{h}_k(P)-\mu)$. 
Next we define $d=l[c\omega_1]+\sum n_i [c\alpha_i]$ and we observe that  $\gamma=((n_i),d)\in \ZZ^{r+1}$ lies in our lattice $\Gamma$. To show that $\gamma$  lies in addition in the ellipsoid $\mathcal E$, we use the mean inequality and linear algebra in order to deduce that $\lambda_k(\sum \abs{n_i})^2\leq r\hat{h}_k(P)$ for  $\lambda_k\in \QQ$ the smallest eigenvalue of the positive definite quadratic form $\hat{h}_k$ on $E(\QQ)\otimes_\ZZ\QQ$. Then we see that \eqref{eq:archifundineq} together with our assumption $M'<\hat{h}_k(P)\leq M$ implies that  $\abs{d}\leq \delta$. Here  $\delta\in \QQ$ is chosen such that $\delta$ has ``small" height in the sense  of Section~\ref{sec:compuaspects} and such that  $\delta\geq \delta_1+\delta_2$ for 
\begin{equation}\label{def:delta12}
\delta_1=2m+2\bigl(r\tfrac{M}{\lambda_k}\bigl)^{1/2} \ \ \ \textnormal{ and } \ \ \ \delta_2=cm\bigl(1+\tfrac{1}{\kappa}\bigl)^2e^{-\tfrac{1}{\sigma}(M'-\mu)}.
\end{equation}
On using again that $\hat{h}_k(P)\leq M$ we obtain $q(\gamma)\leq M+(M/\delta^2) d^2$. This together with  $\abs{d}\leq \delta$ implies that $\gamma\in\mathcal E$ and thus $P$ is determined modulo torsion by $\gamma\in \Gamma\cap \mathcal E$. 
\end{proof}

This lemma provides a sieve for the points $P$ in $\Sigma(x_0)$ with $M'<\hat{h}_k(P)\leq M$. In the following paragraph we discuss the strength of the sieve depending on various parameters.

\paragraph{Strength of the sieve.}  To make the sieve  as efficient as possible, we would like choose the parameter $c$ such that the intersection $\Gamma\cap \mathcal E$ does not contain many points.  In the generic case the cardinality of $\Gamma\cap \mathcal E$ can be approximated (for large $M$) by the euclidean volume of the ellipsoid $\mathcal E_\psi=\{z\in\RR^{r+1}\st q_\psi(z)\leq 2M\}$ inside $\RR^{r+1}$. Here $q_\psi$ denotes the positive definite quadratic form obtained by pulling back $q$ along the linear transformation $\psi$ of $\RR^{r+1}$  which satisfies $\psi (\ZZ^{r+1})=\Gamma$ and which is explicitly given by 
\[ \begin{pmatrix}
1 & & & 0 \\
 & \ddots & & \vdots \\
 & & 1 & 0 \\
[c\alpha_1] & \cdots & [c\alpha_r] & [c\omega_{1}]
\end{pmatrix}.\] 
To compute the euclidean volume $\vol(\mathcal E_\psi)$ of  $\mathcal E_\psi$, we let $R_{E}=2^r\det(\hat{h}_{ij})$ be the regulator of $E(\QQ)$ normalized as in \cite[p.253]{silverman:aoes} and we denote by $V_{r+1}$ the euclidean volume of the unit ball in $\RR^{r+1}$. Then the volume $\vol(\mathcal E_\psi)$ is approximately 
\begin{equation}\label{archivolumecomp}
u\cdot M^{r/2}\tfrac{(\delta_1+\delta_2)}{c}, \ \ \ u=\tfrac{2^{r+1/2}V_{r+1}}{\omega_1R_{E}^{1/2}}.
\end{equation}
We note that $\delta_2/c$ does not depend on $c$. Hence in view of \eqref{archivolumecomp}  we choose $c$ such that $u\cdot M^{r/2}\tfrac{\delta_1}{c}$ is smaller than $u\cdot M^{r/2}\tfrac{\delta_2}{c}$. For example $c$ should always dominate $M^{(r+1)/2}$ if $M$ is large. We next discuss the dependence of the sieve on $M'$ and $M$. For some large $M$, we choose $c$  as indicated above and we assume for a moment that $M'$ dominates 
$
\frac{\sigma r}{2}\log M.
$
Then it follows that $M^{r/2}\tfrac{\delta_2}{c}$ is close to zero and hence  \eqref{archivolumecomp} implies that the volume of $\mathcal E_\psi$ is small. In the generic case this assures that $\Gamma\cap \mathcal E$ has very little points or is even trivial. In particular, we see that the archimedean sieve is very efficient for such $M'$ and $M$. On the other hand, for small $M'$ our sieve is not that efficient in view of \eqref{archivolumecomp}.

\begin{remark}\label{rem:optellinf}
One can replace $\mathcal E$ by the more balanced ellipsoid $\mathcal E^*\subset \RR^{r+1}$ of the form $\mathcal E^*=\{z\in \RR^{r+1}\st q^*(z)\leq M\}$ for $q^*(z)=\tfrac{r}{r+1}\hat{h}_k(z_1,\dotsc,z_r)+\tfrac{1}{r+1}(M/\delta^2) z_{r+1}^2$. Indeed this follows by observing that $\gamma$ appearing in the proof of Lemma~\ref{lem:archicov} satisfies $q^*(\gamma)\leq M$. 
\end{remark}

\paragraph{Archimedean sieve.} In the following  sieve, we use the version (FP) of the Fincke--Pohst algorithm  described in Section~\ref{sec:compuaspects} in order to determine all points in $\Gamma\cap \mathcal E$.

\begin{Algorithm}[Archimedean sieve]\label{algo:archisieve} The inputs are $\kappa\in\ZZ_{\geq 1}$ and $M',M\in\ZZ$ with $\mu\leq M'< M$. The output is the set of points $P\in\Sigma(x_0)$ with $M'<\hat{h}_k(P)\leq M$. 
\begin{itemize}
\item[(i)] First choose the parameter $c\in \ZZ_{\geq 1}$ as explained in the discussion surrounding \eqref{archivolumecomp}. Then compute the lattice $\Gamma\subset \ZZ^{r+1}$,  by determining the period $\omega_1$ and the real elliptic logarithms $\alpha_i=\log(mP_i)$ up to the required precision for all $i\in\{1,\dotsc,r\}$.
\item[(ii)]  Determine $\Gamma\cap \mathcal E$ by using the version of the Fincke--Pohst algorithm in \textnormal{(FP)}. 
\item[(iii)] For each lattice point $\gamma$ in $\Gamma\cap \mathcal E$ and for each torsion point $Q$ in $E(\QQ)_{\textnormal{tor}}$, output the point $P=Q+\sum \gamma_i P_i$ if   $M'<\hat{h}_k(P)\leq M$ and if $P$ is in $\Sigma(x_0)$.
\end{itemize} 
\end{Algorithm}

\paragraph{Correctness.} Suppose that $P\in\Sigma(x_0)$ satisfies $M'<\hat{h}_k(P)\leq M$. Lemma~\ref{lem:archicov} gives that $P$ is determined modulo torsion by some $\gamma\in \Gamma\cap\mathcal E$. In other words, there is $Q\in E(\QQ)_{\textnormal{tor}}$ such that $P=Q+\sum \gamma_i P_i$ and hence step (iii) produces our point $P$ as desired.

\paragraph{Complexity.} We now discuss aspects of Algorithm~\ref{algo:archisieve} which significantly influence the running time. In step (i) the running time of the computation of the lattice $\Gamma=\psi(\ZZ^{r+1})$  crucially depends on the size of $c$. For example if $c$ is approximately $M^{(r+1)/2}$ then we need to know the real logarithms $\omega_1$ and $\log(mP_i)$  up to a number of decimal digits which is approximately $\tfrac{r+1}{2}\log_{10} M$, where $\log_b z=(\log z)/\log b$ for $z,b\in\RR_{>0}$. We shall apply the algorithm with  huge parameters $M$. Therefore we need to compute $(r+1)$ real elliptic logarithms up to a very high precision and this can take a long time.  Step (ii) is essentially always fast in practice. The reason is that the involved Mordell--Weil rank $r$ of $E(\QQ)$ is usually not that large and hence the application of (FP) with the lattice $\Gamma$ of rank $r+1$ is fast. Finally step (iii) needs to compute in particular the coordinates of certain points in $E(\QQ)$ and this can take some time if $\hat{h}(P)$ and $r$ are not small. 

\paragraph{Comparison.} There are important differences between our approach and the known approach. In particular we work with the N\'eron--Tate height $\hat{h}$, while all other authors use the inequality $\lambda\|\cdot\|_\infty^2\leq \hat{h}(\cdot)$ to work with the norm $\|\cdot\|_\infty$. Also we actually determine the intersection $\Gamma\cap \mathcal E$, while the known approach computes a lower bound for the length of the shortest nonzero vector in $\Gamma$ in order to rule out non-trivial points in $\Gamma\cap\mathcal E$. Other, more technical, differences are the following:  The parameter $\kappa$ allows us (up to a certain extent) to adapt  the strength of the sieve to the given situation, and the construction of our ellipsoid $\mathcal E^*$ involving the weights $\tfrac{r}{r+1}$ and $\tfrac{1}{r+1}$ is more balanced in particular for large $r$.  In Sections~\ref{sec:globalsieve} and \ref{sec:comparisonwithelr}, we shall further compare the two approaches and we shall explain in detail the improvements provided by our new ideas.

\subsection{Non-archimedean sieve}\label{sec:nonarchisieve} 

Building on ideas of Smart~\cite{smart:sintegralpoints}, Peth{\H{o}}--Zimmer--Gebel--Herrmann~\cite{pezigehe:sintegralpoints} and Tzanakis~\cite{tzanakis:book}, we construct in this section the non-archimedean sieve. We shall use this sieve to improve inter alia the known reduction process of the elliptic logarithm method at non-archimedean primes, see the discussions in Sections~\ref{sec:globalsieve} and \ref{sec:comparisonwithelr}. 

Throughout this Section~\ref{sec:nonarchisieve} we work with the setup of Section~\ref{sec:setup} and we continue the notation introduced above. Further we fix $p$ in $S$ and we assume that the Weierstrass model \eqref{eq:weieq} of our given elliptic curve $E$ is minimal at $p$. To simplify the notation of this section, we write  $v(\cdot)=\ord_p(\cdot)$ and $\abs{\cdot}=\abs{\cdot}_p$ with $\abs{x}_p=p^{-v(x)}$ for $x\in \QQ_p$.

\paragraph{The $p$-adic elliptic logarithm.} 
Let $E_1(\QQ_p)$ be the subgroup of $E(\QQ_p)$ formed by the points $P$ in $E(\QQ_p)$ with $\pi(P)=0$ for  $\pi:E(\QQ_p)\to E(\mathbb F_p)$ the reduction map. Here $E(\mathbb F_p)$ denotes the set of $\mathbb F_p$-points of the special fiber of the projective closure\footnote{Here we mean $\textnormal{Proj}\bigl(\ZZ_p[x,y,z]/(f)\bigl)$ for $f=y^2z+a_1xyz+a_3yz^2-(x^3+a_2x^2z+a_4xz^2+a_6z^3)$.} of $\eqref{eq:weieq}$ inside $\mathbb P^2_{\ZZ_p}$. Let $\hat{E}$ be the formal group over $\ZZ_p$ associated to $E_{\QQ_p}$. There is an isomorphism $E_1(\QQ_p)\xrightarrow{\sim} \hat{E}(p\ZZ_p)$ of abelian groups, which is given away from zero by $(x,y)\mapsto -\tfrac{x}{y}$.  
Composing this isomorphism  with the formal logarithm of $\hat{E}$ induces a morphism  
\begin{equation}\label{def:plog}
\log:E_1(\QQ_p)\to \mathbb G_a(\QQ_p)
\end{equation}
of abelian groups. 
We call the displayed morphism the $p$-adic elliptic logarithm. Explicitly if $P\in E_1(\QQ_p)$ is nonzero and corresponds to the solution $(x,y)$ of \eqref{eq:weieq}, then it holds that $\log(P)=z+\sum_{n\geq 2} \tfrac{b_{n}}{n} z^{n}$ with $z=-\tfrac{x}{y}$ and $b_n\in\ZZ_p$. 
A priori the $p$-adic elliptic logarithm is only defined on the subgroup $E_1(\QQ_p)$ of $E(\QQ_p)$. One can somehow circumvent this problem by multiplying the points in $E(\QQ_p)$ with a suitable integer. To construct such an integer, let $E_{\textnormal{ns}}(\mathbb F_p)$ be the group formed by the nonsingular points in $E(\mathbb F_p)$ and consider the subgroup $E_0(\QQ_p)=\pi^{-1}(E_{\textnormal{ns}}(\mathbb F_p))$ of $E(\QQ_p)$. 
We denote by $\iota$ the index of $E_0(\QQ_p)$ in $E(\QQ_p)$, and we write $e_t$ and $e_{ns}$ for the exponents of the finite groups $E(\QQ)_{\textnormal{tor}}$ and $E_{\textnormal{ns}}(\mathbb F_p)$ respectively.  The short exact sequence $0\to E_1(\QQ_p)\to E_0(\QQ_p)\overset{\pi}{\to} E_{\textnormal{ns}}(\mathbb F_p)\to 0$ of abelian groups shows that  $(\iota e_{ns})P\in E_1(\QQ_p)$ for all $P\in E(\QQ_p)$. We now define 
\begin{equation}\label{def:mnonarchi}
m=\textnormal{lcm}\bigl(e_t,\iota e_{ns}\bigl).
\end{equation}
Recall that $P_1,\dotsc,P_r$ denotes our given basis of the free part of $E(\QQ)$. Any  $P=Q+\sum n_i P_i$ in $E(\QQ)$, with $Q\in E(\QQ)_{\textnormal{tor}}$ and $n_i\in \ZZ$, satisfies $mP=\sum n_i (mP_i)\in E_1(\QQ_p)$.  
The case distinction in the following lemma takes into account that in general the formal logarithm of $\hat{E}$ is not necessarily an isomorphism of formal groups over the given base. 
\begin{lemma}\label{lem:nonarchiest}
Let $P\in E(\QQ_p)$ be nonzero, and suppose that $(x,y)$ is the solution of \eqref{eq:weieq} corresponding to $P$. Then the following two statements hold.  
\begin{itemize}
\item[(i)] Assume that $p\geq 3$. If $P\notin E_1(\QQ_p)$ with $mP\neq 0$ then $\abs{\log(mP)}^2<\abs{x}^{-1}$, and if $P\in E_1(\QQ_p)$ then  $\abs{\log(nP)}^2=\abs{n}^2\abs{x}^{-1}$ for all $n\in\ZZ$.
\item[(ii)] If $p=2$ and $v(x)<-2$, then any $n\in\ZZ$ satisfies $\abs{\log(nP)}^2=\abs{n}^2\abs{x}^{-1}$.
\end{itemize}
\end{lemma}
\begin{proof}
Our proof given below relies on the classical result that the formal logarithm is compatible with the valuation $v$ in the following sense: 
For any $l\in\ZZ$ with $l>v(p)/(p-1)$, the restriction of the formal logarithm of $\hat{E}$ induces an isomorphism 
\begin{equation}\label{eq:formalgpiso}
\hat{E}((p\ZZ_p)^l)\cong (p\ZZ_p)^l 
\end{equation}
of abelian groups.
Further, we shall use below that if $P\in E_1(\QQ_p)$  then it holds that $3v(x)=2v(y)$, thus $v(x)$ is even and the number $z=-x/y$ satisfies $2v(z)=-v(x)$. 
If $mP\neq 0$ then we denote by $(x_m,y_m)$ the solution of \eqref{eq:weieq}  corresponding to  $mP$.

To prove (i)  we may and do assume that $p\geq 3$.  
Then $p$ satisfies $1>v(p)/(p-1)$ and hence the isomorphism in \eqref{eq:formalgpiso} exists for all $l\geq 1$. This implies that $2v(\log(mP))=-v(x_m)$ since $mP\in E_1(\QQ_p)$ is nonzero by assumption. 
If $P$ is not in $E_1(\QQ_p)$ then $v(x)\geq 0$, and $mP\in E_1(\QQ_p)$ thus shows that $v(x)\geq 0> v(x_m)$. This  together with $2v(\log(mP))=-v(x_m)$ proves the claimed inequality if $P$ is not in $E_1(\QQ_p)$. Suppose now that $P\in E_1(\QQ_p)$. Then we obtain that $n\log(P)=\log(nP)$ for all $n\in \ZZ$ since the formal logarithm is a morphism of abelian groups, and the isomorphisms in \eqref{eq:formalgpiso} provide that $2v(\log(P))=-v(x)$. On combining these two equalities, we deduce the second statement of (i).

To show (ii) we may and do assume that $p=2$ and $v(x)<-2$. The latter assumption implies that $P\in E_1(\QQ_2)$ and $v(x)\leq -4$. We deduce that $v(z)\geq 2$ and hence  $z$ lies in $(2\ZZ_2)^2$. Further, the isomorphism in~\eqref{eq:formalgpiso} exists for all $l\geq 2$ since $2>v(2)/(2-1)$. Thus we obtain that $2v(\log(P))=-v(x)$ and then the equality $n\log(P)=\log(nP)$, which holds for all $n\in \ZZ$ since $P\in E_1(\QQ_2)$, implies (ii).  This completes the proof of the lemma.   
\end{proof}
We remark that the assumptions $P\in E_1(\QQ_p)$ and $v(x)<-2$, in (i) and (ii) respectively, assure in particular that the point $P$ has infinite order in $E(\QQ_p)$.

\paragraph{Construction of $\Gamma$ and $\mathcal E$.}  
Let $\sigma> 0$ be a real number and write $\Sigma$ for the set $\Sigma(v,\sigma)$ defined in \eqref{def:sigmavsigma}. Suppose that we are given $M',M\in\ZZ$ with $\mu\leq M'< M$ for $\mu$ as in \eqref{eq:nthcompa}. We would like to find a lattice $\Gamma\subset \ZZ^r$ such that any $P\in\Sigma$ with $M'<\hat{h}_k(P)\leq M$ is determined modulo torsion by some point in $\Gamma\cap\mathcal E$. Here $\mathcal E\subset \RR^{r}$ is the ellipsoid
$$
\mathcal E=\{z\in\RR^{r}\st \hat{h}_k(z)\leq M\}. 
$$
We identify $\alpha\in \ZZ_p$ with the corresponding element $(\alpha^{(1)},\alpha^{(2)},\dotsc)$ of the inverse limit  $\lim \ZZ/(p^n\ZZ)$ and we set $\alpha_i=\log(mP_i)$ for each $i\in\{1,\dotsc,r\}$. If  $p\geq 3$ then Lemma~\ref{lem:nonarchiest} implies that  $\alpha_i\in \ZZ_p$. To deal with the general case, we choose $i^*\in\{1,\dotsc,r\}$ with $v(\alpha_{i^*})= \min v(\alpha_i)$ and then $\beta_i=\alpha_i/p^{v(\alpha_{i^*})}$ lies in $\ZZ_p$. 
Now we denote by $$\Gamma\subset\ZZ^r$$ the lattice formed by the elements $\gamma\in\ZZ^r$ with $\sum \gamma_i\beta_i^{(c)}=0$ in $\ZZ/(p^c\ZZ)$, where $c\in \ZZ_{\geq 0}$ will be chosen in~\eqref{def:cnonarchi}. 
Further we define the set $\Sigma^*$ by setting $\Sigma^*=\Sigma$ if $p\geq 3$ and $\Sigma^*=\Sigma(4)$ if $p=2$. Here $\Sigma(4)$ denotes the set of points $P$ in $\Sigma$ with $\abs{x}>4$, where $(x,y)$ is the solution of \eqref{eq:weieq}  corresponding to $P$.  We obtain the following lemma.
\begin{lemma}\label{lem:nonarchicov}
Suppose that $P\in \Sigma^*$ satisfies $M'<\hat{h}_k(P)\leq M$. Then the point $P$ is determined modulo torsion by some lattice point in $\Gamma\cap\mathcal E$. 
\end{lemma}
\begin{proof}
Let $(x,y)$ be the solution of \eqref{eq:weieq} which corresponds to $P$, and write $P=Q+\sum n_i P_i$ with $Q\in E(\QQ)_{\textnormal{tor}}$ and $n_i\in\ZZ$. On using that $P$ is in $\Sigma$  and that $\mu\leq M'<\hat{h}_k(P)$, we deduce that $\log\abs{x}>0$ and thus our point $P$ lies in fact in $E_1(\QQ_p)$. Further, if $p=2$ then  our additional assumption $P\in \Sigma(4)$ provides that $v(x)<-2$. Hence on recalling that $P\in \Sigma$, we see that Lemma~\ref{lem:nonarchiest} leads to the inequality
$$
\abs{\sum n_i\alpha_i}\leq \abs{m}e^{-\tfrac{1}{\sigma}(\hat{h}_k(P)-\mu)}.
$$
Here we used that  $v(\log(mP))=v(\sum n_i \alpha_i)$, which in turn follows from $mP=\sum n_i(mP_i)$ and $mP_i\in E_1(\QQ_p)$. The displayed inequality together with  $M'<\hat{h}_k(P)$ shows that $v(\sum n_i\beta_i)\geq c$, where $c$ is the smallest element of $\ZZ_{\geq 0}$ which exceeds 
\begin{equation}\label{def:cnonarchi}
v(m)-v(\alpha_{i^*})+\tfrac{1}{\sigma \log p}(M'-\mu).
\end{equation} 
It follows that $\sum n_i \beta_i^{(c)}=0$ in $\ZZ/(p^c\ZZ)$ and therefore $\gamma=(n_i)$ lies in $\Gamma$. Furthermore, our assumption $\hat{h}_k(P)\leq M$ assures that $\gamma$ lies in $\mathcal E$ and hence $P$ is determined modulo torsion by the lattice point $\gamma\in \Gamma\cap \mathcal E$. This completes the proof of the lemma.
\end{proof}
The above lemma provides a sieve for the points $P$ in $\Sigma^*$ with $M'<\hat{h}_k(P)\leq M$. The discussion of the strength of this sieve, depending on the parameters $M'$, $M$ and $p^c$, is analogous to the discussion of the strength of the archimedean sieve in Section~\ref{sec:archisieve}. However there are some minor differences. For example, in the non-archimedean case we can work entirely in dimension $r$ and the parameter $p^c$ is uniquely determined by \eqref{def:cnonarchi}; note that $p^c$ plays here the role of the parameter $c$ in the archimedean sieve.

\paragraph{Non-archimedean sieve.}  The following  sieve uses the version (FP) of the Fincke--Pohst algorithm  described in Section~\ref{sec:compuaspects} in order to determine all points in $\Gamma\cap \mathcal E$.
 
\begin{Algorithm}[Non-archimedean sieve]\label{algo:nonarchsieve} 
The inputs  are  $M',M\in\ZZ$ with $\mu\leq M'< M$. The output is the set of points $P$ in $\Sigma^*$ with $M'<\hat{h}_k(P)\leq M$. 
\begin{itemize}
\item[(i)] To find the number $m$,  determine $e_{ns}$, $e_t$ and  $\iota $. 
\item[(ii)]Determine the lattice $\Gamma\subset \ZZ^{r}$  by computing  the  $p$-adic elliptic logarithms $\log(mP_i)$ up to the required precision for all $i\in\{1,\dotsc,r\}$. 
\item[(iii)] Use  \textnormal{(FP)} to find all points in $\Gamma\cap\mathcal E$.
\item[(iv)] For each $\gamma$ in $\Gamma\cap \mathcal E$ and for each torsion point $Q\in E(\QQ)_{\textnormal{tor}}$, output the point $P=Q+\sum \gamma_i P_i$ if $M'<\hat{h}_k(P)\leq M$ and if $P\in \Sigma^*$.
\end{itemize}
\end{Algorithm} 

\paragraph{Correctness.} We take a point $P\in \Sigma^*$  which satisfies $M'<\hat{h}_k(P)\leq M$ and we write $P=Q+\sum n_i P_i$ with $n_i\in\ZZ$ and $Q\in E(\QQ)_{\textnormal{tor}}$. Lemma~\ref{lem:nonarchicov} gives that $\gamma=(n_i)$ lies in $\Gamma\cap\mathcal E$ and hence we see that step (iv) produces our point $P$ as desired.

\paragraph{Complexity.} We now discuss the influence of each step on the running time in practice.
 In step (i) standard results and algorithms allow to quickly compute the numbers $e_{ns},e_t$ and $\iota$.  In fact the computation of (a suitable) $m$ is very fast in practice, even if  $p$ is relatively large. Step (ii) needs to compute $r$ distinct $p$-adic elliptic logarithms up to a number of $p$-adic digits which is approximately $c+v(\alpha_{i^*})$. The efficiency of this computation crucially depends on the size of $c+v(\alpha_{i^*})$, which in turn depends in particular on the lower bound $M'$. If the number $M'$ is huge, then this step (ii) can become very slow in practice. Finally we mention that a complexity analysis of steps (iii) and (iv) is contained in the complexity discussions of the analogous steps of the archimedean sieve in  Algorithm~\ref{algo:archisieve}.

\paragraph{Comparison.} Similarly as in the archimedean case, there are important differences between our approach and the known method. For instance, we work with the ellipsoid $\mathcal E$ arising from the N\'eron--Tate height $\hat{h}$ and we actually determine all points in the intersection $\Gamma\cap \mathcal E$. We refer to Sections~\ref{sec:globalsieve} and \ref{sec:comparisonwithelr} for a comparison of the two approaches and for a detailed discussion of the improvements provided by our approach.

\subsection{Height-logarithm sieve}\label{sec:heightlogsieve}

We work with the setup of Section~\ref{sec:setup}. The goal of this section is to construct a  sieve which allows to efficiently determine the set of $S$-integral points in any given finite subset of $E(\QQ)$. The sieve  exploits that the global N\'eron--Tate height is essentially determined by the various local elliptic logarithms and thus we call it the height-logarithm sieve. Throughout this section we use the notation introduced above and we assume that the Weierstrass model \eqref{eq:weieq} of our given elliptic curve $E$ is minimal at all $p\in S$.

\paragraph{Main idea.} To describe the main idea of the sieve, we take $P\in E(\QQ)$. For any finite place $v$ of $\QQ$, we define $\log_v(P)=\tfrac{1}{m_v}\log(m_vP)$ with $\log(\cdot)$ and $m_v=m$ as in Section~\ref{sec:nonarchisieve}.  There are real valued functions $f$ and $f_\infty$ on $E(\QQ)$, with $f$ bounded and $f_\infty$  determined by the real elliptic logarithm, such that any non-exceptional\footnote{Here we exclude the exceptional points (Definition~\ref{def:exceptpoint}) in order to avoid the usual technical problems arising when working with the $v$-adic elliptic logarithm at $v=\infty$ and $v=2$.} point $P\in E(\QQ)$ satisfies 
\begin{equation}\label{eq:heightlogformula}
\hat{h}(P)=f(P)+f_\infty(P)-\log\prod \abs{\log_v(P)}_v
\end{equation} 
with the product taken over certain finite places $v$ of $\QQ$. 
Here if $P$ is an $S$-integral point then the product ranges only over $v$ in $S$, providing a strong condition for points in $E(\QQ)$ to be $S$-integral. Furthermore, one can check this condition requiring only to know the form of $P$ in $E(\QQ)_{\textnormal{tor}}\oplus \ZZ^r$. For most points $P$, this  allows to circumvent the slow process of checking whether the coordinates of $P$ are $S$-integral, that is whether $P\in\Sigma(S)$.     

\paragraph{Construction.} To transform the above idea into an efficient sieve for the set of $S$-integral points $\Sigma(S)$ inside $E(\QQ)$, we shall work with a slightly weaker version of \eqref{eq:heightlogformula} which is suitable for our purpose. More precisely, we shall work with an inequality of the form $\hat{h}_k(P)\leq L(P)$ involving an efficiently computable quantity $L(P)$ which is essentially determined by the right hand side of \eqref{eq:heightlogformula}. We begin to explain how to determine $L(P)$ for any $P\in E(\QQ)$. Suppose that $P=Q+\sum n_i P_i$ with $Q\in E(\QQ)_{\textnormal{tor}}$ and $n_i\in\ZZ$. If $v$ is a finite place of $\QQ$ and if  $\alpha_{i,v}=m_v\log_v(P_i)$, then we define $$l_v(P)=\log\max\bigl(\abs{\tfrac{1}{2}}_v,\abs{\tfrac{1}{m_v}\sum n_i \alpha_{i,v}}_v^{-1}\bigl).$$
To give a similar definition at $v=\infty$, we take $\kappa\in\ZZ_{\geq 1}$ and we let $x_0=x_0(\kappa)$ be as in \eqref{def:x0}. Let $\omega_1$ be the period associated to \eqref{eq:weieq}, see Section~\ref{sec:archisieve}.  If $v=\infty$ then we write $\alpha_{i,v}=\log(m_vP_i)$ with $\log(\cdot)$ and $m_v=m$  as in Section~\ref{sec:archisieve} and for any $l\in \ZZ$ we define
$$l_v(P,l)=\log\max\bigl( x_0^{1/2},(1+\tfrac{1}{\kappa})^2\abs{\tfrac{1}{m_v}\bigl(l\omega_1+\sum n_i\alpha_{i,v}\bigl)}_v^{-1}\bigl).$$
Here we say that $l\in \ZZ$ is admissible for $P$ if $\abs{l}\leq 2m_\infty+\sum \abs{n_i}$. Let $\mu$ be as in \eqref{eq:nthcompa}. For any finite place $v$ of $\QQ$, we denote by $G_v$ the subgroup of $E(\QQ)$ formed by the points whose images in $E(\QQ_v)$ lie in fact in $E_1(\QQ_v)$. We shall use the following lemma. 

\begin{lemma}\label{lem:height-log}
If $P\in \Sigma(S)$ then there is an admissible $l\in\ZZ$ such that $$\hat{h}_k(P)\leq \mu+l_\infty(P,l)+\sum_{v\in S\st P\in G_v} l_v(P).$$
\end{lemma}
\begin{proof}
The statement follows by combining \eqref{eq:nthcompa} with Lemmas~\ref{lem:hk}, \ref{lem:archiest} and \ref{lem:nonarchiest}.
\end{proof}

In the next paragraph we shall explain how to control the quantities $l_\infty(P,l)$ and $\sum l_v(P)$ in order to obtain a suitable upper bound $L(P)$ for the right hand side of the inequality in Lemma~\ref{lem:height-log}. The resulting height-logarithm inequality $\hat{h}_k(P)\leq L(P)$ is the main ingredient of the following algorithm in which we identify  $E(\QQ)$ with $E(\QQ)_{\textnormal{tor}}\oplus \Gamma_E$, where $\Gamma_E=\ZZ^r$ is the image of $E(\QQ)$ inside $E(\QQ)\otimes_\ZZ\RR\cong \RR^r$ as in Section~\ref{sec:heightselllog}.

\begin{Algorithm}[Height-logarithm sieve]\label{algo:heightlogsieve} The inputs are $\kappa\in\ZZ_{\geq 1}$ and a finite subset $\Sigma$ of $E(\QQ)_{\textnormal{tor}}\oplus \Gamma_E$. The output is the set $\Sigma\cap\Sigma(S)$ of $S$-integral points inside $\Sigma$.

Determine the set $S_E$ formed by the places $v\in S$ where the elliptic curve $E$ has bad reduction. Then for each nonzero point $P\in \Sigma$ do the following:
\begin{itemize}
\item[(i)] Use the arguments of \textnormal{(1)} below to determine an upper bound $l_\infty(P)\geq \max l_\infty(P,l)$ with the maximum taken over all $l\in \ZZ$ which are admissible for $P$. 
\item[(ii)] Compute the set $S_P=\{v\in S\st P\in G_v\}\cup S_E$ as described in \textnormal{(2)} below. 
\item[(iii)] For each $v\in S_P$ determine $l_v(P)$ by using the arguments of \textnormal{(3)} below, and then set $L(P)=\mu+l_\infty(P)+\sum_{v\in S_P}l_v(P)$. Output $P$ if  $\hat{h}_k(P)\leq L(P)$ and if $P\in \Sigma(S)$.
\end{itemize}
\end{Algorithm}

\paragraph{Correctness.} Suppose that $P$ lies in $\Sigma\cap\Sigma(S)$. For each $v\in S_E$ we obtain that $l_v(P)\geq 0$ and thus $\sum_{v\in S_P}l_v(P)$ exceeds the sum $\sum l_v(P)$ taken over all $v\in S$ with $P\in G_v$. Hence Lemma~\ref{lem:height-log} implies that $\hat{h}_k(P)\leq L(P)$ and thus (iii) produces our point $P$ as desired.

\paragraph{Computing $L(P)$.}  We consider a nonzero point $P=(Q,(n_i))$ in $E(\QQ)_{\textnormal{tor}}\oplus \Gamma_E$; note that $P=Q+\sum n_i P_i$ in $E(\QQ)$. To compute the quantity $L(P)$ we proceed as follows:

\begin{itemize}
\item[(1)] To control $l_\infty(P,l)$ for any admissible $l\in \ZZ$, we compute the real elliptic logarithms $\omega_1$ and $\alpha_i=\alpha_{i,v}$ up to a certain precision with respect to $\abs{\cdot}=\abs{\cdot}_v$ for $v=\infty$. If the linear form $\Lambda=l\omega_1+\sum n_i\alpha_i$ is nonzero then the required precision can be obtained as follows: After choosing a sufficiently large $n\in\ZZ$, one determines approximations $\alpha_{i}'$ and $\omega_1'$ of $\alpha_{i}$ and $\omega_1$ respectively such that $\Lambda'=l\omega_1'+\sum n_i\alpha_i'$ satisfies $\abs{\Lambda'}>\epsilon=10^{-n}$ and such that the absolute differences $\abs{\alpha_i-\alpha_i'}$ and $\abs{\omega_1-\omega_1'}$ are at most $10^{-c}$ for some fixed integer $c\geq n+\log_{10} (2m_\infty+2\sum \abs{n_i})$.
Then the proof of Lemma~\ref{lem:archicov} gives that  $-\log \abs{\Lambda}\leq -\log(\abs{\Lambda'}-\epsilon)$ and hence we can practically compute an upper bound $l_\infty(P)$ in $\RR\cup\{\infty\}$ for $\max l_\infty(P,l)$ with the maximum taken over all admissible $l\in \ZZ$. Here if $\Lambda$ is close to zero or if $c$ is too large, then we just put $l_\infty(P)=\infty$ to assure that the computation of $l_\infty(P)$ is always efficient.
\item[(2)]  We would like to quickly compute the set $S_P=\{v\in S\st P\in G_v\}\cup S_E$. Here one can directly determine $S_E$, since the Weierstrass model \eqref{eq:weieq} is minimal at all $p\in S$. It remains to deal with the places $p\in S_P-S_E$. 
The elliptic curve $E$ has good reduction at $p$ and the canonical reduction map $E(\QQ)\hookrightarrow E(\QQ_p)\to E(\mathbb F_p)$ is a morphism of abelian groups. We determine the images $\bar{Q}$ and $\bar{P_i}$ in $E(\mathbb F_p)$ of all $Q\in E(\QQ)_{\textnormal{tor}}$ and all $P_i$. It follows that our point $P=Q+\sum n_i P_i$ lies in $G_v$ if and only if the point $\bar{Q}+\sum n_i\bar{P_i}$ is zero in $E(\mathbb F_p)$. Therefore we see that we can quickly determine the set $S_P$ provided we already know all $\bar{Q}$, all $\bar{P_i}$ and the group structure of $E(\mathbb F_p)$. 
\item[(3)] We take $v\in S$ and we now explain how to efficiently determine $l_v(P)$. As already mentioned in Section~\ref{sec:nonarchisieve} one can always quickly compute the number $m_v$  in practice. We write $\alpha_i=m_v\log_v(P_i)$ and we define $\alpha=\sum n_i \alpha_i$.  
To compute $v(\alpha)$ we need to know the $v$-adic elliptic logarithms $\alpha_i$ with a certain precision. In practice it usually suffices here to know  $\alpha_i$ with a small $v$-adic precision. Indeed after computing $v(\alpha_{i^*})=\min v(\alpha_i)$, we consider $\beta=\sum n_i\beta_i$ for $\beta_i\in \ZZ_p$ of the form $\beta_i=\alpha_i/p^{v(\alpha_{i^*})}$. The integer $v(\beta)$ is almost always small in practice. Hence one  can usually compute $v(\beta)$ and  $v(\alpha)$ by knowing only the first coefficients of the $v$-adic power series of $\alpha_i$.
\end{itemize}

\paragraph{Huge and tiny parameters} 
To assure that Algorithm~\ref{algo:heightlogsieve} is still fast for huge parameters, one can slightly weaken the sieve as follows: If one of the steps (except the final check whether $P\in \Sigma(S)$) should take too long for a point $P\in \Sigma$, then abort these steps and directly check whether $P\in \Sigma(S)$. To deal with the case of huge sets $S$ in which step (ii) becomes slow (see Remark~\ref{rem:rank1heightlog}), we can always replace $S_P$ by the usually much larger set $S$.  The resulting sieve is still strong for points $P$ with $\hat{h}(P)\gg\log N_S$. However, replacing $S_P$ by $S$ considerably weakens the sieve for points of small height. We now discuss the case when the rank $r$ is small and the height of the involved point $P$ is tiny. Here one can quickly compute the coordinates of $P$ and the Weil height of these coordinates is not that large. 
Hence in this case it is often faster to skip  steps (i) and (ii) and to directly determine in (iii) whether the coordinates of $P$ are $S$-integers. In our implementation of the height-logarithm sieve we take into account the above observations  to avoid that Algorithm~\ref{algo:heightlogsieve} is unnecessarily slow for huge or tiny parameters.

\paragraph{Complexity.} We now discuss aspects of Algorithm~\ref{algo:heightlogsieve} which considerably influence the running time in practice. In steps (i) and (iii) we need to compute various elliptic logarithms up to a certain precision depending on the height $\hat{h}(P)$ of the points $P\in \Sigma$. In practice we will apply the height-logarithm sieve only in situations in which the heights $\hat{h}(P)$ are not huge and in such situations steps (i) and (iii) are always very fast. 
In step (ii) the running time of the computation of the set $S_P$  crucially depends on the number of primes in $S$. In practice it turned out that this step is fast when $\abs{S}$ is small.  However, if $\abs{S}$ becomes huge then step (ii) can take a long time as explained in the following remark.

\begin{remark}[Rank $r=1$]\label{rem:rank1heightlog}
For huge sets $S$ the computation of $S_P$ takes a long time, since one has to compute with many large groups $E(\mathbb F_p)$. In the case $r=1$ the following observation considerably improves this process. Let $v\in S$ such that $E$ has good reduction at $v$ and write $p=v$. Let $e_v$ be the order of $P_1$ in the finite group $E(\mathbb F_p)$. Consider a point $P\in E(\QQ)$ with $P=Q+n_1P_1$ for $n_1\in\ZZ$ and $Q\in E(\QQ)_{\textnormal{tor}}$, and let $e_Q$ be the order of $Q$ in $ E(\QQ)_{\textnormal{tor}}$. If $P\in G_v$ then the points $-Q$ and $n_1P_1$ coincide in $E(\mathbb F_p)$ and hence $e_v$ divides $n_1e_Q$. In other words, if $e_v$ does not divide $n_1e_Q$ then $v$ is not in $S_P$ and therefore we obtain a sufficient criterion to decide whether $v\in S$ satisfies $v\notin S_P$. 
\end{remark}

\remark[Inequality trick]\label{rem:inequtrick}  In the case when $S$ is empty, one can use the known inequality trick \cite[p.147]{sttz:elllogoverview} which tests whether a given nonzero point $P\in E(\QQ)$ satisfies the inequality $\lambda\|P\|_\infty^2\leq \mu+ l_\infty(P)$. This inequality is weaker than $\hat{h}(P)\leq \mu+l_\infty(P)$ used in our height-logarithm sieve when $S$ is empty, since $\lambda\|P\|_\infty^2\leq \hat{h}(P)$. Hence our height-logarithm sieve is more efficient than the inequality trick, in particular in the case of large rank $r\geq 2$ where the function  $\lambda\|\cdot\|_\infty^2$ is usually much smaller than $\hat{h}(\cdot)$. See also the examples in the next section. In the case when $S$ is nonempty, one could obtain in principle an inequality trick by testing whether $P$ satisfies the inequality $\lambda\|P\|_\infty^2\leq \mu+\sigma l_v(P)$ for some $v\in S^*=S\cup\{\infty\}$ and $\sigma=\abs{S^*}$. However the resulting sieve is not that efficient (and often useless if $\sigma$ is large), since an arbitrary point $P$ usually satisfies at least one of these $\sigma$ different inequalities which are all considerably weakened by the factor $\sigma$.

\subsection{Refined enumeration}\label{sec:refinedenumell}

We work with the setup of Section~\ref{sec:setup}. The goal of this section is to construct a refined enumeration for the set of $S$-integral points $\Sigma(S)\subset E(\QQ)$ of bounded height which improves the standard enumeration.  Throughout this section we assume that the Weierstrass model \eqref{eq:weieq} of $E$ is minimal at all $p\in S$ and we continue the notation introduced above.  

Recall from Section~\ref{sec:heightselllog} that $\Gamma_E=\ZZ^r$ denotes the image of $E(\QQ)$ inside $E(\QQ)\otimes_\ZZ\RR\cong \RR^r$. For any given upper bound $\ellrad\in \RR_{\geq 1}$, consider the ellipsoid  $\mathcal E_\ellrad=\{z\in\RR^r\st \hat{h}_k(z)\leq \ellrad\}$ contained in $\RR^r$. We observe that the following algorithm works correctly.

\begin{Algorithm}[Refined enumeration]\label{algo:refenu} The input consists of $\kappa\in\ZZ_{\geq 1}$ together with  an upper bound $\ellrad\in \RR_{\geq 1}$. The output is the set of points $P\in \Sigma(S)$ with $\hat{h}_k(P)\leq \ellrad.$
\begin{itemize}
\item[(i)] Use \textnormal{(FP)} to determine all points in the intersection $\Gamma_E\cap\mathcal E_\ellrad$. 
\item[(ii)] For each $\gamma\in\Gamma_E\cap \mathcal E_{\ellrad}$ and for any $Q\in E(\QQ)_{\textnormal{tor}}$, output the point $P=Q+\sum \gamma_i P_i$ if $P$ lies in the set obtained by applying Algorithm~\ref{algo:heightlogsieve} with $\kappa=\kappa$ and $\Sigma=\{P\}$.
\end{itemize}
\end{Algorithm}

\paragraph{Complexity.} We now discuss various aspects which influence the running time of Algorithm~\ref{algo:refenu} in practice. As usual, the application of \textnormal{(FP)} in step (i) crucially depends on the rank $r$. 
The running time of step (ii) depends on the cardinality of $\Gamma_E\cap \mathcal E_\ellrad$, which in turn depends on $r$, $\ellrad$ and the regulator of $E(\QQ)$. Here the application of the height-logarithm sieve efficiently throws away most points in $\Gamma_E\cap \mathcal E_\ellrad$, in particular essentially all points in $\Gamma_E\cap \mathcal E_\ellrad$ of large height. This considerably improves the running time. 

\paragraph{Comparison.} We next compare our refined enumeration with the standard enumeration of the points $P\in\Sigma(S)$ with $\|P\|_\infty^2\leq \ellrad'$, where $\ellrad'=\ellrad/\lambda_k$ depends on the smallest eigenvalue $\lambda_k$ of $\hat{h}_k$. 
Recall that the standard enumeration proceeds as follows: For any $\gamma\in\Gamma_E$ with $\max \abs{\gamma_i}^2\leq \ellrad'$ and for each $Q\in E(\QQ)_{\textnormal{tor}}$,  output the point $P=Q+\sum \gamma_i P_i$ if the coordinates of $P$ are $S$-integers. In the case when $S$ is empty, one can use here in addition the known inequality trick explained in Remark~\ref{rem:inequtrick}. In general we observe that our refined sieve working with the ellipsoid $\mathcal E_\ellrad$ is more efficient, in particular for large rank $r$. Indeed the cube $\{\|\cdot\|_\infty^2\leq \ellrad'\}\subset\RR^r$ always contains the ellipsoid $\mathcal E_\ellrad$ and  then on comparing volumes we see that our refined sieve involves much fewer points. Furthermore the application of the height-logarithm sieve in the refined enumeration gives significant running time improvements. See Section~\ref{sec:comparisonwithelr} for examples and tables which illustrate in particular the running time improvements provided by our refined enumeration.

\subsection{Refined sieve}\label{sec:refinedsieveell}

In this section we work out a refinement of the global sieve obtained by patching together the archimedean sieve of Section~\ref{sec:archisieve} with the various non-archimedean sieves of Section~\ref{sec:nonarchisieve}. Throughout this section we work with the setup of Section~\ref{sec:setup} and we continue the notation introduced above. Furthermore we assume that the Weierstrass model \eqref{eq:weieq} of our given elliptic curve $E$ is minimal at all primes $p\in S$.

\paragraph{Main idea.} The main ingredient of the refined sieve is Proposition~\ref{prop:refinedcov}. Therein we construct a refined covering of certain subsets of the set of $S$-integral points $\Sigma(S)$, which allows to improve the global-local passage required to apply the local sieves  obtained in Sections \ref{sec:archisieve} and \ref{sec:nonarchisieve}. The construction of the covering is inspired by the refined sieve for $S$-unit equations developed in Section~\ref{sec:dwsieve+}. However, in the present case of elliptic curves, everything is more complicated. For example, one has to distinguish archimedean and non-archimedean places and one has to take care of certain exceptional points (Definition~\ref{def:exceptpoint}) arising from technical issues of the $v$-adic elliptic logarithm at the places $v=\infty$ and $v=2$. To deal efficiently with the exceptional points, we conducted some effort to work  entirely in the abelian group $E(\QQ)$. This allows here to avoid working with coordinate functions, which in turn is crucial to solve equations \eqref{eq:weieq} with huge parameters.

\paragraph{Construction of the covering.} For any given $M,M'$ in $\ZZ$ with $0\leq M'<M$, we would like to find the set of points $P\in \Sigma(S)$ which satisfy $M'<\hat{h}_k(P)\leq M$. For this purpose we ``cover" this set as follows. 
Let $\kappa,n,\tau$ in $\ZZ_{\geq 1}$ with $\tau\leq n\leq s^*$ for $s^*=\abs{S}+1$ 
and choose an admissible partition $\{S_j\}$ of $S^*=S\cup \{\infty\}$ into  disjoint nonempty parts $S^*= S_1\dotcup\ldots\dotcup S_g$. Here admissible partition means  that $|S_j|\leq n$ for all $j\in\{1,\dotsc,g\}$ and that  $g\leq\ceil{\tfrac{s^*}{n}}+2$. For a motivation of working with admissible partitions of $S^*$, we refer to the efficiency discussion given below. Further we choose ``weights" $w_1,\dotsc,w_\tau$ in $\QQ$ with $w_1=1$ and $w_1\geq\dotsc\geq w_\tau>0$  and for any $t\in \{1,\dotsc,\tau\}$  we put
\begin{equation}\label{def:weightswt}
\sigma_t=\tfrac{w}{w_t}, \ \ \ 
w=\sum_{j=1}^g w(j), \ \ \ w(j)= 
\begin{cases}
(|S_j|-\tau)w_\tau+\sum_{t\leq\tau} w_t & \textnormal{if } \tau\leq \abs{S_j},\\
\sum_{t\leq \abs{S_j}} w_t & \textnormal{if } \tau>\abs{S_j}.
\end{cases}
\end{equation}
Next we take $j\in\{1,\dotsc,g\}$ and we consider a nonempty subset $T$ of $S_j$ with cardinality $\abs{T}$ at most $\tau$. Write $t=\abs{T}$ and suppose that $v\in T$. If $v\in S$ then we denote by $\Gamma_v\subseteq\ZZ^r$ the lattice constructed in Section~\ref{sec:nonarchisieve} with $\sigma=\sigma_t$, and if $v=\infty$ then $\Gamma_v\subseteq\ZZ^{r+1}$ denotes the lattice from Section~\ref{sec:archisieve} with $\sigma=\sigma_t$ and $\kappa=\kappa$. In the case $\mu>M'$, where $\mu$ is as in \eqref{eq:nthcompa}, we set here $\Gamma_v=\ZZ^r$ if $v\in S$ and $\Gamma_v=\ZZ^{r+1}$ if $v=\infty$. Now we define
$$\Gamma_T=\bigcap_{v\in T} \Gamma_v.$$
Here if $T$ contains $\infty$ then for any $v\in S\cap T$ we identify $\Gamma_v\subseteq \ZZ^r$ with the lattice inside $\ZZ^{r+1}$ given by $\phi(\Gamma_v)\oplus e\ZZ$, where $\phi$ denotes the canonical product embedding of $\ZZ^r$ into $\ZZ^r\times\ZZ=\ZZ^{r+1}$ and $e=(0,1)\in\ZZ^{r}\times\ZZ$.   
Next we consider the ellipsoid $$\mathcal E_T$$
which is defined as follows: If $T\subseteq S$ then $\mathcal E_T\subset \RR^r$ is the ellipsoid appearing in Section~\ref{sec:nonarchisieve}, and if $T$ contains $\infty$ then $\mathcal E_T\subset \RR^{r+1}$ is the ellipsoid constructed in Section~\ref{sec:archisieve} with respect to the parameters $\sigma=\sigma_t$ and $\kappa=\kappa$. In the case $\mu>M'$ we define here $\mathcal E_T=\RR^r$ if $T\subseteq S$ and $\mathcal E_T=\RR^{r+1}$ if $T$ contains $\infty$. We next define the exceptional points.

\begin{definition}[Exceptional point]\label{def:exceptpoint}
Consider a point $P\in\Sigma(S)$ and denote by $(x,y)$ the corresponding solution of \eqref{eq:weieq}. We say that $P$ is an exceptional point if $\abs{x}_2\leq 4$ or if $\abs{x}_\infty\leq x_0$, where $x_0$ denotes the number $x_0(\kappa)$ defined in \eqref{def:x0}. 
\end{definition}
To completely ``cover" our set of interest, we need  to take into account the exceptional points.
For this purpose we let $\ellrad$ be the positive real number defined in \eqref{def:ellrad}, which depends inter alia on the parameters $\kappa,\tau,\{S_j\}$ and $w_t$, and we work with the ellipsoid $$\mathcal E_{\ellrad}=\{z\in\RR^r\st \hat{h}_k(z)\leq \ellrad\}.$$
Recall from Section~\ref{sec:heightselllog} that $\Gamma_E=\ZZ^r$ denotes the image of $E(\QQ)$ inside $E(\QQ)\otimes_\ZZ\RR\cong \RR^r$. The following result shows that  our set of interest can be ``covered" by the set $\Gamma_E\cap \mathcal E_{\ellrad}$ together with the sets $\Gamma_T\cap \mathcal E_T$ associated to some $T$ as above.
\begin{proposition}\label{prop:refinedcov}
Suppose that $P$ lies in $\Sigma(S)$ and assume that $M'<\hat{h}_k(P)\leq M$. Then at least one of the following statements holds.
\begin{itemize}
\item[(i)] The point $P$ is determined modulo torsion by some $\gamma\in\Gamma_E\cap\mathcal E_\ellrad$.
\item[(ii)] There is  a nonempty set $T$ with $\abs{T}\leq \tau$ such that $T\subseteq S_j$ for some $S_j$ in $\{S_j\}$ and such that $P$ is determined modulo torsion by an element $\gamma\in \Gamma_T\cap \mathcal E_T$.
\end{itemize}
\end{proposition}
\begin{proof}
If $M'<\mu$ then (ii) holds for example with $T=\{\infty\}$. Thus we may and do assume that $M'\geq \mu$. Let $(x,y)$ be the solution of \eqref{eq:weieq} corresponding to $P$. We claim that there exists a nonempty set  $T$ with $\abs{T}\leq \tau$ such that  $T\subseteq S_j$ for some $S_j$ in $\{S_j\}$ and such that 
\begin{equation}\label{coveringclaim}
P\in \bigcap_{v\in T}\Sigma(v,\sigma_t).
\end{equation}
Here $\Sigma(v,\sigma_t)$ is defined in \eqref{def:sigmavsigma} with $t=\abs{T}$.
To prove this claim by contradiction, we assume that \eqref{coveringclaim} does not hold. Then for each $j$ and for any nonempty subset $T\subseteq S_j$ with $t=\abs{T}\leq \tau$, there exists $v\in T$ such that $P\notin \Sigma(v,\sigma_t)$. In particular, for any $j$ and for each $t\in\ZZ_{\geq 1}$ with $t\leq \min(\tau,\abs{S_j})$, it follows that the $t$-th largest  of the real numbers $\tfrac{1}{2}\log\abs{x}_v$, $v\in S_j,$ is strictly smaller than  $\tfrac{1}{\sigma_t}(\hat{h}_k(P)-\mu)$. 
We deduce that $\tfrac{1}{2}\sum\max(0,\log\abs{x}_v)< \tfrac{w(j)}{w}(\hat{h}_k(P)-\mu)$ with the sum taken over all $v\in S_j$. Here we used that $\mu\leq M'<\hat{h}_k(P)$ and that the weights $w_t$ satisfy  $w_1\geq \dotsc \geq w_\tau>0$. Then our assumption $P\in \Sigma(S)$ together with $S^*=\cup S_j$ implies that $\tfrac{1}{2}h(x)< \hat{h}_k(P)-\mu$. But this contradicts the inequality $\hat{h}_k(P)-\mu\leq \tfrac{1}{2}h(x)$ which follows by combining Lemma~\ref{lem:hk} with \eqref{eq:nthcompa}. Therefore we conclude that our claim  \eqref{coveringclaim} holds as desired.

Let $\mathcal T$ be the nonempty set of all sets $T$ satisfying \eqref{coveringclaim};   put $t=\min\{\abs{T}\st T\in\mathcal T\}$ and define $\mathcal T_{\min}=\{T\in\mathcal T\st \abs{T}=t\}$. Further, on slightly abusing terminology, we write $\Sigma(4)$ for the subset of $\Sigma(2,\sigma_t)$ defined in Section~\ref{sec:nonarchisieve} and we denote by $\Sigma(x_0)$ the subset of $\Sigma(\infty,\sigma_t)$ from Section~\ref{sec:archisieve}. First we consider the case $t=1$. Suppose that we can choose $T\in\mathcal T_{\min}$ with $T\subseteq S-2$. Then we obtain that $T=\{p\}$ with $p\geq 3$. Thus on recalling that the Weierstrass model \eqref{eq:weieq} is minimal at all $p\in S$, we see that the inequalities $\mu\leq M'$ and $M'<\hat{h}_k(P)\leq M$ together with  \eqref{coveringclaim} show that $P$ satisfies the assumptions of Lemma~\ref{lem:nonarchicov}. Hence Lemma~\ref{lem:nonarchicov}  implies (ii). Suppose now that there is no $T\in\mathcal T_{\min}$ with $T\subseteq S-2$. Then  $\{\infty\}$ or $\{2\}$ lies in $\mathcal T_{\min}$ and  any $v\in S-2$ satisfies  
\begin{equation}\label{eq:extrasearch}
\tfrac{1}{2}\log |x|_v< \tfrac{1}{\sigma_1}(\hat{h}_k(P)-\mu).
\end{equation}
To complete the proof for $t=1$, it remains to establish (i) or (ii) in the following  cases (a), (b) and (c). Before we go into these cases we define the number $\ellrad$ appearing in $\mathcal E_\ellrad$: If $t^*=\min(\tau,\max_{v=2,\infty}\abs{S_{j(v)}})$ with $S_{j(v)}$ denoting the set $S_j$ which contains $v$, then 
\begin{equation}\label{def:ellrad}
\ellrad=\mu+\tfrac{1}{2w_{t^*}}(1+s_{t^*})\log \max(x_0,4). 
\end{equation}
Here $x_0=x_0(\kappa)$ is defined in \eqref{def:x0} and $s_{t^*}$ is the number of $p\in S$ with $p^{2w_{t^*}}\leq \max(x_0,4)$. In what follows we shall use that $P\in \Sigma(S)$, that for each finite place $v$ of $\QQ$ it holds $v(x)\leq -2$ if $P\in E_1(\QQ_v)$, that $1=w_1\geq\dotsc\geq w_\tau>0$, and that $\hat{h}_k\leq \hat{h}$ by Lemma~\ref{lem:hk}.

\begin{itemize}
\item[(a)]Case $\{2\}\in\mathcal T_{\min}$ and $\{\infty\}\in\mathcal T_{\min}$. If $\abs{x}_\infty>x_0$ or $\abs{x}_2> 2^2$, then \eqref{coveringclaim} implies that $P\in\Sigma(x_0)$ or $P\in\Sigma(4)$  and thus Lemma~\ref{lem:archicov} or Lemma \ref{lem:nonarchicov} shows (ii) for $T=\{\infty\}$ or $T=\{2\}$ respectively. On the other hand, if  $\abs{x}_\infty\leq x_0$ and $\abs{x}_2\leq 2^2$ then  \eqref{eq:extrasearch} and \eqref{coveringclaim} lead to
an upper bound for $h(x)$ which together with \eqref{eq:nthcompa} proves (i).
\item[(b)] Case $\{2\}\in \mathcal T_{\min}$ and $\{\infty\}\notin\mathcal T_{\min}$. Here inequality \eqref{eq:extrasearch} holds in addition for $v=\infty$, since $\{\infty\}$ is not in $\mathcal T_{\min}$. Therefore, if $\abs{x}_2\leq 2^2$ then we see as above that  \eqref{eq:extrasearch}, \eqref{coveringclaim} and \eqref{eq:nthcompa} imply statement (i). If $\abs{x}_2> 2^ 2$ then \eqref{coveringclaim} gives that $P\in \Sigma(4)$ and hence Lemma~\ref{lem:nonarchicov} shows statement (ii) with $T=\{2\}$.
\item[(c)]Case $\{2\}\notin \mathcal T_{\min}$ and $\{\infty\}\in\mathcal T_{\min}$. Now  \eqref{eq:extrasearch} holds in addition for $v=2$, since $\{2\}$ is not in $\mathcal T_{\min}$. Thus as above we  deduce (i) if $\abs{x}_\infty\leq x_0$. If $\abs{x}_\infty> x_0$ then \eqref{coveringclaim} gives $P\in \Sigma(x_0)$ and hence Lemma~\ref{lem:archicov} proves (ii) with $T=\{\infty\}$. 
\end{itemize} 
We now establish the case $t\geq 2$. If we can choose $T\in\mathcal T_{\min}$ with $T\subseteq S-2$, then \eqref{coveringclaim} together with Lemma~\ref{lem:nonarchicov} implies (ii). Suppose now that there is no $T\in\mathcal T_{\min}$ with $T\subseteq S-2$. 
Then each $T$ in $\mathcal T_{\min}$ contains $\infty$ or $2$. Furthermore any $v\in S^*$ satisfies \eqref{eq:extrasearch}, since  $t\geq 2$. To complete the proof it thus suffices to consider the following cases:
\begin{itemize}
\item[(d)] Case when each $T\in\mathcal T_{\min}$ contains $2$ and $\infty$. If $\abs{x}_\infty > x_0$ and $\abs{x}_2> 2^2$, then \eqref{coveringclaim} gives that $P\in \Sigma(x_0)\cap \Sigma(4)$. Thus on recalling the construction of $\Gamma_T$ and $\mathcal E_T$, we see that  Lemmas~\ref{lem:archicov} and \ref{lem:nonarchicov} together with \eqref{coveringclaim}  show that (ii) holds for any $T\in \mathcal T_{\min}$. 
On the other hand, if  $\abs{x}_\infty\leq x_0$ or $\abs{x}_2\leq 2^2$ then \eqref{eq:extrasearch}, \eqref{coveringclaim} and \eqref{eq:nthcompa} prove (i).
\item[(e)] Case when there is $T\in \mathcal T_{\min}$ with $2\in T$ and $\infty\notin T$. If $\abs{x}_2> 2^2$ then \eqref{coveringclaim} gives that $P\in \Sigma(4)$ and thus we deduce (ii) by using  $\infty\notin T$, \eqref{coveringclaim} and Lemma~\ref{lem:nonarchicov}. On the other hand, if $\abs{x}_2\leq 2^2$ then \eqref{eq:extrasearch}, \eqref{coveringclaim} and \eqref{eq:nthcompa} imply (i).
\item[(f)] Case when there exists $T\in \mathcal T_{\min}$ with $2\notin T$ and $\infty\in T$. If $\abs{x}_\infty\leq x_0$ then \eqref{eq:extrasearch}, \eqref{coveringclaim} and \eqref{eq:nthcompa} imply (i). Finally, if $\abs{x}_\infty> x_0$ then \eqref{coveringclaim} gives $P\in \Sigma(x_0)$. Therefore on using $2\notin T$ and \eqref{coveringclaim}, we see that Lemmas~\ref{lem:archicov} and \ref{lem:nonarchicov} prove (ii).
\end{itemize}
Hence we conclude that in all cases (i) or (ii) holds. This completes the proof.
\end{proof}

The arguments used to prove \eqref{coveringclaim} show in addition that one can further refine the covering in Proposition~\ref{prop:refinedcov} by working with the ellipsoids $\mathcal E_T^*$ discussed in Remark~\ref{rem:refellips}.

\paragraph{Refined Sieve.} A collection of sieve parameters $\mathcal P$ consists of the following data: Parameters $\kappa,\tau,n\in\ZZ_{\geq 1}$ with $\tau\leq n\leq s^*$,  an admissible partition $\{S_j\}$  of $S^*=S_1\dotcup\dotsc \dotcup S_g$ with respect to $n$, and  weights $w_1,\dotsc,w_\tau$ as in \eqref{def:weightswt}. We denote by $\ellrad(\mathcal P)$ the number associated to $\mathcal P$ as in \eqref{def:ellrad} and we obtain the following algorithm.

\begin{Algorithm}[Refined sieve]\label{algo:refinedsieve} The input is a collection of sieve parameters $\mathcal P$ together with bounds $M',M\in\ZZ_{\geq 0}$ satisfying $M'< M$. Put $M'_\ellrad=\max(\ellrad(\mathcal P),M')$. The output is the set of points $P\in \Sigma(S)$ with $M'_\ellrad<\hat{h}_k(P)\leq M$.

For any $j\in \{1,\dotsc,g\}$ and for each $T\subseteq S_j$ with $1\leq \abs{T}\leq \tau$, do the following:
\begin{itemize}
\item[(i)] Determine a basis of $\Gamma_T$ and then compute $\Gamma_T\cap \mathcal E_T$ by using the version of the Fincke--Pohst algorithm in \textnormal{(FP)}. 
\item[(ii)] For each $\gamma\in\Gamma_T\cap \mathcal E_T$ and for any $Q\in E(\QQ)_{\textnormal{tor}}$, output the point $P=Q+\sum \gamma_i P_i$ if $P$ satisfies $M'_\ellrad<\hat{h}_k(P)\leq M$ and if $P$ lies in the set produced by an application of Algorithm~\ref{algo:heightlogsieve} with $\Sigma=\{P\}$ and $\kappa=1$.
\end{itemize}
\end{Algorithm}

\paragraph{Correctness.} Assume that $P\in \Sigma(S)$ satisfies $M'_\ellrad<\hat{h}_k(P)\leq M$. Then it holds that  $\hat{h}_k(P)>\ellrad(\mathcal P)$ and hence there is no lattice point $\gamma$ in $\Gamma_E\cap \mathcal E_{\ellrad(\mathcal P)}$ such that $P$ is determined modulo torsion by $\gamma$. Furthermore, our assumption provides that $M'<\hat{h}_k(P)\leq M$. Therefore Proposition~\ref{prop:refinedcov} shows that step (ii) produces our point $P$ as desired.

\paragraph{Efficiency.} We now discuss the efficiency of the refined sieve and we motivate several concepts appearing therein. First we observe that the case $n=1$ in the refined sieve corresponds to the non-refined sieve obtained by patching together the local sieves at $v\in S^*$ with $\sigma=s^*$. Suppose now  that $n\geq\tau\geq 2$ in the refined sieve. Then the iteration over the sets $T$ ranges in particular over all sets $T=\{v\}$ with $v\in S^*$. However, compared with the non-refined sieve $n=1$, there is the following fundamental difference:  If $T=\{v\}$ then the discussions in Sections~\ref{sec:archisieve} and \ref{sec:nonarchisieve} together with $\sigma_1\leq s^*$ show that the refined sieve involving $\Gamma_T\cap \mathcal E_T$ with $\sigma_1$ is usually much stronger than the non-refined sieve $n=1$ involving the local sieve at $v$ with $\sigma=s^*$. Furthermore, if $\abs{T}\geq 2$ then the intersection $\Gamma_T=\cap \Gamma_v$ is usually considerably smaller than each part $\Gamma_v$. These observations suggest that the improvements coming from $\abs{T}=1$ are significant enough to absorb the additional iterations over sets $T$ with $\abs{T}\geq 2$. In practice this turned out to be correct in many fundamental situations (see Section~\ref{sec:globalsieve}), showing that the refined sieve provides significant running time improvements. To deal with huge sets $S$, we introduced admissible coverings of $S^*$ which allow to control the number of additional iterations over sets $T$ with $\abs{T}\geq 2$. Indeed the conditions $\abs{S_j}\leq n$ and $\abs{T}\leq \tau$ assure that the number of additional iterations are controlled in terms of $n,\tau$. Further we point out that a canonical choice for the weights $w_t$ would be $w_t=\tfrac{1}{t}$.  However in practice it turned out that for $t\geq 2$ it would be better to choose $w_t$ slightly larger than $\tfrac{1}{t}$. In fact this is the reason for working with the more general weights $w_t$ defined above \eqref{def:weightswt}. Finally we mention that the discussion of the influence of the parameters $M'$, $M$ and $\sigma_t$ on the strength of the sieve $\Gamma_T\cap\mathcal E_T$ with $\abs{T}=t$ is similar to the corresponding discussions in Sections~\ref{sec:archisieve} and \ref{sec:nonarchisieve}.

\paragraph{Complexity.} We do not try to analyse the complexity of the refined sieve in general, since it depends on too many parameters. However, in Sections~\ref{sec:globalsieve} and \ref{sec:comparisonwithelr} we shall  discuss aspects influencing the complexity of the refined sieve and we shall illustrate the running improvements provided by the refined sieve in various fundamental cases.

\begin{remark}[Refined ellipsoids]\label{rem:refellips}
For any set $T$ as above with $\abs{T}\geq 2$, the arguments used in the proof of \eqref{coveringclaim}  show in addition the following: Instead of using in Algorithm~\ref{algo:refinedsieve} the ellipsoids $\mathcal E_T$, one can work  with the ellipsoids $\mathcal E_T^*$ obtained by replacing in the definition of $\mathcal E_T$ the bound $M$ by the possibly much smaller number  $$M_T=\min\left(M,w_T(M'-\mu)+\mu\right), \ \ \ w_T=\tfrac{1}{w}\sum_{j=1}^g w^*(j).$$ Here $w^*(j)$ is obtained by replacing in the definition of $w(j)$ the number $\tau$ by  $\abs{T}-1$. Note that $\mathcal E_T^*\subseteq\mathcal E_T$ and $w_T$ does not depend on $M$. Now if $M>M_T$ then  $\mathcal E_T^*$ is strictly contained in $\mathcal E_T$ and hence using $\mathcal E_T^*$ improves Algorithm~\ref{algo:refinedsieve}~(i). In principle further refinements are possible by taking into account the part $S_{j}$ of $S^*$ which contains  $T\subseteq S_{j}$.
\end{remark}

\subsection{Global sieve}\label{sec:globalsieve}

We continue the setup, notation and assumptions of the previous section. After choosing suitable collections of sieve parameters, we combine in this section the refined enumeration with the refined sieve: For any given upper bound $M_1\in\ZZ_{\geq 1}$, we obtain a global sieve which allows to efficiently determine all points $P\in \Sigma(S)$ with $\hat{h}_k(P)\leq M_1$. 

\paragraph{Sieve parameters.} We shall apply our refined sieve with the following collections of sieve parameters. Choose $\kappa\in\ZZ_{\geq 1}$ such that $\abs{\ellrad-10}$ is as small as possible, where $\ellrad$ is defined in \eqref{def:ellrad} with $\tau=1$ and $\kappa=\kappa$. Then let $\mathcal P(1)$ be the collection of sieve parameters determined by $\kappa=\kappa$ and $n=1$. For any $i\in \{2,3,4\}$  we define $\mathcal P(i)$ as follows: We take $\kappa=\kappa$, $\tau=i$ and $n=10$, and we choose weights $w_1,\dotsc,w_\tau$ as in \eqref{def:weightswt} such that each $w_t$ is slightly larger than $\tfrac{1}{t}$ for $t\geq 2$. Further we use here an admissible covering $\{S_j\}$ of $S^*=S\cup\{\infty\}$ with $S_1=\{\infty\}$ and with the following properties: If $2\notin S$ then $\abs{S_j}=n$ for each $j\in\{2,g-1\}$, and if $2\in S$ then $S_2=\{2\}$ and  $\abs{S_j}=n$ for any $j\in\{3,g-1\}$. We shall motivate our choice of sieve parameters in the discussions below.

\begin{remark}\label{rem:equalrad}
For each $i\in\{1,\dotsc,4\}$  the number $\ellrad_i$, associated to $\mathcal P(i)$ in \eqref{def:ellrad}, satisfies $\ellrad=\ellrad_i$. Indeed on using that $\tau=1$ in $\mathcal P(1)$ and that $\max_{v=2,\infty}\abs{S_{j(v)}}=1$ in $\mathcal P(i)$ with $i\geq 2$, we see in all four cases that $t^*=1$  and hence we obtain that $\ellrad=b_i$ as desired. 
\end{remark}

\paragraph{Global sieve.} For any given $M_1\in\ZZ_{\geq 1}$, we would like to efficiently determine the set of points $P\in \Sigma(S)$ with $\hat{h}_k(P)\leq M_1$. For this purpose we enumerate these points  from below  and from above, using  the refined enumeration in Algorithm~\ref{algo:refenu} and the refined sieve in Algorithm~\ref{algo:refinedsieve} respectively. More precisely we proceed as follows:
\begin{itemize} 
\item[(a)] Let $\ellrad$ be as in the above paragraph, and define parameters $\indpar'=\ellrad$, $\indpar=M_1$ and $f(\indpar)=\min(\indpar-1,\lfloor 0.99\indpar\rfloor)$. Then apply the following sieves $(i=1,\dotsc,4)$:
\begin{itemize}
\item[0.] Apply Algorithm~\ref{algo:refenu} with $\ellrad=\indpar'$ and $\kappa=1$, and let $\rho_0$ be the running time divided by the euclidean volume $\vol(\mathcal E_{\indpar'})$. Put $\indpar'= 4^{1/r}\indpar'$.
\item[i.] If $\abs{S^*}\geq i$ and $\indpar>\indpar'/4^{1/r}$ then  apply Algorithm~\ref{algo:refinedsieve} with the parameters $\mathcal P=\mathcal P(i)$,  $M'=f(\indpar)$ and $M=\indpar$, and let $\rho_i$ be the running time divided by the euclidean volume $\vol(\mathcal E_{M}\setminus\mathcal E_{M'})$. Put $\indpar=f(\indpar)$.
 \end{itemize}
\item[(b)] As long as $\indpar>\indpar'/4^{1/r}$ continue with the most efficient sieve in (a), that is the sieve $i^*\in \{0,\dotsc,4\}$  for which the ``efficiency measure" $1/\rho_{i^*}$ is maximal. 
\end{itemize}
This algorithm outputs the set of points $P\in \Sigma(S)$ with $\hat{h}_k(P)\leq M_1$. Indeed Remark~\ref{rem:equalrad} gives that $\ellrad=\ellrad_i$ for each $i\in\{1,\dotsc,4\}$, and therefore we see that the whole space $\{P\in \Sigma(S)\st \hat{h}_k(P)\leq M_1\}$ is covered by an application of steps (a) and (b).

\paragraph{Decomposition.} We now motivate the decomposition of the above algorithm and we explain our choice of the parameters appearing therein. First we discuss the step size functions. In the enumeration from below, we double the volume of the ellipsoid in each step by working with the step size function given by multiplication with $4^{1/r}$. This assures that the repeated enumerations of candidates with tiny height are not that significant for the running time, see also the remark at the end of this paragraph. In the enumeration from above, we work with the step size function $f(\indpar)=\min(\indpar-1,\lfloor 0.99\indpar\rfloor)$. If $\indpar$ is large then the height lower bound $M'=f(\indpar)$ is still large and thus the refined sieve is strong in view of the complexity discussions in Sections~\ref{sec:archisieve} and \ref{sec:nonarchisieve}. Hence, to accelerate the enumeration from above, we work with the relatively big step size $\indpar-f(\indpar)\geq 10^{-2}\indpar$ for large $\indpar$; here the factor $0.99$ turned out to be suitable in practice where usually $M_1\leq 10^9$. On the other hand, if $\indpar\leq 100$ is small then we work with the tiny step size $\indpar-f(\indpar)=1$ to assure that $M'=\indpar-1$ is relatively large  making the sieves stronger.

We next motivate our choice of the sieve parameters $\mathcal P(i)$.  If the parameter $\kappa$ becomes larger then the refined sieve becomes stronger at $v=\infty$. On the other hand, we can not choose $\kappa$ arbitrarily large since the ``square-radius"  $\ellrad$ of the ellipsoid $\mathcal E_\ellrad$ satisfies $\ellrad\geq \log(\kappa)$.  Now, choosing $\kappa$ such that $\abs{\ellrad-10}$ is as small as possible, assures that the refined enumeration via $\Gamma_E\cap \mathcal E_\ellrad$ is efficient in practice where usually $r\leq 12$. To explain our choices for $n$ and $\tau$, we recall that the efficiency of the refined sieve depends inter alia on $\sigma_1=\sigma_1(n,\tau)$ and the number of additional iterations over subsets $T$ with $\abs{T}\leq \tau$. Here it is not clear to us what are the optimal choices for $\tau$ and $n$. In practice it turned out that for $\tau\geq 5$ or $n\geq 11$ there are usually too many additional iterations and thus we only work with $\tau\leq 4$ and $n=10$. The reason for using admissible partitions $\{S_j\}$ of $S^*$ with $S_1=\{\infty\}$ and $S_2=\{2\}$ if $2\in S$, is to assure that $\ellrad=\ellrad_i$ (see Remark~\ref{rem:equalrad}) which means that the  ellipsoids $\mathcal E_{\ellrad_i}$ are not larger than the minimal involved ellipsoid $\mathcal E_{\ellrad}$. In fact  controlling the ellipsoids $\mathcal E_{\ellrad_i}$  is crucial for dealing efficiently with huge parameters.

Finally we mention that in the enumeration from below, the application of (FP) repetitively enumerates candidates. Here it is not clear to us how to avoid these repeated enumerations of candidates, since (FP) is not faster for circular discs than for the whole ellipsoid. In any case  these repeated enumerations of candidates have a small influence on the running time in practice, since (FP) is very fast in our situations of interest where usually the rank $r$ is small. Furthermore, to assure that the height-logarithm sieve is applied at most once for each candidate, we order the candidates with respect to their height $\hat{h}_k$ in the implementation of the enumeration from below.

\paragraph{Main features.} We now discuss the main features of the global sieve. The first steps of the enumeration from above may be viewed as a reduction of $M_1$. Indeed these steps are usually very fast in practice and they often allow to considerably reduce $M_1$. 

Further we point out that each step of the global sieve is more efficient than the standard enumeration. In fact in general it is not clear to us in which situation which sieve of (a) is the most efficient. To overcome this problem, we work with the quantities $1/\rho_i$ in order to ``measure" the efficiency of the sieves in the given situation. In practice this allows our algorithm to choose a suitable sieve in each step. This is very important for the running time, since the efficiency of the involved sieves strongly depends on the given situation. We also mention that Section~\ref{sec:comparisonwithelr}  contains explicit examples which illustrate (up to some extent) the improvements in practice provided by our global sieve.

\subsection{Elliptic logarithm sieve}\label{sec:elllogsievealgo} 
We work with the setup of Section~\ref{sec:setup} and we continue the notation introduced above. In the first part of this section we construct the elliptic logarithm sieve by putting together the sieves obtained in the previous sections. In the second part we compare the elliptic logarithm sieve with the known approach and we explain in detail our improvements. 

We recall that in the setup of Section~\ref{sec:setup} we are given the following information: The coefficients $a_i\in\ZZ$ of a Weierstrass equation \eqref{eq:weieq} of an elliptic curve $E$ over $\QQ$, a finite set of rational primes $S$, a basis $P_1,\dotsc,P_r$ of the free part of $E(\QQ)$ and a number $M_0\in\ZZ$ such that any $P\in\Sigma(S)$ satisfies $\hat{h}(P)\leq M_0$. Given this information, the following algorithm completely determines the set $\Sigma(S)$ formed by the $S$-integral solutions of \eqref{eq:weieq}.

\begin{Algorithm}[Elliptic logarithm sieve]\label{algo:elllogsieve} The inputs are the coefficients $a_i$ of \eqref{eq:weieq}, the set $S$, the basis $P_1,\dotsc,P_r$ and the initial bound $M_0$. The output is the set $\Sigma(S)$.

\begin{itemize}
\item[(i)] Compute the following additional input data.
\begin{itemize}  
\item[(a)] Determine the equation of an affine Weierstrass model $W$ of $E$ over $\ZZ$, which is minimal at all primes in $S$, together with an isomorphism $\varphi$ over $\mathcal O=\ZZ[1/N_S]$ from  $W$ to the affine model defined by \eqref{eq:weieq}. 
\item[(b)] Compute a suitable rational approximation $\hat{h}_k$ of $\hat{h}$.
\item[(c)] Determine the torsion subgroup $E(\QQ)_{\textnormal{tor}}$ of  $E(\QQ)$, and compute the numbers $\ellrad$ and $\kappa$ appearing in the collections of sieve parameters from Section~\ref{sec:globalsieve}.
\end{itemize}  
 \item[(ii)] To reduce locally the initial bound $M_0$,  apply the archimedean sieve and the non-archimedean sieve. More precisely, work with the set $\Sigma'(S)$ formed by the $\mathcal O$-points of $W$ and for each place $v$ in $S^*=S\cup \{\infty\}$ do the following: 
\begin{itemize}
\item[(a)] Find an integer $M_1(v)\geq \ellrad$ with the following property: If $M'=M_1(v)$, $M=M_0$, $\kappa=\kappa$ and $\sigma=\abs{S}+1$, then Algorithm~\ref{algo:archisieve}~(ii) outputs only $0$ when $v=\infty$ or $0$ is the output of Algorithm~\ref{algo:nonarchsieve}~(iii) when $v\in S$. Here first try $M_1(v)=\ellrad$. If this does not work then try a slightly larger number, and so on.
\item[(b)] Having found such an $M_1(v)$, try to reduce it further by repeating (a) with  different parameters $M'$ and $M$. Let $M_1(v)$ be the final reduced bound at $v$.
\end{itemize} 
\item[(iii)] Determine the set $\Sigma'(S)$ by applying the global sieve from Section~\ref{sec:globalsieve} with $M_1=\max_{v\in S^*}M_1(v)$. Then output the set $\varphi(\Sigma'(S))$. 
\end{itemize}
\end{Algorithm}

\paragraph{Correctness.} To prove that this algorithm works correctly, we recall that $\varphi$ is an isomorphism of affine Weierstrass models of $E$ over $\mathcal O$. This shows that $\varphi(\Sigma'(S))=\Sigma(S)$ and hence it remains to verify that $\Sigma'(S)$ is completely determined. Lemma~\ref{lem:hk} gives that $\hat{h}_k\leq\hat{h}$, and any $P\in \Sigma'(S)$ satisfies $\hat{h}(P)\leq M_0$ since $\hat{h}$ is invariant under isomorphisms.  Therefore on using that $\ellrad\leq M_1$, we see that the arguments of Proposition~\ref{prop:refinedcov}, with $n=1$ and $\kappa=\kappa$, prove that $\hat{h}_k(P)\leq M_1$ for all $P\in \Sigma'(S)$. 
It follows that the application of the global sieve with $M_1$ produces the set $\Sigma'(S)$ in step (iii) as desired.

\paragraph{Complexity.} We now discuss various aspects which influence the running time in practice. The computation of the additional input data in step (i) is always very fast. More precisely, in (a) we use Tate's algorithm to transform \eqref{eq:weieq} into a globally minimal Weierstrass model of $E$ over $\ZZ$ which then can be used to directly determine a pair $(W,\varphi)$ with the desired properties. To construct the quadratic form $\hat{h}_k$ in (b) we proceed as described in Section~\ref{sec:heightselllog}; see also Section~\ref{sec:compuaspects} for computational aspects. Finally in (c) the numbers $\ellrad$ and $\kappa$ can be directly determined and the computation of  $E(\QQ)_{\textnormal{tor}}$ is always efficient.

The running time of step (ii) crucially depends on the initial upper bound $M_0$, see the complexity discussions in Sections~\ref{sec:archisieve} and \ref{sec:nonarchisieve} for details. Here step (a) can take a long time for huge $M_0$, while the repetitions in step (b) are then quite fast since at this point we have already computed the involved elliptic logarithms. Step (ii) usually allows to avoid the process of testing  whether candidates of huge height have $S$-integral coordinates. This process is so slow that it is beneficial in (b) to make as many repetitions as required to obtain a reduced global bound $M_1$ which is as small as possible.

The running time of the application of the global sieve in step (iii) crucially depends on the cardinality of $S$ and the rank $r$. For example, as already explained in previous sections, the rank $r$ has a huge influence on the running time of the refined enumeration and on the refined sieve. We also recall that the cardinality of $S$ significantly influences the efficiency of the height-logarithm sieve, which in turn is used in the refined sieve and the refined enumeration. See also the discussions in Sections~\ref{sec:globalsieve} and \ref{sec:comparisonwithelr}.

\paragraph{Bottleneck.} In practice the bottleneck of the elliptic logarithm sieve crucially depends on the situation,  in particular on the Mordell--Weil rank $r$, the cardinality of $S$ and the size of the initial bound $M_0$.  We first suppose that $M_0$ is huge, say $M_0$ is the initial bound coming from the theory of logarithmic forms (see Section~\ref{sec:input}).  If $S$ is empty or when $r\leq 1$, then either the elliptic logarithm sieve is fast or the bottleneck is step (iii).  Assume now that $r\geq 2$ and that $S$ is nonempty.  If in addition $r\leq 4$ then the bottleneck is usually  part (a) of step (ii), in particular when $\abs{S}$ is large.  On the other hand, if in addition $r\geq 5$ then the bottleneck is either step (iii) or part (a) of step (ii).  If $M_0$ is not that large then the bottleneck is usually step (iii). For example, in the case when \eqref{eq:weieq} is a Mordell equation, the initial bounds of Proposition~\ref{prop:mwbound} are strong enough such that either the elliptic logarithm sieve is fast or the bottleneck is step (iii).

\begin{remark}[Generalizations]\label{rem:elllogsievegen}
(i) Algorithm~\ref{algo:elllogsieve} allows to solve Diophantine equations which are a priori more general than \eqref{eq:weieq}. For example, our algorithm can be applied to find all $S$-integral solutions with bounded height of any Weierstrass equation \eqref{eq:weieq} of $E$ with coefficients $a_i$ in $\QQ$. Indeed on multiplying  equation \eqref{eq:weieq} with the sixth power $u^6$ of the least common multiple $u$ of the denominators of the $a_i$, one obtains a Weierstrass equation \eqref{eq:weieq}$^*$ with coefficients $a_i^*=u^ia_i$ in $\ZZ$ and then one checks for each $(x,y)\in \Sigma^*(S)$ whether $u^{-2}x$ and $u^{-3}y$ are in $\mathcal O$. Here $\Sigma^*(S)$ denotes the set of solutions of \eqref{eq:weieq}$^*$ in $\mathcal O\times\mathcal O$ obtained by applying Algorithm~\ref{algo:elllogsieve} with the coefficients $a_i^*$, with the same initial bound $M_0$ and with the transformed coordinates of the same basis $P_i$.

(ii) The above Algorithm~\ref{algo:elllogsieve} works equally well with any initial bound for the usual Weil height $h$ or for the infinity norm $\|\cdot\|_\infty$. Indeed the explicit inequalities  \eqref{eq:nthcompa} and \eqref{eq:nthlowerbound}  translate any initial bound for $h$ or $\|\cdot\|_\infty$ into an initial bound for $\hat{h}$.

(iii) We mention that various authors (including Stroeker, Tzanakis and de Weger)  generalized and modified the known elliptic logarithm approach in order to efficiently solve  more general Diophantine equations defining genus one curves. For an overview we refer to the discussions in  Stroeker--Tzanakis~\cite{sttz:genus1} and Tzanakis~\cite{tzanakis:book}.
\end{remark}

\subsubsection{Comparison with the known approach}\label{sec:comparisonwithelr}

To discuss the improvements provided by the elliptic logarithm sieve, we now compare our sieve with the known approach via elliptic logarithms. We recall that the main steps of the known method are as follows (see for example \cite[Sect 4]{pezigehe:sintegralpoints} or \cite{tzanakis:book}):
\begin{itemize}
\item[(1)] As explained in Section~\ref{sec:elllintroelr}, one tries to obtain a reduced bound $N_1$ which is as small as possible such that any non-exceptional point $P\in \Sigma(S)$ satisfies $\|P\|_\infty\leq N_1$. 
\item[(2)] One goes through all points $P\in E(\QQ)$ with $\|P\|_\infty\leq N_1$ and one tests whether $P$ lies in fact in $\Sigma(S)$. In the case $S^*=\{\infty\}$, one can apply in addition the known inequality trick (see Remark~\ref{rem:inequtrick}) before one tests whether $P$ lies in $\Sigma(S)$.
\item[(3)] One makes a so-called extra search to find all exceptional points.  
\end{itemize}

\paragraph{Reduction.} Steps (i)+(ii) of our elliptic logarithm sieve may be viewed as an analogue of (1).  Here the running times of (1) and (i)+(ii) are essentially equal, since the aspects in which the two approaches differ are irrelevant for the running time. Indeed in both approaches the running time is essentially determined by the computations of the involved elliptic logarithms and these computations are the same in both approaches. On the other hand, in view of the subsequent enumeration, an important difference is that (1) uses the inequality $\lambda\|P\|_\infty^2\leq \hat{h}(P)$ in order to work with $\|\cdot\|_\infty$, while our reduction in (ii) directly works with $\hat{h}$.   Geometrically, this means that (1) uses a cube which always contains the ellipsoids $\mathcal E$ used in (ii).  In fact there exist non-trivial improvements of (1), see Stroeker--Tzanakis~\cite{sttz:elllogoverview} which optimizes the Mordell--Weil basis and see Hajdu--Kov\'acs \cite{hako:elllog} which in the case $S=\emptyset$ intersects cubes containing $\mathcal E$.  Our reduction in (ii) is always as good as these improvements of (1), since we work directly with the ellipsoid $\mathcal E$ which is optimal from a geometric point of view. We define $N_{\textnormal{opt}}=\lfloor(M_1/\lambda)^{1/2}\rfloor$ with $\lambda$ coming from a Mordell--Weil basis which is optimized in the sense of \cite{sttz:elllogoverview}.  For instance, the usual elliptic logarithm reduction can not reduce anymore  (\cite[p.400]{pezigehe:sintegralpoints}) the bound $N_1=17$ in the example of \cite{pezigehe:sintegralpoints} involving an optimized Mordell--Weil basis ($r=4$). On the other hand, our reduction in (ii) gives $M_1\leq 61$ which implies that $N_{\textnormal{opt}}=12$.  Furthermore, in many important situations, the reduction in (ii) is considerably stronger than the known improvements of (1). For example if $r\geq 2$ becomes large then the volume of $\mathcal E$ becomes much smaller than the volume of any cube  containing $\mathcal E$. Hence (ii) is significantly more efficient when  $r\geq 2$ is large. In particular, in the generic case we obtain here a running time improvement by a factor which is exponential in terms of $r$. Furthermore, in the most common nontrivial case (where $r=1$), we obtain huge running time improvements for large $\abs{S}$ by using the following trick: The idea is that we do not need to know the involved $v$-adic elliptic logarithms $(\beta_{1,v}^{(c)}, v\in S)$ appearing in Section~\ref{sec:nonarchisieve}. Indeed on exploiting that the involved lattice has rank $r=1$, one observes that it suffices here to know the orders $(v(\alpha_{1,v}),v\in S)$. These orders can always be efficiently computed in practice.

\paragraph{Enumeration.}   Step (iii) of our sieve plays the role of (2)+(3). In the following comparison, we denote by $t$ and $t^*$  the running times of (iii) and (2)+(3) respectively. 
 Comparing $t$ with $t^*$ is suitable to illustrate our improvements. Indeed it takes into account that the running times of (i)+(ii) and (1) are essentially equal and  it makes the comparison independent of the initial bound.   We denote by $t_2$  the running time of our refined enumeration in Algorithm~\ref{algo:refenu}~(i)  applied with our reduced initial bound $b=M_1$ from (ii).  In practice the running time of (2)  always exceeds $t_2$, which means that $t^*> t_2$.  In fact, in many cases of interest (e.g. when  $r$, $\abs{S}$ or $\max \abs{a_i}$ is not small), the time $t^*$ is often considerably larger than $t_2$. To explain the improvements provided by our sieve, we mention four situations in which the known enumeration (2)+(3) usually becomes very slow:
\begin{itemize}
\item[(S1)] Suppose that  $r$ is not small and $S$ is nonempty. Then the enumeration in (2)  becomes very slow or even hopeless, since one has to consider $(2N_1+1)^r$ points of $E(\QQ)$.
\item[(S2)] Assume that $r$ is large and $S$ is empty. Then the enumeration in (2) becomes often hopeless, since one has to go through $(2N_1+1)^r$ points of $E(\QQ)$. 
\item[(S3)] Suppose that the height $\max \abs{a_i}$ of \eqref{eq:weieq} is large. Then the extra search in (3) has to test many pairs $(x,y)$ of $S$-integers whether they satisfy \eqref{eq:weieq}. Thus (3) becomes very slow when $\max\abs{a_i}$ is large, in particular if $S$ is nonempty.
\item[(S4)] Suppose that $\abs{S}$ and $N_1$ are both large. Then (2) becomes very slow, since one has to compute many rational numbers $x(P)$ of huge height and one has to test whether they are $S$-integral. Here $x(P)$ is the $x$-coordinate of some $P\in E(\QQ)$. 
\end{itemize}
 To deal with situation (S1) we developed the global sieve which combines our refined sieve and our refined enumeration.  We recall from the discussions in previous sections that the first steps of the global sieve are essentially a further reduction of $M_1$, while the subsequent steps of the global sieve are also more efficient than the corresponding enumerations in (2) and (3). To illustrate our improvements in practice, we consider the Mordell curve  $y^2=x^3+1358556$ of rank $r=6$ and three additional examples which shall be further discussed in Section~\ref{sec:largerankexamples} below. These three additional examples are given by Kretschmer's curve  in \cite[Ex.3.3]{siksek:infdescent} with $r=8$, Mestre's curve  in \cite[Ex.3.2]{siksek:infdescent} with (conditional) $r=12$   and the curve of Fermigier with $r=14$.  In the following table, the entries of the first and second row are  rounded up and down respectively.

\begin{table}[h]
\begin{center}
{\small
\begin{tabular}{lccccccccccccccc}
$(r,S)$ & $(6,S(10))$ & $(6,S(20))$ & $(8,S(10))$ & $(12,S(5))$ & $(14,S(2))$\\
\cmidrule(r){2-6}
$t$ & 3.5m & 24m & 3.4h & 12h & 15h\\
$t_2$  & 115m & 1116m & ?  & ?  & ?
\\
$N_{\textnormal{opt}}$ & 15  & 21 & 36  & 36  & 16
\end{tabular}
}
\end{center}
\end{table} 
\noindent  We point out that here the running times of (2) would be significantly larger than the listed  times $t_2$, since $N_1\geq N_{\textnormal{opt}}$ and since  $(2N_{\textnormal{opt}}+1)^r$ is huge in each case. In particular the enumeration (2) would be very slow in the above cases involving $r=6$, while the  cases involving $r=8,12,14$ seem to be completely  out of reach for (2). 

 Furthermore our global sieve leads in addition to a significant improvement in situation (S2) where $r$ is large and $S$ is empty. To illustrate our improvements in practice, we consider again Mestre's curve of rank $r=12$, Fermigier's curve of rank $r=14$ and two curves of Elkies of rank $r=17,19$ respectively.  These four curves shall be  discussed in more detail in Section~\ref{sec:largerankexamples} below. In the case when $r=12,14,17,19$, it turned out that $t$ is less than 2 minutes, 13 minutes, 26 hours, 73 hours  respectively and that $N_{\textnormal{opt}}= 14,10,22,15$  respectively. It seems that all these cases are out of reach for the enumeration (2).

To deal with situation (S3) where  $\max \abs{a_i}$  is large, we constructed the refined covering in Proposition~\ref{prop:refinedcov}. This result allows us to work entirely in the finitely generated abelian group $E(\QQ)$ in order to find the exceptional points. We illustrate our improvements in practice by considering a rather randomly chosen Mordell equation \eqref{eq:weieq} of rank $r=1$ with $a_6=-17817895$. This choice is suitable, since the assumption $r=1$ is satisfied in the nontrivial generic case and since our sieve does not exploit special properties of Mordell equations. Furthermore, if the rank $r$ is large then  $\max \abs{a_i}$ is large and the above table already contains times $t$ for $r\geq 6$.  
To obtain the lower bounds for $t^*$ listed in the following table, we used the running times of the extra search (3); this search was implemented in Sage by Cremona, Mardaus and Nagel following the presentation in \cite[p.393]{pezigehe:sintegralpoints}. 

\begin{table}[h]
\begin{center}
{\small
\begin{tabular}{lccccccccccccccc}
$S$ & $S(1)$ & $S(2)$ & $S(3)$ & $S(4)$ &  $S(5)$ &  $S(10)$ & $S(100)$\\
\cmidrule(r){2-8}
$t $  & 1s & 1.1s & 1.1s & 1.1s & 1.2s & 1.6s & 7s \\
$t^*$ & 1.3s & 34s & 15m & 12.5h & ? & ? & ? 
\end{tabular}
}
\end{center}
\end{table}
\noindent  Here the entries of the first and second row are upper and lower bounds for $t$ and $t^*$ respectively. We note that $t^*$ explodes when $\max\abs{a_i}$ becomes larger. For example, let us consider the Mordell equation \eqref{eq:weieq} of rank $r=1$ with $a_6=- 4211349581402184375$. Here our running times $t$ essentially coincide with the times displayed in the above table, while the extra search (3) did not terminate within 48 hours in the simplest case $S=S(1)$. We mention that our improvements in the case of large $\max\abs{a_i}$ are crucial for the computation (see Section~\ref{sec:shaf}) of elliptic curves over $\QQ$ with good reduction outside $S$. Indeed these computations usually require to find all $S'$-integral solutions of many distinct Mordell equations $y^2=x^3+a$, with $a\in\ZZ$ having huge $\abs{a}\geq 10^{15}$ and $S'=S\cup \{2,3\}$.

 We developed various global constructions to efficiently deal with situation (S4) where $\abs{S}$ and $N_1$ are large. In particular, for many involved points $P\in E(\QQ)$ our height-logarithm sieve allows to avoid  the slow process of testing whether the coordinates of $P$ are $S$-integral.  To demonstrate our improvements in practice, we considered the   curves 37a1, 389a1, 5077a1  in Cremona's database with rank $r=1,2,3$ respectively. They have minimal conductor among all elliptic curves over $\QQ$ of rank $r=1,2,3$ respectively.

\begin{table}[h]
\begin{center}
{\small
\begin{tabular}{lccccccccccccccc}
$(r,S)$ & $(1,S(10^5))$ & $(1,S(2\cdot 10^5))$ & $(2,S(10^4))$ & $(3,S(500))$  \\
\cmidrule(r){2-5}
$t$ & 3.4h & 9.4h & 18d & 7.6h \\
$t_2$ & 150h?  & 900h? & 58d & 85h    \\
$N_{\textnormal{opt}}$ & 11263 & 16015 & 1467 & 246
\end{tabular}
}
\end{center}
\end{table}
\noindent Here the entries of the first row are rounded up. Further,  it is reasonable to expect that the first two entries in the $t_2$ row are larger than 150 hours and 900 hours respectively. Indeed, it took 192 seconds (resp. 777 seconds) to determine whether the coordinates of the point $11000P$ (resp. $16015P$) are $S(10^5)$-integral (resp. $S(2\cdot 10^5)$-integral), where $P\ZZ=E(\QQ)$ and $E$ denotes the rank one curve 37a1 used in the above table. We mention that the ability of Algorithm~\ref{algo:elllogsieve} to solve \eqref{eq:weieq} for large sets $S$ was crucial for obtaining  data motivating our conjectures in Section~\ref{sec:malgoapplications} and in Section~\ref{sec:elllogsieveapp} below.

\subsection{Input data}\label{sec:input}

We continue our notation. The elliptic logarithm sieve requires an initial height bound for the points in $\Sigma(S)$ and a Mordell--Weil basis of $E(\QQ)$. In this section we  recall some results and techniques which allow to compute the required input data in practice.

\subsubsection{Mordell--Weil basis}\label{sec:compmwbasis}

The problem of finding a Mordell--Weil basis of $E(\QQ)$ is difficult in theory and in practice. In fact in the case of an arbitrary elliptic curve $E$ over $\QQ$ there is so far no unconditional method which allows in principle to determine a Mordell--Weil basis of $E(\QQ)$. However, thanks to the work of many authors, it is usually possible to compute such a basis in practice. In fact it turned out in practice that the  methods implemented in Pari, Sage and Magma are remarkably efficient in computing such a basis, even in the case when the height $\max\abs{a_i}$ of \eqref{eq:weieq} is large. Furthermore, Cremona's database contains a Mordell--Weil basis of $E(\QQ)$ for each elliptic curve $E$ over $\QQ$ with conductor at most $350 000$.    Unless mentioned otherwise, we shall use below these bases of Cremona's database.
 
\subsubsection{Initial height bounds}

Starting with the works of Baker~\cite{baker:contributions,baker:mordellequation,baker:elliptic}, there is a long tradition of establishing explicit height bounds for the points in $\Sigma(S)$ using lower bounds for linear forms in logarithms. See for example Baker--W\"ustholz~\cite{bawu:logarithmicforms} for an overview and a discussion of the state of the art.  Furthermore Masser~\cite{masser:ellfunctions} and W\"ustholz~\cite{wustholz:recentprogress} initiated an  approach which provides explicit height bounds for the points in $\Sigma(S)$ using lower bounds for linear forms in elliptic logarithms. Here the actual best lower bounds can be found in the works of David~\cite{david:elllogmemoir} and Hirata-Kohno~\cite{hirata-kohno:p-adicelllogs}.

To obtain our results discussed in the next section, we compute two initial height bounds for the points in $\Sigma(S)$ and then we take the minimum of these two bounds. More precisely, the first initial height bound  is a direct consequence of  the results of Hajdu--Herendi~\cite{hahe:elliptic}  combined with a height comparison in the style of \eqref{eq:nthcompa}, see  Peth{\H{o}}--Zimmer--Gebel--Herrmann~\cite[Thm]{pezigehe:sintegralpoints}. Here the results of Hajdu--Herendi ultimately rely on lower bounds for linear forms in complex and $p$-adic logarithms.  The second initial height bound  depends on the above mentioned lower bounds for linear forms in elliptic logarithms due to David in the archimedean case and due to\footnote{In fact explicit lower bounds for linear forms in two nonarchimedean elliptic logarithms were established in the work of R\'emond--Urfels~\cite{reur:padicelllog}.} Hirata-Kohno in the nonarchimedean case. See for example the proof of Tzanakis~\cite[Thm 11.2.6]{tzanakis:book}; here  one has to take into account that in some cases the normalizations  in \cite{tzanakis:book} do not coincide with our corresponding normalizations  used in Sections~\ref{sec:archisieve} and \ref{sec:nonarchisieve}. In fact, unless mentioned otherwise, we obtained all applications in Section~\ref{sec:elllogsieveapp} by using the first height bound.

\subsection{Applications}\label{sec:elllogsieveapp}

In this section we discuss additional applications of the elliptic logarithm sieve. To obtain the input data required for the  applications of our sieve, we used unless mentioned otherwise the results  described in the previous section. We continue the  notation introduced above  and for any $n\in\ZZ_{\geq 1}$  we denote by $S(n)$ the set of the first $n$ rational primes.

\subsubsection{Elliptic curves database}

 For any elliptic curve $E$ over $\QQ$ of conductor at most 1000, Cremona's database contains in particular a minimal Weierstrass equation \eqref{eq:weieq} of $E$. We used Algorithm~\ref{algo:elllogsieve} to compute the $S$-integral solutions of each of these equations with $S=S(20)$.  Moreover, for any of these minimal equations defining an elliptic curve over $\QQ$ of conductor at most 100, we used Algorithm~\ref{algo:elllogsieve}  to determine its set of $S$-integral solutions with $S=S(10^4)$.

\subsubsection{Elliptic curves of large rank}\label{sec:largerankexamples}

In addition, we used the elliptic logarithm sieve (Algorithm~\ref{algo:elllogsieve}) in order to determine the set of $S$-integral solutions of various Weierstrass equations \eqref{eq:weieq} for which the involved Mordell--Weil rank $r$ of $E(\QQ)$ is relatively large. We now discuss some examples.

\paragraph{Mordell curves.} Recall that \eqref{eq:weieq} is called a Mordell equation if the coefficients $a_1,\dotsc,a_5$ are all zero.   In Section~\ref{sec:malgoapplications} we used the elliptic logarithm sieve to find all $S$-integral solutions of two (resp. four) Mordell equations with $S=S(50)$ and $r=7$ (resp. $S=S(40)$ and $r=8$).  Instead of using initial bounds coming from the theory of logarithmic forms, we applied here our optimized height bound in Proposition~\ref{prop:mwbound} which is based on the method of Faltings (Arakelov, \parshin, Szpiro) \cite{faltings:finiteness} combined with the Shimura--Taniyama conjecture~\cite{wiles:modular,taywil:modular,breuil:modular}. Further, we determined here the required Mordell--Weil bases by using techniques implemented in Pari, Sage and Magma.

\paragraph{Rank eight.} Let $E_{\textnormal{Kr}}$ be the elliptic curve over $\QQ$ considered in \cite[Example 5.3]{siksek:infdescent}, with $E_{\textnormal{Kr}}(\QQ)$ of rank $r=8$  by Kretschmer \cite{kretschmer:largerank}. Siksek~\cite{siksek:infdescent} combined his refined descent techniques with Cremona's ``mwrank" to find a Mordell--Weil basis of $E_{\textnormal{Kr}}(\QQ)$.  On using Siksek's basis as an input for our sieve, we determined all $S$-integral solutions of the minimal Weierstrass equation \eqref{eq:weieq} of $E_{\textnormal{Kr}}$ with $S=S(10)$. This took our sieve less than 35 seconds, 17 hours and 75 hours for $S=\emptyset$, $S=S(8)$ and $S=S(10)$ respectively. Here we notice that our running times for $E_{\textnormal{Kr}}$ are significantly worse than for the four Mordell curves of rank $r=8$ discussed in Section~\ref{sec:malgoapplications}. The reason is that in the case of $E_{\textnormal{Kr}}$ we need to use initial height bounds based on the theory of logarithm forms.  These initial bounds are substantially weaker (see Section~\ref{sec:minitbounds}) than our optimized height bound in Proposition~\ref{prop:mwbound} which is currently only available for Mordell curves.

\paragraph{Rank twelve.}  Mestre~\cite{mestre:rank12} constructed an elliptic curve $E_{\textnormal{Me}}$ over $\QQ$ of analytic rank 12, together with 12 independent points in $E_{\textnormal{Me}}(\QQ)$ of infinite order.  Here again Siksek~\cite[Example 5.2]{siksek:infdescent} applied his refined descent techniques to find a Mordell--Weil basis of $E_{\textnormal{Me}}(\QQ)$. On using Siksek's basis as an input for our sieve, we determined the set of $S$-integral solutions of a minimal Weierstrass equation \eqref{eq:weieq} of $E_{\textnormal{Me}}$ with $S=S(7)$. This took our sieve less than 2 minutes, 4 days and 16 days for $S=\emptyset$, $S=S(6)$ and $S=S(7)$ respectively. We point out that the completeness of our solution sets are here conditional on Siksek's assumption that $E_{\textnormal{Me}}(\QQ)$ has rank $r=12$ which he used in his construction of a basis of $E_{\textnormal{Me}}(\QQ)$. For example, this assumption is satisfied if $r$ is at most the analytic rank of $E_{\textnormal{Me}}$ as predicted by the rank part of the Birch--Swinnerton-Dyer conjecture.   

\paragraph{Rank at most 28.} At the time of writing, we are not aware of an elliptic curve over $\QQ$ of rank $r\geq 13$ for which an explicit Mordell--Weil basis can be computed explicitly.  If $r\leq 28$ and $S$ is empty then the elliptic logarithm sieve would allow to determine $\Sigma(S)$ for such large rank elliptic curves $E$ over $\QQ$, provided that one knows an explicit Mordell--Weil basis of $E(\QQ)$. To demonstrate this feature, we work with independent points $Q_1,\dotsc,Q_r$ generating a rank $r$ subgroup $\Lambda$ of the free part of $E(\QQ)$.  Dujella lists in particular such points (see \url{web.math.pmf.unizg.hr/~duje}) for three elliptic curves over $\QQ$ of rank  $r=14,17,19$ respectively.  Here the curve of rank $r=14$ was constructed by Fermigier, while the other two curves were found by Elkies.  We denote by $\Sigma_\Lambda(S)$  the intersection of $\Sigma(S)$ with $\Lambda\oplus E(\QQ)_{\textnormal{tor}}$.  On using the basis $Q_1,\dotsc,Q_r$ of $\Lambda$ in the input of  Algorithm~\ref{algo:elllogsieve}, we determined the set $\Sigma_\Lambda(S)$.  If $S$ is empty then this took less than 13 minutes, 26 hours, 73 hours for $r=14,17,19$ respectively, and it took less than 18 hours  when $r=14$ and $S=S(2)$. Now, given a Mordell--Weil basis,  one can expect similar running times of Algorithm~\ref{algo:elllogsieve} when computing the full set $\Sigma(S)\supseteq \Sigma_\Lambda(S)$. Indeed  the index of $\Lambda$  in the free part of $E(\QQ)$ is not that large  for these three curves. 
Finally, we consider Elkies' elliptic curve $E_{\textnormal{El}}$ over $\QQ$ of rank $r\geq 28$. In the case when $S$ is empty, we applied Algorithm~\ref{algo:heightlogsieve} with the  28 independent points $Q_i$ constructed by Elkies and we  computed the 
intersection of $\Sigma(S)$ with $(\oplus_{i}Q_i\ZZ)\oplus E_{\textnormal{El}}(\QQ)_{\textnormal{tor}}$ in less than  75 days.

\subsubsection{Conjectures and questions}
In this section we generalize our conjectures and questions for Mordell curves (see Section~\ref{sec:malgoapplications}) to hyperbolic genus one curves. We then provide some motivation by using our data and by generalizing our  constructions for Mordell curves. See also Section~\ref{sec:ia} which contains the initial motivation for our conjectures and questions.

We continue our notation.  As in Section~\ref{sec:elllintroapp}, we let  $Y=(X,D)$ be a hyperbolic genus one curve over some open subscheme $B$ of $\sp(\ZZ)$  and we denote by $r$  the rank of the group formed by the $\QQ$-points of $\textnormal{Pic}^0(X_\QQ)$. Now we recall our conjecture.

\vspace{0.3cm}
\noindent{\bf Conjecture.}
\emph{There are constants $c_{Y}$ and $c_r$, depending only on $Y$ and $r$ respectively, such that any nonempty finite set of rational primes $S$ with $T=\sp(\ZZ)-S$ satisfies} $$\abs{Y(T)}\leq c_{Y} \abs{S}^{c_r}.$$

\noindent If $D$ is given by a section of $X\to B$, then $Y$ identifies with a closed subscheme of $\mathbb A^2_B$ defined by a Weierstrass equation~\eqref{eq:weieq}. 
Hence the family $S[b]$ constructed in Section~\ref{sec:malgoapplications} shows that the exponent $c_r$ has to be at least $\tfrac{r}{r+2}$ when $Y$ is a Weierstrass curve. 

\vspace{0.3cm}
\noindent{\bf Question 1.}
\emph{What is the optimal exponent $c_r$ in the above conjecture?}
\vspace{0.3cm}

\noindent Our data strongly indicates that the exponent $c_r=\tfrac{r}{r+2}$ is still far from optimal for many families of sets $S$ of interest, including the family $S(n)$ with $n\in\ZZ_{\geq 1}$. Further, our data  motivates in addition the following question on the dependence on $q=\max S$.

\vspace{0.3cm}
\noindent{\bf Question 2.}
\emph{Are there constants $c_{Y}$ and $c_r$, depending only on $Y$ and $r$ respectively, such that any nonempty finite set of rational primes $S$ with $T=\sp(\ZZ)-S$ satisfies} $$\abs{Y(T)}\leq c_{Y}(\log \max S)^{c_r} \, \textnormal{?}$$

\noindent   In the case when $S=S(n)$ with $n\in\ZZ_{\geq 2}$, one can replace here $q$ by $n\log n$ without changing the content of the question. However  Question~2 has in general a negative answer when $q$ is replaced by any power of $\max(2,\abs{S})$.  On using again the arguments of Section~\ref{sec:malgoapplications}, we see that  the exponent $c_r$ of  Question~2 has to be at least $r/2$ if $Y$ is a Weierstrass curve. In light of this we ask whether Question~2 has a positive answer with the exponent 
\begin{equation*}
c_r=r/2 \, \textnormal{?}
\end{equation*}
 Our probabilistic model constructed in the discussion surrounding \eqref{refquestbound} predicts that this question has a positive answer when $Y$ is a Weierstrass curve. To conclude  we mention that additional motivation for our conjecture is given by the theory of logarithmic forms \cite{bawu:logarithmicforms}, which was applied in \cite{rvk:ed} to obtain new  bounds for the number of integral points on arbitrary hyperbolic genus one curves over any number field. These bounds establish in particular our Conjecture for certain sets $S$ of interest, including the sets $S(n)$.

\subsection{Computational aspects}\label{sec:compuaspects}

In this final section we discuss various computational aspects of our algorithms. In particular we explain the numerical details of  constructions used in  previous sections.

\paragraph{Interval arithmetic.} Most real numbers $x$ cannot be explicitly presented by a computer. Hence we apply standard interval arithmetic which uses  lower and upper bounds for $x$. This allows us to control numerical errors. In particular, one can detect when the error explodes and in this case one can restart the computation with higher precision.

\paragraph{Rounding real numbers.}   We consider a real number $x\in\RR$. In this paper the symbol $[x]$ denotes an element of $\ZZ$ with $\abs{x-[x]}<1$. There might be two distinct choices for $[x]\in \ZZ$ with the desired property.  However, for our purpose any choice will  be sufficient. In the implementation, real numbers $x$ are only stored up to some chosen precision. Suppose that $x$ is stored with absolute precision $m\in\ZZ_{\geq 1}$, that is the computer stores $x$ as $x'\in\QQ$ with $\abs{x-x'} \leq \tfrac{1}{2^{m}}$.  Then we can compute $[x]$ as a closest integer to $x'$; there may be two choices, in which case we pick one. This $[x]\in\ZZ$ moreover satisfies $\abs{x-[x]}\leq \tfrac{1}{2}+\tfrac{1}{2^m}$.

\paragraph{Construction and eigenvalues of $\hat{h}_k$.} To discuss the numerical details required for the construction of $\hat{h}_k$, we continue the notation and terminology of Section~\ref{sec:heightselllog}. We recall that we need to find a suitable $k\in\ZZ_{\geq 1}$ for which $\hat{h}_k$ is close enough to $\hat{h}$.  To find such a $k$, we start with some $k$ such as for example  $k=10$. We first determine a reasonably good rational approximation $f\in\QQ$ of $2^k/\|\hat{h}_{ij}\|$.  After computing a lower bound  for the smallest eigenvalue $\lambda_k$ of $\hat{h}_k$, we can check whether $\lambda_k>0$ and  $r<\tfrac{1}{100}f\lambda_k$. If both conditions are satisfied, then $k$ is suitable. Otherwise we increase $k$ and we repeat the above procedure as long as required to find a $k$ with the desired properties.   Lemma~\ref{lem:hk} implies that this procedure terminates; in fact it terminates always very quickly in practice.  While the condition $r<\tfrac{1}{100}f\lambda_k$ is in principle not required for the correctness of our algorithms, it assures via Lemma~\ref{lem:hk} that $\hat{h}_k$ is relatively close to $\hat{h}$, which is reasonable in practice.

To explain how we compute a reasonably good lower bound for $\lambda_k$, we recall that the quadratic form $\hat{h}_k$ is given by $A/f$ with $A=([f\hat{h}]-r\cdot \textnormal{id})\in\ZZ^{r\times r}$. First, we apply Newton's method to compute a root of the characteristic polynomial of $A$. Here we start at minus the Cauchy bound for the largest absolute value of its roots.
At each step of Newton's method we stay below all eigenvalues. Indeed
 the characteristic polynomial is a convex or concave function on the interval $(-\infty,f\lambda_k)$, since all $r$ eigenvalues of $A$ are real  as $A$ is symmetric.
In order to remove numerical issues, we use interval arithmetic. After each step of Newton's method, we moreover replace the point by its lower bound. Once this value is smaller than or equal to the one in the previous step, we terminate. An upper bound for the largest eigenvalue of $A$ can be obtained in an analogous way.

\paragraph{Choice of $\delta$.} To explain our choice of the parameter $\delta$ in \eqref{def:delta12}, we continue the notation and terminology of Section~\ref{sec:archisieve}.  We choose $a,b\in\ZZ$  such that $\tfrac{a}{b}$ is a good lower approximation of $M/(\delta_1+\delta_2)^2$, such that $\tfrac{a}{b}=M/\delta^2$ for some positive $\delta\in \QQ$ and such that $\abs{a},\abs{b}$ are not too large.  These are two simultaneous objective functions that we try to optimize: 1. Making the inequality $\tfrac{a}{b} \leq M/(\delta_1+\delta_2)^2$ as sharp as possible assures that the Fincke--Pohst algorithm does not return too many additional candidates. 2. The smaller $|a|$ and $|b|$, the smaller are the entries of the quadratic form on which we run the Fincke--Pohst algorithm.  Further,  this construction of $\delta$ (that is of the integers $a,b$) allows us to apply our implementation of the Fincke--Pohst lattice point enumeration which works with integral quadratic forms.  Indeed on multiplying by the (controlled) denominator $bf$ of the quadratic form $q$ used in the definition of $\mathcal E$, one obtains an integral quadratic form.

\paragraph{Lattice point enumeration.}
To enumerate lattice points in an ellipsoid, we apply the version (FP) of Fincke--Pohst~\cite{fipo:algo} which in turn uses $L^3$ \cite{lelelo:lll}. More precisely, for any $d\in \ZZ_{\geq 1}$ the algorithm (FP) takes as the input a basis of a lattice $\Gamma\subseteq\ZZ^d$ together with an ellipsoid $\mathcal E\subset\RR^d$ centered at the origin, and
it outputs the inter\-section~$\Gamma\cap \mathcal E$. In fact we use here our own implementation of (FP) which is described in Remark~\ref{rem:fp}.

\paragraph{Elliptic logarithms.}
Information on the computation of the elliptic logarithms can be found in Zagier~\cite{zagier:largeintegralpoints} for real elliptic logarithms and in Peth{\H{o}} et al~\cite{pezigehe:sintegralpoints} for $p$-adic elliptic logarithms, see also Tzanakis~\cite{tzanakis:book}.
To compute the elliptic logarithms we use Sage which in turn is based on Pari. Our normalizations of the real and $p$-adic elliptic logarithm in Sections~\ref{sec:archisieve} and \ref{sec:nonarchisieve} respectively coincide with the normalizations of Pari.

\bibliographystyle{alpha} 

{\scriptsize

\bibliography{../../literature}
}

\vspace{0.5cm}

\vspace{0.1cm}
\noindent Rafael von K\"anel, MPIM Bonn,  Vivatsgasse 7, 53111 Bonn, Germany

\vspace{0.1cm}
\noindent Current affiliation: Princeton University, Mathematics, Fine Hall, NJ 08544-1000, USA

\noindent E-mail address: {\sf rkanel@math.princeton.edu}

\vspace{0.5cm}

\noindent Benjamin Matschke, MPIM Bonn,  Vivatsgasse 7, 53111 Bonn, Germany

\vspace{0.1cm}
\noindent Current affiliation: Institut de Math\'ematiques de Bordeaux, Universit\'e de Bordeaux,  
351, cours de la Lib\'eration, 
33405 Talence, France

\noindent E-mail address: {\sf  benjamin.matschke@math.u-bordeaux.fr}

\end{document}